\documentclass[11pt,a4]{article}
\usepackage{indentfirst}
\usepackage{amssymb}
\newtheorem{theo}{THEOREM}[subsection]
\newtheorem{p}[theo]{PROPOSITION}
\newtheorem{lem}[theo]{LEMMA}
\newtheorem{de}[theo]{DEFINITION}
\newtheorem{co}[theo]{COROLLARY}
\newtheorem{e}[theo]{EXAMPLE}
\newcommand{\fr}[1]{\ensuremath{\mathfrak{#1}}}
\newcommand{\h}{Theorem }
\newcommand{\hh}{[C1] Theorem 5.6.}
\newcommand{\pr}{Proposition }
\newcommand{\prr}{[C1] Proposition 5.6.}
\newcommand{\dd}{Definition}
\newcommand{\ee}[1]{Example \ensuremath{\ref{#1}}}
\newcommand{\qedd}{\eqno{\rule{3mm}{3mm}}}
\newcommand{\qed}{\hfill {\rule{3mm}{3mm}}}
\newcommand{\mae}[5]{\ensuremath{#1:#2\longrightarrow #3\,, \quad #4 \longmapsto  #5}}
\newcommand{\mad}[4]{\ensuremath{#1\longrightarrow #2, \quad #3 \longmapsto   #4}}
\newcommand{\mac}[3]{\ensuremath{#1:#2\longrightarrow #3}}
\newcommand{\ri}[1]{Hilbert right $#1$-module}
\newcommand{\ris}[1]{selfdual Hilbert right $#1$-module}
\newcommand{\n}[1]{\ensuremath{\left\| #1 \right\|}}
\newcommand{\me}[2]{\ensuremath{\left\{\,\left. #1\;\right|\; #2 \, \right\}}}
\newcommand{\s}[2]{\ensuremath{\left \langle  \left.\, #1\,\right|\, #2  \, \right \rangle}}
\newcommand{\sa}[2]{\ensuremath{\left \langle \, #1 \,, \, #2 \, \right \rangle}}
\newcommand{\si}[1]{\ensuremath{\sum\limits_{#1}}}
\newcommand{\sii}[2]{\ensuremath{\sum\limits_{#1}^{#2}}}
\newcommand{\siw}[1]{\ensuremath{\widetilde{\sum\limits_{#1}}}}
\newcommand{\pro}[1]{\ensuremath{\prod\limits_{#1}}}
\newcommand{\proo}[2]{\ensuremath{\prod\limits_{#1}^{#2}}}
\newcommand{\ca}{\ensuremath{\mathop{\bigcirc \hspace{-2.7mm} | }}\;} 
\newcommand{\cb}[1]{\ensuremath{\mathop{\bigcirc \hspace{-2.1mm} |}\limits_{ \hspace{2mm} #1}}}
\newcommand{\cw}[1]{\ensuremath{\mathop{\bigcirc \hspace{-2.7mm} |}\limits_{ \hspace{2mm}#1}^{\;\;W}}}
\newcommand{\la}[1]{\ensuremath{\mathcal L(#1)}}
\newcommand{\lb}[2]{\ensuremath{\mathcal L_{#1}(#2)}}
\newcommand{\lc}[3]{\ensuremath{\mathcal L_{#1}(#2,#3)}}

\newcommand{\lh}[2]{\ensuremath{\mathcal L_{#1}(#2)_{\stackrel{...}{#2}}}}
\newcommand{\leh}{\ensuremath{\mathcal L_{E}(H)}}
\newcommand{\ssa}[1]{\ensuremath{\mathcal S(#1)}}
\newcommand{\ssb}[2]{\ensuremath{\mathcal S_{#1}(#2)}}

\newcommand{\ssc}[1]{\ensuremath{\mathcal Cl(#1)}}
\newcommand{\f}[2]{\ensuremath{\mathcal F(#1,#2)}}
\newcommand{\fte}{\ensuremath{\mathcal F(T,E)}}
\newcommand{\cca}[1]{\ensuremath{\mathcal C(#1)}}
\newcommand{\ccb}[2]{\ensuremath{\mathcal C(#1,#2)}}
\newcommand{\ccc}[1]{\ensuremath{\mathcal{#1}}}
\newcommand{\bk}{\ensuremath{\mathrm{I\! K}}}
\newcommand{\bn}{\ensuremath{\mathrm{I\! N}}}
\newcommand{\bnn}[1]{\ensuremath{\mathrm{I\! N_{#1}}}}
\newcommand{\br}{\ensuremath{\mathrm{I\! R}}}
\newcommand{\bc}{\ensuremath{\mathrm{I\!\!\! C}}}
\newcommand{\bz}{\ensuremath{\mathrm{\,Z\hspace{-0.65em} Z\,}}}
\newcommand{\bzz}[1]{\ensuremath{\mathrm{\,Z\hspace{-0.65em} Z_{#1}\,}}}
\newcommand{\bt}{\ensuremath{\mathrm{\,I\hspace{-0.62em}T}}}
\newcommand{\bh}{\ensuremath{\mathrm{\,I\hspace{-0.2em}H}}}
\newcommand{\ab}[4]{\ensuremath{ \left\{ \begin{array}{c@{\quad \mbox{if} \quad}c} #1&#2 \\
#3&#4 \end{array}} \right.}

\newcommand{\ac}[6]{\ensuremath{ \left\{ \begin{array}{c@{\quad \mbox{if} \quad}c} #1&#2 \\
#3&#4 \\ #5&#6 \end{array}} \right.}
\newcommand{\ad}[8]{\ensuremath{ \left\{ \begin{array}{c@{\quad \mbox{if} \quad}c} #1&#2 \\
#3&#4 \\ #5&#6  \\ #7&#8 \end{array}} \right.}
\newcommand{\ti}[3]{\ensuremath{\widetilde{\overbrace{(#1,#2,#3)}}}}
\newcommand{\tia}[3]{\ensuremath{\widetilde{(#1,#2,#3)}}}
\renewcommand{\labelenumi}{\alph{enumi})}
\renewcommand{\labelenumii}{$\alph{enumi}_{\arabic{enumii}}$)}
\newcommand{\mt}[4]{\ensuremath{\left[\begin{array}{cc}#1&#2\\#3&#4\end{array}\right]}}
\newcommand{\un}[1]{\ensuremath{Un\;#1^c}}
\newcommand{\unn}[1]{\ensuremath{Un\;#1}}
\begin{document}
\title{PROJECTIVE REPRESENTATIONS OF GROUPS \\ USING HILBERT RIGHT C*-MODULES}
\author{CORNELIU CONSTANTINESCU}
\maketitle
\begin{center}
Abstract
\end{center}

The projective representation of groups was introduced in 1904 by Issai Schur (1875-1941) in his paper [S]. It differs from the normal representation of groups (introduced by his tutor Ferdinand Georg Frobenius (1849-1917) at the suggestion of Richard Dedekind (1831-1916)) by a twisting factor, which we call Schur function in this paper and which is  called sometimes multipliers or normalized factor set in the literature (other names are also used). It starts with a group $T$ and a Schur function $f$ for $T$. This is a scalar valued function on $T\times T$ satisfying the conditions $f(1,1)=1$ and
$$|f(s,t)|=1\,,\quad\quad f(r,s)f(rs,t)=f(r,st)f(s,t)$$
for all $r,s,t\in T$. The projective representation of $T$ twisted by $f$ is a unital C*-subalgebra of the C*-algebra $\la{l^2(T)}$ of operators on the Hilbert space $l^2(T)$. This representation can be used in order to construct many examples of C*-algebras (see e.g. [C1] Chapter 7). By replacing the scalars $\br$ or $\bc$ with an arbitrary unital (real or complex) C*-algebra $E$ the field of applications is enhanced in an essential way. In this case $l^2(T)$ is replaced by the Hilbert right $E$-module $\cb{t\in T}{E}\approx E\otimes l^2(T)$ and $\la{l^2(T)}$ is replaced by $\lb{E}{E\otimes l^2(T)}$, the C*-algebra of adjointable operators of $\la{E\otimes l^2(T)}$. The projective representation of groups, which we present in this paper, has some similarities with the construction of cross products with discrete groups. It opens the way to create many K-theories.

In a first section we introduce some results which are needed for this construction, which is developed in the second section. In the third section we present examples of C*-algebras obtained by this method. Examples of a special kind (the Clifford algebras) are presented in the last section.

{\it AMS Subject Classification: 22D25 (Primary) 20C25, 46L08 (Secondary)

}

\tableofcontents

\begin{center}
\addtocounter{section}{-1}
\section{Notation and Terminology}
\end{center} 

\fbox{\parbox{12cm}{Throughout this paper we use the following notation: $T$ is a group, $1$ is its neutral element, $K:=l^2(T)$, $1_K:=id_K:=\mbox{identity map of }K$, $E$ is a unital C*-algebra (resp. a W*-algebra), $1_E$ is its unit, $\breve E$ denotes the set $E$ endowed with its canonical structure of a \ri{E} (\prr 1.5),
$$H:=\breve E\otimes K\approx \cb{t\in T}\breve E\,,\quad (\mbox{resp.}\, H:=\breve E\bar \otimes K\approx \cw{t\in T}\breve E)$$
([C3] \pr 2.1, (resp. [C3] Corollary 2.2)). In some examples, in which $T$ is additive, 1 will be replaced by 0.}}

The map
$$\mad{\lb{E}{\breve E}}{E}{u}{\s{u1_E}{1_E}=u1_E}$$
is an isomorphism of C*-algebras with inverse
$$\mad{E}{\lb{E}{\breve E}}{x}{x\cdot }\;.$$
We identify $E$ with $\lb{E}{\breve E}$ using these isomorphisms.

In general we use the notation of [C1]. For tensor products of C*-algebras we use [W], for W*-tensor products of W*-algebras we use [T], for tensor products of \ri{C^*}s we use [L], and for the exterior W*-tensor products of \ris{W^*}s we use [C2] and [C3].

In the sequel we give a list of notation used in this paper.
\renewcommand{\labelenumi}{\arabic{enumi})} 
\begin{enumerate}
\item $\bk$ denotes the field of real numbers ($:=\br$) or the field of complex numbers ($:=\bc$). In general the C*-algebras will be complex or real. $\bh$ denotes the field of quaternions, $\bn$ denotes the set of natural numbers ($0\not\in \bn$), and for every $n\in \bn\cup \{0\}$ we put
$$\bnn{n}:=\me{m\in \bn}{m\leq n}\;.$$
$\bz\;$  denotes the group of integers and for every $n\in \bn$ we put  $\;\bzz{n}:=\bz\,/(n\bz)\;.$
\item For every set $A$, $\fr{P}(A)$ denotes the set of subsets of $A$, $\fr{P}_f(A)$ the set of finite subsets of $A$, and $Card\; A$ denotes the cardinal number of $A$. If $f$ is a function defined on $A$ and $B$ is a subset of $A$ then $f|B$ denotes the restriction of $f$ to $B$.
\item If $A,B$ are sets then $A^B$ denotes the set of maps of $B$ in $A$.
\item For all $i,j$ we denote by $\delta_{i,j}$ Kronecker's symbol:
$$\delta _{i,j}:=\ab{1}{i=j}{0}{i\not=j}\;.$$
\item If $A,B$ are topological spaces then $\ccb{A}{B}$ denotes the set of continuous maps of $A$ into $B$. If $A$ is locally compact space and $E$ is a C*-algebra then $\ccb{A}{E}$ (resp. $\ccc{C}_0(A,E)$) denotes the C*-algebra of continuous maps $A\rightarrow E$, which are bounded (resp. which converge to 0 at the infinity). 
\item For every set $I$ and for every $J\subset I$ we denote by $e_J:=e_J^I$ the characteristic function of $J$ i.e. the function on $I$ equal to $1$ on $J$ and equal to $0$ on $I\setminus J$. For $i\in I$ we put $e_i:=(\delta _{i,j})_{j\in I}\in l^2(I)$.
\item If $F$ is an additive group and $S$ is a set then 
$$F^{(S)}:=\me{x\in F^S}{\me{s\in S}{x_s\not=0} \mbox{is finite}}\;.$$ 
\item If $E,F$ are vector spaces in duality then $E_F$ denotes the vector space $E$ endowed with the locally convex topology of pointwise convergence on $F$, i.e. with the weak topology $\sigma (E,F)$.
\item If $E$ is a normed vector space then $E'$ denotes its dual and $E^{\#}$ denotes its unit ball:
$$E^{\#}:=\me{x\in E}{\n{x}\leq 1}\;.$$
Moreover if $E$ is an ordered Banach space then $E_+$ denotes the convex cone of its positive elements.
If $E$ has a unique predual (up to isomorphisms), then we denote by $\ddot E$ this predual and so by $E_{\ddot E}$ the vector space $E$ endowed with the locally convex topology of pointwise convergence on $\ddot E$. 
\item The expressions of the form "... C*-... (resp. ... W*-...)", which appear often in this paper, will be replaced by expressions of the form "... C**-...". 
\item If $F$ is a unital C*-algebra and $A$ is a subset of $F$ then we denote by $1_F$ the unit of $F$, by $Pr\,F$ the set of orthogonal projections of $F$, by
$$A^c:=\me{x\in F}{y\in A\Rightarrow xy=yx}\,,\quad Re\;F:=\me{x\in F}{x=x^*}\,,$$
and by $\unn{F}$ the set of unitary elements of $F$. If $F$ is a real C*-algebra then $\stackrel{\circ }{F}$ denotes its complexification.
\item If $F$ is a C*-algebra then we denote for every $n\in \bn$ by $F_{n,n}$ the C*-algebra of $n\times n$ matrices with entries in $F$. If $T$ is finite then $F_{T,T}$ has a corresponding signification.
\item Let $F$ be a C*-algebra and $H,K$ \ri{F}s. We denote by $\lc{F}{H}{K}$ the Banach subspace of $\lc{}{H}{K}$ of adjointable operators, by $1_H$ the identity map $H\rightarrow H$ which belongs to
$$\lb{F}{H}:=\lc{F}{H}{H}\;.$$
For $(\xi ,\eta )\in H\times K$ we put 
$$\mae{\eta \s{\cdot }{\xi }}{H}{K}{\zeta }{\eta \s{\zeta }{\xi }}$$
and denote by $\ccc{K}_F(H)$ the closed vector subspace of $\lb{F}{H}$ generated by $\me{\eta \s{\cdot }{\xi }}{\xi ,\eta \in H}$.
\item Let $F$ be a W*-algebra and $H,K$ \ri{F}s. We put for $a\in \ddot F$ and $(\xi ,\eta )\in H\times K$,
$$\mae{\widetilde{(a,\xi )}}{H}{\bk}{\zeta }{\sa{\s{\zeta }{\xi }}{a}}\,,$$
$$\mae{\tia{a}{\xi }{\eta }}{\lc{F}{H}{K}}{\bk}{u}{\sa{\s{u\xi }{\eta }}{a}}$$
and denote by$\ddot H$ the closed vector subspace of the dual $H'$ of $H$ generated by 
$$\me{\widetilde{(a,\xi )}}{a\in \ddot F,\;\xi \in H}$$
and by $\stackrel{...}{H}$ the closed vector subspace of $\lc{F}{H}{K}'$ generated by 
$$\me{\tia{a}{\xi }{\eta }}{(a,\xi ,\eta )\in \ddot F\times H\times K}\;.$$
If $H$ is selfdual then $\stackrel{...}{H}$ is the predual of $\lb{F}{H}$ (\hh 3.5 b)) and $\ddot H$ is the predual of $H$ (\prr 3.3). Moreover a map defined on $F$ is called W*-continuous if it is continuous on $F_{\ddot F}$. If $G$ is a W*-algebra a C*-homomorphism $\varphi :F\rightarrow G$ is called a W*-homomorphism if the map $\varphi :F_{\ddot F}\rightarrow G_{\ddot G}$ is continuous; in this case $\ddot \varphi $ denotes the pretranspose of $\varphi $.
\item If $F$ is a C**-algebra and $(H_i)_{i\in I}$ a family of \ri{F}s then we put
$$\cb{i\in I}H_i:=\me{\xi \in \pro{i\in I}H_i}{\mbox{the family}  \s{\xi _i}{\xi _i}_{i\in I} \mbox{is summable in}\; 
F}$$
respectively
$$\cw{i\in I}H_i:=
\me{\xi \in \pro{i\in I}H_i}{\mbox{the family}  \s{\xi _i}{\xi _i}_{i\in I} \mbox{is summable in}\; F_{ \ddot F}}\;.$$
\item $\odot $ denotes the algebraic tensor product of vector spaces.
\item If $F,G$ are W*-algebras and $H$ (resp. $K$) is a \ris{F} (resp. G-module) then we denote by $H\bar \otimes K$ the W*-tensor product of $H$ and $K$, which is a \ris{F\bar \otimes G} ([C2] Definition 2.3).
\item $\approx $ denotes isomorphic.

\end{enumerate}

If $T$ is finite then (by \hh6.1 f))
$$\lb{E}{H}=E_{T,T}=\bk_{T,T}\otimes E=\ccc{K}_E(H)\;.$$
\renewcommand{\labelenumi}{\alph{enumi})} 
\renewcommand{\labelenumii}{\alph{enumi}_{\arabicenumii}))}
\begin{center}
\section{Preliminaries}

\subsection{Schur functions}
\end{center}
\begin{de}\label{703}
A {\bf Schur $E$-function for $T$ } is a map
$$\mac{f}{T\times T}{\un{E}}$$
such that $f(1,1)=1_E$ and
$$f(r,s)f(rs,t)=f(r,st)f(s,t)$$
for all $r,s,t\in T$. We denote by $\f{T}{E}$ the set of Schur $E$-functions for $T$ and put
$$\mae{\tilde f}{T}{\un{E}}{t}{f(t,t^{-1})^*}\,,$$
$$\mae{\hat f}{T\times T}{\un{E}}{(s,t)}{f(t^{-1},s^{-1})}$$
for every $f\in \f{T}{E}$.
\end{de}

Schur functions are also called normalized factor set or multiplier or two-co-cycle (for $T$ with values in $\un{E}$) in the literature. We present in this subsection only some elementary properties (which will be used in the sequel) in order to fix the notation and the terminology. By the way, $\un{E}$ can be replaced in this subsection by an arbitrary commutative multiplicative group.

\renewcommand{\labelenumi}{\alph{enumi})} 
\renewcommand{\labelenumii}{$\alph{enumi}_{\arabicenumii})$)}
\begin{p}\label{704}
Let $f\in \f{T}{E}$.
\begin{enumerate}
\item For every $t\in T$,
$$f(t,1)=f(1,t)=1_E\,,\quad f(t,t^{-1})=f(t^{-1},t)\,,\quad \tilde f(t)=\tilde f(t^{-1})\;.$$
\item For all $s,t\in T$,
$$f(s,t)\tilde f(s)=f(s^{-1},st)^*\,,\qquad f(s,t)\tilde f(t)=f(st,t^{-1})^*\;.$$
\end{enumerate}
\end{p}

a) Putting $s=1$ in the equation of $f$ we obtain
$$f(r,1)f(r,t)=f(r,t)f(1,t)$$
so
$$f(r,1)=f(1,t)$$
for all $r,t\in T$. Hence
$$f(t,1)=f(1,t)=f(1,1)=1_E\;.$$

Putting $r=t$ and $s=t^{-1}$ in the equation of $f$ we get
$$f(t,t^{-1})f(1,t)=f(t,1)f(t^{-1},t)\,.$$
By the above,
$$f(t,t^{-1})=f(t^{-1},t)\,,\quad \tilde f(t)=\tilde f(t^{-1})\,.$$

b) Putting $r=s^{-1}$ in the equation of $f$, by a),
$$f(s,t)f(s^{-1},st)=f(s^{-1},s)f(1,t)=\tilde f(s)^*\,,$$
$$f(s,t)\tilde f(s)=f(s^{-1},st)^*\;.$$
Putting now $t=s^{-1}$ in the equation of $f$, by a) again,
$$f(r,s)f(rs,s^{-1})=f(r,1)f(s,s^{-1})=\tilde f(s)^*\,,$$
$$f(r,s)\tilde f(s)=f(rs,s^{-1})^*\,,\qquad f(s,t)\tilde f(t)=f(st,t^{-1})^*\;.\qedd$$

\begin{de}\label{705}
We put
$$\Lambda (T,E):=\me{\mac{\lambda }{T}{\un{E}}}{\lambda (1)=1_E}$$
and
$$\mae{\hat \lambda }{T}{\un{E}}{t}{\lambda (t^{-1})}\,,$$
$$\mae{\delta \lambda }{T\times T}{\un{E}}{(s,t)}{\lambda (s)\lambda (t)\lambda (st)^*}$$
for every $\lambda \in \Lambda (T,E)$.
\end{de}

\begin{p}\label{706}
\rule{1em}{0ex}
\begin{enumerate}
\item $\f{T}{E}$ is a subgroup of the commutative multiplicative group $(\un{E})^{T\times T}$ such that $f^*$ is the inverse of $f$ for every $f\in \f{T}{E}$.
\item $\hat f \in \f{T}{E}$ for every $f\in \f{T}{E}$ and the map
$$\mad{\f{T}{E}}{\f{T}{E}}{f}{\hat f}$$
is an involutive group automorphism.
\item $\Lambda (T,E)$ is a subgroup of the commutative multiplicative group $(\un{E})^T$, $\delta \lambda \in \f{T}{E}$ for every $\lambda \in \Lambda (T,E)$, and the map
$$\mae{\delta }{\Lambda (T,E)}{\f{T}{E}}{\lambda }{\delta \lambda }$$
is a group homomorphism with kernel
$$\me{\lambda \in \Lambda (T,E)}{\lambda \;\mbox{\emph{is a group homomorphism}}}$$
such that $\widehat{\delta \lambda }=\delta \hat \lambda $ for every $\lambda \in \Lambda (T,E)$.
\end{enumerate}
\end{p}

a) is obvious.

b) For $r,s,t\in T$,
$$\hat f(r,s)\hat f(rs,t)=f(s^{-1},r^{-1})f(t^{-1},s^{-1}r^{-1})=$$
$$=f(t^{-1},s^{-1})f(t^{-1}s^{-1},r^{-1})=\hat f(r,st)\hat f(s,t)\,,$$
so $\hat f\in \f{T}{E}$. For $f,g\in \f{T}{E}$,
$$\widehat{fg}(s,t)=(fg)(t^{-1},s^{-1})=f(t^{-1},s^{-1})g(t^{-1},s^{-1})=\hat f(s,t)\hat g(s,t)=(\hat f\hat g)(s,t)\,,$$
$$\widehat{fg}=\hat f\hat g\,,$$
$$\hat f^*(s,t)=\hat f(s,t)^*=f(t^{-1},s^{-1})^*=f^*(t^{-1},s^{-1})= \widehat {f^*}(s,t)\,,\qquad (\hat f)^*=\widehat {f^*}\;.$$

c) For $r,s,t\in T$,
$$\delta \lambda (r,s)\delta \lambda (rs,t)=\lambda (r)\lambda (s)\lambda (rs)^*\lambda (rs)\lambda (t)\lambda (rst)^*=\lambda (r)\lambda (s)\lambda (t)\lambda (rst)^*\,,$$
$$\delta \lambda (r,st)\delta \lambda (s,t)=\lambda (r)\lambda (st)\lambda (rst)^*\lambda (s)\lambda (t)\lambda (st)^*=\lambda (r)\lambda (s)\lambda (t)\lambda (rst)^*$$
so $\delta \lambda \in \f{T}{E}$. For $\lambda ,\mu \in \f{T}{E}$ and $s,t\in T$,
$$\delta \lambda (s,t)\delta \mu(s,t)=\lambda (s)\lambda (t)\lambda (st)^*\mu (s)\mu (t)\mu (st)^*=$$
$$=(\lambda \mu )(s)(\lambda \mu )(t)(\lambda \mu )(st)^*=\delta (\lambda \mu)(s,t)\,,$$
$$(\delta \lambda )(\delta \mu )=\delta (\lambda \mu )\,,$$
$$\delta \lambda ^*(s,t)=\lambda ^*(s)\lambda ^*(t)\lambda (st)=(\delta \lambda (s,t))^*=(\delta \lambda )^*(s,t)\,,\qquad \delta \lambda^*=(\delta \lambda )^*\,,$$
so $\delta $ is a group homomorphism. The other assertions are obvious. \qed

\begin{p}\label{6422}
Let $t\in T$, $m,n\in \bz$, and $f\in \f{T}{E}$.
\begin{enumerate}
\item $f(t^m,t^n)=f(t^n,t^m).$
\item $m\in \bn\Longrightarrow f(t^m,t^n)=\left(\proo{j=0}{m-1}f(t^{n+j},t)\right)\left(\proo{k=1}{m-1}f(t^k,t)^*\right)$.
\item We define
$$\mae{\lambda }{\bz}{\un{E}}{n}{\ab{\proo{j=1}{n-1}f(t^j,t)^*}{n\in \bn}{\proo{j=1}{-n}f(t^{-j},t)}{n\not\in \bn}}\;.$$
If $t^p\not=1$ for every $p\in \bn$ then
$$f(t^m,t^n)=\lambda (m)\lambda (n)\lambda (m+n)^*$$
for all $m,n\in \bz$.
\end{enumerate}
\end{p}

a) We may assume $m\in \bn$ because otherwise we can replace $t$ by $t^{-1}$. Put
$$P(m,n):\Longleftrightarrow  f(t^m,t^n)=f(t^n,t^m)\,,$$
$$Q(m):\Longleftrightarrow P(m,n)\;  \mbox{holds for all}\; n\in \bz\;.$$
From
$$f(t^m,t^{n-m})f(t^n,t^m)=f(t^m,t^n)f(t^{n-m},t^m)$$
it follows
$$P(m,n)\Longleftrightarrow P(m,n-m)\Longleftrightarrow P(m,n-km)$$
for all $k\in \bz$.

We prove the assertion by induction. $P(m,0)$ follows from \pr \ref{704} a). By the above
$$P(1,0)\Longleftrightarrow P(1,k)$$
for all $k\in \bz$. Thus $Q(1)$ holds.

Assume $Q(p)$ holds for all $p\in \bnn{m-1}$. Then $P(m,p)$ holds for all $p\in \bnn{m-1}\cup \{0\}$. Let $n\in \bz$. There is a $k\in \bz$ such that
$$p:=n-km\in \bnn{m-1}\cup \{0\}\;.$$
By the above $P(m,n)$ holds. Thus $Q(m)$ holds and this finishes the inductive proof. 

b) We prove the formula by induction with respect to $m$. By a), the formula holds for $m=1$. Assume the formula holds for  an $m\in \bn$. Since
$$f(t^m,t)f(t^{m+1},t^n)=f(t^m,t^{n+1})f(t,t^n)$$
we get by a),
$$f(t^{m+1},t^n)=f(t^m,t^{n+1})f(t,t^n)f(t^m,t)^*=$$
$$=\left(\proo{j=0}{m-1}f(t^{n+1+j},t)\right)\left(\proo{k=1}{m-1}f(t^k,t)^*\right)f(t^n,t)f(t^m,t)^*=$$
$$=\left(\proo{j=0}{m}f(t^{n+j},t)\right)\left(\proo{k=1}{m}f(t^k,t)^*\right)\;.$$
Thus the formula holds also for $m+1$.

c) If $m,n\in \bn$ then by b),
$$\lambda (m)\lambda (n)\lambda (m+n)^*=$$
$$=\left(\proo{k=1}{m-1}f(t^k,t)^*\right)\left(\proo{j=1}{n-1}f(t^j,t)^*\right)\left(\proo{j=1}{m+n-1}f(t^j,t)\right)=$$
$$=\left(\proo{j=0}{m-1}f(t^{n+j},t)\right)\left(\proo{k=1}{m-1}f(t^k,t)^*\right)=f(t^m,t^n)\;.$$

If $m,n\in \bn, n\leq m-1$ then by b),
$$\lambda (m)\lambda (-n)\lambda (m-n)^*=$$
$$=\left(\proo{j=1}{m-1}f(t^j,t)^*\right)\left(\proo{j=1}{n}f(t^{-j},t)\right)\left(\proo{j=1}{m-n-1}f(t^j,t)\right)=$$
$$=\left(\proo{j=0}{m-1}f(t^{-n+j},t)\right)\left(\proo{k=1}{m-1}f(t^k,t)^*\right)=
f(t^m,t^{-n})\;.$$
If $m,n\in \bn, n\geq m$ then by b),
$$\lambda (m)\lambda (-n)\lambda (m-n)^*=$$
$$=\left(\proo{k=1}{m-1}f(t^k,t)^*\right)\left(\proo{j=1}{n}f(t^{-j},t)\right)\left(\proo{j=1}{n-m}f(t^{-j},t)^*\right)=$$
$$=\left(\proo{j=n-m+1}{n}f(t^{-j},t)\right)\left(\proo{k=1}{m-1}f(t^k,t)^*\right)=f(t^m,t^{-n})\;.$$

For all $m,n\in \bn$ put
$$R(m,n):\Longleftrightarrow f(t^{-m},t^{-n})=\lambda (-m)\lambda (-n)\lambda (-m-n)^*\;.$$
By the above and by \pr\ref{704} a),b),
$$\lambda (-1)\lambda (-1)\lambda (-2)^*=f(t^{-1},t)f(t^{-2},t)^*=\tilde f(t^{-1})^*f(t,t^{-2})^*=f(t^{-1},t^{-1})\,,$$
so $R(1,1)$ holds. Let now $m,n\in \bn$ and assume $R(m,n)$ holds. Then
$$\lambda (-m)\lambda (-n-1)\lambda (-m-n-1)^*=$$
$$=\left(\proo{j=1}{m}f(t^{-j},t)\right)\left(\proo{j=1}{n+1}f(t^{-j},t)\right)\left(\proo{j=1}{m+n+1}f(t^{-j},t)^*\right)=$$
$$=f(t^{-m},t^{-n})f(t^{-n-1},t)f(t^{-m-n-1},t)^*=f(t^{-m},t^{-n-1})\,,$$
so $R(m,n)\Rightarrow R(m,n+1)$. By symmetry and a), $R(m,n)$ holds for all $m,n\in \bn$. \qed

\begin{co}\label{719}
The map
$$\mad{\Lambda (\bz,E)}{\f{\bz}{E}}{\lambda }{\delta \lambda }$$
is a surjective group homomorphism with kernel
$$\me{\lambda \in \Lambda (\bz,E)}{n\in \bz\Longrightarrow \lambda (n)=\lambda (1)^n}\;.$$
\end{co}

By \pr\ref{706} c), only the surjectivity of the above map has to be proved and this follows from \pr\ref{6422} c).\qed

\begin{center}
\subsection{$E$-C*-algebras}
\end{center}

By replacing the scalars with the unital C*-algebra $E$ we restrict the category of C*-algebras to the subcategory of those C*-algebras which are connected in a certain way with $E$. The category of unital C*-algebras is replaced by the category of $E$-C*-algebras, while the general category of C*-algebras is replaced by the category of adapted $E$-modules.

\begin{de}\label{32}
We call in this paper {\bf{$E$-module}} a C*-algebra $F$ endowed with the bilinear maps
$$\mad{E\times F}{F}{(\alpha ,x)}{\alpha x}\,,$$
$$\mad{F\times E}{F}{(x,\alpha )}{x\alpha }$$
such that for all $\alpha ,\beta \in E$ and $x,y\in F$,
$$(\alpha \beta )x=\alpha (\beta x)\,,\qquad \alpha (x\beta )=(\alpha x)\beta \,,\qquad x(\alpha \beta )=(x\alpha )\beta \,,$$
$$\alpha (xy)=(\alpha x)y\,,\qquad (xy)\alpha =x(y\alpha )\,,\qquad \alpha \in E^c\Longrightarrow \alpha x=x\alpha \,,$$
$$(\alpha x)^*=x^*\alpha ^*\,,\qquad (x\alpha )^*=\alpha ^*x^*\,,\qquad 1_Ex=x1_E=x \;.$$
If $F,G$ are $E$-modules then a C*-homomorphism $\varphi :F\rightarrow G$ is called {\bf{$E$-linear}} if for all $(\alpha ,x)\in E\times F$,
$$\varphi (\alpha x)=\alpha (\varphi x)\,,\qquad \varphi (x\alpha )=(\varphi x)\alpha \;.$$
\end{de}

For all $(\alpha ,x)\in E\times F$,
$$\n{\alpha x}^2=\n{x^*\alpha ^*\alpha x}\leq \n{x}^2\n{\alpha }^2\,,\quad\quad \n{x\alpha }^2=\n{\alpha ^*x^*x\alpha }\leq \n{\alpha }^2\n{x}^2$$
so
$$\n{\alpha x}\leq \n{\alpha }\n{x}\,,\quad\quad\quad\n{x\alpha }\leq \n{x}\n{\alpha }\;.$$ 

\begin{de}\label{808}
An {\bf $E$-C**-algebra} is a unital C**-algebra $F$ for which $E$ is a {\emph {canonical}} unital C**-subalgebra such that $E^c$ defined with respect to $E$ coincides with $E^c$ defined with respect to $F$ i.e. for every $x\in E$, if $xy=yx$ for all $y\in E$ then $xy=yx$ for all $y\in F$. Every closed ideal of an $E$-C*-algebra is canonically an $E$-module.

Let $F,G$ be $E$-C**-algebras. A map $\mac{\varphi }{F}{G}$ is called an {\bf $E$-C**-homomorphism} if it is an $E$-linear C**-homomorphism . If in addition $\varphi $ is a C*-isomorphism then we say that $\varphi $ is an {\bf $E$-C*-isomorphism} and we use in this case the notation $\approx _E$. A C**-subalgebra $F_0$ of $F$ is called {\bf $E$-C**-subalgebra} of $F$ if $E\subset F_0$.
\end{de}

With the notation of the above Definition $(\alpha -\varphi \alpha )\varphi x=0$ for all $\alpha \in E$ and $x\in F$. Thus $\varphi $ is unital iff $\varphi \alpha =\alpha $ for every $\alpha \in E$. The example
$$\mad{\bk}{\bk\times \bk}{x}{(x,0)}$$
shows that an $E$-C*-homomorphism need not be unital.

If we put $\bt:=\me{z\in \bc}{|z|=1}$, $E:=\ccc{C}(\bt ,\bc)$, and
$$\mae{x}{\bt }{\bc}{z}{z}$$
and if we denote by $\lambda $ the Lebesgue measure on $\bt $ then $L^\infty (\lambda )$ is an $E$-C*-algebra, $x\in \unn{E}$, and $x$ is homotopic to $1_E$ in $\unn{L^\infty (\lambda )}$ but not in $\unn{\ccc{C}(\bt ,\bc)}$.

\begin{de}\label{41}
We denote by $\fr{C}_E$ (resp. by $\fr{C}_E^1$) the category of $E$-C*-algebras for which the morphisms are the $E$-C*-homomorphisms (resp. the unital $E$-C*-homomorphisms). 
\end{de}

\begin{p}\label{33}
Let $F$ be an $E$-module.
\begin{enumerate}
\item We denote by $\check F$ the vector space $E\times F$ endowed with the bilinear map
$$\mad{(E\times F)\times (E\times F)}{E\times F}{((\alpha ,x),(\beta ,y))}{(\alpha \beta ,\alpha y+x\beta +xy)}$$
and with the conjugate linear map
$$\mad{E\times F}{E\times F}{(\alpha ,x)}{(\alpha ^*,x^*)}\;.$$
$\check F$ is an involutive unital algebra with $(1_E,0)$ as unit. 
\item The maps
$$\mae{\pi }{\check F}{E}{(\alpha ,x)}{\alpha }\,,$$
$$\mae{\lambda }{E}{\check F}{\alpha }{(\alpha ,0)}\,,$$
$$\mae{\iota }{F}{\check F}{x}{(0,x)}$$
are involutive algebra homomorphisms such that $\pi \circ \lambda $ is the identity map of $E$, $\lambda $ and $\iota $ are injective, and $\lambda $ and $\pi $ are unital. If there is a norm on $\check F$ with respect to which it is a C*-algebra (in which case such a norm is unique), then we call $ F$ {\bf{adapted}}. We denote by $\fr{M}_E$ the category of adapted $E$-modules for which the morphism are the $E$-linear C*-homomorphisms.
\item If $F$ is adapted then $\check F$ is an $E$-C*-algebra by using canonically the injection $\lambda $ and for all $\alpha \in E$ and $x\in F$,
$$\n{\alpha }\leq \n{(\alpha ,x)}\leq \n{\alpha }+\n{x}\,,\quad \n{(0,x)}=\n{x}\leq 2\n{(\alpha ,x)}\,,$$
$$\n{(\alpha ,0)(0,x)}\leq \n{\alpha }\n{x}\,,\qquad \n{(0,x)(\alpha ,0)}\leq \n{x}\n{\alpha }\;.$$
In particular $F$ (identified with $\iota (F)$) is a closed ideal of $\check F$.
\item If $E$ and $F$ are C*-subalgebras of a C*-algebra $G$ in such a way that the structure of $E$-module of $F$ is inherited from $G$ then
$$\mae{\varphi }{\check{F} }{E\times G}{(\alpha ,x)}{(\alpha ,\alpha +x)}$$
is an injective involutive algebra homomorphism, $\varphi (\check{F} )$ is closed, $F$ is adapted, and for all $\alpha \in E$ and $x\in F$,
$$\n{(\alpha ,x)}_{E\times F}=\sup\{\n{\alpha },\n{\alpha +x}\}\;.$$
In particular every closed ideal of an $E$-C*-algebra is adapted and $\fr{C}_E$ is a full subcategory of $\fr{M}_E$.
\item A closed ideal $G$ of an adapted $E$-module $F$, which is at the same time an $E$-submodule of $F$, is adapted.
\item If $F$ is unital then it is adapted and
$$\mad{\check{F} }{\br_+}{(\alpha ,x)}{\sup\{\n{\alpha },\n{\alpha 1_F+x}\}}$$
is the C*-norm of $\check{F} $.
\item If 
$$\lim_{y,\,\fr{F}}\n{\alpha y-y\alpha }=0$$
for all $\alpha \in E_+$, where $\fr{F}$ denotes the canonical approximate unit of $F$, then $F$ is adapted and 
$$\mad{\check{F} }{\br_+}{(\alpha ,x)}{\sup\left\{\n{\alpha },\,\limsup_{y,\,\fr{F}}\n{\alpha y+x}\right\}}$$
is the C*-norm of $\check{F} $. In particular $F$ is adapted if $E$ is commutative.
\item If $F$ is an adapted $E$-module then (with the notation of b))
$$0\longrightarrow F\stackrel{\iota }{\longrightarrow }\check F{\scriptscriptstyle{\stackrel{\pi }{\longrightarrow }\atop\stackrel{\lambda }{\longleftarrow }}}E\longrightarrow 0$$
is a split exact sequence in the category $\fr{M}_E$.
\end{enumerate}
\end{p}

a) and b) are easy to see.

c) Since $\lambda $ and $\iota $ are injective and $ $,
$$\pi (\alpha ,x)=\alpha\,,\qquad\qquad (\alpha ,x)=(\alpha ,0)+(0,x)\,,$$
$$ (\alpha ,0)(0,x)=(0,\alpha x)\,,\qquad\qquad (0,x)(\alpha ,0)=(0,x\alpha )$$
we get the first and the last two inequalities as well as the identity $\n{(0,x)}=\n{x}$. It follows
$$\n{(0,x)}\leq \n{(\alpha ,x)}+\n{(\alpha ,0)}=\n{(\alpha ,x)}+\n{\lambda \pi (\alpha ,x)}\leq$$
$$\leq  \n{(\alpha ,x)}+\n{(\alpha ,x)}=2\n{(\alpha ,x)}\;.$$

d) It is easy to see that $\varphi $ is an injective involutive algebra homomorphism. Let $(\alpha ,x)\in \overline{\varphi (\check F)}$. There are sequences $(\alpha _n)_{n\in \bn}$ and $(x_n)_{n\in \bn}$ in $E$ and $F$, respectively, such that
$$\lim_{n\rightarrow \infty }(\alpha _n,\alpha _n+x_n)=(\alpha ,x)\;.$$
It follows
$$\alpha  =\lim_{n\rightarrow \infty }\alpha _n\in E\,,\quad x-\alpha =\lim_{n\rightarrow \infty }x_n\in F\,,\quad (\alpha ,x)=\varphi (\alpha ,x-\alpha )\in \varphi (\check F)\;.$$
Thus $\varphi (\check F)$ is closed, which proves the assertion by pulling back the norm of $E\times G$.

e) By c), $F$ is a closed ideal of $\check F$ so $G$ is a closed ideal of $\check F$ (use an approximate unit of $F$). Since $G$ is an $E$-submodule of $F$ its structure of $E$-module is inherited from $\check F$. By d), $G$ is adapted.

f) The map
$$\mad{\check{F} }{E\times F}{(\alpha ,x)}{(\alpha ,\alpha 1_F+x)}$$
is an isomorphism of involutive algebras and so we can pull back the norm of $E\times F$.

g) It is easy to see that the above map is a norm. Since
$$\sup\{\n{\alpha },\frac{1}{2}\n{x}\}\leq \n{(\alpha ,x)}\leq \n{\alpha }+\n{x}$$
for all $(\alpha ,x)\in E\times F$, $\check{F} $ endowed with this norm is complete. For $(\alpha ,x)\in E\times F$,
$$(\alpha ,x)^*(\alpha ,x)=(\alpha ^*\alpha ,\alpha ^*x+x^*\alpha +x^*x)\,,$$
$$\n{(\alpha ,x)^*(\alpha ,x)}=\sup\{\n{\alpha }^2,\,\limsup_{y,\,\fr{F}}\n{\alpha ^*\alpha y+\alpha ^*x+x^*\alpha +x^*x}\}\;.$$
For $y\in F_+^{\#}$,
$$\n{(\alpha y^{\frac{1}{2}}+x)^*(\alpha y^{\frac{1}{2}}+x)-(\alpha ^*\alpha y+\alpha ^*x+x^*\alpha +x^*x)}\leq $$
$$\leq \n{y^{\frac{1}{2}}\alpha ^*\alpha -\alpha ^*\alpha y^{\frac{1}{2}}}+\n{y^{\frac{1}{2}}\alpha ^*x-\alpha ^*x}+\n{x^*\alpha y^{\frac{1}{2}}-x^*\alpha }$$
so
$$\lim_{y,\,\fr{F}}\n{(\alpha y^{\frac{1}{2}}+x)^*(\alpha y^{\frac{1}{2}}+x)-(\alpha ^*\alpha y+\alpha ^*x+x^*\alpha +x^*x)}=0\;.$$
Since the map $F_+\rightarrow F_+,\,y\mapsto y^{\frac{1}{2}}$ maps $\fr{F}$ into itself and
$$\n{\alpha y+x}^2=\n{y\alpha ^*\alpha y+y\alpha ^*x +x^*\alpha y+x^*x}$$
we have by the above,
$$\n{(\alpha ,x)}^2=\sup\left\{\n{\alpha }^2,\,\limsup_{y,\,\fr{F}}\n{\alpha y^{\frac{1}{2}}+x}^2\right\}=$$
$$=\sup\left\{\n{\alpha }^2,\,\limsup_{y,\,\fr{F}}\n{(\alpha y^{\frac{1}{2}}+x)^*(\alpha y^{\frac{1}{2}}+x)}\right\}=$$
$$=\sup\{\n{\alpha }^2,\,\limsup_{y,\,\fr{F}}\n{\alpha ^*\alpha y+\alpha ^*x+x^*\alpha +x^*x}\}=\n{(\alpha ,x)^*(\alpha ,x)}\;.$$
Thus the above norm is a C*-norm and $F$ is adapted.

h) $\iota $ is an injective $E$-C*-homomorphism and its image is equal to $Ker\;\pi $.\qed

\begin{co}\label{51}
Let $F$ an $E$-module, $G$ a C*-algebra, and $\otimes _\sigma $ the spatial tensor product. 
\begin{enumerate}
\item $F\otimes _\sigma G$ is in a natural way an $E$-module the multiplication being given by
$$\alpha (x\otimes y)=(\alpha x)\otimes y\,,\qquad \qquad (x\otimes y)\alpha =(x\alpha )\otimes y$$
for all $\alpha \in E$, $x\in F$, and $y\in G$.
\item If $F$ is an $E$-C*-algebra and $G$ is unital then the map
$$\mad{E}{F\otimes _\sigma G}{\alpha }{\alpha \otimes 1_G}$$
is an injective C*-homomorphism. In particular, the $E$-module $F\otimes _\sigma G$ is an $E$-C*-algebra.
\item If $F$ is an adapted $E$-module then the $E$-module $F\otimes _\sigma G$ is adapted and
$$\n{(\alpha ,z)}=\sup\{\n{\alpha },\n{\alpha +z}\}$$
for all $(\alpha ,z)\in E\times (F\otimes _\sigma G)$.
\item If $F$ is an adapted $E$-module and $G:=\ccc{C}_0(\Omega )$ for a locally compact space $\Omega $ then $\ccc{C}_0(\Omega ,F)$ is adapted and 
$$\n{(\alpha ,x)}=\sup\{\n{\alpha },\n{\alpha e_\Omega +x}\}$$
for all $(\alpha ,x)\in E\times \ccc{C}_0(\Omega ,F)$. 
\end{enumerate}
\end{co} 

a) and b) are easy to see.

c) If $\tilde{G}$ denotes the unitization of $G$ then by b), $\check{F}\otimes _\sigma \tilde{G}$ is an $E$-C*-algebra and $F\otimes _\sigma G$ is a closed ideal of it, so the assertion follows from \pr \ref{33} d),e).

d) follows from c).\qed  

\begin{p}\label{34}
\rule{1ex}{0em}
\begin{enumerate}
\item If $F,G$ are $E$-modules and $\varphi :F\rightarrow G$ is an $E$-linear C*-homomorphism then the map
$$\mae{\check \varphi }{\check F}{\check G}{(\alpha ,x)}{(\alpha ,\varphi x)}$$
is an involutive unital algebra homomorphism, injective or surjective if $\varphi $ is so. If $F=G$ and if $\varphi $ is the identity map then $\check \varphi $ is also the identity map.
\item Let $F_1,F_2,F_3$ be $E$-modules and let $\varphi :F_1\rightarrow F_2$ and $\psi :F_2\rightarrow F_3$ be $E$-linear C*-homomorphisms. Then $\check {\overbrace{\psi \circ \varphi }} =\check \psi \circ \check \varphi $.\qed 
\end{enumerate}
\end{p}

\begin{p}\label{7903}
Let $G$ be an $E$-module, $F$ an $E$-submodule of $G$ which is at the same time an ideal of $G$, and $\varphi :G\rightarrow G/F$ the quotient map.
\begin{enumerate}
\item $G/F$ has a natural structure of $E$-module and $\varphi $ is $E$-linear.
\item If $G$ is adapted then $G/F$ is also adapted. Moreover if $\psi :\check{G}\rightarrow \check{G}/F$ denotes the quotient map (where $F$ is identified to $\{(0,x)|\; x\in F\}$) then there is an $E$-C*-isomorphism  $\theta :\check{\overbrace{G/F}}\rightarrow {\check G}/F$ such that $\psi =\theta \circ \check{\varphi }$.
\end{enumerate}
\end{p}

a) is easy to see.

b) Let $(\alpha ,z)\in \check{\overbrace{G/F}} $ and let $x,y\in \stackrel{-1}{\varphi }(z)$. Then $\psi (\alpha ,x)=\psi (\alpha ,y)$ and we put $\theta (\alpha ,z):=\psi (\alpha ,x)$. It is straightforward to show that $\theta $ is an isomorphism of involutive algebras. By pulling back the norm of $\check{G}/F$ with respect to $\theta $ we see that $G/F$ is adapted.\qed

\begin{lem}\label{879}
Let  $\{(F_i)_{i\in I},(\varphi _{ij})_{i,j\in I})\}$ be an inductive system in the category of C*-algebras, $\{F,(\varphi _i)_{i\in I}\}$ its inductive limit, $G$ a C*-algebra, for every $i\in I$, $\psi _i:F_i\rightarrow G$ a C*-homomorphism such that $\psi _j\circ \varphi _{ji}=\psi _i$ for all $i,j\in I,i\leq j$, and $\psi :F\rightarrow G$ the resulting C*-homomorphism. If  $Ker\,\psi _i\subset Ker\,\varphi _i$ for every $i\in I$ then $\psi $ is injective.
\end{lem}

Let $i\in I$. Since $Ker\,\varphi _i\subset Ker\,\psi _i$ is obvious, we have $Ker\,\varphi _i= Ker\,\psi _i$. Let $\rho :F_i\rightarrow F_i/Ker\,\psi _i$ be the quotient map and 
$$\mac{\varphi _i'}{F_i/Ker\,\psi _i}{F}\,,\qquad \mac{\psi _i'}{F_i/Ker\,\psi _i}{G}$$
the injective C*-homomorphisms with
$$\varphi _i=\varphi _i'\circ \rho \,,\qquad \psi _i=\psi _i'\circ \rho \;.$$
Then
$$\psi _i'\circ \rho =\psi _i=\psi \circ \varphi _i=\psi \circ \varphi _i'\circ \rho \;.$$
For $x\in F_i$, since $\psi _i'$ and $\varphi _i'$ are norm-preserving,
$$\n{\rho x}=\n{\psi _i'\rho x}=\n{\psi \varphi _i'\rho x}\leq \n{\varphi _i'\rho x}=\n{\rho x}\,,$$
$$\n{\psi \varphi _ix}=\n{\psi \varphi _i'\rho x}=\n{\varphi _i'\rho x}=\n{\varphi _ix}\;.$$
Thus $\psi $ preserves the norms on $\cup _{i\in I}\varphi _i(F_i)$. Since this set is dense in $F$, $\psi $ is injective.\qed

\begin{p}\label{36}
Let $\{(F_i)_{i\in I},(\varphi _{ij})_{i,j\in I}\}$ be an inductive system in the category $\fr{M}_E$ and let $(F,(\varphi _i)_{i\in I})$ be its inductive limit in the category of $E$-modules \emph{(\pr \ref{33} c))}.
\begin{enumerate}
\item $F$ is adapted.
\item Let $(G,(\psi _i)_{i\in I})$ be the inductive limit in the category $\fr{C}_E^1$ of the inductive system $\{(\check F_i)_{i\in I},(\check \varphi _{ij})_{i,j\in I}\}$ \emph{(\pr \ref{34} a),b))} and let $\psi :G\rightarrow \check F$ be the unital C*-homomorphism such that $\psi \circ \psi _i=\check \varphi _i$ for every $i\in I$ . Then $\psi $ is an $E$-C*-isomorphism.
\end{enumerate}
 \end{p}

a) Put
$$F_0:=\me{(\alpha ,x)\in \check F}{\alpha \in E,\;x\in \bigcup _{i\in I}\varphi _i(F_i)}\,,$$
$$\mae{p}{F_0}{\br_+}{(\alpha ,x)}{\inf\me{\n{(\alpha ,x_i)}}{i\in I,\;x_i\in F_i,\;\varphi _ix_i=x}}\;.$$
$F_0$ is an involutive unital subalgebra of $\check F$. $p$ is a norm and by \pr \ref{33} c),
$$q(\alpha ,x):=\lim_{(\alpha ,y)\in F_0 \atop y\rightarrow x}p(\alpha ,y)$$
exists and 
$$\n{\alpha }\leq q(\alpha ,x)\leq \n{\alpha }+\n{x}\,,\qquad \n{x}\leq 2q(\alpha ,x)$$
for every $(\alpha ,x)\in \check F$. 

Let $(\alpha ,x)\in F_0$. Let further $i\in I$, $x_i,y_i\in F_i$ with $\varphi _ix_i=x,\,\varphi _iy_i=\alpha ^*x+x^*\alpha +x^*x$. Then
$$(0,\varphi _i(\alpha ^*x_i+x_i^*\alpha +x_i^*x_i-y_i))=\check{\varphi _i}((\alpha ,x_i)^*(\alpha ,x_i)-(\alpha ^*\alpha ,y_i))=0 $$
so
$$\inf_{i\leq j}\n{\varphi _{ji}(\alpha ^*x_i+x_i^*\alpha +x_i^*x_i-y_i)}=0\;.$$
For $\epsilon >0$ there is a $j\in I,\,i\leq j,$ with
$$\n{\varphi _{ji}(\alpha ^*x_i+x_i^*\alpha +x_i^*x_i-y_i)}<\epsilon \;.$$
We get
$$p(\alpha ,x)^2\leq \n{(\alpha ,\varphi _{ji}x_i)}^2=\n{(\alpha ,\varphi _{ji}x_i)^*(\alpha ,\varphi _{ji}x_i)}=$$
$$=\n{(\alpha ^*\alpha , \alpha ^*\varphi _{ji}x_i+(\varphi _{ji}x_i^*)\alpha +\varphi _{ji}(x_i^*x_i))}=$$
$$=\n{(\alpha ^*\alpha ,\varphi _{ji}(\alpha ^*x_i+x_i^*\alpha +x_i^*x_i))}\leq $$
$$\leq \n{(\alpha ^*\alpha ,\varphi _{ji}y_i)}+\n{(0,\varphi _{ji}(\alpha ^*x_i+x_i^*\alpha +x_i^*x_i-y_i))}<\n{(\alpha ^*\alpha ,\varphi _{ji}y_i)}+\epsilon \;.$$
By taking the infimum on the right side it follows, since $\epsilon $ is arbitrary,
$$p(\alpha ,x)^2\leq p(\alpha ^*\alpha ,\alpha ^*x+x^*\alpha +x^*x)=p((\alpha ,x)^*(\alpha ,x))$$
and this shows that $p$ is a C*-norm. It is easy to see that $q$ is a C*-norms. By the above inequalities, $\check F$ endowed with the norm $q$ is complete, i.e. $\check F$ is a C*-algebra and $F$ is adapted.

b) Let $i\in I$ and let $(\alpha ,x)\in Ker\,\check \varphi _i$. Then 
$$0=\check \varphi _i(\alpha ,x)=(\alpha ,\varphi _ix)$$
so
$$\alpha =0\,,\qquad \varphi _ix=0\,,\qquad \inf_{j\in I,\,j\geq i}\n{\varphi _{ji}x}=0\,,$$
$$\n{\check \varphi _{ji}(0,x)}=\n{(0,\varphi _{ji}x)}=\n{\varphi _{ji}x}\,,$$ 
$$\n{\psi _i(\alpha ,x)}=\inf_{j\in I,\,j\geq i}\n{\check \varphi _{ji}(0,x)}=0\,,\qquad (\alpha ,x)\in Ker\,\psi _i\;.$$
By Lemma \ref{879}, $\psi $ is injective.

Let $(\beta ,y)\in \check F$ and let $\varepsilon >0$. There are $i\in I$ and $x\in F_i$ with $\n{\varphi _ix-y}<\varepsilon $.
Then
$$\psi \psi _i(\beta ,x)=\check \varphi _i(\beta ,x)=(\beta ,\varphi _ix)\,,$$
$$\n{\psi \psi _i(\beta ,x)-(\beta ,y)}=\n{\check \varphi _i(\beta ,x)-(\beta ,y)}=\n{\varphi _ix-y}<\varepsilon \;.$$
Thus $\psi (G)$ is dense in $\check F$ and $\psi $ is surjective. Hence $\psi $ is a C*-isomorphism.\qed

\begin{co}\label{42}
We put $\Phi_E(F):=\check F$ for every $E$-module $F$ and similarly
$\Phi_E(\varphi ):=\check \varphi $
for every $E$-linear C*-homomorphism $\varphi $.
\begin{enumerate}
\item $\Phi_E$ is a covariant functor from the category $\fr{M}_E$ in the category $\fr{C}_E^1$.
\item The categories $\fr{C}_E^1$ and $\fr{M}_E$ possess inductive limits and the functor $\Phi_E$ is continuous with respect to the inductive limits.
\end{enumerate}
\end{co}

a) follows from \pr \ref{34}.

b) follows from \pr \ref{36}.\qed

{\it Remark.} The category $\fr{C}_E$ does not possess inductive limits in general. This happens for instance if $\varphi _{ij}=0$ for all $i,j\in I$.

\begin{center}
\subsection{Some topologies}
\fbox{$T$ is only a set in this subsection}
\end{center}

If the group $T$ is infinite then different topologies play a certain role in the construction of the projective representations of $T$. It will be shown that all these topologies conduct to the same construction, but the use of them simplifies the manipulations.  

We introduce the following notation in order to unify the cases of C*-algebras and (resp. W*-algebras).

\begin{de}\label{673}
\rule{1ex}{0em}
$$\widetilde {\ca}:=\ca\qquad\qquad\qquad ( \mbox{\emph{resp.}}\;  \widetilde {\ca}:=\cw{}\,)\,,$$
$$\widetilde {\otimes }:=\otimes \qquad\qquad\qquad (\emph{resp.}\;\widetilde {\otimes }:=\bar {\otimes }\,)\,,$$
$$\widetilde {\sum}:=\sum\qquad\qquad\qquad (\emph{resp.}\;\widetilde {\sum}:=\sii{}{\ddot E}\,)\;.$$
If $\fr{T}$ is a Hausdorff topology on $\lb{E}{H}$ then for every $\ccc{G}\subset \lb{E}{H}$, $\ccc{G}_{\fr{T}}$  denotes the set $\ccc{G}$ endowed with the relative topology $\fr{T}$ and $\stackrel{\fr{T}}{\bar {\ccc{G}}}$ denotes the closure of $\ccc{G}$ in $\lb{E}{H}_{\fr{T}}$. Moreover $\sii{}{\fr{T}}$ denotes the sum with respect to $\fr{T}$.
\end{de}

\begin{lem}\label{18}
For $x\in E$, by the above identification of $E$ with $\lb{E}{\breve E}$,
$$\mae{x\widetilde\otimes 1_K}{H}{H}{\xi }{(x\xi _t)_{t\in T}}$$
is well-defined and belongs to $\lb{E}{H}$.
\begin{enumerate}
\item The map
$$\mae{\varphi }{E}{\lb{E}{H}}{x}{x\widetilde\otimes 1_K}$$
is an injective unital C*-homomorphism.
\item Assume $E$ is a W*-algebra. Then for every $(a,\xi ,\eta )\in \ddot E\times H\times H$, the family $(\xi _t\,a\,\eta _t^*)_{t\in T}$ is summable in $\ddot E_E$ and for every $x\in E$,
$$\sa{\varphi x}{\tia{a}{\xi }{\eta }}=\sa{x}{\sii{t\in T}{E}\xi _t\,a\,\eta _t^*}\;.$$
Thus $\varphi $ is a W*-homomorphism \emph{(\hh 3.5 d))} with
$$\ddot \varphi \tia{a}{\xi }{\eta }=\sii{t\in T}{E}\xi _t\,a\,\eta _t^*\,,$$
where $\ddot \varphi $ denotes the pretranspose of $\varphi $.
\item If we consider $E$ as a canonical unital C**-subalgebra of $\lb{E}{H}$ by using the embedding of a) then $\lb{E}{H}$ is an $E$-C**-algebra. 
\end{enumerate}
\end{lem}

a) follows from [L] page 37 (resp. [C3] \pr 1.4).

b) We have
$$\sa{x\bar \otimes 1_K}{\tia{a}{\xi }{\eta }}=\sa{\s{(x\bar \otimes 1_K)\xi }{\eta }}{a}=\sa{\sii{t\in T}{\ddot E}\eta _t^*\,x\,\xi _t}{a}=$$
$$=\si{t\in T}\sa{\eta _t^*\,x\,\xi _t}{a}=\si{t\in T}\sa{x}{\xi _t\,a\,\eta _t^*}\;.$$
Thus the family $(\xi _t\,a\,\eta _t^*)_{t\in T}$ is summable in $\ddot E_E$ and
$$\sa{\varphi x}{\tia{a}{\xi }{\eta }}=\sa{x}{\sii{t\in T}{E}\xi _t\,a\,\eta _t^*}\;.$$
If $\varphi ':\lb{E}{H}\rightarrow E'$ denotes the transpose of $\varphi $ then
$$\varphi '\tia{a}{\xi }{\eta }=\sii{t\in T}{E}\xi _t\,a\,\eta _t^*\in \ddot E\;.$$
By continuity $\varphi '\left(\overbrace{\lb{E}{H}}^{..}\right)\subset \ddot E$ and $\varphi $ is a unital W*-homomorphism. 

c) Let $x\in E^c$ and $\xi ,\eta \in \lb{E}{H}$. By \prr 3.17 d),
$$\s{(x\widetilde\otimes 1_K)\xi }{\eta }=\siw{t\in T}\eta _t^*((x\widetilde\otimes 1_K)\xi )_t=\siw{t\in T}\eta _t^*x\xi _t=$$
$$=\siw{t\in T}x\eta _t^*\xi _t=x\siw{t\in T}\eta _t^*\xi _t=x\s{\xi }{\eta }\;.$$
Thus for $u\in \lb{E}{H}$,
$$\s{u(x\widetilde\otimes 1_K)\xi }{\eta }=\s{(x\widetilde\otimes 1_K)\xi }{u^*\eta }=x\s{\xi }{u^*\eta }=x\s{u\xi }{\eta }\,,$$
$$u(x\widetilde\otimes 1_K)=(x\widetilde\otimes 1_K)u\,,$$
and so $x\widetilde\otimes 1_K\in \lb{E}{H}^c$.\qed
 
\renewcommand{\labelenumii}{$\alph{enumi}_{\arabic{enumii}}$)}
\begin{de}\label{740}
We put for all $\xi ,\eta \in H$  (resp. and $a\in \ddot E_+$)
$$\mae{p_{\xi ,\eta }}{\lb{E}{H}}{\br_+}{X}{\n{\s{X\xi }{\eta }}}\,,$$
$$(\mbox{\emph{resp.}}\;\mae{p_{\xi ,\eta ,a}}{\lb{E}{H}}{\br_+}{X}{|\sa{\s{X\xi }{\eta }}{a}|})\,,$$
$$\mae{p_\xi }{\lb{E}{H}}{\br_+}{X}{\n{X\xi }=\n{\s{X\xi }{X\xi }}^{1/2}}\,,$$
$$(\mbox{\emph{resp.}}\;\mae{p_{\xi ,a}}{\lb{E}{H}}{\br_+}{X}{\sa{\s{X\xi }{X\xi }}{a}^{1/2}})\,,$$
$$\mae{q_\xi }{\lb{E}{H}}{\br_+}{X}{p_\xi (X^*)}\,,$$
$$(\mbox{\emph{resp.}}\;\mae{q_{\xi ,a}}{\lb{E}{H}}{\br_+}{X}{p_{\xi ,a}(X^*)})\;.$$
and denote, respectively, by $\fr{T}_1\,,\fr{T}_2\,,\fr{T}_3$ the topologies on $\lb{E}{H}$ generated by the set of seminorms
$$\me{p_{\xi ,\eta }}{\xi ,\eta \in H}\,,\qquad \left(\mbox{\emph{resp.}}\;\me{p_{\xi ,\eta ,a}}{\xi ,\eta \in H,\,a\in \ddot E_+}\right)\,,$$
$$\me{p_{\xi }}{\xi \in H}\,,\qquad \left(\mbox{\emph{resp.}}\;\me{p_{\xi ,a}}{\xi \in H,\,a\in \ddot E_+}\right)\,,$$
$$\me{p_\xi }{\xi \in H}\cup \me{q_\xi }{\xi \in H}\,,$$
$$ \left(\mbox{\emph{resp.}}\;\me{p_{\xi ,a}}{\xi \in H,\,a\in \ddot E_+}\cup \me{q_{\xi ,a}}{\xi \in H,\,a\in \ddot E_+}\right)\;.$$
Moreover $\n{\cdot }$ denotes the norm topology on $\lb{E}{H}$.
\end{de}

Of course $\fr{T}_2\subset \fr{T}_3$. In the C*-case, $\fr{T}_2$ is the topology of pointwise convergence. If $E$ is finite-dimensional then the C*-case and the W*-case coincide. 

\begin{p}\label{741} 
Let $X\in \lb{E}{H}$ and $\xi ,\eta \in H$ (resp. and $a\in \ddot E$).
\begin{enumerate}
\item $p_{\xi ,\eta }(X)=p_{\eta ,\xi }(X^*) \qquad ($\emph{resp.} $p_{\xi ,\eta ,|a|}(X)=p_{\eta ,\xi ,|a|}(X^*))$. 
\item $p_{\xi ,\eta }(X)\leq p_\xi (X)\n{\eta }$.
\item If $E$ is a W*-algebra and $a=x|a|$ is the polar representation of $a$ then
$$p_{\xi x,\eta ,|a|}(X)=\left|\sa{X}{\tia{a}{\xi }{\eta }}\right|\leq p_{\xi x,|a|}(X)\sa{\s{\eta }{\eta }}{|a|}^{1/2}\;.$$
\item If $Y,Z\in \lb{E}{H}$ then
$$p_{\xi ,\eta }(YXZ)=p_{Z\xi ,Y^*\eta }(X)\qquad\emph{(resp.}\;p_{\xi ,\eta ,|a|}(YXZ)=p_{Z\xi ,Y^*\eta ,|a|}(X)\emph{)}\,,$$
$$p_\xi (YXZ)\leq \n{Y}p_{Z\xi }(X)\qquad(\emph{resp.} \;p_{\xi ,|a|}(YXZ)\leq \n{Y}p_{Z\xi ,|a|}(X))\;.$$
\end{enumerate}
\end{p}

a) From 
$$\s{X\xi }{\eta }=\s{\xi }{X^*\eta }=\s{X^*\eta }{\xi }^*$$
it follows
$$p_{\xi ,\eta }(X)=\n{\s{X\xi }{\eta }}=\n{\s{X^*\eta }{\xi }}=p_{\eta ,\xi }(X^*)\,,$$
$$(\mbox{resp.}\;p_{\xi ,\eta ,|a|}(X)=|\sa{\s{X^*\eta }{\xi }}{|a|}|=p_{\eta ,\xi ,|a|}(X^*))\;.$$

b) $p_{\xi ,\eta }(X)=\n{\s{X\xi }{\eta }}\leq p_\xi (X)\n{\eta }$.

c) We have
$$p_{\xi x,\eta ,|a|}(X)=|\sa{\s{X(\xi x)}{\eta }}{|a|}|=|\sa{\s{X\xi }{\eta }x}{|a|}|=$$
$$=|\sa{\s{X\xi }{\eta }}{x|a|}|=|\sa{\s{X\xi }{\eta }}{a}|=\left|\sa{X}{\tia{a}{\xi }{\eta }}\right|\;.$$
By Schwarz' inequality ([C1] \pr 2.3.3.9),
$$|\sa{\s{X(\xi x)}{\eta }}{|a|}|^2\leq \sa{\s{X(\xi x)}{X(\xi x)}}{|a|}\sa{\s{\eta }{\eta }}{|a|}\,,$$
so
$$p_{\xi x,\eta ,|a|}(X)\leq p_{\xi x,|a|}(X)\sa{\s{\eta }{\eta }}{|a|}^{1/2}\;.$$

d) The first equation follows from
$$p_{\xi ,\eta}(YXZ)=\n{\s{YXZ\xi }{\eta }}=\n{\s{XZ\xi }{Y^*\eta }}=p_{Z\xi ,Y^*\eta }(X) $$
$$(\mbox{resp.}\;p_{\xi ,\eta ,|a|}(YXZ)=|\sa{\s{YXZ\xi }{\eta }}{|a|}|=$$
$$=|\sa{\s{XZ\xi }{Y^*\eta }}{|a|}|=p_{Z\xi ,Y^*\eta ,|a|}(X))$$
and the second from
$$p_\xi (YXZ)=\n{YXZ\xi }\leq \n{Y}\n{XZ\xi }=\n{Y}p_{Z\xi }(X)$$
$$(\mbox{resp.}\;p_{\xi ,|a|}(YXZ)=\sa{\s{YXZ\xi }{YXZ\xi }}{|a|}^{1/2}\leq $$
$$\leq \n{Y}\sa{\s{XZ\xi }{XZ\xi }}{|a|}^{1/2}=\n{Y}p_{Z\xi ,|a|}(X))\;.\qedd$$

\begin{lem}\label{759} 
Let $n\in \bn$ and $(x_i)_{i\in \bnn{n}}$ a family in $E$. Then
$$\left(\sum_{i\in \bnn{n}}x_i\right)^*\left(\sum_{i\in \bnn{n}}x_i\right)\leq n\sum_{i\in \bnn{n}}x_i^*x_i\;.$$
\end{lem}

We prove the relation by induction with respect to $n$. By [C1] Corollary 4.2.2.4 and by the hypothesis of the induction,
$$\left(\sum_{i\in \bnn{n}}x_i\right)^*\left(\sum_{i\in \bnn{n}}x_i\right)=\left(x_n^*+\sum_{i\in \bnn{n-1}}x_i^*\right)\left(x_n+\sum_{i\in \bnn{n-1}}x_i\right)=$$
$$=x_n^*x_n+\sum_{i\in \bnn{n-1}}(x_n^*x_i+x_i^*x_n)+\left(\sum_{i\in \bnn{n-1}}x_i\right)^*\left(\sum_{i\in \bnn{n-1}}x_i\right)\leq $$
$$\leq x_n^*x_n+\sum_{i\in \bnn{n-1}}(x_n^*x_n+x_i^*x_i)+(n-1)\sum_{i\in \bnn{n-1}}x_i^*x_i=n\sum_{i\in \bnn{n}}x_i^*x_i\;.\qedd$$
\begin{lem}\label{760} 
Let $n\in \bn$, $x\in E_{n,n}$, and for every $j\in \bnn{n}$ put
$$\eta _j:=(\delta_{ji}1_E)_{i\in \bnn{n}}\in \cb{i\in \bnn{n}}\breve E\;.$$  
Then
$$\n{x}\leq \sqrt{n}\sup_{j\in \bnn{n}}\n{x\eta _j}\;.$$
\end{lem}

For $\xi \in \left(\cb{i\in \bnn{n}}\breve E\right)^{\#}$, by Lemma \ref{759}, 
$$\s{x\xi }{x\xi }=\sum_{i\in \bnn{n}}\s{(x\xi )_i}{(x\xi )_i}=\sum_{i\in \bnn{n}}\left(\sum_{j\in \bnn{n}}x_{ij}\xi _j\right)^*\left(\sum_{j\in \bnn{n}}x_{ij\xi _j}\right)\leq $$
$$\leq n\sum_{i\in \bnn{n}}\sum_{j\in \bnn{n}}(x_{ij}\xi _j)^*(x_{ij}\xi _j)=n\sum_{i\in \bnn{n}}\sum_{j\in \bnn{n}}\xi _j^*x_{ij}^*x_{ij}\xi _j=$$
$$=n\sum_{j\in \bnn{n}}\xi _j^*\left(\sum_{i\in \bnn{n}}x_{ij}^*x_{ij}\right)\xi _j\;. $$
For $i,j\in \bnn{n}$,
$$(x\eta _j)_i=\sum_{k\in \bnn{n}}x_{ik}\eta _{jk}=x_{ij}\,,$$
$$\s{x\eta _j}{x\eta _j}=\sum_{i\in \bnn{n}}(x\eta _j)_i^*(x\eta _j)_i=\sum_{i\in \bnn{n}}x_{ij}^*x_{ij}\,,$$
so
$$\s{x\xi }{x\xi }\leq n\sum_{j\in \bnn{n}}\xi _j^*\s{x\eta _j}{x\eta _j}\xi _j\leq n\sum_{j\in \bnn{n}}\n{x\eta _j}^2\xi _j^*\xi _j\leq $$
$$\leq n\sup_{j\in \bnn{n}}\n{x\eta _j}^2\sum_{j\in \bnn{n}}\xi _j^*\xi _j\leq n\sup_{j\in \bnn{n}}\n{x\eta _j}^21_E\,,$$
$$\n{x}^2\leq n\sup_{j\in \bnn{n}}\n{x\eta _j}^2\;.\qedd$$
\begin{co}\label{742} 
\rule{1em}{0ex}
\begin{enumerate}
\item The map
$$\mad{\leh_{\fr{T}_1}}{\leh_{\fr{T}_1}}{X}{X^*}$$
is continuous. In particular $Re\,\leh$ is a closed set of $\leh_{\fr{T}_1}$.
\item $\fr{T}_1\subset \fr{T}_2\subset \fr{T}_3\subset $ norm topology.
\item If $E$ is a W*-algebra then the identity map
$$\lh{E}{H}\longrightarrow \leh_{\fr{T}_1}$$
is continuous so
$$\leh_{\fr{T}_1}^{\#}=\lh{E}{H}^{\#}$$
is compact.
\item For $Y,Z\in \leh$ and $k\in \{1,2\}$, the map
$$\mad{\leh_{\fr{T}_k}}{\leh_{\fr{T}_k}}{X}{YXZ}$$
is continuous.
\item $\leh_{\fr{T}_3}$ is complete in the C*-case.
\item If $T$ is finite then $\fr{T}_2$ is the norm topology in the C*-case.
\item $ \ccc{K}_{E}(H)$ is dense in $\leh_{\fr{T}_3}$.
\end{enumerate}
\end{co}

a) follows from \pr\ref{741} a).

b) $\fr{T}_1\subset \fr{T}_2$ follows from \pr\ref{741} b),c). $\fr{T}_2\subset \fr{T}_3\subset $ norm topology is trivial.

c) follows from \pr\ref{741} c) (and \hh 3.5 a)).

d) follows from \pr\ref{741} d).

e) Let $\fr{F}$ be a Cauchy filter on $\leh_{\fr{T}_3}$. Put
$$\mae{Y}{H}{H}{\xi }{\lim_{X,\fr{F}}(X\xi )}\,,$$
$$\mae{Z}{H}{H}{\xi }{\lim_{X,\fr{F}}(X^*\xi )}\,,$$
where the limits are considered in the norm topology of $H$. For $\xi ,\eta \in H$,
$$\s{Y\xi }{\eta }=\lim_{X,\fr{F}}\s{X\xi }{\eta }=\lim_{X,\fr{F}}\s{\xi }{X^*\eta }=\s{\xi }{Z\eta }\,,$$
so $Y,Z\in \leh$ and $Z=Y^*$. Thus $\fr{F}$ converges to $Y$ in $\leh_{\fr{T}_3}$ and $\leh_{\fr{T}_3}$ is complete.

f) follows from b) and Lemma \ref{760}.

g) Let $X\in \leh$ and $\xi \in H$. For every $S\in \fr{P}_f(T)$ put
$$P_S:=\si{s\in S}e_s\s{\cdot }{e_s}\in Pr\;\ccc{K}_E(H)$$
and let $\fr{F}_T$ be the upper section filter or $\fr{P}_f(T)$. Then $P_SX\in \ccc{K}_E(H)$ for every $S\in \fr{P}_f(T)$ and
$$\lim_{S,\fr{F}_T}P_SX\xi =X\xi $$
in $H$ (resp. in $H_{\ddot H}$) (\prr 4.1 e) (resp. \prr 4.6 c))). Thus
$$\lim_{S,\fr{F}_T}P_SX=X$$
with respect to the topology $\fr{T}_2$. Since the same holds for $X^*$, it follows that $X$ belongs to the closure of $\ccc{K}_E(H)$ in $\leh_{\fr{T}_3}$.\qed

{\it{Remark.}} The inclusions in b) can be strict as it is known from the case $E:=\bk$. 

\begin{lem}\label{766} 
Let $G$ be a W*-algebra and $F$ a C*-subalgebra of $G$. Then the following are equivalent.
\begin{enumerate}
\item $F$ generates $G$ as a W*-algebra.
\item $F^{\#}$ is dense in $G_{\ddot G}^{\#}$.
\item $F$ is dense in $G_{\ddot G}$.
\end{enumerate}
\end{lem}

$a\Longrightarrow b$ follows from [C1] Corollary 6.3.8.7.

$b\Longrightarrow c$ is trivial.

$c\Longrightarrow a$ follows from [C1] Corollary 4.4.4.12 a).\qed

\begin{p}\label{764} 
Let $G$ be a W*-algebra, $F$ a C*-subalgebra of $G$ generating it as W*-algebra, $I$ a set, and
$$L:=\cb{i\in I}\breve F\,,\qquad M:=\cw{i\in I}\breve G\;.$$
\begin{enumerate}
\item $M$ is the extension of $L$ to a \ris{G} \emph{([C2] \pr 1.3 f))} and $L^{\#}$ is dense in $M_{\ddot M}^{\#}$.
\item If we denote for every $X\in \lb{F}{L}$ by $\bar X\in \lb{G}{M}$ its unique extension \emph{([C3] \pr 1.4 a))} then the map
$$\mad{\lb{F}{L}}{\lb{G}{M}}{X}{\bar X}$$
is an injective C*-homomorphism and its image is dense in $\lh{G}{M}$.
\item The map
$$\mad{\lb{F}{L}_{\fr{T}_2}^{\#}}{\lb{G}{M}_{\fr{T}_1}^{\#}}{X}{\bar X}$$
is continuous.
\end{enumerate}
\end{p}

a) By Lemma \ref{766} $a\Rightarrow b$, $F^{\#}$ is dense in $G_{\ddot G}^{\#}$ so $\breve F^{\#}$ is dense in $\breve G_{\ddot {\breve G}}^{\#}$ and $\breve G$ is the extension of $\breve F$ to a \ris{G} ([C3] Corollary 1.5 $a_2\Rightarrow a_1$). By [C3] \pr 1.8, $M$ is the extension of $L$ to a \ris{G}. By [C3] Corollary 1.5 $a_1\Rightarrow a_2$, $L^{\#}$ is dense in $M_{\ddot M}^{\#}$.

b) By a) and [C3] \pr 1.4 e), the map
$$\mad{\lb{F}{L}}{\lb{G}{M}}{X}{\bar X}$$
is an injective C*-homomorphism. By [C3] \pr 1.9 b), its image is dense in $\lh{G}{M}$.

c) Denote by $N$ the vector subspace of $\stackrel{...}{M}$ generated by
$$\me{\tia{a}{\xi }{\eta }}{(a,\xi ,\eta )\in \ddot G\times L\times L}\;.$$
By a) and [C3] \pr 1.9 a), $N$ is dense in $\stackrel{...}{M}$ so by Corollary \ref{742} c),
$$\lb{G}{M}_{\fr{T}_1}^{\#}=\lb{G}{M}_N^{\#}\;.$$
For $(a,\xi ,\eta )\in \ddot G_+\times L\times L$ and $X\in \lb{F}{L}$, by \pr \ref{741} c),
$$p_{\xi ,\eta ,a}(\bar X)=|\sa{\s{\bar X\xi }{\eta }}{a}|=$$
$$=|\sa{\s{X\xi }{\eta }}{a}|\leq p_{\xi x,|a|}(X)\sa{\s{\eta }{\eta }}{|a|}^{\frac{1}{2}}\,,$$
where $a=x|a|$ is the polar representation of $a$, so the map
$$\mad{\lb{F}{L}_{\fr{T}_2}^{\#}}{\lb{G}{M}_{\fr{T}_1}^{\#}}{X}{\bar X}$$
is continuous.\qed

\begin{lem}\label{0}
Let $n\in \bn$, $\xi \in \cb{i\in \bnn{n}}\breve E$, and
$$x:=[\xi _i\delta _{j,1}]_{i,j\in \bnn{n}}\in E_{n,n}\;.$$
Then $\n{x}=\n{\xi }$.
\end{lem}

For $\eta \in \cb{i\in \bnn{n}}\breve E$ and $i\in \bnn{n}$,
$$(x\eta )_i=\si{j\in \bnn{n}}x_{ij}\eta _j=\si{j\in \bnn{n}}\xi _i\delta _{j,1}\eta _j=\xi _i\eta _1\,,$$
$$\s{x\eta }{x\eta }=\si{i\in \bnn{n}}\s{(x\eta )_i}{(x\eta )_i}=
\si{i\in \bnn{n}}\s{\xi _i\eta _1}{\xi _i\eta _1}=
\si{i\in \bnn{n}}\eta _1^*\xi _i^*\xi _i\eta _1=$$
$$=\eta _1^*\left(\si{i\in \bnn{n}}\xi _i^*\xi _i\right)\eta _1=
\eta _1^*\s{\xi }{\xi }\eta _1\leq \n{\xi }^2\eta _1^*\eta _1\,,$$
$$\n{x\eta }^2\leq \n{\xi }^2\n{\eta _1}^2\leq \n{\xi }^2\n{\eta }^2\,,\qquad \n{x}\leq \n{\xi }\;.$$
On the other hand if we put $\zeta :=(\delta _{i,1}1_E)_{i\in \bnn{n}}$ then for $i\in \bnn{n}$,
$$(x\zeta )_i=\si{j\in \bnn{n}}x_{ij}\zeta _j=\si{j\in \bnn{n}}\xi _i\delta _{j,1}1_E=\xi _i\,,$$
$$\s{x\zeta }{x\zeta }=\si{i\in \bnn{n}}(x\zeta )_i^*(x\zeta )_i=\si{i\in \bnn{n}}\xi _i^*\xi _i=\s{\xi }{\xi }\,,$$
$$\n{x}\geq \n{x\zeta }=\n{\xi }\,,\quad \n{x}=\n{\xi }\;.\qedd$$

\begin{lem}\label{755}
Let $F,G$ be unital C**-algebras, $\varphi :F\rightarrow G$ a surjective C**-homomorphism, $I$ a set,
$$L:=\widetilde{\cb{i\in I}}\breve F\approx \breve F\widetilde\otimes l^2(I)\,,\quad M:=\widetilde{\cb{i\in I}}\breve G\approx \breve G\widetilde\otimes l^2(I)\,,$$
and for every $\xi \in L$ put $\tilde \xi :=(\varphi \xi _i)_{i\in I}$.
\begin{enumerate}
\item If $\xi ,\eta \in L$ and $x\in F$ then
$$\tilde \xi \in M\,,\quad \n{\tilde \xi }\leq \n{\xi }\,,\quad \widetilde{(\xi x)}=(\tilde \xi )\varphi x\,,\quad \s{\tilde \xi }{\tilde \eta }=\varphi \s{\xi }{\eta }\;.$$
\item For every $\eta \in M$ there is a $\xi \in L$ with $\tilde \xi =\eta ,\;\n{\xi }=\n{\eta }$. 
\item In the W*-case the map 
$$\mad{L_{\ddot L}}{M_{\ddot M}}{\xi }{\tilde \xi }$$
is continuous.
\end{enumerate}
\end{lem}

a) For $J\in \fr{P}_f(I)$,
$$\si{i\in J}\s{\varphi \xi _i}{\varphi \eta _i}=\si{i\in J}(\varphi \eta _i)^*(\varphi \xi _i)=\varphi \si{i\in J}\eta _i^*\xi _i\;.$$
It follows $\tilde \xi \in M,\;\n{\tilde \xi }\leq \n{\xi },\;\s{\tilde \xi }{\tilde \eta }=\varphi \s{\xi }{\eta }$. Moreover for $i\in I$,
$$(\widetilde{\xi x})_i=\varphi (\xi x)_i=\varphi (\xi _ix)=(\varphi \xi _i)(\varphi x)=\tilde \xi_i(\varphi x),\quad \widetilde{\xi x}=\tilde \xi (\varphi x)\;. $$

b)\begin{center}
Case 1 $\me{i\in I}{\eta _i\not=0}$ is finite
\end{center}

For simplicity we assume $\me{i\in I}{\eta _i\not=0}=\bnn{n}$ for some $n\in \bn$. We put
$$\mae{\theta }{F_{n,n}}{G_{n,n}}{[x_{ij}]_{i,j\in \bnn{n}}}{[\varphi x_{ij}]_{i,j\in \bnn{n}}}\;.$$
$\theta $ is obviously a surjective C*-homomorphism. So if we put
$$y:=[\eta _i\delta _{j,1}]_{i,j\in \bnn{n}}\in G_{n,n}\,,$$
then there is an $x\in F_{n,n}$ with $\theta x=y,\,\n{x}=\n{y}$ ([K] \h 10.1.7). If we put
$$\mae{\xi }{I}{\breve F}{i}{\ab{x_{i1}}{i\in \bnn{n}}{0}{i\in I\setminus \bnn{n}}}$$
and $z:=[x_{ij}\delta _{j1}]_{i,j\in \bnn{n}}\in F_{n,n}$ then
$$\theta z=[\varphi (x_{ij}\delta _{j1})]_{i,j\in \bnn{n}}=[y_{ij}\delta _{j1}]_{i,j\in \bnn{n}}=y$$
and by \hh 6.1 a), $\n{z}\leq \n{x}$. We get for $i\in \bnn{n}$, 
$$\tilde \xi _i=\varphi \xi _i=\varphi x_{i1}=y_{i1}=\eta _i\;.$$
By a) and Lemma \ref{0},
$$\n{\xi }=\n{z}\leq \n{x}=\n{y}=\n{\eta }=\n{\tilde \xi }\leq \n{\xi }\,,\qquad \n{\xi }=\n{\eta }\;.$$

\begin{center}
Case 2 $\eta $ arbitrary in the W*-case
\end{center}

We may assume $\n{\eta }=1$. We put for every $J\in \fr{P}_f(I)$,
$$\mae{\eta _J}{I}{G}{i}{\ab{\eta _i}{i\in J}{0}{i\in I\setminus J}}\;.$$
By Case 1, for every $J\in \fr{P}_f(I)$ there is a $\xi _J\in L$ with $\tilde \xi _J=\eta _J$ and $\n{\xi _J}=\n{\eta _J}\leq 1$. Let $\fr{F}$ be an ultrafilter on $\fr{P}_f(I)$ finer than the upper section filter of $\fr{P}_f(I)$. By \prr 3.3 $a\Rightarrow b$,
$$\xi :=\lim_{J,\fr{F}}\xi _J$$
exists in $L^{\#}_{\ddot L}$. For $i\in I$,
$$\tilde \xi _i=\varphi \xi _i=\varphi \lim_{J,\fr{F}}(\xi _J)_i=\lim_{J,\fr{F}}\varphi (\xi _J)_i=\eta _i$$
so $\tilde \xi =\eta $. By a), $1=\n{\eta }=\n{\tilde \xi }\leq \n{\xi }\leq 1$, so $\n{\xi }=\n{\eta }$.

\begin{center}
Case 3 $\eta $ arbitrary in the C*-case
\end{center}

We put for every $J\in \fr{P}_f(I)$ and every $\zeta \in M$,
$$\mae{\zeta _J}{I}{G}{i}{\ab{\zeta _i}{i\in J}{0}{i\in I\setminus J}}\;.$$
Moreover we denote by $\fr{F}_I$ the upper section filter of $\fr{P}_f(I)$, set 
$$M_0:=\me{\zeta \in M}{\me{i\in I}{\zeta  _i\not=0}\;\mbox{is finite}}\,,$$
and denote by $\ccc{M}$ the vector subspace of $\ccc{K}_G(M)$ generated by the set 
$$\me{\zeta _1 \s{\cdot }{\zeta _2}}{\zeta _1,\zeta _2\in M_0}\;.$$

Let $\ccc{G}$ be the vector subspace of $\ccc{K}_F(L)$ generated by the set
$$\me{\alpha  \s{\cdot }{\beta }}{\alpha ,\beta \in L}\;.$$
$\ccc{G}$ is an involutive subalgebra of $\ccc{K}_F(L)$. Let $(\alpha _q)_{q\in Q}$, $(\beta  _q)_{q\in Q}$ be finite families in $L$ such that
$$\si{q\in Q}\alpha _q\s{\cdot }{\beta _q}=0\;.$$
Let further $\alpha ',\beta '\in M_0$. By Case 1, there are $\alpha ,\beta \in L$ with $\tilde \alpha =\alpha '\,,\tilde \beta =\beta '$ and we get by a),
$$\s{\si{q\in Q}\tilde \alpha _q\s{\beta '}{\tilde \beta _q}}{\alpha '}=\si{q\in Q}\s{\tilde \alpha _q}{\alpha '}\s{\beta '}{\tilde \beta _q}=$$
$$=\si{q\in Q}\s{\tilde \alpha _q}{\tilde \alpha }\s{\tilde \beta }{\tilde \beta _q}
=\varphi \left(\si{q\in Q}\s{\alpha _q}{\alpha }\s{\beta }{\beta _q}\right)=$$
$$=\varphi \left(\s{\left(\si{q\in Q}\alpha _q\s{\cdot }{\beta _q}\right)\beta }{\alpha }\right)=0\;.$$
It follows (\prr 4.1 e))
$$\si{q\in Q}\tilde \alpha _q\s{\cdot }{\tilde \beta _q}=0\;.$$
Thus the linear map
$$\mae{\psi }{\ccc{G}}{\ccc{K}_G(M)}{\si{q\in Q}\alpha _q\s{\cdot }{\beta _q}}{\si{q\in Q}\tilde \alpha _q\s{\cdot }{\tilde \beta _q}}$$
is well-defined and it is easy to see (by a)) that $\psi $ is an involutive algebra homomorphism.

\begin{center}
Step 1 $\n{\psi }\leq 1$; we extend $\psi $ by continuity to a map $\psi :\ccc{K}_F(L)\rightarrow \ccc{K}_G(M)$
\end{center}

Let 
$$u:=\si{q\in Q}\alpha _q\s{\cdot }{\beta _q}\in \ccc{G}$$
and let $\zeta \in M_0^{\#}$. By Case 1, there is an $\alpha \in L^{\#}$ with $\tilde \alpha =\zeta $. By a),
$$(\psi u)\zeta =\si{q\in Q}\tilde \alpha _q\s{\tilde \alpha }{\tilde \beta _q}=\si{q\in Q}\tilde \alpha _q\varphi \s{\alpha }{\beta _q}=\si{q\in Q}\widetilde{\overbrace{\alpha _q\s{\alpha }{\beta _q}}}=\widetilde{u\alpha }\,,$$
$$\n{(\psi u)\zeta }=\n{\widetilde{u\alpha }}\leq \n{u\alpha }\leq \n{u}\;.$$
Since $M_0$ is dense in $M$ (\prr 4.1 e)), it follows
$$\n{\psi u}\leq \n{u},\qquad \n{\psi }\leq 1\;.$$

\begin{center}
Step 2 $\ccc{M}$ is dense in $\ccc{K}_G(M)$
\end{center}

Let $\alpha ,\beta \in M$. By \prr 4.1 e),
$$\alpha =\lim_{J,\fr{F}_I}\alpha _J,\qquad \beta =\lim_{J,\fr{F}_I}\beta _J$$
so by \prr 5.2 a),
$$\alpha \s{\cdot }{\beta }=\lim_{J,\fr{F}_I}\alpha _J\s{\cdot }{\beta _J}\,,$$
which proves the assertion.

\begin{center}
Step 3 $\psi $ is a surjective C*-homomorphism
\end{center}

By Step 1, $\psi $ is a C*-homomorphism. Since its image contains $\ccc{M}$ (by Case 1) it is surjective by Step 2.

\begin{center}
Step 4 The assertion
\end{center}

Let $j\in I$. By Step 3 and [K] \h 10.1.7 (and \prr 5.2 a)), there is a $u\in \ccc{K}_F(L)$ with 
$$\psi u=\eta \s{\cdot }{1_G\otimes e_j},\qquad \n{u}=\n{\eta \s{\cdot }{1_G\otimes e_j}}=\n{\eta }\;.$$
From
$$\psi (u((1_F\otimes e_j)\s{\cdot }{1_F\otimes e_j}))=(\eta \s{\cdot }{1_G\otimes e_j})((1_G\otimes e_j)\s{\cdot }{1_G\otimes e_j})=$$
$$=\eta \s{\cdot }{1_G\otimes e_j}\,,$$
$$\n{\eta }=\n{\eta \s{\cdot }{1_G\otimes e_j}}\leq \n{u((1_F\otimes e_j)\s{\cdot }{1_F\otimes e_j})}\leq $$
$$\leq \n{u}\n{(1_F\otimes e_j)\s{\cdot }{1_F\otimes e_j}}=
\n{u}=\n{\eta }\,,$$
$$\n{u((1_F\otimes e_j)\s{\cdot }{1_F\otimes e_j})}=\n{\eta }$$
we see that we may assume
$$u=u((1_F\otimes e_j)\s{\cdot }{1_F\otimes e_j})\;.$$
Then
$$u=(u(1_F\otimes e_j))\s{\cdot }{1_F\otimes e_j}\;.$$
If we put $\xi :=u(1_F\otimes e_j)\in L$ then $u=\xi \s{\cdot }{1_F\otimes e_j}$, $\n{\eta }=\n{u}=\n{\xi }$,
$$\eta \s{\cdot }{1_G\otimes e_j}=\psi u=\tilde \xi \s{\cdot }{1_G\otimes e_j})\,,$$
$$\eta =\eta \s{1_G\otimes e_j}{1_G\otimes e_j}=\tilde \xi \s{1_G\otimes e_j}{1_G\otimes e_j}=\tilde \xi \;.$$

c) Let $(a,\eta _0)\in \ddot G\times M$. By b), there is a $\xi _0\in L$ with $\tilde \xi _0=\eta _0$. By a), for $\xi \in L$,
$$\sa{\tilde \xi }{\widetilde{(a,\eta _0)}}=\sa{\s{\tilde \xi }{\eta _0}}{a}=\sa{\s{\tilde \xi }{\tilde \xi _0}}{a}=$$
$$=\sa{\varphi \s{\xi }{\xi _0}}{a}=\sa{\s{\xi }{\xi _0}}{\ddot \varphi a}=\sa{\xi }{\widetilde{(\ddot \varphi a,\xi _0)}}\;.$$
We put
$$\mae{\theta }{L}{M}{\xi }{\tilde \xi }$$
and denote by $\theta ':M'\rightarrow L'$ its transpose. By the above, $\theta '\widetilde{(a,\eta _0)}\in \ddot L$. Since $\theta '$ is continuous, $\theta '(\ddot M)\subset \ddot L$ and this
 proves the assertion.\qed

\begin{p}\label{1}
We use the notation of \emph{Lemma \ref{755}}.
\begin{enumerate}
\item If $X\in \lb{F}{L}$ and $\xi \in L$ with $\tilde \xi =0$ then $\widetilde{X\xi }=0$; we define
$$\mae{\tilde X}{M}{M}{\eta }{\widetilde{X\xi }}\,,$$
where $\xi \in L$ with $\tilde \xi =\eta $ \emph{(Lemma \ref{755} b))}.
\item For every $X\in \lb{F}{L}$, $\tilde X$ belongs to $\lb{G}{M}$ and the map
$$\mad{\lb{F}{L}}{\lb{G}{M}}{X}{\tilde X}$$
is a surjective C**-homomorphism continuous with respect to the topologies $\fr{T}_k$ with $k\in \{1,2,3\}$. 
\item For $\xi ,\eta \in L$,
$$\widetilde{\overbrace{\eta \s{\cdot }{\xi }}}=\tilde \eta \s{\cdot }{\tilde \xi }$$
and
$$\ccc{K}_G(M)=\me{\tilde X}{X\in \ccc{K}_F(L)}\;.$$
\end{enumerate} 
\end{p}

a) For $i\in I$, $\varphi \xi _i=\tilde \xi _i=0$ so by Lemma \ref{755} a),
$$\widetilde{X(e_i\xi _i)}=\widetilde{(Xe_i)\xi _i}=\widetilde{(Xe_i)}\varphi \xi _i=0\;.$$
By \prr 4.1 e) (resp. \prr  4.6 c) and \prr 3.4 c)),
$$X\xi =X\left(\si{i\in I}e_i\xi _i\right)=\si{i\in I}X(e_i\xi _i)$$
$$\left(\mbox{resp.}\;X\xi =X\left(\sii{i\in I}{\ddot L}e_i\xi _i\right)=\sii{i\in I}{\ddot L}X(e_i\xi _i)\right)$$
so by Lemma \ref{755} a) (resp. c)),
$$\widetilde{X\xi }=\widetilde{\overbrace{\si{i\in I}X(e_i\xi _i)}}=\si{i\in I}\widetilde{X(e_i\xi _i)}=0$$
$$\left(\mbox{resp.}\;\widetilde{X\xi }=\widetilde{\overbrace{\sii{i\in I}{\ddot L}X(e_i\xi _i)}}=\sii{i\in I}{\ddot M}\widetilde{X(e_i\xi _i)}=0\right)\;.$$

b) For $X,Y\in \lb{F}{L}$ and $\xi ,\eta \in L$, by Lemma \ref{755} a),
$$\s{\tilde X\tilde \xi }{\tilde \eta }=\s{\widetilde{X\xi }}{\tilde \eta }=\varphi \s{X\xi }{\eta }=$$
$$=\varphi \s{\xi }{X^*\eta }=\s{\tilde \xi }{\widetilde{X^*\eta }}=\s{\tilde \xi }{\widetilde {X^*}\tilde \eta }\,,$$
$$\tilde X\tilde Y\tilde \xi =\tilde X\widetilde{Y\xi }=\widetilde{X(Y\xi )}=\widetilde{(XY)\xi }=\widetilde{XY}\tilde \xi \;.$$
By Lemma \ref{755} b), $\tilde X\in \lb{G}{M}$, $(\tilde X)^*=\tilde X^*$, and $\tilde X\tilde Y=\widetilde{XY}$, i.e. the map is a C*-homomorphism.

For $X\in \lb{F}{L}$ and $\xi ,\eta \in L$ (resp. and $a\in \ddot M_+$), by Lemma \ref{755} a),
$$p_{\tilde \xi ,\tilde \eta }(\tilde X)=\n{\s{\tilde X\tilde \xi }{\tilde \eta }}=\n{\s{\widetilde{X\xi }}{\tilde \eta }}=\n{\varphi \s{X\xi }{\eta }}\leq p_{\xi ,\eta }(X)$$
$$(\mbox{resp.}\;p_{\tilde \xi ,\tilde \eta ,a}(X)=\left|\sa{\s{\tilde X\tilde \xi }{\tilde \eta }}{a}\right|=|\sa{\varphi \s{X\xi }{\eta }}{a}|=$$
$$=|\sa{\s{X\xi }{\eta }}{\ddot \varphi a}|=p_{\xi ,\eta ,\ddot \varphi a}(X))\,,$$
so by Lemma \ref{755} b), the map is continuous with respect to the topology $\fr{T}_1$. The proof for the other topologies is similar.

c) For $\zeta \in L$, by Lemma \ref{755} a),
$$\widetilde{\overbrace{\eta \s{\cdot }{\xi }}}\,\tilde \zeta =\widetilde{\overbrace{(\eta \s{\cdot }{\xi })\zeta }}=\widetilde{\overbrace{\eta \s{\zeta }{\xi }}}=$$
$$=\tilde \eta \,\varphi \s{\zeta }{\xi }=\tilde \eta \s{\tilde \zeta }{\tilde \xi }=\left(\tilde \eta \s{\cdot }{\tilde \xi }\right)\tilde \zeta $$
so by Lemma \ref{755} b), 
$$\widetilde{\overbrace{\eta \s{\cdot }{\xi }}}=\tilde \eta \s{\cdot }{\tilde \xi }\;.$$ 
The last assertion follows now from b).\qed

\begin{center}
\section{Main Part}
\end{center}
\begin{center}
\fbox{\parbox[s]{7.5cm}{Throughout this section we fix $f\in \fte$}}
\subsection{The representations}
\end{center}

We present here the projective representation of the groups and its main properties.

\begin{de}\label{671} 
We put for every $t\in T$ and $\xi \in H$,
$$\mae{u_t}{\breve E}{H}{\zeta }{\zeta \otimes e_t}\,,$$
$$\mae{V_t\xi }{T}{\breve E}{s}{f(t,t^{-1}s)}\xi (t^{-1}s)\;.$$
\end{de}

If we want to emphasize the role of $f$ then we put $V_t^f$ instead of $V_t$. For $x\in E$,
$$\mae{(x\widetilde \otimes 1_K)V_t\xi }{T}{\breve E}{s}{f(t,t^{-1}s)x\xi (t^{-1}s)}\;.$$

\begin{p}\label{674} 
Let $s,t\in T$, $x\in E$, $\zeta \in \breve E$, and $\xi \in H$.
\begin{enumerate}
\item $V_t\xi \in H$.
\item $V_sV_t=(f(s,t)\widetilde \otimes 1_K)V_{st}$.
\item $V_t(\zeta \otimes e_s)=(f(t,s)\zeta )\otimes e_{ts}$.
\item $V_t(x\widetilde \otimes 1_K)=(x\widetilde \otimes 1_K)V_t$.
\item $V_t\in \unn\leh\,,\quad V_t^*=(\tilde f(t)\widetilde \otimes 1_K)V_{t^{-1}}$.
\item $(x\widetilde \otimes 1_K)V_t(\zeta \otimes e_s)=(f(t,s)x\zeta )\otimes e_{ts}$.
\item If $T$ is infinite and $\fr{F}$ denotes the filter on $T$ of cofinite subsets, i.e. 
$$\fr{F}:=\me{S}{S\in \fr{P}(T)\,,\;T\setminus S\in \fr{P}_f(T)}\,,$$
then
$$\lim_{t,\fr{F}}V_t=0$$
in $\lb{E}{H}_{\fr{T}_1}$.
\end{enumerate}
\end{p}

a) For $R\in \fr{P}_f(T)$,
$$\sum_{r\in R}\s{(V_t\xi )_r}{(V_t\xi )_r}=\sum_{r\in R}\s{f(t,t^{-1}r)\xi _{t^{-1}r}}{f(t,t^{-1}r)\xi _{t^{-1}r}}=$$
$$=\sum_{r\in R}\s{\xi _{t^{-1}r}}{\xi _{t^{-1}r}}=\sum_{r\in R}\s{\xi _r}{\xi _r}\leq \s{\xi }{\xi }$$
so $V_t\xi \in H$.

b) For $r\in T$,
$$(V_sV_t\xi )_r=f(s,s^{-1}r)(V_t\xi )_{s^{-1}r}=f(s,s^{-1}r)f(t,t^{-1}s^{-1}r)\xi _{t^{-1}s^{-1}r}=$$
$$=f(s,t)f(st,t^{-1}s^{-1}r)\xi _{t^{-1}s^{-1}r}=f(s,t)(V_{st}\xi )_r=((f(s,t)\widetilde\otimes 1_K)V_{st}\xi )_r$$
so
$$V_sV_t=(f(s,t)\widetilde\otimes 1_K)V_{st}\;.$$

c) For $r\in T$,
$$(V_t(\zeta \otimes e_s))_r=f(t,t^{-1}r)(\zeta \otimes e_s)_{t^{-1}r}=$$
$$=\delta _{s,t^{-1}r}f(t,t^{-1}r)\zeta =
\delta _{r,ts}f(t,s)\zeta =((f(t,s)\zeta )\otimes e_{ts})_r$$
so
$$V_t(\zeta \otimes e_s)=(f(t,s)\zeta )\otimes e_{ts}\;.$$

d) We have
$$(V_t(x\widetilde\otimes 1_K)\xi )_s=f(t,t^{-1}s)((x\widetilde\otimes 1_K)\xi )_{t^{-1}s}=f(t,t^{-1}s)x\xi _{t^{-1}s}=
((x\widetilde\otimes 1_K)V_t\xi )_s$$
so
$$V_t(x\widetilde\otimes 1_K)=(x\widetilde\otimes 1_K)V_t\;.$$

e) For $\eta \in H$, by \pr \ref{704} a),b),
$$\s{V_t\xi }{\eta }=\siw{s\in T}\s{(V_t\xi )_s}{\eta _s}=\siw{s\in T}\s{f(t,t^{-1}s)\xi _{t^{-1}s}}{\eta _s}=$$
$$=\siw{r\in T}\s{f(t,r)\xi _r}{\eta _{tr}}=\siw{r\in T}\s{\xi _r}{\tilde f(t)f(t^{-1},tr)\eta _{tr}}=$$
$$=\siw{r\in T}\s{\xi _r}{(((\tilde f(t)\widetilde\otimes 1_K)V_{t^{-1}})\eta )_r}
=\s{\xi }{((\tilde f(t)\widetilde\otimes 1_K)V_{t^{-1}})\eta }$$
so $V_t\in \lb{E}{H}$ with $V_t^*=(\tilde f(t)\widetilde\otimes 1_K)V_{t^{-1}}$. By b) and d),
$$V_t^*V_t=(\tilde f(t)\widetilde\otimes 1_K)V_{t^{-1}}V_t=(\tilde f(t)\widetilde\otimes 1_K)(f(t^{-1},t)\widetilde\otimes 1_K)V_{t^{-1}t}=id_H\,,$$
$$V_tV_t^*=V_t(\tilde f(t)\widetilde\otimes 1_K)V_{t^{-1}}=(\tilde f(t)\widetilde\otimes 1_K)V_tV_{t^{-1}}=$$
$$=(\tilde f(t)\widetilde\otimes 1_K)(f(t,t^{-1})\widetilde\otimes 1_K)V_{tt^{-1}}=id_H\;.$$

f) follows from c).

g) Let us consider first the C*-case. Let $\xi ,\eta \in H$, $t\in T$, and $\varepsilon >0$. There is an $S\in \fr{P}_f(T)$ such that $\n{\eta e_{T\setminus S}}<\varepsilon $. By e),
$$|\s{V_t\xi }{\eta e_{T\setminus S}}|\leq \n{V_t\xi }\n{\eta e_{T\setminus S}}\leq \varepsilon \n{\xi }$$
so
$$p_{\xi ,\eta }(V_t)=|\s{V_t\xi }{\eta }|\leq |\s{V_t\xi }{\eta e_S}|+|\s{V_t\xi }{\eta e_{T\setminus S}}|<|\s{V_t\xi }{\eta e_S}|+\varepsilon \;.$$
From
$$\s{V_t\xi }{\eta e_S}=\si{s\in S}\eta _s^*f(t,t^{-1}s)\xi _{t^{-1}s}$$
it follows
$$\lim_{t,\fr{F}}\s{V_t\xi }{\eta e_S}=0\,,\qquad\lim_{t,\fr{F}}p_{\xi ,\eta }(V_t)=0\;.$$

The W*-case can be proved similarly.\qed

{\it Remark.} By e), $\fr{T}_1$ cannot be replaced by $\fr{T}_2$ in g).

\begin{p}\label{676}
Let $s,t\in T$.\begin{enumerate}
\item $u_t\in \lc{E}{\breve E}{H}\,,\quad u_t^*=\s{\;\cdot\;}{1_E\otimes e_t}$.
\item $u_s^*u_t=\delta _{s,t}1_E$.
\item $u_su_t^*=1_E\widetilde\otimes (\s{\cdot }{e_t}e_s)$.
\item $\sii{r\in T}{\fr{T_2}}u_ru_r^*=id_H$.
\end{enumerate}
\end{p}

a) For $\zeta \in \breve E$ and $\xi \in H$,
$$\s{u_t\zeta }{\xi }=\s{\zeta \otimes e_t}{\xi }=\siw{s\in T}\xi _s^*(\zeta \otimes e_t)_s=\xi _t^*\zeta =\s{\zeta }{\xi _t}$$
so
$$u_t\in \lc{E}{\breve E}{H}\,,\qquad u_t^*\xi =\xi _t=\s{\xi }{1_E\otimes e_t}\;.$$

b) For $\zeta \in \breve E$, by a),
$$u_s^*u_t\zeta =u_s^*(\zeta \otimes e_t)=\s{\zeta \otimes e_t}{1_E\otimes e_s}=\delta _{s,t}\zeta $$
so $u_s^*u_t=\delta _{s,t}1_E$.

c) For $\zeta \in \breve E$ and $r\in T$, by a),
$$u_su_t^*(\zeta \otimes e_r)=u_s\delta _{r,t}\zeta =\delta _{r,t}(\zeta \otimes e_s)=$$
$$=\zeta \otimes \s{e_r}{e_t}e_s=(1_E\widetilde\otimes (\s{\cdot }{e_t}e_s))(\zeta \otimes e_r)\,,$$
so (by a) and \prr 4.1 e) (resp. and \prr 4.6 c), \prr 3.4 c))) $u_su_t^*=1_E\widetilde\otimes (\s{\cdot }{e_t}e_s)$.

d) For $\xi \in H$ (resp. and $a\in \ddot E_+$) and $S\in \fr{P}_f(T)$, by c),
$$p_\xi \left(\si{t\in S}u_tu_t^*-id_H\right)=\n{\si{t\in T\setminus S}\s{\xi }{\xi }}^{1/2}$$
$$\Bigg{(}\mbox{resp.}\; p_{\xi ,a}\left(\si{t\in S}u_tu_t^*-id_H\right)=$$
$$=\sa{\s{\si{t\in S}(u_tu_t^*-id_H)\xi }{\si{t\in S}(u_tu_t^*-id_H)\xi }}{a}^{1/2}=$$
$$=(\si{t\in T\setminus S}\sa{\s{\xi }{\xi }}{a})^{1/2}\Bigg{)}$$
and the assertion follows.\qed

\begin{p}\label{677}
Let $s,t\in T$ and $x\in E$.
\begin{enumerate}
\item $V_su_t=u_{st}f(s,t)$.
\item $u_s^*V_t=f(t,t^{-1}s)u_{t^{-1}s}^*$.
\item $(x\widetilde\otimes 1_K)u_t=u_tx$.
\item $xu_t^*=u_t^*(x\widetilde\otimes 1_K)$.
\end{enumerate}
\end{p}

a) For $\zeta \in \breve E$, by \pr \ref{674} c),
$$V_su_t\zeta =V_s(\zeta \otimes e_t)=(f(s,t)\zeta )\otimes e_{st}=u_{st}f(s,t)\zeta $$
so $V_su_t=u_{st}f(s,t)\;.$

b) For $\zeta \in \breve E$ and $r\in T$, by \pr \ref{674} c) and \pr \ref{676} a),
$$u_s^*V_t(\zeta \otimes e_r)=u_s^*((f(t,r)\zeta )\otimes e_{tr})=\delta _{s,tr}f(t,r)\zeta =$$
$$=\delta _{t^{-1}s,r}f(t,t^{-1}s)\zeta=f(t,t^{-1}s)u_{t^{-1}s}^*(\zeta \otimes e_r)$$
so $u_s^*V_t=f(t,t^{-1}s)u_{t^{-1}s}^*$.

c) For $\zeta \in \breve E$,
$$(x\widetilde\otimes 1_K)u_t\zeta =(x\widetilde\otimes 1_K)(\zeta \otimes e_t)=(x\zeta )\otimes e_t=u_tx\zeta $$
so $(x\widetilde\otimes 1_K)u_t=u_tx$.

d) follows from c).\qed

\begin{de}\label{680}
We put for all $s,t\in T$ \emph{(\pr \ref{676} a))}
$$\mae{\varphi _{s,t}}{\lb{E}{H}}{\lb{E}{\breve E}\approx E}{X}{u_s^*Xu_t}$$
and set $X_t:=\varphi _{t,1}X$ for every $X\in \lb{E}{H}$.
\end{de}

\begin{p}\label{681}
Let $s,t\in T$.
\begin{enumerate}
\item $\varphi _{s,t}$ is linear with $\n{\varphi _{s,t}}=1$. 
\item For $X\in \lb{E}{H}$ and $x,y\in \breve E$,
$$\s{(\varphi _{s,t}X)x}{y}=\s{X(x\otimes e_t)}{y\otimes e_s}\;.$$
\item The map
$$\mac{\varphi _{s,t}}{\lb{E}{H}_{\fr{T}_1}}{E\;(\mbox{\emph{resp.}}\;E_{\ddot E})}$$
is continuous.
\item $\varphi _{t,t}$ is involutive and completely positive.
\item For $r\in T$ and $x\in E$,
$$\varphi _{s,t}((x\widetilde\otimes 1_K)V_r)=\delta _{s,rt}f(r,t)x\;.$$
\item If $(x_r)_{r\in T}\in E^{(T)}$ and
$$X:=\si{r\in T}(x_r\widetilde\otimes 1_K)V_r$$
then
$$\varphi _{s,t}X=f(st^{-1},t)x_{st^{-1}}\,,\qquad X_t=x_t\;.$$
\item For $X\in \lb{E}{H}$ and $x,y\in E$,
$$\varphi _{s,t}((x\widetilde\otimes 1_K)X(y\widetilde\otimes 1_K))=x(\varphi _{s,t}X)y\,,$$
$$((x\widetilde\otimes 1_K)X(y\widetilde\otimes 1_K))_t=xX_ty\;.$$
\end{enumerate}
\end{p}

a) follows from \pr \ref{676} a),b).

b) We have
$$\s{(\varphi _{s,t}X)x}{y}=\s{u_s^*Xu_tx}{y}=\s{Xu_tx}{u_sy}=\s{X(x\otimes e_t)}{y\otimes e_s}\;.$$

c) 
\begin{center}
\emph{The C*-case}
\end{center}

By b), for $X\in \lb{E}{H}$,
$$\n{\varphi _{s,t}X}=\n{\s{(\varphi _{s,t}X)1_E}{1_E}}=$$
$$=\n{\s{X(1_E\otimes e_t)}{1_E\otimes e_s}}=p_{1_E\otimes e_t,1_E\otimes e_s}(X)\;.$$

\begin{center}
\emph{The W*-case}
\end{center}

Let $a\in \ddot E$ and let $a=x|a|$ be its polar representation. By b), for $X\in \lb{E}{H}$,
$$|\sa{\varphi _{s,t}X}{a}|=|\sa{\s{(\varphi _{s,t}X)1_E}{1_E}}{x|a|}|=|\sa{\s{(\varphi _{s,t}X)x}{1_E}}{|a|}|=$$
$$=|\sa{\s{X(x\otimes e_t)}{1_E\otimes e_s}}{|a|}|=p_{x\otimes e_t,1_E\otimes e_s,|a|}(X)\;.$$

d) For $X\in \lb{E}{H}$,
$$(\varphi _{t,t}X)^*=(u_t^*Xu_t)^*=u_t^*X^*u_t=\varphi _{t,t}(X^*)$$
so $\varphi _{t,t}$ is involutive. For $n\in \bn$, $X\in ((\lb{E}{H})_{n,n})_+$, and $\zeta \in \breve E^n$,
$$\si{i\in \bnn{n}}\s{\si{j\in \bnn{n}}((\varphi _{t,t}X_{ij})\zeta _j)}{\zeta _i}=\si{i,j\in \bnn{n}}\s{u_t^*X_{ij}u_t\zeta _j}{\zeta _i}=$$
$$=\si{i,j\in \bnn{n}}\s{X_{ij}u_t\zeta _j}{u_t\zeta _i}\geq 0$$
(\hh 6.1 f) and \hh 1.11 $c_1\Rightarrow c_2$) so $\varphi _{t,t}$ is completely positive (\hh 6.1 f) and \hh 1.11 $c_2\Rightarrow c_1$).

e) By \pr \ref{677} a),d) and \pr \ref{676} b),
$$\varphi _{s,t}((x\widetilde\otimes 1_K)V_r)=u_s^*(x\widetilde\otimes 1_K)V_ru_t=xu_s^*V_ru_t=xu_s^*u_{rt}f(r,t)=\delta _{s,rt}f(r,t)x\;.$$

f) By e) (and \pr \ref{704} a)),
$$\varphi _{s,t}X=\si{r\in T}\varphi _{s,t}((x_r\widetilde\otimes 1_K)V_r)=\si{r\in T}\delta _{s,rt}f(r,t)x_r=f(st^{-1},t)x_{st^{-1}}\,,$$
$$X_t=\varphi _{t,1}X=f(t,1)x_t=x_t\;.$$

g) By \pr \ref{677} c),d),
$$\varphi _{s,t}((x\widetilde\otimes 1_K)X(y\widetilde\otimes 1_K))=u_s^*(x\widetilde\otimes 1_K)X(y\widetilde\otimes 1_K)u_t=$$
$$=xu_s^*Xu_ty=x(\varphi _{s,t}X)y\;.\qedd$$

\begin{de}\label{682}
We put
$$\ccc{R}(f):=\me{\si{t\in T}(x_t\widetilde\otimes 1_K)V_t}{(x_t)_{t\in T}\in E^{(T)}}\,,$$
$$\ssa{f}:=\stackrel{\fr{T}_3}{\overline{\ccc{R}(f)}}\,,\qquad \qquad \ssb{\n{\cdot }}{f}:=\stackrel{\n{\cdot }}{\overline{\ccc{R}(f)}}\;.$$
Moreover we put $\ssb{C}{f}:=\ssa{f}$ in the C*-case and $\ssb{W}{f}:=\ssa{f}$ in the W*-case. If $F$ is a subset of $E$ then we put
$$\ssa{f,F}:=\me{X\in \ssa{f}}{t\in T\Longrightarrow X_t\in F}$$
and use similar notation for the other $\ccc{S}$.
\end{de}

By \pr \ref{674} b),d),e), $\ccc{R}(f)$ is an involutive unital $E$-subalgebra of $\lb{E}{H}$ (with $V_1$ as unit). In particular $\ssb{\n{\cdot }}{f}$ is an $E$-C*-subalgebra of $\lb{E}{H}$. If $T$ is finite then $\ccc{R}(f)=\ssa{f}$. By Corollary \ref{742} e), $\ssb{C}{f}_{\fr{T}_3}$ is complete.

\begin{p}\label{745a}
For $X\in \stackrel{\fr{T}_1}{\overline{\ccc{R}(f)}}$ and $s,t\in T$,
$$\varphi _{s,t}X=f(st^{-1},t)X_{st^{-1}}\;.$$
\end{p}

Let $\fr{F}$ be a filter on $\ccc{R}(f)$ converging to $X$ in the $\fr{T}_1$-topology. By \pr \ref{681} c),f) (and Corollary \ref{742} d)), 
$$\varphi _{s,t}X=\lim_{Y,\fr{F}}\varphi _{s,t}Y=\lim_{Y,\fr{F}}f(st^{-1},t)Y_{st^{-1}}=f(st^{-1},t)\lim_{Y,\fr{F}}Y_{st^{-1}}=$$
$$=f(st^{-1},t)\lim_{Y,\fr{F}}\varphi _{st^{-1},1}Y=f(st^{-1},t)\varphi _{st^{-1},1}X=f(st^{-1},t)X_{st^{-1}}\;.\qedd$$ 

\begin{theo}\label{745}
Let $X\in\stackrel{\fr{T}_1}{\overline {\ccc{R}(f)}}$.
\begin{enumerate}
\item If $(x_t)_{t\in T}$ is a family in $E$ such that
$$X=\sii{t\in T}{\fr{T}_1}(x_t\widetilde\otimes 1_K)V_t$$
then $X_t=x_t$ for every $t\in T$. In particular, if $T$ is finite then the map
$$\mad{E^T}{\ssa{f}}{x}{\si{t\in T}(x_t\otimes 1_K)V_t}$$
is bijective and $E$-linear \emph{(\pr \ref{674} d))}.
\item We have
$$X=\sii{t\in T}{\fr{T}_3}(X_t\widetilde\otimes 1_K)V_t\in \ssa{f}\;.$$
\item $(X^*)_t=\tilde f(t)(X_{t^{-1}})^*$ for every $t\in T$ and
$$X^*=\sii{t\in T}{\fr{T}_3}((X_t)^*\widetilde\otimes  1_K)V_t^*\in \stackrel{\fr{T}_3}{\overline{\ccc{R}(f)}}\;.$$
\item $\ssa{f}=\stackrel{\fr{T}_1}{\overline{\ccc{R}(f)}}=\stackrel{\fr{T}_2}{\overline{\ccc{R}(f)}}$.
\item For $\xi \in H$ and $t\in T$,
$$(X\xi )_t=\siw{s\in T}f(s,s^{-1}t)X_s\xi _{s^{-1}t}\;.$$
\item If $T$ is finite and if we identify $\lb{E}{H}$ with $E_{T,T}$ then $X$ is identified with the matrix
$$[f(st^{-1},t)X_{st^{-1}}]_{s,t\in T}\,,$$
and for every $r\in T$, $V_r$ is identified with the matrix
$$[f(st^{-1},t)\delta _{s,rt}]_{s,t\in T}\;.$$ 
\item If $X,Y\in \ssa{f}$ and $t\in T$ then $XY\in \ssa{f}$ and 
$$(XY)_t=\siw{s\in T}f(s,s^{-1}t)X_sY_{s^{-1}t}\,,$$ 
$$(X^*Y)_t=\siw{s\in T}f(s,t)^*X_s^*Y_{st}\,,\qquad (XY^*)_t=\siw{s\in T}f(t,s)^*X_{ts}Y_s^*\,,$$
$$ (X^*Y)_1=\siw{s\in T}X_s^*Y_s\,,\qquad (XY^*)_1=\siw{s\in T}X_sY_s^*\;.$$
\item The map
$$\mad{E}{\ssa{f}}{x}{x\widetilde\otimes 1_K}$$
is an injective unital C**-homomorphism and so $\ssa{f}$ is an $E$-C**-subalge-bra of $\lb{E}{H}$ and $Re\,\ssa{f}$ is closed in $\ssa{f}_{\fr{T}_1}$. In the W*-case, $\ssb{W}{f}$ is the W*-subalgebra of $\lb{E}{H}$ generated by $\ccc{R}(f)$ and $\ccc{R}(f)^{\#}$ is dense in $\ssb{W}{f}_{\fr{T}_1}^{\#}=\ssb{W}{f}_{\stackrel{...}{H}}^{\#}$, which is compact. 
\item If $E$ is a W*-algebra then $\ssb{C}{f}$ may be identified canonically with a unital C*-subalgebra of $\ssb{W}{f}$ by using the map of \pr \ref{764} b). By this identification $\ssb{C}{f}$ generates $\ssb{W}{f}$ as W*-algebra.
\item If $F$ is a closed ideal of $E$ (resp. of $E_{\ddot E}$ ) then $\ssa{f,F}$
is a closed ideal of $\ssa{f}$ {\LARGE(}resp. of $\ssa{f}_{\stackrel{..}{\overbrace{\ssa{f}}}}${\LARGE)}.
\item If $F$ is a unital C**-subalgebra of $E$ such that $f(s,t)\in F$ for all $s,t\in T$ then $\ssa{f,F}$ is a unital C**-subalgebra of $\ssa{f}$ and the map
$$\mad{\ssa{f,F}}{\ssa{g}}{X}{\sii{t\in T}{\fr{T}_3}(X_t\widetilde\otimes 1_K)V_t^g}$$
is an injective C**-homomorphism, where
$$\mae{g}{T\times T}{\un{F}}{(s,t)}{f(s,t)}\;.$$
This map induces a C*-isomorphism $\ssb{\n{\cdot }}{f,F}\rightarrow \ssb{\n{\cdot }}{g}$.
\item $(X,Y)\in \stackrel{\circ }{\overbrace{\ssa{f}}}_+\Longrightarrow (X_1,Y_1)\in \stackrel{\circ }{E}_+$.
\end{enumerate}
\end{theo}

a) By \pr \ref{681} c),e),
$$X_t=\varphi _{t,1}X=\siw{s\in T}\varphi _{t,1}((x_s\widetilde\otimes 1_K)V_s)=\siw{s\in T}\delta _{t,s}f(s,1)x_s=x_t\;.$$

b\&c\&d
\begin{center}
Step 1 $X=\sii{t\in T}{\fr{T}_2}(X_t\widetilde\otimes 1_K)V_t$
\end{center}

By \pr \ref{676} d), Corollary \ref{742} d), \pr \ref{745a}, and \pr \ref{677} b),d),
$$X=\left(\sii{s\in T}{\fr{T}_2}u_su_s^*\right)X\left(\sii{t\in T}{\fr{T}_2}u_tu_t^*\right)=\sii{s\in T}{\fr{T}_2}\sii{t\in T}{\fr{T}_2}u_su_s^*Xu_tu_t^*=$$
$$=\sii{s\in T}{\fr{T}_2}\sii{t\in T}{\fr{T}_2}u_s(\varphi _{s,t}X)u_t^*=\sii{s\in T}{\fr{T}_2}\sii{t\in T}{\fr{T}_2}u_sf(st^{-1},t)X_{st^{-1}}u_t^*=$$
$$=\sii{s\in T}{\fr{T}_2}\sii{r\in T}{\fr{T}_2}u_sX_rf(r,r^{-1}s)u_{r^{-1}s}^*=\sii{s\in T}{\fr{T}_2}\sii{r\in T}{\fr{T}_2}u_sX_ru_s^*V_r=$$
$$=\sii{s\in T}{\fr{T}_2}\sii{r\in T}{\fr{T}_2}u_su_s^*(X_r\widetilde\otimes 1_K)V_r=\sii{s\in T}{\fr{T}_2}u_su_s^*\left(\sii{t\in T}{\fr{T}_2}(X_t\widetilde\otimes 1_K)V_t\right)=
\sii{t\in T}{\fr{T}_2}(X_t\widetilde\otimes 1_K)V_t\;.$$

\begin{center}
Step 2  b\&c\&d
\end{center}

By Step 1, Corollary \ref{742} a), and \pr \ref{674} d),e) (and \pr \ref{704} a)),
$$X^*=\left(\sii{s\in T}{\fr{T}_1}(X_s\widetilde\otimes 1_K)V_s\right)^*=\sii{s\in T}{\fr{T}_1}(X_s^*\widetilde\otimes 1_K)V_s^*=$$
$$=\sii{s\in T}{\fr{T}_1}(X_s^*\widetilde\otimes 1_K)(\tilde f(s)\widetilde\otimes 1_K)V_{s^{-1}}=\sii{r\in T}{\fr{T}_1}((\tilde f(r)X_{r^{-1}}^*)\widetilde\otimes 1_K)V_r\in \stackrel{\fr{T}_1}{\overline{\ccc{R}(f)}}\;.$$
By a),
$$(X^*)_t=\tilde f(t)(X_{t^{-1}})^*\;.$$
By Step 1 and \pr \ref{674} e) (and \pr \ref{704} a)),
$$X^*=\sii{t\in T}{\fr{T}_2}((X^*)_t\widetilde\otimes 1_K)V_t=\sii{t\in T}{\fr{T}_2}((X_{t^{-1}})^*\widetilde\otimes 1_K)(\tilde f(t)\widetilde\otimes 1_K)V_t=$$
$$=\sii{t\in T}{\fr{T}_2}((X_{t^{-1}})^*\widetilde\otimes 1_K)V_{t^{-1}}^*=\sii{t\in T}{\fr{T}_2}((X_t)^*\widetilde\otimes 1_K)V_t^*\;.$$
Together with Step 1 this proves
$$X=\sii{t\in T}{\fr{T}_3}(X_t\widetilde\otimes 1_K)V_t\in \ssa{f}\,,\qquad X^*=\sii{t\in T}{\fr{T}_3}((X_t)^*\widetilde\otimes 1_K)V_t^*\in \ssa{f}\;.$$
In particular $\ssa{f}=\stackrel{\fr{T}_1}{\overline{\ccc{R}(f)}}=\stackrel{\fr{T}_2}{\overline{\ccc{R}(f)}}$.

e) By b) and Corollary \ref{742} b), in the C*-case,
$$(X\xi )_t=\s{\left(\sii{s\in T}{\fr{T}_1}(X_s\widetilde\otimes 1_K)V_s\right)\xi }{1_E\otimes e_t}=\si{s\in T}\s{(X_s\widetilde\otimes 1_K)V_s\xi }{1_E\otimes e_t}=$$
$$=\si{s\in T}X_sf(s,s^{-1}t)\xi _{s^{-1}t}=\si{s\in T}f(s,s^{-1}t)X_s\xi _{s^{-1}t}\;.$$
The proof is similar in the W*-case.

f) For $\xi \in H$ and $s\in T$, by e),
$$(X\xi )_s=\si{t\in T}f(t,t^{-1}s)X_t\xi _{t^{-1}s}=\si{r\in T}f(sr^{-1},r)X_{sr^{-1}}\xi _r\;.$$

g) By b), Corollary \ref{742} b),d), and \pr \ref{674} b),d),
$$XY=\left(\sii{s\in T}{\fr{T}_2}(X_s\widetilde\otimes 1_K)V_s\right)\left(\sii{t\in T}{\fr{T}_2}(X_t\widetilde\otimes 1_K)V_t\right)=$$
$$=\sii{s\in T}{\fr{T}_2}\sii{t\in T}{\fr{T}_2}(X_s\widetilde\otimes 1_K)V_s(Y_t\widetilde\otimes 1_K)V_t=\sii{s\in T}{\fr{T}_2}\sii{t\in T}{\fr{T}_2}(X_s\widetilde\otimes 1_K)(Y_t\widetilde\otimes 1_K)V_sV_t=$$
$$=\sii{s\in T}{\fr{T}_2}\sii{t\in T}{\fr{T}_2}(X_s\widetilde\otimes 1_K)(Y_t\widetilde\otimes 1_K)(f(s,t)\widetilde\otimes 1_K)V_{st}=$$
$$=\sii{s\in T}{\fr{T}_2}\sii{r\in T}{\fr{T}_2}((f(s,s^{-1}r)X_sY_{s^{-1}r})\widetilde\otimes 1_K)V_r\;.$$
Since by d), 
$$\sii{r\in T}{\fr{T}_2}((f(s,s^{-1}r)X_sY_{s^{-1}r})\widetilde\otimes 1_K)V_r\in \ssa{f}$$
for every $s\in T$ we get $XY\in \ssa{f}$, again by d). By Corollary \ref{742} b) and \pr \ref{681} c),e),
$$(XY)_t=\varphi _{t,1}(XY)=\siw{s\in T}\,\siw{r\in T}\varphi _{t,1}((f(s,s^{-1}r)X_sY_{s^{-1}r})\widetilde\otimes 1_K)V_r=$$
$$=\siw{s\in T}\,\siw{r\in T}\delta _{t,r}f(r,1)f(s,s^{-1}r)X_sY_{s^{-1}r}=\siw{s\in T}f(s,s^{-1}t)X_sY_{s^{-1}t}\;.$$
By the above, c), and \pr \ref{704} b),
$$(X^*Y)_t=\siw{s\in T}f(s,s^{-1}t)(X^*)_sY_{s^{-1}t}=\siw{s\in T}f(s,s^{-1}t)\tilde f(s)(X_{s^{-1}})^*Y_{s^{-1}t}=$$
$$=\siw{s\in T}f(s^{-1},t)^*(X_{s^{-1}})^*Y_{s^{-1}t}=\siw{s\in T}f(s,t)^*X_s^*Y_{st}\,,$$
$$(XY^*)_t=\siw{s\in T}f(s,s^{-1}t)X_s(Y^*)_{s^{-1}t}=\siw{s\in T}f(s,s^{-1}t)X_s\tilde f(s^{-1}t)(Y_{t^{-1}s})^*=$$
$$=\siw{s\in T}f(t,t^{-1}s)^*X_s(Y_{t^{-1}s})^*=\siw{s\in T}f(t,s)^*X_{ts}Y_s^*\;.$$
It follows by \pr \ref{704} a),
$$(X^*Y)_1=\siw{s\in T}X_s^*Y_s\,,\qquad (XY^*)_1=\siw{s\in T}X_sY_s^*\;.$$

h) By c) and g), $\ssa{f}$ is an involutive unital subalgebra of $\lb{E}{H}$. Being closed (resp. closed in $\lh{E}{H}$ (d) and Corollary \ref{742} c))) it is a C**-subalgebra of $\lb{E}{H}$ (resp. generated by $\ccc{R}(f)$ \hh 3.5 b) and [C1] Corollary 4.4.4.12 a) and by [C1] Corollary 6.3.8.7 $\ccc{R}(f)^{\#}$ is dense in $\ssb{W}{f}_{\fr{T}_1}^{\#}$, which is compact by Corollary \ref{742} c)). The assertion concerning $E$ follows from \pr \ref{674} d) and Lemma \ref{18} c). By Corollary \ref{742} a), $Re\,\ssa{f}$ is a closed set of $\ssa{f}_{\fr{T}_1}$.

i) The assertion follows from h), \pr \ref{764} b), and Lemma \ref{766} $c)\Rightarrow a)$.

j) For $X\in \ccc{S}(f,F)$, $Y\in \ssa{f}$, and $t\in T$, by g), $(XY)_t,(YX)_t\in \ccc{S}(f,F)$ so $\ccc{S}(f,F)$ is an ideal of $\ssa{f}$. The closure properties follow from \pr \ref{681} c).

k) By c) and g), $\ssa{f,F}$ is a unital involutive subalgebra of $\ssa{f}$ and by \pr \ref{681} c), $\ssa{f,F}$ is a C**-subalgebra of $\ssa{f}$. The last assertion follows from the fact that the image of the map contains $\ccc{R}(g)$.

l) There are $U,V\in \ssa{f}$ with
$$(X,Y)=(U,V)^*(U,V)=(U^*,-V^*)(U,V)=(U^*U+V^*V,U^*V-V^*U)\;.$$
For $t\in T$, 
$$0\leq (U_t,V_t)^*(U_t,V_t)=(U_t^*,-V_t^*)(U_t,V_t)=(U_t^*U_t+V_t^*V_t,U_t^*V_t-V_t^*U_t)\;.$$
By g),
$$X_1=(U^*U+V^*V)_1=\siw{t\in T}(U_t^*U_t+V_t^*V_t)\,,$$
$$Y_1=(U^*V-V^*U)_1=\siw{t\in T}(U_t^*V_t-V_t^*U_t)$$
so
$$(X_1,Y_1)=\siw{t\in T}(U_t^*U_t+V_t^*V_t,U_t^*V_t-V_t^*U_t)\in \stackrel{\circ }{E}_+\;.\qedd$$

{\it Remark.} It may happen that by the identification of i), $\ssb{C}{f}\not=\ssb{W}{f}$ (Remark of \pr \ref{750}).

\begin{co}\label{775}
\rule{0em}{1ex}
\begin{enumerate}
\item If $(x_t)_{t\in T}$ is a family in $E$ such that $(\n{x_t})_{t\in T}$ is summable then 
$$((x_t\widetilde\otimes 1_K)V_t)_{t\in T}$$
 is norm summable in $\lb{E}{H}$ and
$$\n{\si{t\in T}(x_t\widetilde\otimes 1_K)V_t}\leq \si{t\in T}\n{x_t}\;.$$
\item The set
$$\ccc{A}:=\me{X\in \ssa{f}}{\si{t\in T}\n{X_t}<\infty }$$
is a dense involutive unital subalgebra of $\ssb{\n{\cdot }}{f}$ with
$$\si{t\in T}\n{(X^*)_t}=\si{t\in T}\n{X_t}\,,$$
$$\si{t\in T}\n{(XY)_t}\leq \left(\si{t\in T}\n{X_t}\right)\left(\si{t\in T}\n{Y_t}\right)$$
for all $X,Y\in \ccc{A}$.
\item $\ccc{A}$ endowed with the norm
$$\mad{\ccc{A}}{\br_+}{X}{\si{t\in T}\n{X_t}}$$
is an involutive Banach algebra and $\ssb{\n{\cdot }}{f}$ is its C*-hull.
\end{enumerate}
\end{co}

a) For $S\in \fr{P}_f(T)$, by \pr \ref{674} e),
$$\n{\si{t\in S}(x_t\widetilde\otimes 1_K)V_t}\leq \si{t\in S}\n{x_t\widetilde\otimes 1_K}\n{V_t}=\si{t\in S}\n{x_t}$$
and the assertion follows.

b) By \h \ref{745} c), $X^*\in \ssa{f}$ and
$$\n{(X^*)_t}=\n{(X_{t^{-1}})^*}=\n{X_{t^{-1}}}$$
for all $t\in T$ so
$$\si{t\in T}\n{(X^*)_t}=\si{t\in T}\n{X_{t^{-1}}}=\si{t\in T}\n{X_t}\;.$$
By \h \ref{745} g), $XY\in \ssa{f}$ and
$$\n{(XY)_t}=\n{\siw{s\in T}f(s,s^{-1}t)X_sY_{s^{-1}t}}\leq \si{s\in T}\n{X_s}\n{Y_{s^{-1}t}}$$
for every $t\in T$ so
$$\si{t\in T}\n{(XY)_t}\leq \si{t\in T}\si{s\in T}\n{X_s}\n{Y_{s^{-1}t}}=\si{s\in T}\n{X_s}\left(\si{t\in T}\n{Y_{s^{-1}t}}\right)=$$
$$=\si{s\in T}\n{X_s}\left(\si{t\in T}\n{Y_t}\right)=\left(\si{t\in T}\n{X_t}\right)\left(\si{t\in T}\n{Y_t}\right)\;.$$

c) is easy to see.\qed

{\it{Remark.}} There may exist $X\in \ssb{\n{\cdot }}{f}$ for which  $((X_t\widetilde\otimes 1_K)V_t)_{t\in T}$ is not norm summable, as it is known from the theory of trigonometric series (see Proposition \ref{819}). In particular the inclusion $\ccc{A}\subset \ssb{\n{\cdot }}{f}$ may be strict.

\begin{co}\label{711}
Let $F$ be a unital C**-algebra and $\tau :E\rightarrow F$ a positive continuous (resp. W*-continuous) unital trace.
\begin{enumerate}
\item $\tau \circ \varphi _{1,1}$ is a positive continuous (resp. W*-continuous) unital trace.
\item If $\tau $ is faithful then $\tau \circ \varphi _{1,1}$ is faithful and $V_1$ is finite.
\item In the W*-case, $\ssb{W}{f}$ is finite iff $E$ is finite.
\end{enumerate}
\end{co}

a) Let $X,Y\in \ssa{f}$. By \h \ref{745} g) (and \pr \ref{704} a)),
$$\tau \varphi _{1,1}(XY)=\tau \left(\siw{t\in T}f(t,t^{-1})X_tY_{t^{-1}}\right)=\tau \left(\siw{t\in T}f(t,t^{-1})X_{t^{-1}}Y_t\right)=$$
$$=\siw{t\in T}\tau (f(t,t^{-1})X_{t^{-1}}Y_t)=\siw{t\in T}\tau (f(t,t^{-1})Y_tX_{t^{-1}})=\tau \left(\siw{t\in T}f(t,t^{-1})Y_tX_{t^{-1}}\right)=$$
$$=\tau \varphi _{1,1}(YX)\;.$$
Thus $\tau \circ \varphi _{1,1}$ is a trace which is obviously positive, continuous (resp. W*-continuous), and unital (\pr \ref{681} c),d)).

b) By \h \ref{745} g), $\varphi _{1,1}$ is faithful, so $\tau \circ \varphi $ is also faithful. Let $X\in \ssa{f}$ with $X^*X=V_1$. By a),
$$\tau \varphi _{1,1}(XX^*)=\tau \circ \varphi _{1,1}(X^*X)=\tau \varphi _{1,1}V_1=1_F$$
so
$$\tau \varphi _{1,1}(V_1-XX^*)=1_F-1_F=0\,,\qquad V_1=XX^*\,,$$
and $V_1$ is finite.

c) By b), if $E$ is finite then $\ssb{W}{f}$ is also finite. The reverse implication follows from the fact that $E\bar \otimes 1_K$ is a unital W*-subalgebra of $\ssb{W}{f}$ (\h \ref{745} h)).\qed 

\begin{co}\label{885}
Assume $T$ finite and for every $x'\in (E')^T$ put
$$\mae{\widetilde{x'}}{\ssa{f}}{\bk}{X}{\si{t\in T}}\sa{X_t}{x'_t}\;.$$\begin{enumerate}
\item $\widetilde{x'}\in \ssa{f}'$ and
$$\sup_{t\in T}\n{x'_t}\leq \n{\widetilde{x'}}\leq \si{t\in T}\n{x'_t}$$
for every $x'\in (E')^T$ and the map
$$\mae{\varphi }{(E')^T}{\ssa{f}'}{x'}{\widetilde{x'}}$$
is an isomorphism of involutive vector spaces such that
$$\varphi (xx')=(x\otimes 1_K)(\varphi x')\,,\qquad \varphi (x'x)=(\varphi x')(x\otimes 1_K)$$
(\emph{[C1] \pr 2.2.7.2}) for every $x\in E$ and $x'\in (E')^T$.
\item If $E$ is a W*-algebra then the map
$$\mae{\psi }{(\ddot E)^T}{\overbrace{\ssa{f}}^{..}}{(a_t)_{t\in T}}{(\tilde a_t)_{t\in T}}$$
is an isomorphism of involutive vector spaces such that
$$\psi (xa)=(x\otimes 1_K)(\psi a)\,,\qquad \psi (ax)=(\psi a)(x\otimes 1_K)$$
for every $x\in E$ and $a\in (\ddot E)^T$.\qed
\end{enumerate}
\end{co}

\begin{co}\label{867}
Assume $T$ finite and let $M$ be a \ri{\ssa{f}}. $M$ endowed with the right multiplication
$$\mad{M\times E}{M}{(\xi ,x)}{\xi (x\tilde \otimes 1_K)}$$
and with the inner-product
$$\mad{M\times M}{E}{(\xi ,\eta )}{\s{\xi }{\eta }_1}$$
is a \ri{E} denoted by $\widetilde M$, $\lb{\ssa{f}}{M}$ is a unital C*-subalgebra of $\lb{E}{\widetilde M}$, and $M$ is selfdual if $\widetilde M$ is so.
\end{co}

By \pr \ref{681} d),g) and \h \ref{745} g),l), for $X,Y\in \ssa{f}$ and $x\in E$,
$$\varphi _{1,1}(X(x\tilde \otimes 1_K))=(\varphi _{1,1}X)x\,,\qquad X\geq 0\Longrightarrow \varphi _{1,1}X\geq 0\,,$$
$$(X,Y)\in \stackrel{\circ }{\overbrace{\ssa{f}}}_+\Longrightarrow (\varphi _{1,1}X,\varphi _{1,1}Y)\in \stackrel{\circ }{E}_+\,,$$
$$\inf\me{\n{\varphi _{1,1}X}}{X\in \ssa{f}_+\,,\n{X}=1}>0$$
and the assertion follows from \pr \ref{681} a),c),d) and \prr 2.5 a),c),d).\qed

\begin{co}\label{9}
Let $n\in \bn$ and let $\varphi :\ssa{f}\rightarrow E_{n,n}$ be an $E$-C*-homomorphism. Then $(\varphi V_t)_{i,j}\in E^c$ for all $t\in T$ and all $i,j\in \bnn{n}$.
\end{co}

For $x\in E$,  by \pr \ref{674} d) and \h \ref{745} h),
$$x(\varphi V_t)=\varphi (x\widetilde \otimes 1_K)(\varphi V_t)=\varphi ((x\widetilde \otimes 1_K)V_t)=$$
$$=\varphi (V_t(x\widetilde\otimes 1_K))=(\varphi V_t)\varphi (x\widetilde\otimes 1_K)=(\varphi V_t)x$$
so $(\varphi V_t)_{i,j}\in E^c$.\qed

\begin{co}\label{5}
Let $S$ be a group and $g\in \ccc{F}(S,\ccc{S}(f))$. If we put
$$\mae{h}{(T\times S)\times (T\times S)}{\un{\ccc{S}(f)}}{((t_1,s_1),(t_2,s_2))}$$
$${(f(t_1,t_2)\widetilde\otimes 1_K)g(s_1,s_2)}$$
then $h\in \ccc{F}(T\times S,\ccc{S}(f))$.
\end{co}

The assertion follows from \h \ref{745} h).\qed

\begin{co}\label{40}
Let $X\in \ssa{f}\; (resp. \;X\in \ssb{\n{\cdot }}{f})$. 
\begin{enumerate}
\item For every $S\subset T$,
$$\sii{s\in S}{\fr{T}_3}(X_s\widetilde\otimes 1_K)V_s\in \ssa{f}\qquad (\emph{resp.}\;\sii{s\in S}{\n{\cdot }}(X_s\widetilde\otimes 1_K)V_s\in \ssb{\n{\cdot }}{f})$$
and 
$$\gamma :=\sup\me{\n{\si{t\in S}(X_t\widetilde\otimes 1_K)V_t}}{S\in \fr{P}_f(T)}<\infty \;.$$
\item We put for every $\alpha \in l^\infty (T)$
$$\mae{\alpha X}{T}{E}{t}{\alpha _tX_t}\;.$$
Then $\alpha X\in \ssa{f} \;(\mbox{resp.} \;\alpha X\in \ssb{\n{\cdot }}{f})$ for every $\alpha \in l^\infty (T)$ and the map
$$\mad{l^\infty (T)}{\ssa{f}\; (\emph{resp.}\; \ssb{\n{\cdot }}{f})}{\alpha }{\alpha X}$$
is norm-continuous.
\item Assume $E$ is a W*-algebra and let $l^\infty (T,E)$ be the C*-direct product of the family $(E)_{t\in T}$, which is a W*-algebra \emph{([C1] \pr 4.4.4.21 a))}. We put for every $\alpha \in l^\infty (T,E)$,
$$\mae{\alpha X}{T}{E}{t}{\alpha _tX_t}\;.$$
Then $\alpha X\in \ssb{W}{f}$ for every $\alpha \in l^\infty (T,E)$ and the map
$$\mad{l^\infty (T,E)}{\ssb{W}{f}}{\alpha }{\alpha X}$$
is continuous and W*-continuous.
\end{enumerate}
\end{co}

a) In the C*-case the family $((X_s\otimes 1_K)V_s)_{s\in S}$ is summable since $\ssb{C}{f}_{\fr{T}_3}$ is complete. By Banach-Steinhaus Theorem, $\gamma $ is finite. In the W*-case the summability follows now from Corollary \ref{742} b),c) and \h \ref{745} b).  

b) Let $G$ be the vector subspace $\me{\alpha \in l^\infty (T)}{\alpha (T)\;\mbox{is finite}}$ of $l^\infty (T)$. By a), the map
$$\mad{G}{\ssa{f} \;(\mbox{resp.}\; \ssb{\n{\cdot }}{f})}{\alpha }{\alpha X}$$
is well-defined, linear, and continuous. The assertion follows by continuity.

c) Let $x\in E$, $S\subset T$, and $\alpha :=xe_S$. For $\xi ,\eta \in H$ and $a\in \ddot E$, by a) and Lemma \ref{18} b) (and \h \ref{745} b)),
$$\sa{\alpha X}{\ti{a}{\xi }{\eta }}=\sa{\s{\alpha X\xi }{\eta }}{a}=\sa{\sii{t\in T}{\ddot E}\eta _t^*x((e_SX)\xi )_t}{a}=$$
$$=\si{t\in T}\sa{x}{((e_SX)\xi )_ta\eta _t^*}=\sa{x}{\sii{t\in T}{E}((e_SX)\xi )_ta\eta _t^*}\;.$$
Let $G$ be the involutive subalgebra 
$\me{\alpha \in l^\infty (T,E)}{\alpha (T)\;\mbox{is finite}}$
 of $l^\infty (T,E)$ and let $\bar G$ be its norm-closure in $l^\infty (T,E)$, which is a C*-subalgebra of $l^\infty (T,E)$. By [C1] \pr 4.4.4.21 a), $G$ is dense in $l^\infty (T,E)_{\ddot F}$, where $F:=l^\infty (T,E)$. 

Let $\alpha \in l^\infty (T,E)^{\#}$ and let $\fr{F}$ be a filter on $G^{\#}$ converging to $\alpha $ in $l^\infty (T,E)_{\ddot F}$ ([C1] Corollary 6.3.8.7). By the above (and by \h \ref{745} h)),
$$\lim_{\beta ,\fr{F}}\beta X=\alpha X$$
in $\ssb{W}{f}_{\overbrace{\ssb{W}{f}}^{..}}$ and so $\alpha X\in \ssb{W}{f}$. The assertion follows.\qed

\begin{co}\label{776}
Let $S$ be a subgroup of $T$. Put
$$f_S:=f|(S\times S)\,,\quad K_S:=l^2(S)\,,\quad \ccc{G}:=\me{X\in \ssa{f}}{t\in T\setminus S\Longrightarrow X_t=0}\;.$$
\begin{enumerate}
\item $f_S\in \f{S}{E}$.
\item $\ccc{G}$ is an $E$-C**-subalgebra of $\ssa{f}$.
\item For every $X\in \ccc{G}$, the family $((X_s\widetilde\otimes 1_{K_S})V_s^{f_S})_{s\in S}$ is summable in $\lb{E}{K_S}_{\fr{T}_3}$ and the map
$$\mae{\varphi }{\ccc{G}}{\ssa{f_S}}{X}{\sii{s\in S}{\fr{T}_3}(X_s\widetilde\otimes 1_{K_S})V_s^{f_S}}$$
is an injective $E$-C**-homomorphism.
\item If $X\in \ccc{G}\cap \ssb{\n{\cdot }}{f}$ then $\varphi X\in \ssb{\n{\cdot }}{f_S}$ and the map 
$$\mad{\ccc{G}\cap \ssb{\n{\cdot }}{f}}{\ssb{\n{\cdot }}{f_S}}{X}{\varphi X}$$
is an $E$-C*-isomorphism.
\item If $S$ is finite then the map
$$\mad{\ccc{G}}{\ssa{f_S}}{X}{\si{t\in S}(X_t\otimes 1_{K_S})}V_t^{f_S}$$
is an $E$-C*-isomorphism. 
\end{enumerate}
\end{co}

a) is obvious.

b) By \h \ref{745} c),g), $\ccc{G}$ is an involutive unital subalgebra of $\ssa{f}$ and by \pr \ref{681} a) (resp. \pr \ref{681} c) and Corollary \ref{742} c)) and \h \ref{745} h), it is an $E$-C**-subalgebra of $\ssa{f}$.

c) follows from \h \ref{745} b) and Corollary \ref{40} a).

d) follows from c).

e) is contained in d).\qed

\begin{de}\label{779}
We denote by $\fr{S}_T$ the set of finite subgroups of $T$ and call $T$ {\bf locally finite} if $\fr{S}_T$ is upward directed and
$$\bigcup _{S\in {\fr{S}_T}}S=T\;.$$
\end{de}

T is locally finite iff the subgroups of $T$ generated by finite subsets of $T$ are finite.

\begin{co}\label{780}
Assume $T$ locally finite. We put $f_S:=f|(S\times S)$ for every $S\in \fr{S}_T$ and identify $\ssa{f_S}$ with $\me{X\in \ssa{f}}{t\in T\setminus S\Rightarrow X_t=0}$ \emph{(Corollary \ref{776} e))}.
\begin{enumerate}
\item For every $X\in \ssb{\n{\cdot }}{f}$ and $\varepsilon >0$ there is an $S\in \fr{S}_T$ such that
$$\n{\si{t\in R}(X_t\otimes 1_K)V_t-X}<\varepsilon $$
for every $R\in \fr{S}_T$ with $S\subset R$. 
\item $\ssb{\n{\cdot }}{f}$ is the norm closure of $\cup_{s\in \fr{S}_T}\ssa{f_S}$ and so it is canonically isomorphic to the inductive limit of the inductive system $\me{\ssa{f_S}}{S\in \fr{S}_T}$ and for every $S\in \fr{S}_T$ the inclusion map $\ssa{f_S}\rightarrow \ssb{\n{\cdot }}{f}$ is the associated canonical morphism.
\end{enumerate}
\end{co}

a) There is a $Y\in \ccc{R}(f)$ with $\n{X-Y}<\frac{\varepsilon}{2} $. Let $S\in \fr{S}_T$ with $Y\in \ssa{f_S}$. By Corollary \ref{776} b), for $R\in \fr{S}_T$ with $S\subset R$,
$$\n{\si{t\in R}((X_t-Y_t)\widetilde\otimes 1_K)V_t}\leq \n{X-Y}<\frac{\varepsilon }{2}$$
so
$$\n{\si{t\in R}(X_t\widetilde\otimes 1_K)V_t-X}\leq \n{\si{t\in R}((X_t-Y_t)\widetilde\otimes 1_K)V_t}+\n{Y-X}<\frac{\varepsilon }{2}+\frac{\varepsilon }{2}=\varepsilon \;.$$

b) follows from a).\qed

{\it Remark.} The C*-algebras of the form $\ssb{\n{\cdot }}{f}$ with $T$ locally finite can be seen as a kind of AF-$E$-C*-algebras.

\begin{p}\label{25}
The following are equivalent for all $t\in T$ with $t^2=1$ and $\alpha \in \unn{E}$.
\begin{enumerate}
\item $\frac{1}{2}(V_1+(\alpha \widetilde\otimes 1_K)V_t)\in Pr\,\ssa{f}$.
\item $\alpha ^2=\tilde f(t)$.
\end{enumerate}
\end{p}

By \pr \ref{674} b),d),e),
$$(V_t)^*=(\tilde f(t)\widetilde\otimes 1_K)V_t\,,\qquad (V_t)^2=(\tilde f(t)^*\widetilde\otimes 1_K)V_1$$
so
$$\frac{1}{2}(V_1+(\alpha \widetilde\otimes 1_K)V_t)^*=\frac{1}{2}(V_1+((\alpha ^*\tilde f(t))\widetilde\otimes 1_K)V_t)\,,$$
$$\left(\frac{1}{2}(V_1+(\alpha \widetilde\otimes 1_K)V_t)\right)^2=\frac{1}{4}((1_E+\alpha ^2\tilde f(t)^*)\widetilde\otimes 1_K)V_1+\frac{1}{2}(\alpha \widetilde\otimes 1_K)V_t\;.$$
Thus a) is equivalent to $\alpha ^*\tilde f(t)=\alpha $ and $\alpha ^2\tilde f(t)^*=1_E$, which is equivalent to b).\qed 

\begin{co}\label{869}
Let $t\in T$ such that $t^2=1$ and $\tilde f(t)=1_E$. Then
$$\frac{1}{2}(V_1\pm V_t)\in Pr\;\ssa{f}\,,\qquad (V_1+V_t)(V_1-V_t)=0\;.$$
\end{co}

The assertion follows from \pr \ref{25}.\qed

\begin{co}\label{26}
Let $\alpha ,\beta \in \unn{E}$, $s,t\in T$ with $s^2=t^2=1$, $st=ts$,
$$\gamma :=\frac{1}{2}(\alpha ^*\beta f(s,st)^*+ \beta ^*\alpha f(t,st)^*)\,,\quad \gamma ':=\frac{1}{2}(\alpha \beta ^*f(st,t)^*+\beta \alpha ^*f(st,s)^*)\,,$$
and
$$X:=\frac{1}{2}((\alpha \widetilde\otimes 1_K)V_s+(\beta \widetilde\otimes 1_K)V_t)\;.$$
\begin{enumerate}
\item $f(s,st)f(t,st)=f(st,t)f(st,s)=\tilde f(st)^*$.
\item $f(st,t)f(s,st)=f(st,s)f(t,st)$.
\item $X^*X=\frac{1}{2}(V_1+(\gamma \widetilde\otimes 1_K)V_{st})\,,\quad XX^*=\frac{1}{2}(V_1+(\gamma '\widetilde\otimes 1_K)V_{st})$.
\item The following are equivalent.
\begin{enumerate}
\item $X^*X\in Pr\,\ssa{f}$.
\item $XX^*\in Pr\,\ssa{f}$.
\item $\alpha ^*\beta f(t,st)=\beta ^*\alpha f(s,st)$.
\item $\alpha ^*\beta f(st,t)=\beta ^*\alpha f(st,s)$.
\end{enumerate}
\end{enumerate}
\end{co}

a) and b) follow from the equation of Schur functions (\dd \ref{703}) and \pr \ref{704} a).

c) By \pr \ref{674} b),e) and \pr \ref{704} b),
$$X^*=\frac{1}{2}(((\alpha ^*\tilde f(s))\widetilde\otimes 1_K)V_s+((\beta ^*\tilde f(t))\widetilde\otimes 1_K)V_t)\,,$$
$$X^*X=\frac{1}{2}V_1+\frac{1}{4}((\alpha ^*\beta \tilde f(s)f(s,t)+\beta ^*\alpha \tilde f(t)f(t,s))\widetilde\otimes 1_K)V_{st}=$$
$$=\frac{1}{2}V_1+\frac{1}{4}((\alpha ^*\beta f(s,st)^*+\beta ^*\alpha f(t,st)^*)\widetilde\otimes 1_K)V_{st}=\frac{1}{2}(V_1+(\gamma \widetilde\otimes 1_K)V_{st})\,,$$
$$XX^*=\frac{1}{2}V_1+\frac{1}{4}((\alpha \beta ^*\tilde f(t)f(s,t)+\beta \alpha ^*\tilde f(s)f(t,s))\widetilde\otimes 1_K)V_{st}=$$
$$=\frac{1}{2}V_1+\frac{1}{4}((\alpha \beta ^*f(st,t)^*+\beta \alpha ^*f(st,s)^*)\widetilde\otimes 1_K)V_{st}=\frac{1}{2}(V_1+(\gamma '\widetilde\otimes 1_K)V_{st})\;.$$

$d_1\Leftrightarrow  d_2$ is known.

$d_1\Leftrightarrow d_3$. By a),
$$\gamma ^2-\tilde f(st)=
\frac{1}{4}(\alpha ^*\beta \alpha ^*\beta f(s,st)^{*2}+\beta ^*\alpha \beta ^*\alpha f(t,st)^{*2}+2f(s,st)^*f(t,st)^*)-$$
$$-f(s,st)^*f(t,st)^*
=\frac{1}{4}(\alpha ^*\beta f(s,st)^*-\beta ^*\alpha f(t,st)^*)^2\;.$$
By \pr \ref{25} $d_1$) is equivalent to $\gamma ^2=\tilde f(st)$ so, by the above, since $\alpha ^*\beta f(s,st)^*-\beta ^*\alpha f(t,st)^*$ is normal, it is equivalent to
$$\alpha ^*\beta f(s,st)^*=\beta ^*\alpha f(t,st)^*\quad\mbox{or to}\quad \beta ^*\alpha f(s,st)=\alpha ^*\beta f(t,st)\;.$$

$d_3\Leftrightarrow d_4$ follows from b).\qed  

\begin{p}\label{750}
Let $X\in \ssa{f}$.
\begin{enumerate}
\item $\siw{t\in T}X_t^*X_t=(X^*X)_1\,,\quad \siw{t\in T}(X_tX_t^*)=(XX^*)_1$.
\item $(X_t)_{t\in T}\,,(X_t^*)_{t\in T}\in \widetilde{\cb{t\in T}}\breve E$,
$$\n{(X_t)_{t\in T}}\leq \n{X}\,,\quad \n{(X_t^*)_{t\in T}}\leq \n{X}\;.$$
\item If $T$ is finite and $f$ is constant then there is an $X\in \ssa{f}$ with
$$\n{X}\geq \sqrt{Card\;T}\n{(X_t)_{t\in T}}\,,\quad \n{X}\geq \sqrt{Card \;T}\n{(X_t^*)_{t\in T}}\;.$$
\item If $T$ is infinite and locally finite and $f$ is constant then the map
$$\mad{\ssa{f}}{\widetilde{\cb{t\in T}}\breve E}{X}{(X_t)_{t\in T}}$$
is not surjective.
\end{enumerate}
\end{p}

a) follows from  \h \ref{745} g).

b) By a),
$$(X_t)_{t\in T}\,,(X_t^*)_{t\in T}\in \widetilde{\cb{t\in T}}\breve E$$
and by \pr \ref{681} a),
$$\n{(X_t)_{t\in T}}^2=\n{\varphi _{1,1}(X^*X)}\leq \n{X^*X}=\n{X}^2\,,$$
$$\n{(X_t^*)_{t\in T}}^2=\n{\varphi _{1,1}(XX^*)}\leq \n{XX^*}=\n{X}^2\;.$$

c) Let $n:=Card \;T$ and for every $t\in T$ put $X_t:=1_E$, $\xi _t:=1_E$. Then
$$\n{(X_t)_{t\in T}}^2=\n{(X_t^*)_{t\in T}}^2=n\,,\quad\n{(\xi _t)_{t\in T}}^2=n\;.$$
For $t\in T$, by \h \ref{745} e),
$$(X\xi )_t=\si{s\in T}f(s,s^{-1}t)X_s\xi _{s^{-1}t}=n1_E$$
so
$$\s{X\xi }{X\xi }=n^31_E\,,\quad n\n{X}^2=\n{X}^2\n{\xi }^2\geq \n{X\xi }^2=n^3\,,$$
$$\n{X}^2\geq n\n{(X_t)_{t\in T}}^2\,,\quad \n{X}\geq \sqrt{n}\n{(X_t)_{t\in T}}\;.$$

d) follows from c), \h \ref{745} a), and the Principle of Inverse Operator.\qed

{\it Remark.} If $E$ is a W*-algebra then it may exist a family $(x_t)_{t\in T}$ in $E$ such that the family $((x_t\widetilde\otimes 1_K)V_t)_{t\in T}$ is summable in $\lb{E}{H}_{\fr{T}_2}$ in the W*-case but not in the C*-case as the following example shows. Take $T:=\bz$, $f$ constant, $E:=l^\infty (\bz)$, and $x_t:=(\delta _{t,s})_{s\in T}\in E$ for every $t\in T$. By \pr \ref{750} b), $((x_t\otimes 1_K)V_t)_{t\in T}$ is not summable in $\lb{E}{H}_{\fr{T}_2}$ in the C*-case. In the W*-case for $\xi \in H$ and $s,t\in T$, 
$$\s{((x_t\bar \otimes 1_K)V_t\xi )_s}{((x_t\bar \otimes 1_K)V_t\xi )_s}=e_t|\xi _{s-t}|^2\,,$$
$$\s{(x_t\bar \otimes 1_K)V_t\xi }{(x_t\bar \otimes 1_K)V_t\xi }=e_t\n{\xi }^2\;.$$
Thus
$$X:=\sii{t\in T}{\fr{T}_2}(x_t\bar  \otimes 1_K)V_t\in \ssb{W}{f}\;.$$
Using the identification of \h \ref{745} i), we get $X\in \ssb{W}{f}\setminus \ssb{C}{f}$.

\begin{co}\label{752}
Let $X\in \ssa{f}$.
\begin{enumerate}
\item $X\in \me{x\widetilde\otimes 1_K}{x\in E}^c$ iff $X_t\in E^c$ for all $t\in T$.
\item $X\in \me{V_t}{t\in T}^c$ iff
$$X_{s^{-1}ts}=f(s,s^{-1}ts)^*f(t,s)X_t=f(s^{-1},ts)f(t,s)\tilde f(s)X_t$$
for all $s,t\in T$.
\item $X\in \ssa{f}^c$ iff for all $s,t\in T$
$$X_t\in E^c\,,\quad X_{s^{-1}ts}=f(s,s^{-1}ts)^*f(t,s)X_t=f(s^{-1},ts)f(t,s)\tilde f(s)X_t\;.$$
In particular if $f(s,t)=f(t,s)$ for all $s,t\in T$ then $X\in \ssa{f}^c$ iff $X_t\in E^c$ for all $t\in T$.
\item $\varphi _{1,1}(\ssa{f}^c)=E^c$.
\item If the conjugacy class of $t\in T$ (i.e. the set $\me{s^{-1}ts}{s\in T}$) is infinite and $X\in \me{V_t}{t\in T}^c$ then $X_t=0$.
\item If the conjugacy class of every $t\in T\setminus \{1\}$ is infinite then
$$\me{V_t}{t\in T}^c=\me{x\widetilde\otimes 1_K}{x\in E}\,,\quad \ssa{f}^c=\me{x\widetilde\otimes 1_K}{x\in E^c}\;.$$
Thus in this case $\ssa{f}$ is a kind of $E$-factor.
\item The following are equivalent:
\begin{enumerate}
\item $\ssa{f}$ is commutative.
\item $T$ and $E$ are commutative and $f(s,t)=f(t,s)$ for all $s,t\in T$.
\end{enumerate}
\end{enumerate}
\end{co}

For $s,t\in T\,,\,x\in E$, and $Y:=(x\widetilde\otimes 1_K)V_s$, by \h \ref{745} g),
$$(XY)_t=\siw{r\in T}f(r,r^{-1}t)X_rY_{r^{-1}t}=\siw{r\in T}f(r,r^{-1}t)X_r\delta _{s,r^{-1}t}x=f(ts^{-1},s)X_{ts^{-1}}x\,,$$
$$(YX)_t=\siw{r\in T}f(r,r^{-1}t)Y_rX_{r^{-1}t}=\siw{r\in T}f(r,r^{-1}t)\delta _{r,s}xX_{r^{-1}t}=f(s,s^{-1}t)xX_{s^{-1}t}\;.$$

a) follows from the above by putting $s:=1$ (\pr \ref{704} a)).

b) follows from the above by putting $x:=1_E$ and $t:=rs$ (\pr \ref{704}).

c) follows from a),b), and Corollary \ref{742} d). The last assertion follows using \pr \ref{6422} a).

d) follows from c) (and \pr \ref{704} a)).

e) follows from b) and \pr \ref{750} b).

f) follows from c), e), and \pr \ref{674} d).

$g_1\Rightarrow g_2$. By a), $E$ is commutative. By \pr \ref{674} b),
$$f(s,t)V_{st}=V_sV_t=V_tV_s=f(t,s)V_tV_s=f(t,s)V_{ts}$$
and so by \h \ref{745} a), $st=ts$ and $f(s,t)=f(t,s)$.

$g_2\Rightarrow g_1$ follows from c).\qed

\begin{co}\label{713}
If $\bk=\br$ then the following are equivalent:
\begin{enumerate}
\item $\ssa{f}^c=\ssa{f}=Re\; \ssa{f}$.
\item $T$ is commutative, $E^c=E=Re\;E$, and 
$$f(s,t)=f(t,s)\,,\quad \tilde f(t)=1_E\,,\quad t^2=1$$
for all $s,t\in T$.
\end{enumerate}
\end{co}

$a\Rightarrow b$. By Corollary \ref{752} $g_1\Rightarrow g_2$, $T$ is commutative, $E=E^c$, and $f(s,t)=f(t,s)$ for all $s,t\in T$. Since $E$ is isomorphic with a C*-subalgebra of $\ssa{f}$ (\h \ref{745} h)), $E=Re\;E$. By \pr \ref{674} e),
$$V_t=V_t^*=(\tilde f(t)\widetilde\otimes 1_K)V_{t^{-1}}$$
so by \h \ref{745} a), $t=t^{-1}\,,\,\tilde f(t)=1_E$, so $t^2=1$.

$b\Rightarrow a$. By Corollary \ref{752} $g_2\Rightarrow g_1$, $\ssa{f}^c=\ssa{f}$. For $X\in \ssa{f}$ and $t\in T$, by \h \ref{745} c),
$$(X^*)_t=\tilde f(t)(X_{t^{-1}})^*=(X_t)^*=X_t$$
so $X^*=X$ (\h \ref{745} a)).\qed

\begin{p}\label{14}
Let $(E_i)_{i\in I}$ be a family of unital C**-algebras such that $E$ is the C*-direct product of this family. For every $i\in I$, we identify $E_i$ with the corresponding closed ideal of $E$ (resp. of $E_{\ddot E}$) and put
$$\mae{f_i}{T\times T}{\un{E_i}}{(s,t)}{f(s,t)_i}\;.$$
\begin{enumerate}
\item For every $i\in I$, $f_i\in \f{T}{E_i}$. We put \emph{(by \h \ref{745} b))}
$$\mae{\varphi _i}{\ssa{f}}{\ssa{f_i}}{X}{\sii{t\in T}{\fr{T}_2}((X_t)_i\widetilde\otimes 1_K)V_t^{f_i}}\;.$$
$\varphi _i$ is a surjective C**-homomorphism.
\item In the C*-case, if $T$ is finite then $\ccc{R}(f)=\ssb{\n{\cdot }}{f}=\ssb{C}{f}$ is isomorphic to the C*-direct product of the family 
$$(\ccc{R}(f_i)=\ssb{\n{\cdot }}{f_i}=\ssb{C}{f_i})_{i\in I}\;.$$
\item In the C*-case, if $I$ is finite then $\ssb{C}{f}$ (resp. $\ssb{\n{\cdot }}{f}$) is isomorphic to $\pro{i\in I}\ssb{C}{f_i}$ (resp. $\pro{i\in I}\ssb{\n{\cdot }}{f_i}$).  
\item In the W*-case, $\ssb{W}{f}$ is isomorphic to the C*-direct product of the family $(\ssb{W}{f_i})_{i\in I}$.\qed
\end{enumerate}  
\end{p}

{\it Remark.} The C*-isomorphisms of b) and c) cease to be surjective in general if $T$ and $I$ are both infinite. Take $T:=(\bzz{2})^{\bn}$, $I:=\bn$, $E_i:=\bk$ for every $i\in I$, and $E:=l^\infty $ (i.e. $E$ is the C*-direct product of the family $(E_i)_{i\in I}$). For every $n\in \bn$ put $t_n:=(\delta _{m,n})_{m\in \bn}\in T$. 
Assume there is an $X\in \ssb{C}{f}$ (resp. $X\in \ssb{\n{\cdot }}{f}$) with $\psi X=(V_{t_i}^{f_i})_{i\in I}$ (resp. $\varphi X=(V_{t_i}^{f_i})_{i\in I}$), where $\psi $ and $\varphi $ are the maps of b) and c), respectively. Then $(X_{t_n})_i=\delta _{i,n}$ for all $i,n\in \bn$ and this implies $(X_t)_{t\in T}\not\in \cb{t\in T}\breve E$, which contradicts \pr \ref{750} b).  

\begin{p}\label{868}
Let $S$ be a finite group, $K':=l^2(S)$, $K'':=l^2(S\times T)$, and $g\in \f{S}{\ssa{f}}$ such that $g(s_1,s_2)\in \un{E}$ (where $\un{E}$ is identified with $(\un{E})\widetilde \otimes 1_{K}\subset \un{\ssa{f}}$) for all $s_1,s_2\in S$ and put 
$$\mae{h}{(S\times T)\times (S\times T)}{\un{E}}{((s_1,t_1),(s_2,t_2))}{g(s_1,s_2)f(t_1,t_2)}\;.$$
\begin{enumerate}
\item $h\in \f{S\times T}{E}$; for every $X\in \ssa{g}$ put
$$\varphi X:=\si{s\in S}\sii{t\in T}{\fr{T}_3}((X_s)_t\widetilde \otimes 1_{K''})V_{(s,t)}^h\in \ssa{h}\;.$$
\item $\mac{\varphi }{\ssa{g}}{\ssa{h}}$ is an $E$-C*-isomorphism.
\end{enumerate}
\end{p}

a) is obvious.

b) For $X,Y\in \ssa{g}$ and $(s,t)\in S\times T$, by \h \ref{745} c),g) and \pr \ref{681} g),
$$(\varphi X^*)_{(s,t)}=((X^*)_s)_t=\tilde g(s)((X_{s^{-1}})^*)_t=$$
$$=\tilde g(s)\tilde f(t)((X_{s^{-1}})_{t^{-1}})^*=\tilde h(s,t)(X_{(s,t)^{-1}})^*=((\varphi X)^*)_{(s,t)}\,,$$
$$(\varphi (XY))_{(s,t)}=((XY)_s)_t=\si{r\in S}g(r,r^{-1}s)(X_rY_{r^{-1}s})_t=$$
$$=\si{r\in S}g(r,r^{-1}s)\siw{q\in T}f(q,q^{-1}t)(X_r)_q(Y_{r^{-1}s})_{q^{-1}t}=$$
$$=\siw{(r,q)\in S\times T}h((r,q),(r,q)^{-1}(s,t))X_{(r,q)}Y_{(r,q)^{-1}(s,t)}=((\varphi X)(\varphi Y))_{(s,t)}\,,$$
so $\varphi $ is a C*-homomorphism. If $\varphi X=0$ then $X_{(s,t)}=0$ for all $(s,t)\in S\times T$, so $X=0$ and $\varphi $ is injective. Let $Z\in \ssa{h}$. For every $s\in S$ put
$$X_s:=\sii{t\in T}{\fr{T}_3}(Z_{(s,t)}\tilde \otimes 1_K)V_t^f\in \ssa{f}\,,$$
$$X:=\si{s\in S}(X_s\otimes 1_{K'})V_s^g\in \ssa{g}\;.$$
Then $\varphi X=Z$ and $\varphi $ is surjective.\qed 

\begin{p}\label{714}
If $T$ is infinite and $X\in \ssa{f}\setminus \{0\}$ then $X(H^{\#})$ is not precompact. 
\end{p}

Let $t\in T$ with $X_t\not=0$. There is an $x'\in E'_+$ (resp. $x'\in \ddot E_+$) with $\sa{X_t^*X_t}{x'}>0$. We put $t_1:=1$ and construct a sequence $(t_n)_{n\in \bn}$ recursively in $T$ such that for all $m,n\in \bn,\,m<n$,
$$\left|\sa{f(t,t_m)^*f(tt_mt_n^{-1},t_n)X_t^*X_{tt_mt_n^{-1}}}{x'}\right|<\frac{1}{2}\sa{X_t^*X_t}{x'}\;.$$ 
Let $n\in \bn\setminus \{1\}$ and assume the sequence was constructed up to $n-1$. Since (\pr \ref{750} a))
$$\si{s\in T}\sa{X_{tt_ms^{-1}}^*X_{tt_ms^{-1}}}{x'}<\infty $$
for all $m\in \bnn{n-1}$ there is a $t_n\in T$ with 
$$\sa{X_{tt_mt_n^{-1}}^*X_{tt_mt_n^{-1}}}{x'}<\frac{1}{4}\sa{X_t^*X_t}{x'}$$
for all $m\in \bnn{n-1}$. By Schwarz' inequality ([C1] \pr 2.3.4.6 c)) for $m\in \bnn{n-1}$,
$$\left|\sa{f(t,t_m)^*f(tt_mt_n^{-1},t_n)X_t^*X_{tt_mt_n^{-1}}}{x'}\right|^2\leq$$
$$\leq  \sa{X_t^*X_t}{x'}\sa{X_{tt_mt_n^{-1}}^*X_{tt_mt_n^{-1}}}{x'}<\frac{1}{4}\sa{X_t^*X_t}{x'}^2\;.$$
This finishes the recursive construction.

For $r,s\in T$, by \h \ref{745} e),
$$(X(1_E\otimes e_r))_s=\siw{q\in T}f(q,q^{-1}s)X_q\delta _{r,q^{-1}s}=f(sr^{-1},r)X_{sr^{-1}}\,,$$
$$\s{X(1_E\otimes e_r)}{X_t\otimes e_s}=f(sr^{-1},r)X_t^*X_{sr^{-1}}\;.$$
For $m,n\in \bn,\,m<n$, it follows 
$$\s{X(1_E\otimes e_{t_m})}{X_t\otimes e_{tt_m}}=f(t,t_m)X_t^*X_t\,,$$
$$\sa{\s{X(1_E\otimes e_{t_m})}{X_t\otimes e_{tt_m}}}{x'f(t,t_m)^*}=\sa{X_t^*X_t}{x'}\,,$$
$$\s{X(1_E\otimes e_{t_n})}{X_t\otimes e_{tt_m}}=f(tt_mt_n^{-1},t_n)X_t^*X_{tt_mt_n^{-1}}\,,$$
$$|\sa{\s{X(1_E\otimes e_{t_n})}{X_t\otimes e_{tt_m}}}{x'f(t,t_m)^*}|=$$
$$\left|\sa{f(t,t_m)^*f(tt_mt_n^{-1},t_n)X_t^*X_{tt_mt_n^{-1}}}{x'}\right|<
\frac{1}{2}\sa{X_t^*X_t}{x'}\,,$$
$$\n{x'}\n{X(1_E\otimes e_{t_m})-X(1_E\otimes e_{t_n})}
\n{X_t}\geq$$
$$\geq  \left|\sa{\s{X(1_E\otimes e_{t_m})-X(1_E\otimes e_{t_n})}{X_t\otimes e_{tt_m}}}{x'f(t,t_m)^*}\right|\geq $$
$$\geq \left|\sa{\s{X(1_E\otimes e_{t_m})}{X_t\otimes e_{tt_m}}}{x'f(t,t_m)^*}\right|-$$
$$-\left|\sa{\s{X(1_E\otimes e_{t_n})}{X_t\otimes e_{tt_m}}}{x'f(t,t_m)^*}\right|>$$
$$>\sa{X_t^*X_t}{x'}-\frac{1}{2}\sa{X_t^*X_t}{x'}=\frac{1}{2}\sa{X_T^*X_T}{x'}\;.$$
Thus the sequence $(X(1_E\otimes e_{t_n}))_{n\in \bn}$ has no Cauchy subsequence and therefore $X(H^{\#})$ is not precompact.\qed

\begin{p}\label{aw}
Assume $T$ finite and let $\Omega $ be a compact space, $\omega _0\in \Omega $,
$$\mae{g}{T\times T}{\unn{\ccc{C}(\Omega ,E)}}{(s,t)}{f(s,t)1_\Omega }\,,$$
$$A:=\me{X\in \ssa{g}}{t\in T,t\not=1\Longrightarrow X_t(\omega _0)=0}\,,$$
$$B:=\me{Y\in \ccc{C}(\Omega ,\ssa{f})}{t\in T,t\not=1\Longrightarrow Y(\omega _0)_t=0 }\;.$$
Then $g\in \f{T}{\ccb{\Omega }{E}}$ and we define for every $X\in A$ and $Y\in B$,
$$\mae{\varphi X}{\Omega }{\ssa{f}}{\omega }{\si{t\in T}(X_t(\omega )\otimes 1_K)V_t^f}\,,$$
$$\psi Y:=\si{t\in T}(Y(\cdot )_t\otimes 1_K)V_t^g\;.$$
Then $A$ (resp. $B$) is a unital C*-subalgebra of $\ssa{g}$ (resp. of $\ccc{C}(\Omega ,\ssa{f})$)
$$\mac{\varphi }{A}{B}\,,\qquad \qquad \mac{\psi }{B}{A}$$
are C*-isomorphisms, and $\varphi =\psi ^{-1}$.
\end{p} 

It is easy to see that $A$ (resp. $B$) is a unital C*-subalgebra of $\ssa{g}$ (resp. of $\ccc{C}(\Omega ,\ssa{f})$) and that $\varphi $ and $\psi $ are well-defined. For $X,X'\in A$, $t\in T$, and $\omega \in \Omega $, by \h \ref{745} c),g) and \pr \ref{674} e),
$$(((\varphi X)(\varphi X'))(\omega ))_t=\si{s\in T}f(s,s^{-1}t)((\varphi X)(\omega ))_s((\varphi X')(\omega ))_{s^{-1}t}=$$
$$=\si{s\in T}f(s,s^{-1}t)X_s(\omega )X'_{s^{-1}t}(\omega )=
\si{s\in T}(f(s,s^{-1}t)X_sX'_{s^{-1}t})(\omega )=$$
$$=(XX')_t(\omega )=
(\varphi (XX')(\omega ))_t\,,$$
$$(\varphi X^*)(\omega )=\si{s\in T}(((X^*)_s(\omega ))\otimes 1_K)V_s^f=\si{s\in T}((\tilde{f}(s)((X_{s^{-1}})^*(\omega )))\otimes 1_K) V_s^f=$$
$$=\si{s\in T}((X_{s^{-1}})(\omega )^*\otimes 1_K)(V_{s^{-1}}^f)^*=\si{s\in T}(X_s(\omega )^*\otimes 1_K)(V_s^f)^*=(\varphi X)^*(\omega )$$
so $\varphi $ is a C*-homomorphism and we have
$$(\psi \varphi X)_t=(\varphi X)_t=X_t\;.$$
Moreover for $Y\in B$,
$$(\varphi \psi Y)_t(\omega )=((\psi Y)(\omega ))_t=Y_t(\omega )$$
which proves the assertion.\qed

\begin{center}
\subsection{Variation of the parameters}
\end{center}

In this subsection we examine the changes produced by the replacement of the groups and of the Schur functions.

\begin{de}\label{720}
We put for every $\lambda \in \Lambda (T,E)$ \emph{(\dd$\,$\ref{705})}
$$\mae{U_\lambda }{H}{H}{\xi }{(\lambda (t)\xi _t)_{t\in T}}\;.$$
\end{de}

It is easy to see that $U_\lambda $ is well-defined, $U_\lambda \in Un\;\lb{E}{H}$, and the map
$$\mad{\Lambda (T,E)}{Un\;\lb{E}{H}}{\lambda }{U_\lambda }$$
is an injective group homomorphism with $U_\lambda ^*=U_{\lambda ^*}$ (\pr \ref{706} c)). Moreover 
$$\n{U_\lambda -U_\mu }\leq \n{\lambda -\mu }_\infty$$ 
for all $\lambda ,\mu \in \Lambda (T,E)$.

\begin{p}\label{716}
Let $f,g\in \fte$ and $\lambda \in \Lambda (T,E)$.
\begin{enumerate}
\item The following are equivalent:
\begin{enumerate}
\item $g=f\delta \lambda $.
\item There is a (unique) $E$-C*-isomorphism 
$$\mac{\varphi }{\ssa{f}}{\ssa{g}}$$
continuous with respect to the $\fr{T}_2$-topologies such that for all $t\in T$ and $x\in E$,
$$\varphi V_t^f=(\lambda (t)^*\widetilde\otimes 1_K)V_t^g$$
(we call such an isomorphism an {\bf $\ccc{S}$-isomorphism} and denote it by $\approx _{\ccc{S}}$)
\end{enumerate}
\item If the above equivalent assertions are fulfilled then for $X\in \ssa{f}$ and $t\in T$,
$$\varphi X=U_\lambda ^*XU_\lambda \,,\quad (\varphi X)_t=\lambda (t)^*X_t\;.$$
\item There is a natural bijection
$$\me{\ssa{f}}{f\in \fte}/\approx _{\ccc{S}}\longrightarrow \fte/\me{\delta \lambda }{\lambda \in \Lambda (T,E)}\;.$$
\end{enumerate}
\end{p}

By \pr \ref{706} c), $\delta \lambda \in \fte$ for every $\lambda \in \Lambda (T,E)$.

$a_1\Rightarrow a_2\&b$. For $s,t\in T$ and $\zeta \in \breve E$, by \pr \ref{674} c),
$$U_\lambda ^*V_t^fU_\lambda (\zeta \otimes e_s)=U_\lambda ^*V_t^f((\lambda (s)\zeta )\otimes e_s)=U_\lambda ^*((f(t,s)\lambda (s)\zeta )\otimes e_{ts})=$$
$$=(\lambda (ts)^*f(t,s)\lambda (s)\zeta )\otimes e_{ts}=(\lambda (t)^*g(t,s)\zeta )\otimes e_{ts}=(\lambda (t)^*\widetilde\otimes 1_K)V_t^g(\zeta \otimes e_s)$$
so (by \pr \ref{674} e))
$$U_\lambda ^*V_t^fU_\lambda =(\lambda (t)^*\widetilde\otimes 1_K)V_t^g\;.$$
Thus the map
$$\mae{\varphi }{\ssa{f}}{\ssa{g}}{X}{U_\lambda ^*XU_\lambda }$$
is well-defined. It is obvious that it has the properties described in $a_2$). The uniqueness follows from \h \ref{745} b).

We have
$$\varphi ((X_t\widetilde\otimes 1_K)V_t^f)=(X_t\widetilde\otimes 1_K)(\lambda (t)^*\widetilde\otimes 1_K)V_t^g=((\lambda (t)^*X_t)\widetilde\otimes 1_K)V_t^g$$
so $(\varphi X)_t=\lambda (t)^*X_t$.

$a_2\Rightarrow a_1$. Put $h:=f\delta \lambda $. By the above, for $t\in T$,
$$(\lambda (t)^*\widetilde\otimes 1_K)V_t^g=\varphi V_t^f=(\lambda (t)^*\widetilde\otimes 1_K)V_t^h$$
so $V_t^g=V_t^h$ and this implies $g=h$.

c) follows from a).\qed

{\it Remark.} Not every $E$-C*-isomorphism $\ssa{f}\rightarrow \ssa{g}$ is an $\ccc{S}$ isomorphism  (see Remark of \pr \ref{10}). 
\begin{co}\label{718}
Let
$$\Lambda _0(T,E):=\me{\lambda \in \Lambda (T,E)}{\lambda  \;\mbox{is a group homomorphism}}$$
and for every $\lambda \in \Lambda _0(T,E)$ put
$$\mae{\varphi _\lambda }{\ssa{f}}{\ssa{f}}{X}{U_\lambda ^*XU_\lambda }\;.$$
Then the map $\lambda \mapsto \varphi _\lambda $ is an injective group homomorphism.
\end{co}

By \pr \ref{706} c), $\Lambda _0(T,E)$ is the kernel of the map
$$\mad{\Lambda (T,E)}{\fte}{\lambda }{\delta \lambda }$$
so by \pr \ref{716}, $\varphi _\lambda $ is well-defined. Thus only the injectivity of the map has to be proved. For $t\in T$ and $\zeta \in \breve E$, by \pr \ref{674} c),
$$U_\lambda ^*V_tU_\lambda (\zeta \otimes e_1)=U_\lambda ^*V_t(\zeta \otimes e_1)=U_\lambda ^*(\zeta \otimes e_t)=$$
$$=(\lambda (t)^*\zeta )\otimes e_t=(\lambda (t)^*\widetilde\otimes 1_K)V_t(\zeta \otimes e_1)\;.$$
So if $\varphi _\lambda $ is the identity map then $\lambda (t)=1_E$ for every $t\in T$.\qed

\begin{p}\label{733}
Let $F$ be a unital C**-algebra, $\varphi :E\rightarrow F$ a surjective C**-homomorphism, $g:=\varphi \circ f\in \ccc{F}(T,F)$, and  $L:=\widetilde{\cb{t\in T}}\breve F$. We put for all $\xi \in H,\;\eta \in L$, and $X\in \leh$, 
$$\tilde \xi :=(\varphi \xi _i)_{i\in I}\in L\,,\qquad \tilde X\eta :=\widetilde{X\zeta }\in L\,,$$
where $\zeta \in H$ with $\tilde \zeta =\eta $ \emph{(Lemma \ref{755} a),b) and \pr \ref{1} a))}. Then 
$$\tilde X=\sii{t\in T}{\fr{T}_3}((\varphi X_t)\widetilde\otimes 1_K)V_t^g\in \ssa{g}$$
for every $X\in \ssa{f}$ and the map
$$\mae{\tilde \varphi }{\ssa{f}}{\ssa{g}}{X}{\tilde X}$$
is a surjective C**-homomorphism, continuous with respect to the topologies $\fr{T}_k$, $k\in \{1,2,3\}$ such that
$$Ker \,\tilde \varphi =\me{X\in \ssa{f}}{t\in T\Longrightarrow X_t\in Ker\, \varphi }\;.$$
\end{p}   

For $s,t\in T$ and $\xi \in H$,
$$\left(\widetilde{\overbrace{(X_t\widetilde\otimes 1_K)V_t^f}}\tilde \xi \right)\hspace{-0.2em}_s=(\widetilde{\overbrace{(X_t\widetilde\otimes 1_K)V_t^f\xi }})_s=\varphi ((X_t\widetilde\otimes 1_K)V_t^f\xi )_s=$$
$$=\varphi (f(t,t^{-1}s)X_t\xi _{t^{-1}s})=g(t,t^{-1}s)(\varphi X_t)\tilde \xi _{t^{-1}s}=(((\varphi X_t)\widetilde\otimes 1_K)V_t^g\tilde \xi )_s$$
so by Lemma \ref{755} b),
$$\widetilde{\overbrace{(X_t\widetilde\otimes 1_K)V_t^f}}=((\varphi X_t)\widetilde\otimes 1_K)V_t^g\;.$$

By \h \ref{745} b),
$$X=\sii{t\in T}{\fr{T}_3}(X_t\widetilde\otimes 1_K)V_t^f$$
so by the above and by \pr \ref{1} b),
$$\tilde X=\sii{t\in T}{\fr{T}_3}((\varphi X_t)\widetilde\otimes 1_K)V_t^g\in \ssa{g}\;.$$
By \pr \ref{1} b), $\tilde \varphi $ is a surjective C**-homomorphism, continuous with respect to the topologies ${\fr{T}_k}\;(k\in \{1,2,3\})$. The last assertion is easy to see.\qed

\begin{co}\label{35}
Let $F$ be a unital C*-algebra, $\varphi :E\rightarrow F$ a unital C*-homomorphism such that $\varphi (\un{E})\subset F^c$, $g:=\varphi \circ f\in \f{T}{F}$, and $L:=\cb{t\in T}\breve F$. Then the map
$$\mae{\tilde \varphi }{\ssb{\n{\cdot }}{f}}{\ssb{\n{\cdot }}{g}}{X}{\sii{t\in T}{\n{\cdot }}((\varphi X_t)\otimes 1_L)V_t^g}$$
is C*-homomorphism.
\end{co}

Put $G:=E/Ker\,\varphi $ and denote by $\varphi _1:E\rightarrow G$ the quotient map and by $\varphi _2:G\rightarrow F$ the corresponding injective C*-homomorphism. By \pr \ref{733}, the corresponding map
$$\mac{\tilde \varphi _1}{\ssb{\n{\cdot }}{f}}{\ssb{\n{\cdot }}{\varphi _1\circ f}}$$
is a C*-homomorphism and by \h \ref{745} k), the corresponding map
$$\mac{\tilde \varphi _2}{\ssb{\n{\cdot }}{\varphi _1\circ f}}{\ssb{\n{\cdot }}{g}}$$
is also a C*-homomorphism. The assertion follows from $\tilde \varphi =\tilde \varphi _2\circ \tilde \varphi _1$.\qed

\begin{p}\label{721}
Let $T'$ be a group, $K':=l^2(T')$, $H':=\breve E\widetilde\otimes K'$, $\psi :T\rightarrow T'$ a surjective group homomorphism such that
$$\sup_{t'\in T'}Card\stackrel{-1}{\psi }\hspace{-0.4em}(t')\in \bn\,,$$
and $f'\in \ccc{F}(T',E)$ such that $f'\circ (\psi \times \psi )=f$. If we put 
$$X'_{t'}:=\si{t\in \stackrel{-1}{\psi }(t')}X_t$$
for every $X\in \ssa{f}$ and $t'\in T'$ then the family $((X'_{t'}\widetilde\otimes 1_{K'})V_{t'}^{f'})_{t'\in T'}$ is summable in $\lb{E}{H'}_{\fr{T}_2}$ for every $X\in \ssa{f}$ and the map 
$$\mae{\tilde \psi }{\ssa{f}}{\ssa{f'}}{X}{X':=\sii{t'\in T'}{\fr{T}_1}(X'_{t'}\widetilde\otimes 1_{K'})V_{t'}^{f'}}$$
is a surjective $E$-C**-homomorphism.  

We may drop the hypothesis that $\psi $ is surjective if we replace $\ccc{S}$ by $\ccc{S}_{\n{\cdot }}$.
\end{p}

Let $X\in \ssa{f}$. By Corollary \ref{40} a), since $\psi $ is surjective and
$$\sup_{t'\in T'}Card\stackrel{-1}{\psi }\hspace{-0.4em}(t')\in \bn\,,$$
it follows that the family $((X'_{t'}\widetilde\otimes 1_{K'})V_{t'}^{f'})_{t'\in T'}$ is summable in $\lb{E}{H'}_{\fr{T}_2}$ and therefore $X'\in \ssa{f'}$.

Let $X,Y\in \ssa{f}$. By \h \ref{745} c),g), for $t'\in T'$,
$$(X'^*)_{t'}=\widetilde{f'}(t')(X_{t'^{-1}})^*=\widetilde{f'}(t')\left(\si{t\in \stackrel{-1}{\psi }(t'^{-1})}X_t\right)^*=\widetilde{f'}(t')\si{s\in \stackrel{-1}{\psi }(t')}(X_{s^{-1}})^*=$$
$$=\si{s\in \stackrel{-1}{\psi }(t')}\tilde f(s)(X_{s^{-1}})^*=\si{s\in \stackrel{-1}{\psi }(t')}(X^*)_s=(X^*)'_{t'}\,,$$
$$(X'Y')_{t'}=\siw{s'\in T'}f'(s',s'^{-1}t')X'_{s'}Y'_{s'^{-1}t'}=$$
$$=\siw{s'\in T'}f'(s',s'^{-1}t')\left(\si{s\in \stackrel{-1}{\psi }(s')}X_s\right)\left(\si{r\in \stackrel{-1}{\psi }(s'^{-1}t')}Y_r\right)=$$
$$=\siw{s'\in T'}f'(s',s'^{-1}t')\left(\si{s\in \stackrel{-1}{\psi }(s')}\si{t\in \stackrel{-1}{\psi }(t')}X_sY_{s^{-1}t}\right)=$$
$$=\siw{s'\in T'}\left(\si{s\in \stackrel{-1}{\psi }(s')}\si{t\in \stackrel{-1}{\psi }(t')}f(s,s^{-1}t)X_sY_{s^{-1}t}\right)=$$
$$=\si{t\in \stackrel{-1}{\psi }(t')}\siw{s\in T}f(s,s^{-1}t)X_sY_{s^{-1}t}=\si{t\in \stackrel{-1}{\psi }(t')}(XY)_t=(XY)'_{t'}\;.$$
Thus $\psi $ is a C*-homomorphism. The other assertions are easy to see.

The last assertion follows from Corollary \ref{776} d).\qed

\begin{co}\label{763}
If we use the notation of \emph{\pr \ref{721}} and \emph{Corollary \ref{35}} and define $\widetilde {\varphi '}$ and $\widetilde {\psi '}$ in an obvious way then $\widetilde {\varphi '}\circ \tilde \psi =\widetilde {\psi '}\circ \tilde \varphi $.
\end{co}

For $X\in \ssa{f}$ and $t'\in T'$,
$$(\widetilde {\varphi '}\tilde \psi X)_{t'}=\varphi ((\tilde \psi X)_{t'})=\varphi \si{t\in \stackrel{-1}{\psi }(t')}X_t=\si{t\in \stackrel{-1}{\psi }(t')}\varphi X_t\,,$$
$$(\widetilde {\psi '}\tilde \varphi X)_{t'}=\si{t\in \stackrel{-1}{\psi }(t')}(\tilde \varphi X)_t=\si{t\in \stackrel{-1}{\psi }(t')}\varphi X_t\,,$$
so
$$\widetilde {\varphi '}\circ \tilde \psi =\widetilde {\psi '}\circ \tilde \varphi \;.\qedd$$ 

\begin{p}\label{726}
Let $F$ be a unital C*-subalgebra of $E$ such that $f(s,t)\in F$ for all $s,t\in T$. We denote by $\psi :F\rightarrow E$ the inclusion map and put
$$\mae{f^F}{T\times T}{Un\;F^c}{(s,t)}{f(s,t)}\,,$$
$$H^F:=\cb{t\in T}\breve F\approx \breve F\otimes K\,,$$
$$\mae{\tilde \psi }{H^F}{H}{\xi }{(\psi \xi _t)_{t\in T}}\;.$$
Moreover we denote for all $s,t\in T$ by $u_t^F,\;V_t^F$, and $\varphi _{s,t}^F$ the corresponding operators associated with $F$ $(f^F\in \f{T}{F})$. Let $X\in \ssb{C}{f}$ such that $X(\tilde \psi \xi )\in \tilde \psi (H^F)$ for every $\xi \in H^F$ and put
$$\mae{X^F}{H^F}{H^F}{\xi }{\xi '}\,,$$
where $\xi '\in H^F$ with $\tilde \psi \xi '=X(\tilde \psi \xi )$,  and $X^F_t:=(u_1^F)^*X^Fu_t^F\in F$ (by the canonical  identification of $F$ with $\lb{F}{\breve F}$) for every $t\in T$.
\begin{enumerate}
\item $\xi ,\eta \in H^F\Rightarrow \s{\tilde \psi \xi }{\tilde \psi \eta }=\psi \s{\xi }{\eta }$.
\item $\tilde \psi $ is linear and continuous with $\n{\tilde \psi }=1$.
\item $X^F$ is linear and continuous with $\n{X^F}=\n{X}$.
\item For $s,t\in T$,
$$\psi \varphi _{s,t}^FX^F=\varphi _{s,t}X\,,\qquad  \psi X_t^F=X_t\,,\qquad 
 \varphi _{s,t}^FX^F=f^F(st^{-1},t)X_{st^{-1}}^F\;.$$
\item $X^F\in \ssa{f^F}$.
\item $\xi \in H^F\Rightarrow X(\tilde \psi \xi )=\sii{t\in T}{\n{\cdot }}(X_t\otimes 1_K)V_t(\tilde \psi \xi )$.
\end{enumerate}
\end{p}

a\&b\&c are easy to see.

d) By a) and \pr \ref{681} b),
$$\varphi _{s,t}^FX^F=\s{X^F(1_F\otimes e_t)}{1_F\otimes e_s}\,,$$
$$\psi \varphi _{s,t}^FX^F=\psi \s{X^F(1_F\otimes e_t)}{1_F\otimes  e_s}=$$
$$=\s{\tilde \psi (X^F(1_F\otimes e_t))}{\tilde \psi (1_F\otimes e_s)}=
\s{X(1_E\otimes e_t)}{1_E\otimes e_s}=\varphi _{s,t}X\;.$$
In particular
$$\psi X_t^F=\psi \varphi _{1,t}^FX^F=\varphi _{1,t}X=X_t$$
and by \pr \ref{745a},
$$\psi \varphi _{s,t}^FX^F=\varphi _{s,t}X=f(st^{-1},t)X_{st^{-1}}=\psi (f^F(st^{-1},t)X_{st^{-1}}^F)\,,$$
$$\varphi _{s,t}^FX^F=f^F(st^{-1},t)X_{st^{-1}}^F\;.$$

e) By c) and \pr \ref{676} d), for $\xi \in H^F$,
$$\sii{t\in T}{\n{\cdot }}u_t^F(u_t^F)^*\xi =\xi \,,$$
$$X^F\xi =X^F\sii{t\in T}{\n{\cdot }}u_t^F(u_t^F)^*\xi =\sii{t\in T}{\n{\cdot }}X^Fu_t^F(u_t^F)^*\xi \,,$$
$$X^F\xi =\sii{s\in T}{\n{\cdot }}u_s^F(u_s^F)^*X^F\xi =\sii{s\in T}{\n{\cdot }}\sii{t\in T}{\n{\cdot }}u_s^F((u_s^F)^*X^Fu_t^F)(u_t^F)^*\xi \;.$$
By d) and \pr \ref{677} b),d),
$$X^F\xi =\sii{s\in T}{\n{\cdot }}\sii{t\in T}{\n{\cdot }}u_s^Ff^F(st^{-1},t)X_{st^{-1}}^F(u_t^F)^*\xi =\sii{s\in T}{\n{\cdot }}\sii{t\in T}{\n{\cdot }}u_s^FX_{st^{-1}}^F(u_s^F)^*V_{st^{-1}}^F\xi =$$
$$=\sii{s\in T}{\n{\cdot }}\sii{r\in T}{\n{\cdot }}u_s^FX_r^F(u_s^F)^*V_r^F\xi =\sii{s\in T}{\n{\cdot }}\sii{r\in T}{\n{\cdot }}u_s^F(u_s^F)^*(X_r^F\otimes 1_F)V_r^F\xi =$$
$$=\sii{s\in T}{\n{\cdot }}u_s^F(u_s^F)^*\sii{t\in T}{\n{\cdot }}(X_t^F\otimes 1_K)V_t^F\xi =\sii{t\in T}{\n{\cdot }}(X_t^F\otimes 1_K)V_t^F\xi $$
by \pr \ref{676} d), again. Thus
$$X^F=\sii{t\in T}{\fr{T}_2}(X_t^F\otimes 1_K)V_t^F\in \ssb{C}{f^F}\;.$$

f) For $s,t\in T$, by d),
$$(\tilde \psi ((X_t^F\otimes 1_K)V_t^F\xi ))_s=\psi ((X_t^F\otimes 1_K)V_t^F\xi )_s=\psi (f^F(t,t^{-1}s)X_t^F\xi _{t^{-1}s})=$$
$$=f(t,t^{-1}s)X_t(\tilde \psi \xi )_{t^{-1}s}=((X_t\otimes 1_K)V_t\tilde \psi \xi )_s\,,$$
$$\tilde \psi ((X_t^F\otimes 1_K)V_t^F\xi )=(X_t\otimes 1_K)V_t\tilde \psi \xi $$
so by b) and e),
$$X(\tilde \psi \xi )=\tilde \psi (X^F\xi )=\tilde \psi \left(\sii{t\in T}{\n{\cdot }}(X_t^F\otimes 1_K)V_t^F\xi \right)=$$
$$=\sii{t\in T}{\n{\cdot }}\tilde \psi ((X_t^F\otimes 1_K)V_t^F\xi )=\sii{t\in T}{\n{\cdot }}(X_t\otimes 1_K)V_t(\tilde \psi \xi )\;.\qedd$$
 
 \begin{p}\label{729}
Let $F$ be a W*-algebra such that $E$ is a unital C*-subalgebra of $F$ generating it as W*-algebra, $\varphi :E\rightarrow F$ the inclusion map, and $\tilde \xi :=(\varphi \xi _t)_{t\in T}\in L$ for every $\xi \in H$, where
$$L:=\cw{t\in T}\breve F\approx \breve F\bar \otimes K\;.$$
\begin{enumerate}
\item $\varphi (Un\;E^c)\subset Un\;F^c$ and $g:=\varphi \circ f\in \ccc{F}(T,F)$.
\item If
$$\mae{\psi }{\leh}{\lb{F}{L}}{X}{\bar X}$$
is the injective C*-homomorphism defined in \pr \ref{764} b), then $\psi (\ssb{C}{f})\subset \ssb{W}{g}$, $\psi (\ssb{C}{f})$ generates $\ssb{W}{g}$ as W*-algebra, and for every $X\in \ssb{C}{f}$ and $t\in T$ we have $(\bar X)_t=\varphi X_t$.
\item The following are equivalent for every $Y\in \ssb{W}{g}$:
\begin{enumerate}
\item $Y\in \psi (\ssb{C}{f})$.
\item $\xi \in H\Rightarrow Y\tilde \xi \in H$.

If these conditions are fulfilled then
\item $(Y_t)_{t\in T}\in H$.
\item $(Y^*_t)_{t\in T}\in H$.
\item $\xi \in H\Rightarrow Y\tilde \xi  =\sii{t\in T}{\n{\cdot }}(Y_t\bar \otimes 1_K)V_t^g\tilde \xi  \in H$.
\end{enumerate}
\end{enumerate}
\end{p}

a) follows from the density of $\varphi (E)$ in $F_{\ddot F}$ (Lemma \ref{766} $a\Rightarrow c$).

b) For $x\in E$, $t\in T$, and $\xi \in H$, 
$$(((\varphi x)\bar \otimes 1_K)V_t^g\tilde \xi )_s=g(t,t^{-1}s)(\varphi x)\tilde \xi _{t^{-1}s}=$$
$$=\varphi (f(t,t^{-1}s)x\xi _{t^{-1}s})=\varphi ((x\otimes 1_K)V_t\xi _s)$$
so
$$((\varphi x)\bar \otimes 1_K)V_t^g=\overline{(x\otimes 1_K)V_t^f}\;.$$

Let now $X\in \ssa{f}$. By \h \ref{745} b),
$$X=\sii{t\in T}{\fr{T}_2}(X_t\otimes 1_K)V_t^f$$
so by the above and by \pr \ref{764} c) (and \h \ref{745} d)),
$$\bar X=\sii{t\in T}{\fr{T}_1}\overline{(X_t\otimes 1_K)V_t^f}=\sii{t\in T}{\fr{T}_1}((\varphi X_t)\bar \otimes 1_K)V_t^g\in \ssb{W}{g}$$
so $\psi (\ssb{C}{f})\subset \ssb{W}{f}$. By \h \ref{745} a), $(\bar X)_t=\varphi X_t$ for every $t\in T$.

Since $\varphi (E)$ is dense in $F_{\ddot F}$ (Lemma \ref{766} $a)\Rightarrow c)$) it follows that 
$$\ccc{R}(g)\subset \stackrel{\fr{T}_1}{\overline{\varphi (\ccc{R}(f))}}$$
so $\psi  (\ssa{f})$ is dense in $\ssa{g}_{\overbrace{\ssa{g}}^{..}}$ and therefore generates $\ssa{g}$ as W*-algebra (Lemma \ref{766} $c\Rightarrow a$).

$c_1\Rightarrow c_2$ follows from the definition of $\psi $.

$c_2\Rightarrow c_1$ follows from \pr \ref{726} e).

$c_2\Rightarrow c_3\&c_4$ follows from  \pr \ref{750} b).

$c_2\Rightarrow c_5$ follows from \pr \ref{726} f).\qed

\begin{lem}\label{771}
Let $E,F$ be W*-algebras, $G:=E\bar \otimes F$, and 
$$L:=\cw{t\in T}\breve G\approx \breve G\bar \otimes K$$.
\begin{enumerate}
\item If $z\in G^{\#}$ then $z\bar \otimes 1_K$ belongs to the closure of 
$$\me{w\bar \otimes 1_K}{w\in E\odot F,\;\n{w}\leq 1}$$
 in $\lb{G}{L}_{\stackrel{...}{L}}$ .
 \item For every $y\in F$, the map
 $$\mad{E_{\ddot E}^{\#}}{G_{\ddot G}}{x}{x\otimes y}$$
 is continuous.
\end{enumerate}
\end{lem}

a) By [C1] Corollary 6.3.8.7, there is a filter $\fr{F}$ on $\me{w\in E\odot F}{\n{w}\leq 1}$ converging to $z$ in $G_{\ddot G}^{\#}$. By Lemma \ref{18} b), for $(a,\xi ,\eta )\in \ddot G\times L\times L$,
$$\sa{z\bar \otimes 1_K}{\tia{a}{\xi }{\eta }}=\sa{z}{\sii{t\in T}{G}\xi _t\,a\,\eta _t^*}=$$
$$=\lim_{w,\fr{F}}\sa{w}{\sii{t\in T}{G}\xi _t\,a\,\eta _t^*}=\lim_{w,\fr{F}}\sa{w\bar \otimes 1_K}{\tia{a}{\xi }{\eta }}$$
which proves the assertion.

b) Let $(a_i,b_i)_{i\in I}$ be a finite family in $\ddot E\times \ddot F$. For $x\in E$,
$$\sa{x\otimes y}{\si{i\in I}a_i\otimes b_i}=\si{i\in I}\sa{x}{a_i}\sa{y}{b_i}=\sa{x}{\si{i\in I}\sa{y}{b_i}a_i}\;.$$
Since $\me{x\otimes y}{x\in E^{\#}}$ is a bounded set of $G$, the above identity proves the continuity.\qed

\begin{p}\label{781}
Let $F$ be a unital C**-algebra, $S$ a group, and $g\in \ccc{F}(S,F)$. We denote by $\otimes _\sigma $ the spatial tensor product and put
$$G:=E \otimes _\sigma F\qquad (resp.\; G:=E\bar \otimes F)\,,$$
$$ L:=\widetilde{\cb{s\in S}}\breve F\approx \breve F\widetilde\otimes l^2(S),\qquad M:=\widetilde{\cb{(t,s)\in T\times S}}\breve G\approx \breve G\widetilde\otimes l^2(T\times S)\,,$$
$$\mae{h}{(T\times S)\times (T\times S)}{\un{G}}{((t_1,s_1),(t_2,s_2))}{f(t_1,t_2)\otimes g(s_1,s_2)}\;.$$
\begin{enumerate}
\item $h\in \ccc{F}(T\times S,G),\qquad M\approx H\widetilde\otimes L,$

$\leh\otimes _\sigma \lb{F}{L}\subset  \lb{G}{M}\;\mbox{in the C*-case},$

$\leh\bar\otimes \lb{F}{L}\approx \lb{G}{M}\;\mbox{in the W*-case}\;.$
\item For $t\in T$, $s\in S$, $x\in E$, $y\in F$,
$$((x\widetilde \otimes 1_{l^2(T)})V_t^f) \otimes ((y\widetilde \otimes 1_{l^2(S)})V_s^g)=((x\otimes y)\widetilde \otimes 1_{l^2(T\times S)})V_{(t,s)}^h\;.$$
\item In the C*-case, $\ssb{\n{\cdot }}{f}\otimes _\sigma \ssb{\n{\cdot }}{g}\approx \ssb{\n{\cdot }}{h}$ and $\ssb{C}{f}\otimes _\sigma \ssb{C}{g}\approx \ssb{C}{h}$.
\item In the W*-case, if $z\in G^{\#}$ and $(t,s)\in T\times S$ then $(z\bar \otimes 1_{l^2(T\times S)})V_{(t,s)}^h$ belongs to the closure of $\me{(w\bar \otimes 1_{l^2(T\times S)})V_{(t,s)}^h}{w\in (E\odot F)^{\#}}$ in $\lb{G}{M}_{\stackrel{...}{M}}$
\item In the W*-case, $\ssb{W}{f}\bar \otimes \ssb{W}{g}\approx \ssb{W}{h}$.
\end{enumerate}
\end{p}

a) $h\in \ccc{F}(T\times S,G)$ is obvious. 

Let us treat the C*-case first. For $\xi,\xi ' \in H$ and $\eta,\eta ' \in L$,
$$\s{\xi' \otimes \eta' }{\xi \otimes \eta }=\s{\xi' }{\xi }\otimes \s{\eta' }{\eta }=\left(\si{t\in T}\xi _t^*\xi' _t\right)\otimes \left(\si{s\in S}\eta _s^*\eta' _s\right)=$$
$$=\si{(t,s)\in T\times S}((\xi _t^*\xi' _t)\otimes (\eta _s^*\eta' _s))=\si{(t,s)\in T\times S}(\xi _t^*\otimes \eta _s^*)(\xi' _t\otimes \eta' _s)=$$
$$=\si{(t,s)\in T\times S}(\xi _t\otimes \eta _s)^*(\xi' _t\otimes \eta' _s)\,,$$
so the linear map
$$\mad{H\odot L}{M}{\xi \otimes \eta }{(\xi _t\otimes \eta _s)_{(t,s)\in T\times S}}$$
preserves the scalar products and it may be extended to a linear map $\varphi :H\otimes L\rightarrow M$ preserving the scalar products.

Let $z\in G$, $(t,s)\in T\times S$, and $\varepsilon >0$. There is a finite family $(x_i,y_i)_{i\in I}$ in $E\times F$ such that
$$\n{\si{i\in I}x_i\otimes y_i-z}<\varepsilon \;.$$
Then
$$\n{\si{i\in I}(x_i\otimes e_t)\otimes (y_i\otimes e_s)-z\otimes e_{(t,s)}}<\varepsilon $$
so $z\otimes e_{(t,s)}\in \overline{\varphi (H\otimes L)}=\varphi (H\otimes L)$. It follows that $\varphi $ is surjective and so $H\otimes L\approx M$.

The proof for the inclusion $\leh\otimes _\sigma \lb{F}{L}\subset \lb{G}{M}$  can be found in [L] page 37.

Let us now discus the W*-case. $\breve E\bar \otimes \breve F\approx \breve G$ follows from [C2] \pr 1.3 e), $M\approx H\bar \otimes L$ follows from [C3] Corollary 2.2, and $\leh\bar \otimes \lb{F}{L}\approx \lb{G}{M}$ follows from [C2] \h 2.4 d) or [C3] \h 2.4.

b) For $t_1,t_2\in T$, $s_1,s_2\in S$, $\xi \in \breve E$, and $\eta \in \breve F$, by \pr \ref{674} f) and [C3] Corollary 2.11,
$$(((x\widetilde \otimes 1_{l^2(T)})V_{t_1}^f)\widetilde \otimes ((y\widetilde \otimes 1_{l^2(S)})V_{s_1}^g))((\xi \otimes e_{t_2})\otimes (\eta \otimes e_{s_2}))=$$
$$=(((x\widetilde \otimes 1_{l^2(T)})V_{t_1}^f)(\xi \otimes e_{t_2}))\widetilde \otimes (((y\widetilde \otimes 1_{l^2(S)})V_{s_1}^g)(\eta \otimes e_{s_2}))\,,$$
$$((((x\otimes y)\widetilde \otimes 1_{l^2(T\times S)}))V_{(t_1,s_1)}^h)((\xi \otimes \eta )\otimes e_{(t_2,s_2)})=$$
$$=(h((t_1,s_1),(t_2,s_2))(x\otimes y)(\xi \otimes \eta ))\otimes e_{(t_1t_2,s_1s_2)}=$$
$$=((f(t_1,t_2)x\xi )\otimes (g(s_1,s_2)y\eta ))\otimes e_{t_1t_2}\otimes e_{s_1s_2}=$$
$$=(((x\widetilde \otimes 1_{l^2(T)})V_{t_1}^f)(\xi \otimes e_{t_2}))\widetilde \otimes (((y\widetilde \otimes 1_{l^2(S)})V_{s_1}^g)(\eta \otimes e_{s_2}))\;.$$
We put 
$$u:=((x\widetilde \otimes 1_{l^2(T)})V_t^f)\widetilde \otimes ((y\widetilde \otimes 1_{l^2(S)})V_s^g)-((x\otimes y)\widetilde \otimes 1_{l^2(T\times S)})V_{t,s}^h\in \lb{G}{M}\;.$$
By the above, $u(\zeta \otimes e_r)=0$ for all $\zeta \in \breve E\odot \breve F$ and $r\in T\times S$.

Let us consider the C*-case first. Since $\breve E\odot \breve F$ is dense in $\breve G$ , we get $u(z \otimes e_r)=0$ for all $z \in \breve G$ and $r\in T\times S$. For $\zeta \in M$, by \prr 4.1 e),
$$u\zeta =u\left(\si{r\in T\times S}(\zeta _r\otimes e_r)\right)=
\si{r\in T\times S}u(\zeta _r\otimes e_r)=0\,,$$
which proves the assertion in this case.

Let us consider now the W*-case. Let $z\in G^{\#}$ and $r\in T\times S$ and let $\fr{F}$ be a filter on $(E\odot F)^{\#}$ converging to $z$ in $G_{\ddot G}$ ([C1] Corollary 6.3.8.7). For $\eta \in M$, $a\in \ddot G$, and $r\in T\times S$,
$$\sa{z\otimes e_r}{\widetilde{(a,\eta )}}=\sa{\s{z\otimes e_r}{\eta }}{a}=\sa{\eta _r^*z}{a}=\sa{z}{a\eta _r^*}=$$
$$=\lim_{w,\fr{F}}\sa{w}{a\eta _r^*}=\lim_{w,\fr{F}}\sa{w\otimes e_r}{\widetilde{(a,\eta )}}\,,$$
$$\lim_{w,\fr{F}}w\otimes e_r=z\otimes e_r$$
in $M_{\ddot M}$. Since $u:M_{\ddot M}\rightarrow M_{\ddot M}$ is continuous (\prr 3.4 c)), we get by the above $u(z\otimes e_r)=0$. For $\zeta \in M$ it follows by \prr 4.6 c),
$$u\zeta =u\left(\sii{r\in T\times S}{\ddot M}(\zeta _r\otimes e_r)\right)=\sii{r\in T\times S}{\ddot M}u(\zeta _r\otimes e_r)=0$$
which proves the assertion in the W*-case.

c) By b), $\ccc{R}(f)\odot \ccc{R}(g)\subset \ccc{R}(h)$ so by a),
$$\ssb{\n{\cdot }}{f} \odot \ssb{\n{\cdot }}{g}\subset  \ssb{\n{\cdot }}{h}\,,\qquad \ssb{C}{f}\odot \ssb{C}{g}\subset \ssb{C}{h}\,,$$
$$\ssb{\n{\cdot }}{f}\otimes _\sigma  \ssb{\n{\cdot }}{g}\subset  \ssb{\n{\cdot }}{h}\,,\qquad \ssb{C}{f}\otimes _\sigma  \ssb{C}{g}\subset \ssb{C}{h}\;.$$
 Let $z\in G^{\#}$, $(t,s)\in T\times S$, and $\varepsilon >0$. There is a finite family $(x_i,y_i)_{i\in I}$ in $E\times F$ such that
$$\n{\si{i\in I}(x_i\otimes y_i)}<1\,,\qquad \n{\si{i\in I}(x_i\otimes y_i)-z}<\varepsilon \;.$$
By b),
$$\n{\si{i\in I}(((x_i\otimes 1_{l^2(T)})V_t^f)\otimes ((y_i\otimes 1_{l^2(S)})V_s^g))-(z\otimes 1_{l^2(T\times S)})V_{(t,s)}^h}<\varepsilon $$
and so by a),
$$\ccc{R}(h)\subset \;\stackrel{\n{\cdot }}{\overline{\ccc{R}(f)\odot \ccc{R}(g)}}\;\subset \;\stackrel{\fr{T}_2}{\overline{\ccc{R}(f)\odot \ccc{R}(g)}}\,,$$
$$\ssb{\n{\cdot }}{h}\subset \ssb{\n{\cdot }}{f}\otimes _\sigma \ssb{\n{\cdot }}{g}\,,\qquad \ssb{C}{h}\subset \ssb{C}{f}\otimes _\sigma \ssb{C}{g}\;.$$

d) By a) and Lemma \ref{771} a), there is a filter $\fr{F}$ on 
$$\me{w\bar \otimes 1_{l^2(T\times S)}}{w\in (E\odot F)^{\#}}$$
converging to $z\bar \otimes 1_{l^2(T\times S)}$ in $\lb{G}{M}_{\stackrel{...}{M}}$. For $\xi ,\eta \in M$ and $a\in \ddot G$,
$$\sa{(z\bar \otimes 1_{l^2(T\times S)})V_{(t,s)}^h}{\tia{a}{\xi }{\eta }}=\sa{z\bar \otimes 1_{l^2(T\times S)}}{V_{(t,s)}^h\tia{a}{\xi }{\eta }}=$$
$$=\lim_{w,\fr{F}}\sa{w\bar \otimes 1_{l^2(T\times S)}V_{(t,s)}^h}{{\tia{a}{\xi }{\eta }}}=\lim_{w,\fr{F}}\sa{(w\bar \otimes 1_{l^2(T\times S)})V_{(t,s)}^h}{\tia{a}{\xi }{\eta }}\,,$$
which proves the assertion.

e) By \h \ref{745} h), 
$$\left(\stackrel{\stackrel{...}{H}}{\overline{\ccc{R}(f)}}\right)^{\#}\hspace{-0.5em}=\ssb{W}{f}^{\#}\subset \leh,\qquad 
 \left(\stackrel{\stackrel{...}{L}}{\overline{\ccc{R}(g)}}\right)^{\#}\hspace{-0.5em}=\ssb{W}{g}^{\#}\subset \lb{F}{L}\;.$$
 By b), $\ccc{R}(f)\odot \ccc{R}(g)\subset \ccc{R}(h)$, so by Lemma \ref{771} b),
$$\ssb{W}{f}^{\#}\odot \ccc{R}(g)^{\#}\subset \ssb{W}{h}^{\#},\qquad
 \ssb{W}{f}^{\#}\otimes \ssb{W}{g}^{\#}\subset \ssb{W}{h}^{\#}\;.$$
 $$\ssb{W}{f}\otimes \ssb{W}{g}\subset \ssb{W}{h}$$
By [C3] \pr 2.5, 
  $$\ssb{W}{f}\bar \otimes \ssb{W}{g}\approx\; \stackrel{\stackrel{...}{M}}{\overline{\ssb{W}{f}\otimes \ssb{W}{g}}}\subset \ssb{W}{h}\;.$$
  
For $x\in E$, $y\in F$, and $(t,s)\in T\times S$, by b),
  $$((x\otimes y)\bar \otimes 1_{l^2(T\times S)})V_{(t,s)}^h
=((x\bar \otimes 1_{l^2(T)})V_t^f)\bar \otimes ((y\bar \otimes 1_{l^2(S)})V_s^g)\in 
 \ssb{W}{f}\bar \otimes \ssb{W}{g}\;.$$

Let $z\in G^{\#}$. By d), there is a filter $\fr{F}$ on 
$$\me{(w\bar \otimes 1_{l^2(T\times S)})V_{(t,s)}^h}{w\in (E\odot F)^{\#}}$$
 converging to 
$(z\bar \otimes 1_{l^2(T\times S)})V_{(t,s)}^h$
 in $\lb{G}{M}_{\stackrel{...}{M}}$, so by the above 
$$(z\bar \otimes 1_{l^2(T\times S)})V_{(t,s)}^h\in \ssb{W}{f}\bar \otimes \ssb{W}{g}\;.$$
We get
$$\ccc{R}(h)\subset \ssb{W}{f}\bar \otimes \ssb{W}{g},\qquad 
\ssb{W}{h}\subset \ssb{W}{f}\bar \otimes \ssb{W}{g}\,,$$
$$ \ssb{W}{h}= \ssb{W}{f}\bar \otimes \ssb{W}{g}\;.\qedd $$

\begin{co}\label{783}
Let $n\in \bn$ and
$$\mae{g}{T\times T}{Un\;(E_{n,n})^c}{(s,t)}{[\delta _{i,j}f(s,t)]_{i,j\in \bnn{n}}}\;.$$
\begin{enumerate}
\item $(\ssa{f})_{n,n}\approx \ssa{g},\quad (\ssb{\n{\cdot }}{f})_{n,n}\approx  \ssb{\n{\cdot }}{g}\;$.
\item Let us denote by $\rho :\ssa{g}\rightarrow (\ssa{f})_{n,n}$ the isomorphism of a). For $X\in \ssa{g}$, $t\in T$, and $i,j\in \bnn{n}$,
$$((\rho X)_{i,j})_t=(X_t)_{i,j}\;.$$
\end{enumerate}
\end{co}

a) Take $F:=\bk_{n,n}$ and $S:=\{1\}$ in \pr \ref{781}. Then $G\approx E_{n,n}$ and 
$$\mae{g}{T\times T}{Un\;G^c}{(s,t)}{f(s,t)\otimes 1_F}\;.$$
By \pr \ref{781} c),e),
$$\ssa{g}\approx \ssa{f} \otimes \bk_{n,n}\approx (\ssa{f})_{n,n}\,$$
$$\ssb{\n{\cdot }}{g}\approx \ssb{\n{\cdot }}{f}\otimes {\bk_{n,n}}\approx (\ssb{\n{\cdot }}{f})_{n,n}\;.$$

b) By \h \ref{745} b),
$$X=\sii{s\in T}{\fr{T}_3}(X_s\widetilde\otimes 1_K)V_s^g$$
so
$$(\rho X)_{i,j}=\sii{s\in t}{\fr{T}_3}((X_s)_{i,j}\widetilde\otimes 1_K)V_s^f\,,$$
$$((\rho X)_{i,j})_t=(X_t)_{i,j}$$
by \h \ref{745} a).\qed

\begin{co}\label{782}
Let $n\in \bn$. If $\bk=\bc$ (resp. if $n=4^m$ for some $m\in \bn$) then there is an $f\in \ccc{F}(\;\bzz{n}\times \bzz{n},E)$ (resp. $f\in \ccc{F}(\;(\bzz{2})^{2m},E)$) such that
$$\ccc{R}(f)=\ssa{f}\approx E_{n,n}\;.$$
\end{co}

By [C1] \pr 7.1.4.9 b),d) (resp. [C1] \h 7.2.2.7 i),k)) there is a $g\in \ccc{F}(\;\bzz{n}\times \bzz{n},\bc)$ (resp. $g\in \ccc{F}(\;(\bzz{2})^{2m},\bk)$) such that
$$\ssa{g}\approx \bc_{n,n}\qquad (\mbox{resp.} \;\ssa{g}\approx \bk_{n,n})\;.$$
If we put
$$\mae{f}{(\;\bzz{n}\times \bzz{n})\times (\;\bzz{n}\times \bzz{n})}{Un\;E^c}{(s,t)}{g(s,t)\otimes 1_E}$$
$$(\mbox{resp.}\;\mae{f}{(\bzz{2})^{2m}\times (\bzz{2})^{2m}}{Un\;E^c}{(s,t)}{g(s,t)\otimes 1_E})$$
then by \pr \ref{781} a),e), $f\in \ccc{F}(\;\bzz{n}\times \bzz{n},E)$ (resp. $f\in \ccc{F}(\;(\bzz{2})^{2m},E)$) and
$$\ssa{f}\approx \ssa{g}\otimes E\approx \bk_{n,n}\otimes E\approx E_{n,n}\;.\qedd$$

\begin{co}\label{12}
Let $F$ be a unital C**-algebra, $G:=E\widetilde\otimes F$, and
$$\mae{h}{T\times T}{\un{G}}{(s,t)}{f(s,t)\otimes 1_F}\;.$$
Then $h\in \f{T}{G}$ and
$$\ssb{\n{\cdot }}{h}\approx \ssb{\n{\cdot }}{f}\otimes F\,,\qquad \ssa{h}\approx \ssa{f}\widetilde\otimes F\;.\qedd$$
\end{co}

\begin{co}\label{15}
If $E$ is a W*-algebra then the following are equivalent:
\begin{enumerate}
\item $E$ is semifinite.
\item $\ssb{W}{f}$ is semifinite.
\end{enumerate}
\end{co}

$a\Rightarrow b.$ Assume first that there are a finite W*-algebra $F$ and a Hilbert space $L$ such that $E\approx F\bar \otimes \la{L}$. Put
$$\mae{g}{T\times T}{\un{F}}{(s,t)}{f(s,t)}\;.$$
By Corollary \ref{12},
$$\ssb{W}{f}\approx \ssb{W}{g}\bar \otimes \la{L}\;.$$
By Corollary \ref{711} c), $\ssb{W}{g}$ is finite and so $\ssb{W}{f}$ is semifinite.

The general case follows from the fact that $E$ is the C*-direct product of W*-algebras of the above form ([T] \pr V.1.40).

$b\Rightarrow a.$ $E$ is isomorphic to a W*-subalgebra of $\ssb{W}{f}$ (\h \ref{745} h)) and the assertion follows from [T] \h V.2.15.\qed

\begin{p}\label{neu}
Let $S,T$ be finite groups and $g\in \f{S}{\ssa{f}}$ and put $L:=l^2(S)$, $M:=l^2(S\times T)$, and
$$\mae{h}{(S\times T)\times (S\times T)}{\un{\ssa{f}}}{((s_1,t_1),(s_2,t_2))}{f(t_1,t_2)g(s_1,s_2)}\;.$$
Then $h\in \f{S\times T}{\ssa{f}}$ and the map
$$\mae{\varphi }{\ssa{g}}{\ssa{h}}{X}{\si{(s,t)\in S\times T}((X_s)_t\otimes 1_M)V_{(s,t)}^h}$$
is an $\ssa{f}$-C*-isomorphism.
\end{p}

For $X,Y\in \ssa{g}$, $Z\in \ssa{f}$, and $(s,t)\in S\times T$, by \h \ref{745} c),g),
$$(\varphi (X^*))_{(s,t)}=((X^*)_s)_t=(\tilde{g}(s)(X_{s^{-1}})^*)_t=((\tilde{g}(s)^*X_{s^{-1}} )^*)_t=$$
$$=\tilde{f}(t)((\tilde{g}(s)^*X_{s^{-1}} )_{t^{-1}})^*
=\tilde{f}(t)\tilde{g}(s)((X_{s^{-1}})_{t^{-1}})^*=$$
$$=\tilde{h}(s,t)((\varphi X)_{(s^{-1},t^{-1})})^*=\tilde{h}(s,t)((\varphi X)_{(s,t)^{-1}} )^*=((\varphi X)^*)_{(s,t)}\,,$$
$$((\varphi X)(\varphi Y))_{(s,t)}=\si{(r,u)\in S\times T}h((r,u),(r,u)^{-1}(s,t))(\varphi X)_{(r,u)}(\varphi Y)_{(r,u)^{-1}(s,t)}=$$
$$=\si{(r,u)\in S\times T}g(r,r^{-1}s)f(u,u^{-1}t)(X_r)_u(Y_{r^{-1}s})_{u^{-1}t}=\si{r\in S}g(r,r^{-1}s)(X_rY_{r^{-1}s})_t=$$
$$=\left(\si{r\in S}g(r,r^{-1}s)X_rY_{r^{-1}s}\right)_t=((XY)_s)_t=(\varphi (XY))_{(s,t)}\,,$$
$$(\varphi (ZX))_{(s,t)}=((ZX)_s)_t=((ZX)_s)_t=(ZX_s)_t=Z(X_s)_t=Z(\varphi X)_{(s,t)}$$
so
$$\varphi (X^*)=(\varphi X)^*\,,\qquad \varphi (XY)=(\varphi X)(\varphi Y),\qquad \varphi (ZX)=Z\varphi (X)$$
and $\varphi $ is an $\ssa{f}$-C*-homomorphism.

If $X\in \ssa{g}$ with $\varphi X=0$ then for $(s,t)\in S\times T$,
$$(X_s)_t=(\varphi X)_{(s,t)}=0\,,\qquad X_s=0\,,\qquad X=0$$
so $\varphi $ is injective.

Let $x\in E$ and $(s,t)\in S\times T$. Put
$$Z:=(x\otimes 1_K)V_t^f\in \ssa{f}\,,\qquad X:=(Z\otimes 1_L)V_s^g\in \ssa{g}\;.$$
Then for $(r,u)\in S\times T$,
$$(\varphi X)_{(r,u)}=(X_r)_u=\delta _{r,s}Z_u=\delta _{r,s}\delta _{u,t}x$$
so
$$\varphi X=(x\otimes 1_M)V_{(s,t)}^h$$
and $\varphi $ is surjective.\qed

\begin{p}\label{20}
Let $S$ be a finite subgroup of $T$ and $g:=f|(S\times S)$. We identify $\ssa{g}$ with the $E$-C**-subalgebra $\me{Z\in \ssa{f}}{t\in T\setminus S\Rightarrow Z_t=0}$ of $\ssa{f}$ \emph{(Corollary \ref{776} e))}. Let $X\in \ssa{f}\cap \ssa{g}^c$, $P_+:=X^*X$, and $P_-:=XX^*$ and assume $P_\pm \in Pr\;\ssa{f}$.
\begin{enumerate}
\item $P_\pm \in \ssa{g}^c$.
\item The map
$$\mae{\varphi _\pm }{\ssa{g}}{P_\pm \ssa{f}P_\pm }{Y}{P_\pm YP_\pm }$$
is a unital C**-homomorphism.
\item For every $Z\in \varphi _+(\ssa{g})$, $XZX^*\in \varphi _-(\ssa{g})$ and the map
$$\mae{\psi }{\varphi _+(\ssa{g})}{\varphi _-(\ssa{g})}{Z}{XZX^*}$$
is a C*-isomorphism with inverse
$$\mad{\varphi _-(\ssa{g})}{\varphi _+(\ssa{g})}{Z}{X^*ZX}$$
such that $\varphi _-=\psi \circ \varphi _+$.
\item If $p\in Pr\,\ssa{g}$ then
$$(X(\varphi _+p))^*(X(\varphi _+p))=\varphi _+p\,,\qquad (X((\varphi _+p))(X(\varphi _+p))^*=\varphi _-p\;.$$
\item If $\varphi _+$ is injective then $\varphi _-$ is also injective, the map
$$\mad{E}{P_\pm \ssa{f}P_\pm }{x}{P_\pm (x\widetilde\otimes 1_K)P_\pm }$$
is an injective unital C**-homomorphism, $P_\pm \ssa{f}P_\pm $ is an $E$-C**-algebra, $\varphi _\pm (\ssa{g})$ is an $E$-C**-subalgebra of it, and $\varphi _\pm $ and $\psi $ are $E$-C**-homo-morphisms.
\item The above results still hold for an arbitrary subgroup $S$ of $T$ if we replace $\ccc{S}$ by $\ccc{S}_{\n{\cdot }}$. 
\end{enumerate}
\end{p}

a) follows from the hypothesis on $X$.

b) follows from a).

c) Let $Y\in \ssa{g}$ with $Z=P_+YP_+$. By the hypotheses of the \pr,
$$XZX^*=XP_+YP_+X^*=XX^*XYX^*XX^*=$$
$$=XX^*YXX^*XX^*=P_-YP_-\in \varphi _-(\ssa{g})$$
and $\psi $ is a C*-homomorphism. The other assertions follow from
$$X^*(XZX^*)X=P_+ZP_+=P_+YP_+\;.$$

d) By b) and c),
$$(X(\varphi _+p))^*(X(\varphi _+p))=(\varphi _+p)X^*X(\varphi _+p)=(\varphi _+p)P_+(\varphi _+p)=\varphi _+p\,,$$
$$(X(\varphi _+p))(X(\varphi _+p))^*=X(\varphi _+p)(\varphi _+p)^*X^*=X(\varphi _+p)X^*=\psi \varphi _+p=\varphi _-p\;.$$

e) follows from b), c), and Lemma \ref{18}.

f) follows from Corollary \ref{776} d).\qed

{\it Remark.} Even if $\varphi _\pm $ is injective $P_\pm \ssa{f}P_\pm $ is not an $E$-C*-subalgebra of $\ssa{f}$.

\begin{theo}\label{895}
Let $S$ be a finite subgroup of $T$, $L:=l^2(S)$, $g:=f|(S\times S)$, $\omega :\bzz{2}\times \bzz{2}\rightarrow T$ an injective group homomorphism such that $S\cap \omega (\bzz{2}\times \bzz{2})=\{1\}$,
$$a:=\omega (1,0)\,,\quad b:=\omega (0,1)\,,\quad c:=\omega (1,1)\,,\quad \alpha _1:=f(a,a)\,,\quad \alpha _2:=f(b,b)\,,$$
$\beta _1,\beta _2\in \un{E}$ such that $\alpha _1\beta _1^2+\alpha _2\beta _2^2=0$,
$$\gamma :=\frac{1}{2}(\alpha _1^*\beta _1^*\beta _2-\alpha _2^*\beta _1\beta _2^*)=\alpha _1^*\beta _1^*\beta _2=-\alpha _2^*\beta _1\beta _2^*\,,$$
$$X:=\frac{1}{2}((\beta _1\widetilde\otimes 1_K)V_a^f+(\beta _2\widetilde\otimes 1_K)V_b^f)\,,\qquad P_+:=X^*X\,,\quad P_-:=XX^*\;.$$
We assume $f(s,c)=f(c,s)$ and $cs=sc$ for every $s\in S$, and $f(a,b)=-f(b,a)=1_E$. Moreover we consider $\ssa{g}$ as an $E$-C**-subalgebra of $\ssa{f}$ \emph{
(Corollary \ref{776} e))}.
\begin{enumerate}
\item We have
$$f(a,c)=-f(c,a)=\alpha _1\,,\quad f(b,c)=-f(c,b)=-\alpha _2\,,\quad f(c,c)=-\alpha _1\alpha _2\,,$$
$$\gamma ^2=-\alpha _1^*\alpha _2^* \,,\qquad \qquad\qquad V_c^f\in \ssa{g}^c\;.$$
\item We have
$$P_\pm =\frac{1}{2}(V_1^f\pm (\gamma \widetilde\otimes 1_K)V_c^f)\in \ssa{g}^c\cap Pr\,\ssa{f},\; P_++P_-=V_1^f,\; P_+P_-=0\,,$$
 $$X^2=0,\; XP_+=X,\;P_-X=X,\;P_+X=XP_-=0,\;X+X^*\in \unn{\ssa{f}}\,,$$
 $$Y\in \ssa{g}\Longrightarrow XYX=0\;.$$
\item The map
$$\mad{E}{P_\pm \ssa{f}P_\pm }{x}{(x\widetilde\otimes 1_K)P_\pm }$$
is a unital injective C**-homomorphism; we shall consider $P_\pm \ssa{f}P_\pm $ as an $E$-C**-algebra using this map. 
\item The maps
$$\mae{\varphi _+ }{\ssa{g}}{P_+ \ssa{f}P_+ }{Y}{P_+ YP_+ }\,,$$
$$\mae{\varphi _-}{\ssa{g}}{P_-\ssa{f}P_-}{Y}{XYX^*}$$
are orthogonal injective $E$-C**-homomorphisms and $\varphi _++\varphi _-$ is an injective $E$-C*-homomorphism. If $Y_1,Y_2\in \unn{\ssa{g}}$ (resp. $Y_1,Y_2\in Pr\;\ssa{g}$) then $\varphi _+Y_1+\varphi _-Y_2\in \unn{\ssa{f}}$  (resp. $\varphi _+Y_1+\varphi _-Y_2\in Pr\;\ssa{f}$). Moreover the map
$$\mae{\psi }{\ssa{f}}{\ssa{f}}{Z}{(X+X^*)Z(X+X^*)}$$
is an $E$-C**-isomorphism such that
$$\psi ^{-1}=\psi \,,\qquad \psi (P_+\ssa{f}P_+)=P_-\ssa{f}P_-\,,\qquad \psi \circ \varphi _+=\varphi _-\;.$$
If $\bk=\bc$ then $X+X^*$ is homotopic to $V_1^f$ in $\unn{\ssa{f}}$ and $\psi $ is homotopic to the identity map of $\ssa{f}$. Using this homotopy we find that $\varphi _+Y$ is homotopic in the above sense to $\varphi _-Y$ for every $Y\in \ssa{g}$ and $\varphi _+Y_1+\varphi _-Y_2$, $\varphi _-Y_1+\varphi _+Y_2$, $\varphi _+(Y_1Y_2)+P_-$, and $\varphi _+(Y_2Y_1+P_-$ are homotopic in the above sense for all $Y_1,Y_2\in \ssa{g}$.
\item Let $s\in S$ such that $sa=as$. Then
$$sb=bs\,,\qquad f(sc,c)f(s,c)=-\alpha _1\alpha _2\,,$$
$$f(sa,c)f(c,sa)^*=-1_E\,,\qquad f(a,s)f(s,a)^*=f(b,s)f(s,b)^*\;.$$
\item If $sa=as$ for every $s\in S$ then the map
$$\mad{S\times (\bzz{2}\times \bzz{2})}{T}{(s,r)}{s(\omega r)}$$
is an injective group homomorphism.
\item If $T$ is generated by $S\cup \omega (\bzz{2}\times \bzz{2})$ and $sa=as$ for every $s\in S$ then $\varphi _+ $ and $\psi _-$ are $E$-C*-isomorphisms with inverse
$$\mad{P_\pm \ssa{f}P_\pm }{\ssa{g}}{Z}{2\si{s\in S}(Z_s\widetilde\otimes 1_L)V_s^g}\,,$$
where
$$\mae{\psi _-}{\ssa{g}}{P_-\ssa{f}P_-}{Y}{P_-YP_-}\;.$$
\item If $sa=as$ and $f(a,s)=f(s,a)$ for every $s\in S$ then $X\in \ssa{g}^c$, $\varphi _-Y=P_-Y$ for every $Y\in \ssa{g}$, and there is a unique $\ssa{g}$-C**-homomorphism $\phi :\ssa{g}_{2,2}\rightarrow \ssa{f}$ such that
$$\phi \mt{0}{0}{(\alpha _1\beta _1^2)\otimes 1_L}{0}=X\;.$$
$\phi $ is injective and
$$\phi \mt{V_1^g}{0}{0}{0}=P_+\,,\qquad \phi \mt{0}{0}{0}{V_1^g}=P_-\;.$$
\item If $sa=as$ and $f(a,s)=f(s,a)$ for all $s\in S$ and if $T$ is generated by $S\cup \omega (\bzz{2}\times \bzz{2})$ then $\phi $ is an $\ssa{g}$-C*-isomorphism and
$$\phi ^{-1}V_1^f=\mt{1_E\otimes 1_L}{0}{0}{1_E\otimes 1_L}\,,\; \phi ^{-1}V_c^f=\mt{\gamma ^*\otimes 1_L}{0}{0}{-\gamma ^*\otimes 1_L}\,,$$
$$\phi ^{-1}V_a^f=\mt{0}{-\beta _1^*\otimes 1_L}{(\beta _2\gamma ^*)\otimes 1_L}{0}\,,$$
$$\phi ^{-1}V_b^f=\mt{0}{-\beta _2^*\otimes 1_L}{(\beta _1\gamma ^*)\otimes 1_L}{0}\,,$$
$$\phi ^{-1}P_+=\mt{V_1^g}{0}{0}{0}\,,\qquad \phi ^{-1}P_-=\mt{0}{0}{0}{V_1^g}\,,$$
and for every $s\in S$
$$\phi ^{-1}V_s^f=\mt{V_s^g}{0}{0}{V_s^g}\;.$$
\item The above results still hold for an arbitrary subgroup $S$ of $T$ if we replace $\ccc{S}$ with $\ccc{S}_{\n{\cdot }}$. \end{enumerate}
\end{theo}

a) By the equation of the Schur functions,
$$f(a,a)=f(a,c)f(a,b)\,,\; f(a,b)f(c,a)=f(a,c)f(b,a)\,,\; f(a,b)f(c,b)=f(b,b)\,,$$
$$f(b,a)f(c,b)=f(b,c)f(a,b)\,,\qquad f(a,b)f(c,c)=f(a,a)f(b,c)$$
and so
$$\alpha _1=f(a,c)\,,\qquad f(c,a)=-f(a,c)=-\alpha _1\,,\qquad f(c,b)=\alpha _2\,,$$
$$-\alpha _2=-f(c,b)=f(b,c)\,,\qquad f(c,c)=\alpha _1f(b,c)=-\alpha _1\alpha _2\;.$$
For $s\in S$, by \pr \ref{674} b),
$$V_c^fV_s^f=(f(c,s)\widetilde\otimes 1_K)V_{cs}^f=(f(s,c)\widetilde\otimes 1_K)V_{sc}^f=V_s^fV_c^f$$
and so $V_c^f\in \ssa{g}^c$ (by \pr \ref{674} d)). 

b) By \pr \ref{674} b),d),e) (and Corollary \ref{26} c)),
$$X^*=\frac{1}{2}(((\alpha _1^*\beta _1^*)\widetilde\otimes 1_K)V_a^f+((\alpha _2^*\beta _2^*)\widetilde\otimes 1_K)V_b^f)\,,$$
$$P_+=\frac{1}{4}(2V_1^f+((\alpha _1^*\beta _1^*\beta _2)\widetilde\otimes 1_K)V_c^f-((\alpha _2^*\beta _2^*\beta _1)\widetilde\otimes 1_K)V_c^f)=\frac{1}{2}(V_1^f+(\gamma\widetilde \otimes 1_K)V_c^f)\,,$$
$$P_-=\frac{1}{4}(2V_1^f+((\beta _1\alpha _2^*\beta _2^*)\widetilde\otimes 1_K)V_c^f-((\beta _2\alpha _1^*\beta _1^*)\widetilde\otimes 1_K)V_c^f)=\frac{1}{2}(V_1^f-(\gamma \widetilde\otimes 1_K)V_c^f)\;.$$
By a),
$$P_\pm ^*=\frac{1}{2}(V_1^f\pm (\gamma ^*\widetilde\otimes 1_K)((-\alpha _1^*\alpha _2^*)\widetilde\otimes 1_K)V_c^f)=P_\pm \,,$$
$$P_\pm ^2=\frac{1}{4}(V_1^f\pm 2(\gamma\widetilde \otimes 1_K)V_c^f+(\gamma ^2\widetilde\otimes 1_K)((-\alpha _1\alpha _2)\widetilde\otimes 1_K)V_1^f)=$$
$$=\frac{1}{2}(V_1^f\pm (\gamma \widetilde\otimes 1_K)V_c^f)=P_\pm \,,$$
so, by a) again, $P_\pm \in \ssa{g}^c\cap Pr\,\ssa{f}$. By \pr \ref{674} b),d),
$$X^2=\frac{1}{4}(((\beta _1^2\alpha _1+\beta _2^2\alpha _2)\widetilde\otimes 1_K)V_1^f+((\beta _1\beta _2)\widetilde\otimes 1_K)(V_a^fV_b^f+V_b^fV_a^f))=0\,,$$
$$(X+X^*)^2=X^2+XX^*+X^*X+X^{*2}=P_++P_-=V_1^f\;.$$
For the last relation we remark that by the above,
$$XYX=X(P_++P_-)YX=XP_+YX=XYP_+X=0\;.$$

c) follows from b) and Lemma \ref{18}.

d) By b) and c), the map $\varphi _\pm $ is an $E$-C**-homomorphism. Let $Y\in \ssa{g}$ with $\varphi _\pm Y=0$. By b), $Y=\mp Y(\gamma \widetilde\otimes 1_K)V_c^f$ so by \pr \ref{674} b),d) and \h \ref{745} b),
$$\si{s\in S}(Y_s\widetilde\otimes 1_K)V_s^f=\mp Y(\gamma \widetilde\otimes 1_K)V_c^f=\mp \si{s\in S}((Y_s\gamma f(s,c))\widetilde\otimes 1_K)V_{sc}^f\,,$$
which implies $Y_s=0$ for every $s\in S$ (\h \ref{745} a)). Thus $\varphi _\pm $ is injective. It follows that $\varphi _++\varphi _-$ is also injective. 

Assume first $Y_1,Y_2\in \unn{\ssa{g}}$. By b),
$$(\varphi _+Y_1+\varphi _-Y_2)^*(\varphi _+Y_1+\varphi _-Y_2)=(\varphi _+Y_1^*+\varphi _-Y_2^*)(\varphi _+Y_1+\varphi _-Y_2)=$$
$$=\varphi _+(Y_1^*Y_1)+\varphi _-(Y_2^*Y_2)=P_++P_-=V_1^f\;.$$
Similarly $(\varphi _+Y_1+\varphi _-Y_2)(\varphi _+Y_1+\varphi _-Y_2)^*=V_1^f$. The case $Y_1,Y_2\in Pr\;\ssa{g}$ is easy to see.

By b), $\psi $ is an $E$-C**-isomorphism with
$$\psi ^{-1}=\psi \,,\qquad \psi P_+=(X+X^*)X^*X(X+X^*)=XX^*XX^*=P_-\;.$$
Moreover for $Y\in \ssa{g}$,
$$\psi \varphi _+Y=(X+X^*)P_+YP_+(X+X^*)=XYX^*=\varphi _-Y\;.$$

Assume now $\bk=\bc$. By b), $X+X^*\in \unn{\ssa{f}}$. Being selfadjoint its spectrum is contained in $\{-1,+1\}$ and so it is homotopic to $V_1^f$ in $\unn{\ssa{f}}$.

e) We have $sb=sac=asc=acs=bs$. By a),
$$f(s,c)f(sc,c)=f(s,1)f(c,c)=-\alpha _1\alpha _2\,,$$
$$f(s,a)f(sa,c)=
f(s,b)f(a,c)=\alpha _1f(s,b)\,,$$ 
$$f(c,as)f(a,s)=f(c,a)f(b,s)=-\alpha _1f(b,s)\,,$$
$$f(c,bs)f(b,s)=f(c,b)f(a,s)=\alpha _2f(a,s)\,,$$
$$f(s,c)f(sc,b)=f(s,a)f(c,b)=\alpha _2f(s,a)\,,$$
$$f(c,s)f(cs,b)=f(c,sb)f(s,b)$$
so
$$f(sa,c)f(c,as)^*=-f(s,b)f(s,a)^*f(b,s)^*f(a,s)=$$
$$=-f(c,s)f(cs,b)f(c,sb)^*\alpha _2f(s,c)^*f(sc,b)^*\alpha _2^*f(c,bs)=-1_E\;.$$
From
$$f(s,c)f(sc,a)=f(s,b)f(c,a)\,,\qquad f(c,a)f(b,s)=f(c,as)f(a,s)\,,$$
$$ f(c,s)f(cs,a)=f(c,sa)f(s,a)$$
we get
$$f(a,s)f(s,a)^*=f(b,s)f(s,b)^*\;.$$

f) Since $S$ and $\omega (\bzz{2}\times \bzz{2})$ commute, the map is a group homomorphism. If $s(\omega r)=1$ for $(s,r)\in S\times (\bzz{2}\times \bzz{2})$ then $\omega r=s^{-1}\in S\cap \omega (\bzz{2}\times \bzz{2})$, which implies $s=1$ and $r=(0,0)$. Thus this group homomorphism is injective.

g) By e) and the hypothesis of f), for every $t\in T$ there are uniquely $s\in S$ and $d\in \{1,a,b,c\}$ with $t=sd$. Let $Z\in P_\pm \ssa{f}P_\pm $. By b) and \h \ref{745} b) (and Corollary \ref{742} d)),
$$Z=\pm (\gamma \widetilde\otimes 1_K)ZV_c^f=\pm (\gamma \widetilde\otimes 1_K)V_c^fZ$$
By \pr \ref{674} b),
$$ZV_c^f=\si{s\in S}((Z_sf(s,c))\widetilde\otimes 1_K)V_{sc}^f+\si{s\in S}((Z_{sa}f(sa,c))\widetilde\otimes 1_K)V_{sb}^f+$$
$$+\si{s\in S}((Z_{sb}f(sb,c))\widetilde\otimes 1_K)V_{sa}^f+\si{s\in S}((Z_{sc}f(sc,c))\widetilde\otimes 1_K)V_s^f\,,$$
$$V_c^fZ=\si{s\in S}((f(c,s)Z_s)\widetilde\otimes 1_K)V_{sc}^f+\si{s\in S}((f(c,sa)Z_{sa})\widetilde\otimes 1_K)V_{sb}^f+$$
$$+\si{s\in S}((f(c,sb)Z_{sb})\widetilde\otimes 1_K)V_{sa}^f+\si{s\in S}((f(c,sc)Z_{sc})\widetilde\otimes 1_K)V_s^f$$
and so by \h \ref{745} a),
$$Z_s=\pm \gamma f(sc,c)Z_{sc}=\pm \gamma f(c,sc)Z_{sc}\,,$$
$$Z_{sc}=\pm \gamma f(s,c)Z_s=\pm \gamma f(c,s)Z_s\,,$$
$$Z_{sa}=\pm \gamma f(sb,c)Z_{sb}=\pm \gamma f(c,sb)Z_{sb}\,,$$
$$Z_{sb}=\pm \gamma f(sa,c)Z_{sa}=\pm \gamma f(c,sa)Z_{sa}\;.$$
By e), $Z_{sa}=Z_{sb}=0$ for every $s\in S$. We get (by a), d), and \pr \ref{674} b))
$$\varphi _\pm (2\si{s\in S}(Z_s\widetilde\otimes 1_L)V_s^g)=\si{s\in S}(Z_s\widetilde\otimes 1_K)V_s^f\pm (\gamma \widetilde\otimes 1_K)V_c^f\si{s\in S}(Z_s\widetilde\otimes 1_K)V_s^f=$$
$$=\si{s\in S}(Z_s\widetilde\otimes 1_K)V_s^f\pm \si{s\in S}((\gamma f(c,s)Z_s)\widetilde\otimes 1_K)V_{sc}^f=$$
$$=\si{s\in S}(Z_s\widetilde\otimes 1_K)V_s^f+\si{s\in S}(Z_{sc}\widetilde\otimes 1_K)V_{sc}^f=Z\;.$$
Thus $\varphi _\pm $ is an $E$-C*-isomorphism with the mentioned inverse.

h) is a long calculation using e).

i) follows from h). 

j) follows from Corollary \ref{776} d).\qed

{\it Remark.} An example in which the above hypotheses are fulfilled is given in \h \ref{27}.

\begin{center}
\subsection{The functor $\ccc{S}$}

\fbox{Throughout this subsection we assume $T$ finite}
\end{center}

In this subsection we present the construction in the frame of category theory. Some of the results still hold for $T$ locally finite.

\begin{de}\label{46}
The above construction of $\ssa{f}$ can be done for an arbitrary $E$-module $F$, in which case we shall denote the result by $\ssa{F}$. Moreover we shall write $V_t^F$ instead of $V_t^f$ in this case. 
\end{de}

If $F$ is an $E$-module then $\ssa{F}$ is canonically an $E$-module. If in addition $F$ is adapted then $\ssa{F}$ is adapted and isomorphic to $\ssa{\check F,F}$. If $F$ is an $E$-C*-algebra then $\ssa{F}$ is also an $E$-C*-algebra.

\begin{p}\label{47}
If $F,G$ are $E$-modules and $\varphi :F\rightarrow G$ is an $E$-linear C*-homomorphism then the map
$$\mae{\ssa{\varphi }}{\ssa{F}}{\ssa{G}}{X}{\si{t\in S}((\varphi X_t)\otimes 1_K)V_t^G}$$
is an $E$-linear C*-homomorphism, injective or surjective if $\varphi $ is so.
\end{p}

The assertion follows from \h \ref{745} a),c),g).\qed
 
\begin{co}\label{48}
Let $F_1,\:F_2,\:F_3$ be $E$-modules and let $\varphi :F_1\rightarrow F_2$, $\psi :F_2\rightarrow F_3$ be $E$-linear C*-homomorphisms.
\begin{enumerate}
\item $\ssa{\psi }\circ \ssa{\varphi }=\ssa{\psi \circ \varphi }$.
\item If the sequence
$$0\longrightarrow F_1\stackrel{\varphi }{\longrightarrow }F_2\stackrel{\psi }{\longrightarrow }F_3$$
is exact then the sequence
$$0\longrightarrow \ssa{F_1}\stackrel{\ssa{\varphi }}{\longrightarrow }\ssa{F_2}\stackrel{\ssa{\psi }}{\longrightarrow }\ssa{F_3}$$
is also exact.
\item The covariant functor $\ccc{S}:\fr{M}_E\rightarrow \fr{M}_E$ is exact.
\end{enumerate}
\end{co}

a) is obvious.

b) Let $Y\in Ker\,\ssa{\psi }$. For every $t\in T$, $Y_t\in Ker\,\psi =Im\,\varphi $. If we identify $F_1$ with $Im\,\varphi $ then $Y_t\in F_1$. It follows $Y\in Im\,\ssa{\varphi }$, $Ker\,\ssa{\psi} =Im\,\ssa{\varphi}$.

c) follows from b) and \pr \ref{47}.\qed

\begin{co}\label{49}
Let $F$ be an adapted $E$-module and put
$$\mae{\iota }{F}{\check F}{x}{(0,x)}\,,$$
$$\mae{\pi }{\check F}{E}{(\alpha ,x)}{\alpha }\,,$$
$$\mae{\lambda }{E}{\check F}{\alpha }{(\alpha ,0)}\;.$$
Then the sequence
$$0\longrightarrow \ssa{F}\stackrel{\ssa{\iota }}{\longrightarrow }\ssa{\check F}{\scriptscriptstyle{\stackrel{\ssa{\pi }}{\longrightarrow }\atop\stackrel{\ssa{\lambda }}{\longleftarrow }}}\ssa{E}\longrightarrow 0$$
is split exact.\qed
\end{co}

\begin{p}\label{50}
The covariant functor $\ccc{S}:\fr{M}_E\rightarrow \fr{M}_E$ (resp. $\ccc{S}:\fr{C}_E^1\rightarrow \fr{C}_E^1$) \emph{(\pr \ref{47}, Corollary \ref{48} a))} is continuous with respect to the inductive limits \emph{(\pr \ref{36} a),b))}.
\end{p}

Let $\{(F_i)_{i\in I},\;(\varphi _{ij})_{i,j\in I}\}$ be an inductive system in the category $\fr{M}_E$ (resp. $\fr{C}_E^1$) and let $\{F,\;(\varphi _i)_{i\in I}\}$ be its limit in the category $\fr{M}_E$ (resp. $\fr{C}_E^1$). Then $\{(\ssa{F_i})_{i\in I},\;(\ssa{\varphi _{ij}}_{i,j\in I})\}$ is an inductive system in the category  $\fr{M}_E$ (resp. $\fr{C}_E^1$). Let $\{G,\;(\psi _i)_{i\in I}\}$ be its limit in this category and let $\psi :G\rightarrow \ssa{F}$ be the $E$-linear  C*-homomorphism such that $\psi \circ \psi _i=\ssa{\varphi _i}$ for every $i\in I$. In the $\fr{C}_E^1$ case, for $\alpha \in E$ and $i\in I$,
$$\psi (\alpha \otimes 1_K)=\psi \circ \psi _i(\alpha \otimes 1_K)=(\ssa{\varphi _i})(\alpha \otimes 1_K)=\alpha \otimes 1_K$$
so that $\psi $ is an $E$-C*-homomorphism.

Let $i\in I$ and let $X\in Ker\,\ssa{\varphi _i}$. Then $\varphi _iX_t=0$ for every $t\in T$. Since $T$ is finite, for every $\varepsilon >0$ there is a $j\in I$, $j\geq i$, with
$$\n{\varphi _{ji}X_t}<\frac{\varepsilon }{Card\;T}$$
for every $t\in T$. Then
$$\n{(\ssa{\varphi _{ji}})X}=\n{\si{t\in T}((\varphi _{ji}X_t)\otimes 1_K)V_t^{F_j}}<\varepsilon \;.$$
It follows
$$\n{\psi _iX}=\inf_{j\in I,j\geq i}\n{(\ssa{\varphi _{ji}})X}=0\,,$$
$$\psi _iX=0\,,\qquad X\in Ker\,\psi _i\,,\qquad Ker\,\ssa{\varphi _i}\subset Ker\,\psi _i\;.$$
By Lemma \ref{879}, $\psi $ is injective. Since
$$\bigcup _{i\in I}Im\,\ssa{\varphi _i}\subset Im\,\psi \,,$$
$Im\,\psi $ is dense in $\ssa{F}$. Thus $\psi $ is surjective and so an $E$-C*-isomorphism.\qed

\begin{p}\label{52}
Let $\theta :F\rightarrow G$ be a surjective morphism in the category $\fr{C}_E^1$. We use the notation of \emph{\h \ref{895}} and mark with an exponent if this notation is used with respect to $F$ or to $G$. For every $Y\in \unn{\ssa{g^G}}$,
 there is a $Z\in \ssa{g^F}$ such that
$$Z^*Z=P_+^F\,,\qquad\qquad \ssa{\theta }Z=\varphi _+^GY\;.$$
\end{p}

By \pr \ref{47} c), $\ssa{\theta }$ is surjective and so there is a $Z_0\in \ssa{g^F}$ with $\n{Z_0}=1$ and $\ssa{\theta }Z_0=Y$. Put
$$Z:=P_+^FZ_0+X^F(1-Z_0^*Z_0)^{\frac{1}{2}}\;.$$
By \h \ref{895} b),
$$Z^*Z=P_+^FZ_0^*Z_0+(1-Z_0^*Z_0)^{\frac{1}{2}}(X^F)^*X^F(1-Z_0^*Z_0)^{\frac{1}{2}}=$$
$$=P_+^FZ_0^*Z_0+P_+^F(1-Z_0^*Z_0)=P_+^F\;.$$
Since
$$\ssa{\theta }(1-Z_0^*Z_0)=1-Y^*Y=0$$
we get
$$\ssa{\theta }(1-Z_0^*Z_0)^{\frac{1}{2}}=0\,,\qquad\qquad \ssa{\theta }Z=P_+^GY=\varphi _+^GY\;.\qedd$$

\begin{p}\label{7904}
Let $F$ be an adapted $E$-module and $\Omega $ a locally compact space. We define for $X\in \ssa{\ccc{C}_0(\Omega ,F)}$ \emph{(see Corollary \ref{51} d))}  and $Y\in \ccc{C}_0(\Omega ,\ssa{F})$,
$$\mae{\varphi X}{\Omega }{\ssa{F}}{\omega }{\si{t\in T}(X_t(\omega )\otimes 1_K)V_t^F}\,,$$
$$\psi Y:=\si{t\in T}(Y(\cdot )_t\otimes 1_K)V_t^{\ccc{C}_0(\Omega ,F)}\;.$$
Then
$$\mac{\varphi }{\ssa{\ccc{C}_0(\Omega ,F)}}{\ccc{C}_0(\Omega ,\ssa{F})}\,,$$
$$\mac{\psi }{\ccc{C}_0(\Omega ,\ssa{F})}{\ssa{\ccc{C}_0(\Omega ,F)}}$$
are $E$-linear C*-isomorphisms and $\varphi =\psi ^{-1}$.

Let $\omega _0\in \Omega $ and assume $F$ is an $E$-C*-algebra. Then the above maps $\varphi $ and $\psi $ induce the following $E$-C*-isomorphisms
$$\ssa{\me{X\in \ccc{C}_0(\Omega ,F)}{X(\omega _0)\in E}}{\longrightarrow \atop \longleftarrow }\me{Y\in \ccc{C}_0(\Omega ,\ssa{F})}{Y(\omega _0)\in \ssa{E}}\;.$$
\end{p}

Let $X,X'\in \ssa{\ccc{C}_0(\Omega ,F)}$ and $Y,Y'\in \ccc{C}_0(\Omega ,\ssa{F})$. By Proposition \ref{750} b) and Corollary \ref{775} a),
$$\varphi X\in \ccc{C}_0(\Omega ,\ssa{F})\,,\qquad \psi Y\in \ssa{\ccc{C}_0(\Omega ,F)}$$
and it is easy to see that $\varphi $ and $\psi $ are $E$-linear. By \h \ref{745} c),g), for $t\in T$ and $\omega \in \Omega $,
$$((\varphi X)^*(\omega ))_t=\tilde{f}(t)(((\varphi X)(\omega )_{t^{-1}}))^*=$$
$$=\tilde{f}(t)X_{t^{-1}}(\omega )^*=(X^*(\omega ))_t=((\varphi X^*)(\omega ))_t\,,  $$
$$(((\varphi X)(\varphi X'))(\omega ))_t=\si{s\in T}f(s,s^{-1}t)((\varphi X)(\omega ))_s((\varphi X')(\omega ))_{s^{-1}t}=$$
$$=\si{s\in T}f(s,s^{-1}t)X_s(\omega )X'_{s^{-1}t}(\omega )=\left(\si{s\in T}f(s,s^{-1}t)X_sX'_{s^{-1}t}\right)(\omega )=$$
$$=(XX')_t(\omega )=((\varphi (XX'))(\omega ))_t\,,$$
so
$$(\varphi X)^*=\varphi X^*\,,\qquad (\varphi X)(\varphi X')=\varphi (XX')$$
and $\varphi $ is a C*-homomorphism. Similarly
$$(\psi Y^*)_t(\omega )=(Y^*(\omega ))_t=\tilde{f}(t)(Y(\omega )_{t^{-1}})^*=\tilde{f}(t)((\psi Y)_{t^{-1}}(\omega ))^*=((\psi Y)^*)_t(\omega )\,,  $$
$$((\psi Y)(\psi Y'))_t(\omega )=\left(\si{s\in T}f(s,s^{-1}t)(\psi Y)_s(\psi Y')_{s^{-1}t}\right)(\omega )=$$
$$=\si{s\in T}f(s,s^{-1}t)(\psi Y)_s(\omega )(\psi Y')_{s^{-1}t}(\omega )=\si{s\in T}f(s,s^{-1}t)Y(\omega )_sY'(\omega )_{s^{-1}t}=$$
$$=(Y(\omega )Y'(\omega ))_t=(\psi (YY')_t)(\omega )$$
so
$$\psi Y^*=(\psi Y)^*\,,\qquad (\psi Y)(\psi Y')=\psi (YY')$$
and $\psi $ is a C*-homomorphism. Moreover
$$(\psi \varphi X)_t(\omega )=((\varphi X)(\omega ))_t=X_t(\omega )\,,\qquad ((\varphi \psi Y)(\omega ))_t=(\psi Y)_t(\omega )=(Y(\omega ))_t\,,$$
so $\psi \varphi X=X$ and $\varphi \psi Y=Y$ which proves the assertion.

The last assertion is easy to see.\qed

\begin{p}\label{7906}
Let $F$ be an adapted $E$-module,
$$0\longrightarrow F\stackrel{\iota }{\longrightarrow }\check{F}\stackrel{\pi }{\longrightarrow }E\longrightarrow 0\,, $$
$$0\longrightarrow \ssa{F}\stackrel{\iota _0}{\longrightarrow }\check{\overbrace{\ssa{F}}}\stackrel{\pi _0}{\longrightarrow } E\longrightarrow 0 $$
the associated exact sequences \emph{(\pr \ref{33} h))}, and
$$\mae{j}{E}{\ssa{E}}{\alpha }{(\alpha \otimes 1_K)V_1^E}\,,$$
$$\mae{\varphi }{\check{\overbrace{\ssa{F}}} }{\ssa{\check{F} }}{(\alpha ,X)}{\ssa{\iota }X+(\alpha \otimes 1_K)V_1^{\check{F} }}\;.$$
Then $\varphi $ is an injective $E$-C*-homomorphism and $\ssa{\pi }\circ \varphi =j\circ \pi _0$.\qed
\end{p}

\begin{p}\label{7907}
If $E$ is commutative and $F$ is an $E$-module then the map
$$\mae{\varphi }{\ssa{E}\otimes F}{\ssa{F}}{X\otimes x}{\si{t\in T}((X_tx)\otimes 1_K)V_t^F}$$
is a surjective C*-homomorphism. If in addition $E=\bk$ then $\varphi $ is a C*-isomorphism with inverse
$$\mae{\psi }{\ssa{F}}{\ssa{E}\otimes F}{Y}{\si{t\in T}(V_t^E\otimes Y_t)}\;.$$
\end{p}

It is obvious that $\varphi $ is surjective. For $X,Y\in \ssa{E}$ and $x,y\in F$, by \h \ref{745} c),g) and \pr \ref{674} b),d),e),
$$\varphi ((X\otimes x)^*)=\varphi (X^*\otimes x^*)=\si{t\in T}(((X^*)_tx^*)\otimes 1_K)V_t^F=$$
$$=\si{t\in T}((\tilde{f}(t)(X_{t^{-1}})^*x^* )\otimes 1_K)V_t^F
=\si{t\in T}(((X_{t^{-1}})^*x^*)\otimes 1_K)(V_{t^{-1}}^F)^*=$$
$$=\si{t\in T}((x^*(X_t)^*)\otimes 1_K)(V_t^F)^*=(\varphi (X\otimes x))^*\,,$$
$$\varphi (X\otimes x)\varphi (Y\otimes y)=\si{s,t\in T}((X_sxY_ty)\otimes 1_K)V_s^FV_t^F=$$
$$=\si{s,t\in T}((f(s,t)X_sxY_ty)\otimes 1_K)V_{st}^F=
\si{r\in T}\si{s\in T}((f(s,s^{-1}r)X_sY_{s^{-1}r}xy)\otimes 1_K)V_r^F=$$
$$=\si{r\in T}(((XY)_rxy)\otimes 1_K)V_r^F=\varphi ((X\otimes x)(Y\otimes y))$$
so $\varphi $ is a C*-homomorphism.

Assume now $E=\bk$ and let $X\in \ssa{E}$ and $x\in F$. Then
$$\psi \varphi (X\otimes x)=\psi \si{t\in T}((X_tx)\otimes 1_K)V_t^F=\si{t\in T}V_t^E\otimes (X_tx)=$$
$$=\left(\si{t\in T}X_tV_t^E\right)\otimes x=X\otimes x$$
which proves the last assertion (by using the first assertion).\qed

\begin{center}
\section{Examples}

We draw the reader's attention to the fact that in additive groups the neutral element is denoted by 0 and not by 1.

\subsection{$T:=\bzz{2}$}
\end{center}

\begin{p}\label{785}
\rule{1em}{0ex}
\begin{enumerate}
\item The map
$$\mae{\psi }{\f{\bzz{2}}{E}}{\un{E}}{f}{f(1,1)}$$
is a group isomorphism.
\item $\psi (\me{\delta \lambda }{\lambda \in \Lambda (\bzz{2},E)})=\me{x^2}{x\in \un{E}}$.
\item If there is an $x\in E^c$ with $x^2=f(1,1)$ (in which case $x\in \un{E}$) then the map
$$\mae{\varphi }{\ssa{f}}{E\times E}{X}{(X_0+xX_1, X_0-xX_1)}$$
is an $E$-C*-isomorphism.
\item If $\bk=\bc$ and if $A$ is a connected and simply connected compact space or a totally disconnected compact space then for every $x\in \unn{\cca{A}}$ there is a $y\in \ccb{A}{\bc}$ with $x=e^y$.
\item Assume $\bk=\br$.
\begin{enumerate}
\item There are uniquely $p,q\in \Pr\;E^c$ with
$$p+q=1_E,\qquad \qquad pf(1,1)=p,\qquad\qquad qf(1,1)=-q\;.$$
\item The map
$$\mae{\varphi }{\ssa{f}}{(pE)\times (pE)\times \overbrace{qE}^{\circ }}{X}{\tilde X}\,,$$
where $\overbrace{qE}^{\circ }$ denotes the complexification of the C*-algebra $qE$ and 
$$\tilde X:=(p(X_0+X_1),p(X_0-X_1),(qX_0,qX_1))$$ 
for every $X\in \ssa{f}$, is an $E$-C*-isomorphism. In particular if $f(1,1)=-1_E$ then $\ssa{f}$ is isomorphic to the complexification of $E$.
\end{enumerate}
\item Assume $\bk=\bc$, let $\sigma (E^c)$ be the spectrum of $E^c$, and let $\widehat{f_{11}}$ be the function of $\ccc{C}(\sigma (E^c),\bc)$ corresponding to $f_{11}$ by the Gelfand transform. Then
$$\me{e^{i\theta }}{\theta \in \br,\;e^{2i\theta }\in \widehat{f_{11}}(\sigma (E^c))}$$
is the spectrum of $V_1$.
\end{enumerate}
\end{p}

a) follows from \pr \ref{704} a) (and \pr \ref{706} a) ).

b) follows from Definition \ref{705}.

c) For $X,Y\in \ssa{f}$, by \h \ref{745} c),g) (and \pr \ref{704} a)),
$$(X^*)_0=(X_0)^*\,,\qquad (X^*)_1=(x^*)^2(X_1)^*\,,$$
$$(XY)_0=X_0Y_0+x^2X_1Y_1\,,\qquad (XY)_1=X_0Y_1+X_1Y_0\,,$$
so
$$\varphi (X^*)=((X_0)^*+x(x^*)^2(X_1)^*, (X_0)^*-x(x^*)^2(X_1)^*)=$$
$$=((X_0)^*+x^*(X_1)^*\,, (X_0)^*-x^*(X_1)^*)=(\varphi X)^*\,,$$
$$(\varphi X)(\varphi Y)=((X_0+xX_1)(Y_0+xY_1),(X_0-xX_1)(Y_0-x Y_1))=$$
$$=(X_0Y_0+xX_0Y_1+xX_1Y_0+x^2X_1Y_1,X_0Y_0-xX_0Y_1-xX_1Y_0+x^2X_1Y_1)=$$
$$=((XY)_0+x(XY)_1,(XY)_0-x(XY)_1=\varphi (XY)$$
i.e. $\varphi $ is an $E$-C*-homomorphism. $\varphi $ is obviously injective.

Let $(y,z)\in E\times E$. If we take $X\in \ssa{f}$ with
$$X_0:=\frac{1}{2}(y+z),\qquad X_1:=\frac{1}{2}x^*(y-z)$$
then $\varphi X=(y,z)$, i.e. $\varphi $ is surjective.

d) is known.

$e_1)$ follows by using the spectrum of $E^c$.

$e_2)$ Put
$$\mae{\psi }{\ssa{f}}{\overbrace{qE}^{\circ }}{X}{(qX_0,qX_1)}\;.$$
For $X,Y\in \ssa{f}$, by \h \ref{745} c),g),
$$\psi (X^*)=(q(X^*)_0,q(X^*)_1)=(q(X_0)^*,qf(1,1)^*(X_1)^*)=$$
$$=((qX_0)^*,-(qX_1)^*)=(\psi X)^*\,,$$
$$(\psi X)(\psi Y)=(qX_0,qX_1)(qY_0,qY_1)=$$
$$=(q(X_0Y_0-X_1Y_1),(q(X_0Y_1+X_1Y_0)))=\psi (XY)$$
so $\psi $ is an $E$-C*-homomorphism. Thus by c), $\varphi $ is an $E$-C*-homomorphism. The bijectivity of $\varphi $ is easy to see.

f) By \pr \ref{674} e), $V_1$ is unitary so its spectrum is contained in $\me{e^{i\theta }}{\theta \in \br}$. For $\theta \in \br$ and $X\in \ssa{f}$,
$$(e^{i\theta }V_0-V_1)X=X(e^{i\theta }-V_1)=$$
$$=((e^{i\theta }X_0)\otimes 1_K)V_0+((e^{i\theta }X_1)\otimes 1_K)V_1-(X_0\otimes 1_K)V_1-((f_{11}X_1)\otimes 1_K)V_1=$$
$$=((e^{i\theta }X _0-f_{11}X_1)\otimes 1_K)V_0+((e^{i\theta }X_1-X_0)\otimes 1_K)V_1\;.$$
Thus $X$ is the inverse of $e^{i\theta }V_0-V_1$ iff $X_0=e^{i\theta }X_1$ and $e^{i\theta }X_0-f_{11}X_1=1_E$, i.e. $(e^{2i\theta }-f_{11})X_1=1_E$. Therefore $e^{i\theta }V_0-V_1$ is invertible iff $e^{2i\theta }-\widehat{f_{11}}$ does not vanish on $\sigma (E^c)$.\qed

\begin{co}\label{886}
Assume $\bk:=\br$ and let $S$ be a group, $F$ a unital C*-algebra, $g\in \f{S}{F}$, and
$$\mae{h}{(S\times \bzz{2})\times (S\times \bzz{2})}{\un{F}}{((s_1,t_1),(s_2,t_2))}$$
$${\ab{-g(s_1,s_2)}{(t_1,t_2)=(1,1)}{g(s_1,s_2)}{(t_1,t_2)\not=(1,1)}}\;.$$
\begin{enumerate}
\item $h\in \f{S\times \bzz{2}}{F}$. 
\item $\ssa{h}\approx \overbrace{\ssa{g}}^\circ \,,\qquad \ssb{\n{\cdot }}{h}\approx \overbrace{\ssb{\n{\cdot }}{g}}^\circ $. 
\end{enumerate}
\end{co}

Put $E:=\br$ in the above \pr and define $f\in \f{\bzz{2}}{\br}$ by $f(1,1)=-1$ (\pr \ref{785} a)). By this \pr $e_2)$, $\ssa{f}\approx \bc$. Thus by \pr \ref{781} c),e),
$$\ssa{h}\approx \ssa{g}\otimes \ssa{f}\approx \overbrace{\ssa{g}}^\circ \,,\qquad \ssb{\n{\cdot }}{h}\approx \ssb{\n{\cdot }}{g}\otimes \ssb{\n{\cdot }}{f}\approx \overbrace{\ssb{\n{\cdot }}{g}}^\circ \;.\qedd$$

\begin{de}\label{787}
We put
$$\bt:=\me{z\in \bc}{|z|=1}\;.$$
\end{de}

\begin{e}\label{789}
Let $E:=\ccb{\bt}{\bc}$ and $f\in \f{\bzz{2}}{E}$ with
$$\mae{f(1,1)}{\bt}{\unn{\bc}}{z}{z}\;.$$
If we put
$$\mae{\tilde X}{\bt}{\bc}{z}{X_0(z^2)+zX_1(z^2)}$$
for every $X\in \ssa{f}$ then the map
$$\mae{\varphi }{\ssa{f}}{E}{X}{\tilde X}$$
is an isomorphism of C*-algebras (but not an $E$-C*-isomorphism).
\end{e}

For $X,Y\in \ssa{f}$, by \h \ref{745} c),g),
$$(X^*)_0=(X_0)^*,\qquad (X^*)_1=\overline{f(1,1)}(X_1)^*\,,$$
$$(XY)_0=X_0Y_0+f(1,1)X_1Y_1,\qquad (XY)_1=X_0Y_1+X_1Y_0$$
so for $z\in \bt$,
$$\widetilde{X^*}(z)=X_0^*(z^2)+z\bar z^2X_1^*(z^2)=\overline{X_0(z^2)+zX_1(z^2)}=\tilde X^*(z)\,,$$
$$(\tilde X(z))(\tilde Y(z))=(X_0(z^2)+zX_1(z^2))(Y_0(z^2)+zY_1(z^2))=$$
$$=X_0(z^2)Y_0(z^2)+zX_0(z^2)Y_1(z^2)+zX_1(z^2)Y_0(z^2)+z^2X_1(z^2)Y_1(z^2)=$$
$$=(XY)_0(z^2)+z(XY)_1(z^2)=\widetilde{XY}(z)\,,$$
$$\widetilde{X^*}=\tilde X^*\,,\qquad\tilde X\tilde Y=\widetilde{XY}\,,$$
i.e. $\varphi $ is a C*-homomorphism. If $\varphi X=0$ then for $z\in \bt$,
$$X_0(z^2)+zX_1(z^2)=0$$
so, successively,
$$X_0(z^2)-zX_1(z^2)=0\,,\quad X_0(z^2)=X_1(z^2)=0\,,\quad X_0=X_1=0\,,\quad X=0$$
and $\varphi $ is injective.

Put
$$\ccc{G}:=\me{\si{k\in \bz}c_kz^k}{(c_k)_{k\in \bz}\in \bc^{(\bz)}}\subset E\;.$$
Let
$$x:=\si{k\in \bz}c_kz^k\in \ccc{G}$$
and take $X\in \ssa{f}$ with 
$$X_0:=\si{k\in \bz}c_{2k}z^k,\qquad X_1:=\si{k\in \bz}c_{2k+1}z^k\;.$$
Then
$$\tilde X=\si{k\in \bz}c_{2k}z^{2k}+z\si{k\in \bz}c_{2k+1}z^{2k}=x$$
so $\ccc{G}\subset \varphi (\ssa{f})$. Since $\ccc{G}$ is dense in $E$, $\varphi (\ssa{f})=E$ and $\varphi $ is surjective.\qed

\begin{de}\label{787'}
For every $x\in \ccb{\bt}{\bc}$ which does not take the value $0$ we put
$$w(x):={\bf{winding \;number \;of\; x}}:=\frac{1}{2\pi i}\int_x\frac{\mathrm{d}z}{z}=\frac{1}{2\pi i}[log\;x(e^{i\theta })]_{\theta =0}^{\theta =2\pi }\in \bz\;.$$
\end{de}

If $A$ is a connected compact space and $\gamma $ is a cycle in $A$ (i.e. a continuous map of $\bt$ in $A$), which is homologous to 0 (or more generally, if a multiple of $\gamma $ is homologous to 0), then for every $x\in \ccb{A}{\unn{\bc}}$ we have $w(x\circ \gamma )=0$. If $A$ is a compact space and $x\in \ccb{A}{\unn{\bc}}$ such that $w(x\circ \gamma )=0$ for every cycle $\gamma $ in $A$ then there is a $y\in \ccb{A}{\bc}$ with $x=e^y$.

\begin{e}\label{791}
Let $E:=\ccb{\bt}{\bc}$, $f\in \f{\bzz{2}}{E}$, and $n:=w(f(1,1))$.
\begin{enumerate}
\item If $n$ is even then there is an $x\in Un\;E$ with winding number equal to $\frac{n}{2}$ such that the map
$$\mad{\ssa{f}}{E\times E}{X}{(X_0+xX_1,X_0-xX_1)}$$
is an $E$-C*-isomorphism.
\item If $n$ is odd then $\ssa{f}$ is isomorphic to $E$.
\item The group $\f{\bzz{2}}{E}/\!\Lambda (\bzz{2},E)$ is isomorphic to $\bzz{2}$ and
$$Card\,(\me{\ssa{g}}{g\in \f{\bzz{2}}{E}}\hspace{-0.2em}/\hspace{-0.3em}\approx _\ccc{S})=2\;.$$
\item There is a complex unital C*-algebra $E$ and a family $(f_\beta )_{\beta \in \fr{P}(\bn)}$ in $\f{\bzz{2}}{E}$ such that for distinct $\beta,\gamma  \in \fr{P}(\bn)$, $\ssa{f_\beta }\not\approx \ssa{f_\gamma }$.
\end{enumerate}
\end{e}

Put
$$\mae{\alpha }{\bt}{Un\;\bc}{z}{z}\;.$$
Since $w(f(1,1)\alpha ^{-n})=0$, there is a $y\in Un\;E$ with $w(y)=0$ and $f(1,1)\alpha ^{-n}=y^2$.

a) If we put $x:=y\alpha ^{\frac{n}{2}}$ then $w(x)=\frac{n}{2}$ and $f(1,1)=x^2$ and the assertion follows from \pr \ref{785} c).

b) We put $x:=y\alpha ^{\frac{n-1}{2}}$. Then $f(1,1)=\alpha x^2$. Take $g\in \f{\bzz{2}}{E}$ with $g(1,1)=\alpha $ and $\lambda \in \Lambda (\bzz{2},E)$ with $(\delta \lambda )(1,1)=x^2$ (\pr \ref{785} a),b)). Then $f=g\delta \lambda $. By \ee{789}, $\ssa{g}$ is isomorphic to $E$ and by \pr \ref{716} $a_1\Rightarrow a_2$, $\ssa{f}$ is also isomorphic to $E$.

c) follows from \pr \ref{785} b) and \pr \ref{716} a),c).

d) Denote by $E$ the C*-direct product of the sequence $(\ccb{\bt}{\bc_{n,n}})_{n\in \bn}$ and for every $\beta \in \{0,1\}^{\bn}$ define $f_\beta \in \f{\bzz{2}}{E}$ by
$$\mae{f_\beta (1,1)}{\bn}{\un{E}}{n}{\alpha ^{\beta (n)}1_{\bc_{n,n}}}\;.$$
By a) and b), for distinct $\beta ,\gamma \in \{0,1\}^{\bn}$, $\ssa{f_\beta }\not\approx \ssa{f_\gamma }$ (\pr \ref{14} a)).\qed 

\begin{e}\label{797}
Let $I,J$ be finite disjoint sets and for all $i\in I\cup J$ and $j\in J$ put
$A_i:=B_j:=\bt$. We define the compact spaces $A$ and $B$ in the following way. For $A$ we take first the disjoint union of the spaces $A_i$ for all $i\in I\cup J$ and identify then the points $1\in A_i$ for all $i\in I\cup J$. For $B$ we take first the disjoint union of all the spaces $A_i$ for all $i\in I\cup J$ and of the spaces $B_j$ for all $j\in J$ and identify first the points $1\in A_i$ for all $i\in I\cup J$ and identify then also the points $-1\in A_i$ for all $i\in I$ and $1\in B_j$ for all $j\in J$.

Let $E:=\ccb{A}{\bc}$ and $f\in \f{\bzz{2}}{E}$ with
$$\mae{f(1,1)}{A}{\unn{\bc}}{z}{\ab{z}{z\in A_i\;\mbox{with}\;i\in I }{1}{z\in A_i\;with\;i\in J}}\;.$$
For every $X\in \ssa{f}$ define $\tilde X\in \ccb{B}{\bc}$ by
$$\mae{\tilde X}{B}{\bc}{z}{\ac{X_0(z^2)+zX_1(z^2)}{z\in A_i\;with\;i\in I}{X_0(z)+X_1(z)}{z\in A_i\;with\;i\in J}{X_0(z)-X_1(z)}{z\in B_j\;with\;j\in J}}\;.$$
Then the map
$$\mae{\varphi }{\ssa{f}}{\ccb{B}{\bc}}{X}{\tilde X}$$
is an isomorphism of C*-algebras.
\end{e}

Let $X,Y\in \ssa{f}$. By \h \ref{745} c),g),
$$(X^*)_0=(X_0)^*\,,\qquad (X^*)_1=\overline{f(1,1)}(X_1)^*\,,$$
$$(XY)_0=X_0Y_0+f(1,1)X_1Y_1\,,\qquad (XY)_1=X_0Y_1+X_1Y_0\;.$$
For $z\in A_i$ with $i\in I$,
$$\widetilde{X^*}(z)=(X^*)_0(z^2)+z(X^*)_1(z^2)=\overline{X_0(z^2)}+z\bar z^2\overline{X_1(z^2)}=$$
$$=\overline{X_0(z^2)+zX_1(z^2)}=(\tilde X)^*(z)\,,$$
$$\tilde X(z)\tilde Y(z)=(X_0(z^2)+zX_1(z^2))(Y_0(z^2)+zY_1(z^2))=$$
$$=X_0(z^2)Y_0(z^2)+zX_0(z^2)Y_1(z^2)+zX_1(z^2)Y_0(z^2)+z^2X_1(z^2)Y_1(z^2)=$$
$$=(XY)_0(z^2)+z(XY)_1(z^2)=\widetilde{XY}(z)\;.$$
For $z\in A_j$ or $z\in B_j$ with $j\in J$,
$$\widetilde{X^*}(z)=(X^*)_0(z)\pm (X^*)_1(z)=\overline{X_0(z)}\pm \overline{X_1(z)}=(\tilde X)^*(z)\,,$$
$$\tilde X(z)\tilde Y(z)=(X_0(z)\pm X_1(z))(Y_0(z)\pm Y_1(z))=$$
$$=X_0(z)Y_0(z)\pm X_0(z)Y_1(z)\pm X_1(z)Y_0(z)+X_1(z)Y_1(z)=$$
$$=(XY)_0(z)\pm (XY)_1(z)=\widetilde{XY}(z)\;.$$
Thus $\varphi $ is a C*-homomorphism. Assume $\tilde X=0$. For $z\in A_i$ with $i\in I$,
$$X_0(z^2)+zX_1(z^2)=0$$
so, successively,
$$X_0(z^2)-zX_1(z^2)=0\,,\qquad X_0(z^2)=X_1(z^2)=0\,,\qquad X(z)=0\;.$$
For $z\in A_j$ with $j\in J$,
$$\left\{\begin{array}{r@{=0}}
X_0(z)+X_1(z)\\X_0(z)-X_1(z)
\end{array}\right.\,,$$
so
$$X_0(z)=X_1(z)=0\,,\qquad X(z)=0\;.$$
Thus $\varphi $ is injective.

Let $x\in \ccb{B}{\bc}$ such that for every $i\in I$ there is a family $(c_{i,k})_{k\in \bz}\in \bc^{(\bz)}$ with
$$x(z)=\si{k\in \bz}c_{i,k}z^k$$
for all $z\in A_i$. Define $X_0,X_1\in E$ in the following way. If $z\in A_i$ with $i\in I$ we put
$$X_0(z):=\si{k\in \bz}c_{i,2k}z^k\,,\qquad X_1(z):=\si{k\in \bz}c_{i,2k+1}z^k\;.$$
If $z\in A_j$ with $j\in J$ then we put $z':=z\in B_j$,
$$X_0(z):=\frac{1}{2}(x(z)+x(z'))\,,\qquad X_1(z):=\frac{1}{2}(x(z)-x(z'))\;.$$
It is easy to see that $X_0$ and $X_1$ are well defined. Then
$$\tilde X(z)=\si{k\in \bz}c_{i,2k}z^{2k}+z\si{k\in \bz}c_{i,2k+1}z^{2k}=x(z)$$
for all $z\in A_i$ with $i\in I$ and $\tilde X(z)=x(z)$ for all $z\in A_j\cup B_j$ with $j\in J$. Since the elements $x$ of the above form are dense in $\ccb{B}{\bc}$, $\varphi $ is surjective.\qed

\begin{e}\label{792}
Let $E:=\ccb{\bt^2}{\bc}$ and $f,g\in \f{\bzz{2}}{E}$ with
$$\left\{\begin{array}{r}
\mae{f(1,1)}{\bt^2}{\unn{\bc}}{(z_1,z_2)}{z_1}\\
\mae{g(1,1)}{\bt^2}{\unn{\bc}}{(z_1,z_2)}{z_2}
\end{array}\right.\;.$$
Then the maps
$$\left\{\begin{array}{r}
\mad{\ssa{f}}{E}{X}{X_0(z_1^2,z_2)+z_1X_1(z_1^2,z_2)}\\
\mad{\ssa{g}}{E}{X}{X_0(z_1,z_2^2)+z_2X_1(z_1,z_2^2)}
\end{array}\right.$$
are isomorphisms of C*-algebras.\qed
\end{e}

{\it Remark.} $\ssa{f}$ and $\ssa{g}$ are isomorphic but not $E$-C*-isomorphic. 

\begin{e}\label{793}
Let $E:=\ccb{\bt^2}{\bc}$ and $f\in \f{\bzz{2}}{E}$ with
$$\mae{f(1,1)}{\bt^2}{Un\;\bc}{(z_1,z_2)}{z_1z_2}\;.$$
If we put
$$\mae{\tilde X}{\bt^2}{\bc}{(z_1,z_2)}{X_0(z_1^2,z_2^2)+z_1z_2X_1(z_1^2,z_2^2)}$$
for every $X\in\ssa{f}$ then the map 
$$\mae{\varphi }{\ssa{f}}{E}{X}{\tilde X}$$
is an injective unital C*-homomorphism with
$$\varphi (\ssa{f})=\ccc{G}:=\me{x\in E}{(z_1,z_2)\in \bt^2\Longrightarrow x(z_1,z_2)=x(-z_1,-z_2)}\;.$$
In particular $\ssa{f}$ is isomorphic to $E$.
\end{e}

Let $X,Y\in \ssa{f}$. By \h \ref{745} c),g),
$$(X^*)_0=(X_0)^*\,,\qquad (X^*)_1=\overline{f(1,1)}(X_1)^*\,,$$
$$(XY)_0=X_0Y_0+f(1,1)X_1Y_1\,,\qquad (XY)_1=X_0Y_1+X_1Y_0$$
so for $(z_1,z_2)\in \bt^2$,
$$\widetilde{X^*}(z_1,z_2)=X_0^*(z_1^2,z_2^2)+z_1z_2\bar z_1^2\bar z_2^2X_1^*(z_1^2,z_2^2)=$$
$$=\overline{X_0(z_1^2,z_2^2)+z_1z_2X_1(z_1^2,z_2^2)}=\overline{\tilde X(z_1,z_2)}\,,$$
$$(\tilde X(z_1,z_2))(\tilde Y(z_1,z_2))=$$
$$=(X_0(z_1^2,z_2^2)+z_1z_2X_1(z_1^2,z_2^2))(Y_0(z_1^2,z_2^2)+z_1z_2Y_1(z_1^2,z_2^2))=$$
$$=X_0(z_1^2,z_2^2)Y_0(z_1^2,z_2^2)+z_1z_2X_0(z_1^2,z_2^2)Y_1(z_1^2,z_2^2)+$$
$$+z_1z_2X_1(z_1^2,z_2^2)Y_0(z_1^2,z_2^2)+z_1^2z_2^2X_1(z_1^2,z_2^2)Y_1(z_1^2,z_2^2)=$$
$$=(XY)_0(z_1^2,z_2^2)+z_1z_2(XY)_1(z_1^2,z_2^2)=\widetilde{XY}(z_1,z_2)\,,$$
i.e. $\varphi $ is a unital C*-homomorphism. If $\tilde X=0$ then for $(z_1,z_2)\in \bt^2$,
$$X_0(z_1^2,z_2^2)+z_1z_2X_1(z_1^2,z_2^2)=0$$
so, successively,
$$X_0(z_1^2,z_2^2)-z_1z_2X_1(z_1^2,z_2^2)=0\,,\quad X_0(z_1^2,z_2^2)=X_1(z_1^2,z_2^2)=0\,,$$
$$X_0=X_1=0\,,\quad X=0$$
and $\varphi $ is injective.

The inclusion $\ssa{f}\subset \ccc{G}$ is obvious. Let $(a_{j,k})_{j,k\in \bz},\;(b_{j,k})_{j,k\in \bz}\in \bc^{(\bz\times \bz)}$ and
$$x=\si{j,k\in \bz}a_{j,k}z_1^{2j}z_2^{2k}+\si{j,k\in \bz}b_{j,k}z_1^{2j+1}z_2^{2k+1}\in \ccc{G}\;.$$
Define
$$X_0:=\si{j,k\in \bz}a_{j,k}z_1^jz_2^k\,,\qquad X_1:=\si{j,k\in \bz}b_{j,k}z_1^jz_2^k\;.$$
Then $\tilde X=x$. Since the elements of the above form are dense in $\ccc{G}$, $\varphi (\ssa{f})=\ccc{G}$.

If we consider the equivalence relation $\sim $ on $\bt^2$ defined by 
$$(z_1,z_2)\sim (w_1,w_2):\Longleftrightarrow z_1=-w_1,\; z_2=-w_2$$ 
then the quotient space $\bt^2/\!\sim $ is homeomorphic to $\bt^2$. Thus $\ssa{f}$ is isomorphic to $E$.\qed

\begin{e}\label{795}
Let $E:=\ccb{\bt^2}{\bc}$. 
\begin{enumerate}
\item For $x\in \unn{E}$ and $z\in \bt$, $w(x(\,\cdot\, ,z))$ and $w(x(z,\,\cdot \,))$ do not depend on $z$, where $w$ denotes the winding number \emph{(\dd \  \ref{787'})}.
\item If $x\in \unn{E}$ and if 
$$w(x(\,\cdot\, ,1))=w(x(1,\,\cdot\, ))=0$$
then there is a $y\in \unn{E}$ with $x=y^2$.
\item Let $f\in \f{\bzz{2}}{E}$ and put
$$\mae{\alpha }{\bt}{\bt^2}{z}{(z,1)}\,,\qquad\qquad \mae{\beta }{\bt}{\bt^2}{z}{(1,z)}\,,$$
$$m:=w(f(1,1)\circ \alpha )\,,\qquad n:=w(f(1,1)\circ \beta )\;.$$
\begin{enumerate}
\item If $m+n$ is odd then $\ssa{f}$ is isomorphic to $E$.	
\item If $m$ and $n$ are even then $\ssa{f}$ is isomorphic to $E\times E$.
\item If $m$ and $n$ are odd then $\ssa{f}$ is isomorphic to $E$.
\end{enumerate}
\item The group $\f{\bzz{2}}{E}/\!\Lambda (\bzz{2},E)$ is isomorphic to $\bzz{2}\times \bzz{2}$ and $$Card\,(\me{\ssa{f}}{f\in \ccc{F}(\bzz{2},E)}/\approx _{\ccc{S}})=4\;.$$
\end{enumerate}
\end{e}

a) follows by continuity.

b) follows from a).

c) Let $g\in \f{\bzz{2}}{E}$ with
$$\mae{g(1,1)}{\bt^2}{\unn{\bc}}{(z_1,z_2)}{z_1^mz_2^n}\;.$$
Then
$$w(g(1,1)\circ \alpha )=m\,,\qquad w(g(1,1)\circ \beta )=n\;.$$
By b), there is an $x\in \unn{E}$ with $f(1,1)=x^2g(1,1)$. By \pr \ref{785} b) and \pr \ref{716} $a_1\Rightarrow a_2$, $\ssa{f}\approx \ssa{g}$.

$c_1)$ Assume $m$ even and put
$$\mae{y}{\bt^2}{\unn{\bc}}{(z_1,z_2)}{z_1^{\frac{m}{2}}z_2^{\frac{n-1}{2}}}\,.$$
If $h\in \f{\bzz{2}}{E}$ with 
$$\mae{h(1,1)}{\bt^2}{\unn\bc}{(z_1,z_2)}{z_2}$$
then $g(1,1)=y^2h(1,1)$. By \pr \ref{785} b) and \pr \ref{716} $a_1\Rightarrow a_2$, $\ssa{g}\approx \ssa{h}$ and by \ee{792} $a_1\Rightarrow a_2$, $\ssa{h}\approx E$. Thus $\ssa{f}\approx E$.

$c_2)$ If we put 
$$\mae{y}{\bt^2}{\unn{\bc}}{(z_1,z_2)}{z_1^{\frac{m}{2}}z_2^{\frac{n}{2}}}$$
then $g(1,1)=y^2$ and the assertion follows from \pr \ref{785} c).

$c_3)$ We put
$$\mae{y}{\bt^2}{\unn{\bc}}{(z_1,z_2)}{z_1^{\frac{m-1}{2}}z_2^{\frac{n-1}{2}}}$$
and take $h\in \f{\bzz{2}}{E}$ with 
$$\mae{h(1,1)}{\bt^2}{\unn\bc}{(z_1,z_2)}{z_1z_2}$$
then $g(1,1)=y^2h(1,1)$ so by \pr \ref{785} b) and \pr \ref{716} $a_1\Rightarrow a_2$, $\ssa{g}\approx \ssa{h} $. By \ee{793} $\ssa{h}\approx E$, so $\ssa{f}\approx E$.

d) follows from b), \pr \ref{785} b), and \pr \ref{716} a),c).\qed

{\it Remark.} In a similar way it is possible to show that for every $n\in \bn$, $\f{\bzz{2}}{\bt^n}/\Lambda (\bzz{2},\bt^n)$ is isomorphic to $(\bzz{2})^n$ and
$$Card\,(\me{\ssa{f}}{f\in \f{\bzz{2}}{\bt^n}}/\approx _{\ccc{S}})=2^n\;.$$

\begin{e}\label{800}
Let $I,J,K$ be finite pairwise disjoint sets and for every $i\in I\cup J\cup K$ and $k\in K$ put $A_i:=B_k:=\bt^2$. We define the compact spaces $A$ and $B$ in the following way. For $A$ we take first the disjoint union of the spaces $A_i$ with $i\in I\cup J\cup K$ and then identify the points $(1,1)\in A_i$ for all $i\in I\cup J\cup K$. For $B$ we take first the disjoint union of the spaces $A_i$ with $i\in I\cup J\cup K$ and of the spaces $B_k$ with $k\in K$. Then we identify the points $(1,1)\in A_i$ for all $i\in I\cup J\cup K$ and then we identify for every $j\in J$ the points $(z_1,z_2)\in A_j$ with the points $(-z_1,-z_2)\in A_j$ and finally we identify the points $(-1,1)\in A_i$ for all $i\in I\cup J$ with the points $(1,1)\in B_k$ for all $k\in K$.

Let $E:=\ccb{A}{\bc}$ and $f\in \f{\bzz{2}}{A}$ such that
$$\mae{f(1,1)}{A}{\unn{\bc}}{(z_1,z_2)}{\ac{z_1}{(z_1,z_2)\in A_i\;with\;i\in I}{z_1z_2}{(z_1,z_2)\in A_i\;with\;i\in J}{1}{(z_1,z_2)\in A_i\;with\;i\in K}}\;.$$
We define for every $X\in \ssa{f}$ a map $\tilde X:B\rightarrow \bc$ by
$$(z_1,z_2)\mapsto {\ad{X_0(z_1^2,z_2)+z_1X_1(z_1^2,z_2)}{(z_1,z_2)\in A_i\;with\;i\in I}{X_0(z_1^2,z_2^2)+z_1z_2X_1(z_1^2,z_2^2)}{(z_1,z_2)\in A_i\;with\;i\in J}{X_0(z_1,z_2)+X_1(z_1,z_2)}{(z_1,z_2)\in A_i\;with\;i\in K}{X_0(z_1,z_2)-X_1(z_1,z_2)}{(z_1,z_2)\in B_k\;with\;k\in K}}\;.$$ 
Then the map
$$\mad{\ssa{f}}{\ccb{B}{\bc}}{X}{\tilde X}$$
is an isomorphism of C*-algebras.
\end{e}

The proof is similar to the proof of \ee{797}.\qed

\begin{e}\label{8}
If $n\in \bn$, $E:=\ccb{\bt^n}{\bc}$, and $f\in \f{\bzz{2}}{\ccb{\bt^n}{\bc}}$ then $\ssa{f}$ is isomorphic either to $\ccb{\bt^n}{\bc}$ or to $\ccb{\bt^n}{\bc}\times \ccb{\bt^n}{\bc}$.\qed
\end{e}

\begin{e}\label{826}
Assume $E:=\ccb{A}{\bc}$, where $A$ denotes Moebius's band (resp. Klein's bottle), i.e. the topological space obtained from $[0,2\pi ]\times [-\pi ,\pi ]$ by identifying the points $(0,\alpha )$ and $(2\pi ,-\alpha )$ for all $\alpha \in [-\pi ,\pi ]$ (resp. and the points $(\theta ,-\pi )$ and $(\theta ,\pi )$ for all $\theta \in [0,2\pi ]$). We put $B:=\bt\times [-\pi ,\pi ]$ (resp. $B:=\bt^2$) and
$$\mae{\tilde x}{[0,2\pi ]\times [-\pi ,\pi ]}{\bc}{(\theta ,\alpha )}{\ab{x(2\theta ,\alpha )}{\theta \in [0,\pi ]}{x(2(\theta -\pi ),-\alpha )}{\theta \in [\pi ,2\pi ]}}$$
for every $x\in E$.
\begin{enumerate}
\item $\tilde x$ is well-defined and belongs to $\ccb{B}{\bc}$ for every $x\in E$.
\item If $f_{1,1}(\theta ,\alpha )=e^{i\theta }$ for all $(\theta ,\alpha )\in [0,2\pi ]\times [-\pi ,\pi ]$ then the map
$$\mae{\varphi }{\ssa{f}}{\ccb{B}{\bc}}{X}{\widetilde{ X_0}+e^{i\theta }\widetilde {X_1}}$$
is a C*-isomorphism.
\item Let $x\in \unn{E}$. If $w(x(\,\cdot \,,0))=0$ (where $w$ denotes the winding number) then there is a $y\in E$ with $e^y=x$.
\item Let $x\in \unn{E}$ and put $n:=w(x(\,\cdot \,,0))$. Then there is a $y\in E$ with $e^y=e^{-in\theta }x$.
\item The group $\f{\bzz{2}}{A}/\Lambda (\bzz{2},A)$ is isomorphic to $\bzz{2}$.
\item If $w(f_{1,1}(\,\cdot \,,0))$ is even (resp. odd) then $\ssa{f}$ is isomorphic to $E\times E$ (resp. to $\ccb{B}{\bc}$).
\end{enumerate}
\end{e}

a) For $\alpha \in [-\pi ,\pi ]$,
$$\tilde x(\pi ,\alpha )=x(2\pi ,\alpha )=x(0,-\alpha )=\tilde x(\pi ,\alpha )$$
so $\tilde x$ is well-defined. Moreover
$$\tilde x(0,\alpha )=x(0,\alpha )=x(2\pi ,-\alpha )=\tilde x(2\pi ,\alpha )$$
and in the case of Klein's bottle
$$\ab{\tilde x(\theta ,-\pi )=x(2\theta ,-\pi )=x(2\theta ,\pi )=\tilde x(\theta ,\pi )}{\theta \in [0,\pi ]}{\tilde x(\theta ,-\pi )=x(2(\theta -\pi ),\pi )=x(2(\theta -\pi ),-\pi )=\tilde x(\theta ,\pi )}{\theta \in [\pi ,2\pi ]}$$
i.e. $\tilde x\in \ccb{B}{\bc}$.

b) For $X,Y\in \ssa{f}$ and $(\theta ,\alpha )\in [0,2\pi ]\times [-\pi ,\pi ]$, by \h \ref{745} c),g),
$$(\varphi X^*)(\theta ,\alpha )=\widetilde{(X^*)_0}(\theta ,\alpha )+e^{i\theta }\widetilde{(X^*)_1}(\theta ,\alpha )=$$
$$=\widetilde{(X_0)^*}(\theta ,\alpha )+e^{i\theta }\widetilde{\overbrace{(e^{-i\theta }(X_1)^*)}}(\theta ,\alpha )=$$
$$=\ab{\overline{X_0(2\theta ,\alpha )}+e^{i\theta }(e^{-2i\theta }\overline{X_1(2\theta ,\alpha )})}{\theta \in [0,\pi ]}{\overline{X_0(2(\theta -\pi ),-\alpha )}+e^{i\theta }(e^{-2i(\theta -\pi )}\overline{X_1(2(\theta -\pi ),-\alpha )})}{\theta \in [\pi ,2\pi ]}=$$
$$=\ab{\overline{X_0(2\theta ,\alpha )+e^{i\theta }X_1(2\theta ,\alpha )}}{\theta \in [0,\pi ]}{\overline{X_0(2(\theta -\pi ),-\alpha )+e^{i\theta }X_1(2(\theta -\pi ),-\alpha )}}{\theta \in [\pi ,2\pi ]}=\overline{\varphi X}(\theta ,\alpha )\,,$$
$$(\varphi X)(\varphi Y)=(\widetilde{X_0}+e^{i\theta }\widetilde{X_1})(\widetilde{Y_0}+e^{i\theta }\widetilde{Y_1})=\widetilde{X_0}\widetilde{Y_0}+e^{i\theta }\widetilde{X_0}\widetilde{Y_1}+e^{i\theta }\widetilde{X_1}\widetilde{Y_0}+e^{2i\theta }\widetilde{X_1}\widetilde{Y_1}\,,$$
$$\varphi (XY)=\widetilde{(XY)_0}+e^{i\theta }\widetilde{(XY)_1}=$$
$$=\widetilde{X_0}\widetilde{Y_0}+e^{2i\theta }\widetilde{X_1}\widetilde{Y_1}+e^{i\theta }(\widetilde{X_0}\widetilde{Y_1}+\widetilde{X_1}\widetilde{Y_0})=(\varphi X)(\varphi Y)\,,$$
i.e. $\varphi $ is a C*-homomorphism. If $\varphi X=0$ then for $\alpha \in [-\pi ,\pi ]$,
$$\ab{X_0(2\theta ,\alpha )+e^{i\theta }X_1(2\theta ,\alpha )=0}{\theta \in [0,\pi ]}{X_0(2(\theta -\pi ),-\alpha )+e^{i\theta }X_1(2(\theta -\pi ),-\alpha )=0}{\theta \in [\pi ,2\pi ]}$$
so for $\theta \in [0,\pi ]$, replacing $\theta $ by $\theta +\pi $ and $\alpha $ by $-\alpha $ in the second relation,
$$X_0(2\theta ,\alpha )-e^{i\theta }X_1(2\theta ,\alpha )=0\;.$$
It follows successively
$$X_0(2\theta ,\alpha )=X_1(2\theta ,\alpha )=0\,,$$
$$X_0=X_1=0\,,\qquad X=0\;.$$
Thus $\varphi $ is injective.

Let $y\in \ccb{B}{\bc}$. Put
$$\left\{
\begin{array}{c}
\mae{X_0}{[0,2\pi ]\times [-\pi ,\pi ]}{\bc}{(\theta ,\alpha )}{\frac{1}{2}(y(\frac{\theta }{2},\alpha )+y(\frac{\theta }{2}+\pi ,-\alpha ))}\\
\mae{X_1}{[0,2\pi ]\times [-\pi ,\pi ]}{\bc}{(\theta ,\alpha )}{\frac{1}{2}e^{-i\frac{\theta }{2}}(y(\frac{\theta }{2},\alpha )-y(\frac{\theta }{2}+\pi ,-\alpha ))}
\end{array}\;.
\right.$$
For $\alpha \in [-\pi ,\pi ]$,
$$\left\{
\begin{array}{c}
X_0(0,\alpha )=\frac{1}{2}(y(0,\alpha )+y(\pi ,-\alpha ))\\
X_0(2\pi ,-\alpha )=\frac{1}{2}(y(\pi ,-\alpha )+y(2\pi ,\alpha ))
\end{array}
\right.$$
$$\left\{
\begin{array}{c}
X_1(0,\alpha )=\frac{1}{2}(y(0,\alpha )-y(\pi ,-\alpha ))\\
X_1(2\pi ,-\alpha )=-\frac{1}{2}(y(\pi ,-\alpha )-y(2\pi ,\alpha ))
\end{array}
\right.$$
so $X_0,X_1\in E$. Moreover for $(\theta ,\alpha )\in [0,2\pi ]\times [-\pi ,\pi ]$,
$$\widetilde{X_0}(\theta ,\alpha )+e^{i\theta }\widetilde{X_1}(\theta ,\alpha )=$$
$$=\ab{X_0(2\theta ,\alpha )+e^{i\theta }X_1(2\theta ,\alpha )}{\theta \in [0,\pi ]}{X_0(2(\theta -\pi ),-\alpha )+e^{i\theta }X_1(2(\theta -\pi ),-\alpha )}{\theta \in [\pi ,2\pi ]}=$$
$$=\ab{\frac{1}{2}(y(\theta ,\alpha )+y(\theta +\pi ,-\alpha )+y(\theta ,\alpha )-y(\theta +\pi ,-\alpha ))}{\theta \in [0,\pi ]}{\frac{1}{2}(y(\theta -\pi ,-\alpha )+y(\theta ,\alpha )-y(\theta -\pi ,-\alpha )+y(\theta ,\alpha ))}{\theta \in [\pi ,2\pi ]}=$$
$$=y(\theta ,\alpha )$$
i.e. $\varphi $ is surjective.

c) If $A$ is Moebius's band then the assertion is obvious so assume $A$ is Klein's bottle. The winding numbers of 
$$\left\{\begin{array}{r}
\mad{[0,2\pi ]}{\bc}{\alpha }{x(0,\alpha )}\\\mad{[0,2\pi ]}{\bc}{\alpha }{x(2\pi ,\alpha )}
\end{array}\right.$$
are equal by homotopy, but their sum is equal to 0. Thus these winding numbers are equal to 0. The paths $\theta $ and $\alpha $ on $A$ generate the homotopy group of $A$. Thus the winding number of $x$ on any path of $A$ is 0 and the assertion follows.

d) The winding number of
$$\mad{[0,2\pi ]}{\bc}{\theta }{e^{-in\theta }x(\theta ,0)}$$
is 0 and the assertion follows from c).

e) The assertion follows from d) and \pr \ref{785} b).

f) The assertion follows from b), d), \pr \ref{716} $a_1\Rightarrow a_2$, and \pr \ref{785} c).\qed
\begin{center}
\subsection{$T:=\bzz{2}\times \bzz{2}$}
\end{center}
\begin{p}\label{807}
Let $E$ be a unital C*-algebra and let $a,b,c$ be the three elements of $(\bzz{2}\times \bzz{2})\setminus \{(0,0)\}$. Put
$$A:=\me{(\alpha ,\beta ,\gamma ,\varepsilon )\in (\un{E})^4}{\varepsilon ^2=1_E}$$
and for every $\varrho  \in A$ and $\sigma \in (\un{E})^3$ denote by $f_\varrho $ and $g_\sigma $ the functions defined by the following tables:

$$\begin{array}{||c||c|c|c||}\hline \hline
f_\varrho &a&b&c\\\hline\hline a&\beta \gamma &\gamma &\beta \\\hline b&\varepsilon \gamma &\varepsilon \alpha \gamma &\alpha \\\hline c&\varepsilon \beta &\varepsilon \alpha &\alpha \beta \\\hline\hline
\end{array}\qquad\qquad
\begin{array}{||c||c|c|c||}\hline\hline
g_\sigma &a&b&c\\\hline\hline a&\alpha ^2&\alpha \beta \gamma ^*&\alpha \gamma \beta ^*\\\hline b&\alpha \beta \gamma ^*&\beta ^2&\beta \gamma \alpha ^*\\\hline c&\alpha \gamma \beta ^*&\beta \gamma \alpha ^*&\gamma ^2\\\hline\hline
\end{array}$$
\begin{enumerate}
\item $f_\varrho \in \f{\bzz{2}\times \bzz{2}}{E}$ for every $\varrho \in A$ and the map
$$\mad{A}{\f{\bzz{2}\times \bzz{2}}{E}}{\varrho }{f_\varrho }$$
is bijective.
\item $g_\sigma \in \me{\delta \lambda }{\lambda \in \Lambda (\bzz{2}\times \bzz{2},E)}$ for every $\sigma \in (\un{E})^3$ and the map
$$\mad{(\un{E})^3}{\me{\delta \lambda }{\lambda \in \Lambda (\bzz{2}\times \bzz{2},E)}}{\sigma }{g_\sigma }$$
is bijective.
\item The following are equivalent for all $\varrho :=(\alpha ,\beta ,\gamma ,\epsilon )\in A$ and $\varrho ':=(\alpha ',\beta ',\gamma',\epsilon ')\in A$:
\begin{enumerate}
\item $\ssa{f_\varrho }\approx _\ccc{S}\ssa{f_{\varrho '}}$.
\item $\varepsilon =\varepsilon '$ and there are $x,y,z\in \un{E}$ with
$$x^2=\beta \beta '^*\gamma \gamma '^*\,,\qquad y^2=\alpha \alpha '^*\gamma \gamma '^*\,,\qquad z^2=\alpha \alpha '^*\beta \beta '^*\;.$$
\item $\varepsilon =\varepsilon '$ and there are $x,y\in \un{E}$ with
$$x^2=\beta \beta '^*\gamma \gamma '^*\,,\qquad y^2=\alpha \alpha '^*\gamma \gamma '^*\;.$$
\end{enumerate}
\item The following are equivalent for all $\varrho :=(\alpha ,\beta ,\gamma ,\varepsilon \in A)$ and $X\in \ssa{f_\varrho }$:
\begin{enumerate}
\item $X\in \me{V_t^{f_\varrho }}{t\in \bzz{2}\times \bzz{2}}^c$.
\item $t\in \bzz{2}\times \bzz{2}\Longrightarrow \varepsilon X_t=X_t$.
\end{enumerate}
\item The following are equivalent for all $\varrho :=(\alpha ,\beta ,\gamma ,\varepsilon \in A)$ and $X\in \ssa{f_\varrho }$:
\begin{enumerate}
\item $X\in \ssa{f_\varrho }^c$.
\item $t\in \bzz{2}\times \bzz{2}\Longrightarrow \varepsilon X_t=X_t\in E^c$.
\end{enumerate}
\item For $\varrho :=(\alpha ,\beta ,\gamma ,\varepsilon )\in A$ and $X,Y\in \ssa{f_\varrho }$,
$$(X^*)_0=X_0^*\,,\; (X^*)_a=\beta ^*\gamma ^*X_a^*\,,\; (X^*)_b=\varepsilon \alpha ^*\gamma ^*X_b^*\,,\; (X^*)_c=\alpha ^*\beta ^*X_c^*\,,$$
$$(XY)_0=X_0Y_0+\beta \gamma X_aY_a+\varepsilon \alpha \gamma X_bY_b+\alpha \beta X_cY_c\,,$$
$$(XY)_a=X_0Y_a+X_aY_0+\alpha X_bY_c+\varepsilon \alpha X_cY_b\,,$$
$$(XY)_b=X_0Y_b+\beta X_aY_c+X_bY_0+\varepsilon \beta X_cY_a\,,$$
$$(XY)_c=X_0Y_c+\gamma X_aY_b+\varepsilon \gamma X_bY_a+X_cY_0\;.$$
\item Assume $\bk=\bc$, let $\sigma (E^c)$ be the spectrum of $E^c$, and for every $\delta \in E^c$ let $\hat \delta $ be its Gelfand transform. Then
$$\sigma (V_a)=\me{e^{i\theta }}{\theta \in \br,\;e^{2i\theta }\in \widehat{\beta \gamma }(\sigma (E^c))}\,,$$
$$\sigma (V_b)=\me{e^{i\theta }}{\theta \in \br,\;e^{2i\theta }\in \widehat{\alpha \gamma }(\sigma (E^c))}\,,$$
$$\sigma (V_c)=\me{e^{i\theta }}{\theta \in \br,\;e^{2i\theta }\in \widehat{\alpha \beta}(\sigma (E^c))}\;.$$
\end{enumerate}
\end{p}

a) is a long calculation.

b) is easy to verify.

$c_1\Rightarrow c_2$ By \pr \ref{716} $a_2\Rightarrow a_1$ there is a $\lambda \in \Lambda (\bzz{2}\times \bzz{2},E)$ with $f_\varrho =f_{\varrho '}\delta \lambda $. By b), there is a $\sigma :=(x,y,z)\in (\un{E})^3$ with $f_\varrho =f_{\varrho '}g_\sigma $. We get $\varepsilon =\varepsilon '$ and 
$$\alpha \alpha '^*=x^*yz\,,\qquad \beta \beta '^*=xy^*z\,,\qquad \gamma \gamma '^*=xyz^*\;.$$
It follows $xyz=\alpha \alpha '^*\beta \beta '^*\gamma \gamma '^*$ so
$$x^2=\beta \beta '^*\gamma \gamma '^*\,,\qquad y^2=\alpha \alpha '^*\gamma \gamma '^*\,,\qquad z^2=\alpha \alpha '^*\beta \beta '^*\;.$$

$c_2\Rightarrow c_3$ is trivial.

$c_3\Rightarrow c_2$ If we put $z:=xy\gamma ^*\gamma '$ then 
$$z^2=\beta \beta '^*\gamma \gamma '^*\alpha \alpha '^*\gamma \gamma '^*\gamma ^{*2}\gamma '^2=\alpha \alpha '^*\beta \beta '^*\;.$$

$c_2\Rightarrow c_1$ follows from b) and \pr \ref{716} $a_1\Rightarrow a_2$.

d) follows from Corollary \ref{752} b).

e) follows from Corollary \ref {752} c).

f) follows from \h \ref{745} c),g).

g) follows from f).\qed

\begin{co}\label{809}
We use the notation of \emph{\pr \ref{807}} and take $\varrho :=(\alpha ,\beta ,\gamma ,\varepsilon )\in A$.
\begin{enumerate}
\item Assume $\varepsilon =1_E$ and there are $x,y\in \unn{E}$ with $x^2=\beta \gamma \,,y^2=\alpha \gamma $. Put $z:=xy\gamma^*$.
\begin{enumerate}
\item $x,y,z\in \un{E}\,,\;z^2=\alpha \beta $.
\item For every $\lambda ,\mu \in \{-1,1\}$ the map
$$\mae{\varphi _{\lambda ,\mu }}{\ssa{f_\varrho }}{E}{X}{X_0+\lambda xX_a+\mu yX_b+\lambda \mu zX_c}$$
is an $E$-C*-homomorphism.
\item The map
$$\mad{\ssa{f_\varrho }}{E^4}{X}{(\varphi _{1,1}X,\varphi _{1,-1}X,\varphi _{-1,1}X.\varphi _{-1,-1}X)}$$
is an $E$-C*-isomorphism.
\end{enumerate}
\item Assume $\bk:=\br\,,\varepsilon =1_E$, and there are $x,y\in \unn{E}$ with 
$$x^2=-\beta \gamma \,,\quad y^2=\alpha \gamma \,,\quad (resp.\;y^2=-\alpha \gamma  )\;.$$
Put $z:=xy\gamma ^*$. Then $x,y,z\in \un{E}$, $z^2=-\alpha \beta $ (resp. $z^2=\alpha \beta $), and the maps
$$\mad{\ssa{f_\varrho }}{(\stackrel{\circ }{E})^2}{X}{(X_0+ixX_a+yX_b+izX_c,X_0+ixX_a-yX_b-izX_c)}$$
$$\mad{\ssa{f_\varrho }}{(\stackrel{\circ }{E})^2}{X}{(X_0+ixX_a+iyX_b-zX_c,X_0+ixX_a-iyX_b+zX_c)}$$
are respectively $E$-C*-isomorphisms (where $\stackrel{\circ }{E}$ denotes the complexification of $E$).
\item Assume $\bk:=\br$, $\varepsilon =-1_E$, and there are $x,y\in E^c$ with $x^2=-\beta \gamma \,,\;y^2=\alpha \gamma $. Put $z:=xy\gamma ^*$. Then $x,y,z\in \un{E}$, $z^2=-\alpha \beta $, and the map
$$\mad{\ssa{f_\varrho }}{\bh\otimes E}{X}{X_0+ixX_a+jyX_b+kzX_c}\,,$$
where $i,j,k$ are the canonical units of $\bh$, is an $E$-C*-isomorphism.
\item If $\varepsilon =-1_E$ and there is an $x\in \un{E}$ with $x^2=\alpha \beta $ then for every $\delta \in \un{E}$ the map
$$\mad{\ssa{f_\varrho }}{E_{2,2}}{X}{\mt{X_0+xX_c}{\gamma \delta ^*(\beta X_a-xX_b)}{\delta (X_a+x\beta ^* X_b)}{X_0-xX_c}}$$
is an $E$-C*-isomorphism.
\end{enumerate}
\end{co}

The proof is a long calculation using \pr \ref{807} f).\qed

{\it Remarks.} d) is contained in \pr \ref{10} c). An example with $\varepsilon =1_E$ but different from a) is presented in \pr \ref{803}. 

\begin{p}\label{10}
We use the notation of \emph{\pr \ref{807}} and take $\varrho :=(\alpha ,\beta ,\gamma ,\varepsilon )\in A$. 
\begin{enumerate}
\item Let $\varphi :\ssa{f_\varrho }\rightarrow E_{2,2}$ be an $E$-C*-isomorphism and put
$$\mt{A_t}{B_t}{C_t}{D_t}:=\varphi V_t$$
for every $t\in \bzz{2}\times \bzz{2}\setminus \{(0,0)\}$. Then $\varepsilon =-1_E$, $A_t,B_t,C_t,D_t\in E^c$ and $A_t+D_t=0$ for every $t\in \bzz{2}\times \bzz{2}\setminus \{(0,0)\}$, and
$$A_a^*=\beta ^*\gamma ^*A_a\,,\qquad A_b^*=-\alpha ^*\gamma ^*A_b\,,\qquad A_c^*=\alpha ^*\beta ^*A_c\,,$$
$$B_a^*=\beta ^*\gamma ^*C_a\,,\qquad B_b^*=-\alpha^*\gamma ^*C_b\,,\qquad B_c^*=\alpha ^*\beta ^*C_c\,,$$
$$A_a^2+B_aC_a=\beta \gamma \,,\qquad A_b^2+B_bC_b=-\alpha \gamma \,,\qquad A_c^2+B_cC_c=\alpha \beta \,,$$ 
$$A_a^2=\beta \gamma (1_E-|B_a|^2)\,,\quad A_b^2=-\alpha \gamma (1_E-|B_b|^2)\,,\quad A_c^2=\alpha \beta (1_E-|B_c|^2)\,,$$
$$2A_aA_b+B_aC_b+B_bC_a=0\,,\qquad 2A_bA_c+B_bC_c+B_cC_b=0\,,$$
$$2A_cA_a+B_cC_a+B_aC_c=0\,,$$
$$\alpha A_a=A_bA_c+B_bC_c\,,\qquad \alpha B_a=A_bB_c-A_cB_b\,,\qquad \alpha C_a=A_cC_b-A_bC_c\,,$$
$$\beta A_b=A_aA_c+B_aC_c\,,\qquad \beta B_b=A_aB_c-A_cB_a\,,\qquad \beta C_b=A_cC_a-A_aC_c\,,$$
$$\gamma A_c=A_aA_b+B_aC_b\,,\qquad \gamma B_c=A_aB_b-A_bB_a\,,\qquad \gamma C_c=A_bC_a-A_aC_b\,,$$
$$|A_a|+|A_b|+|A_c|\not=0\,,\qquad |B_a|+|B_b|+|B_c|\not=3.1_E\;.$$
\item Let $(A_t)_{t\in T}$, $(B_t)_{t\in T}$, $(C_t)_{t\in T}$, $(D_t)_{t\in T}$ be families in $E^c$ satisfying the above conditions and put 
$$X':=A_aX_a+A_bX_b+A_cX_c\,,\qquad X'':=B_aX_a+B_bX_b+B_cX_c\,,$$
$$X''':=C_aX_a+C_bX_b+C_cX_c$$
for every $X\in \ssa{f_\varrho }$. If $\varepsilon =-1_E$ then the map
$$\mad{\ssa{f_\varrho }}{E_{2,2}}{X}{\mt{X_0+X'}{X''}{X'''}{X_0-X'}}$$
is an $E$-C*-isomorphism.
\item Let $\varepsilon =-1_E$ and assume there is an $x\in E^c$ with $x^2=\beta \gamma $. Let $y\in \un{E}$ and put $z:=\gamma ^*xy$. Then $x,y,z\in \un{E}$ and the map
$$\mae{\varphi }{\ssa{f_\varrho }}{E_{2,2}}{X}{\mt{X_0+xX_a}{\alpha (yX_b+zX_c)}{-\gamma y^*X_b+\beta z^*X_c}{X_0-xX_a}}$$
is an $E$-C*-isomorphism such that 
$$\varphi( \frac{1}{2}(V_0+(x^*\otimes 1_K)V_a))=\mt{1}{0}{0}{0}\;.$$
In particular (by the symmetry of a,b,c), if $\varepsilon =-1_E$ and if there is an $x\in E^c$ with $x^2=\beta \gamma $, or $x^2=-\alpha \gamma $, or $x^2=\alpha \beta $ then $ \ssa{f_\varrho }\approx _E E_{2,2}$.\qed
\end{enumerate}
\end{p}

{\it Remark.} Take $\varrho :=(1_E,1_E,1_E,-1_E)$, $\varrho ':=(1_E,1_E,\gamma ',-1_E)$. By c), $\ssa{f_\varrho }\approx _E\ssa{f_{\varrho '}}$ and by \pr \ref{807} $c_1\Rightarrow c_2$, $\ssa{f_\varrho }\approx _{\ccc{S}}\ssa{f_{\varrho '}}$ implies the existence of an $x\in \un{E}$ with $x^2=\gamma '$.

\begin{co}\label{887}
We use the notation of \emph{\pr \ref{10}} and take $E:=\bk$, $\alpha =1$, and $\beta =\gamma =\varepsilon =-1$. Let $S$ be a group, $F$ a unital C*-algebra, $g\in \f{S}{F}$, and
$$\mae{h}{((S\times (\bzz{2})^2)\times (S\times (\bzz{2})^2))}{\un{F}}
{((s_1,t_1),(s_2,t_2))}$$
$${f_\varrho (t_1,t_2)g(s_1,s_2)}\;.$$
\begin{enumerate}
\item $h\in \f{S\times (\bzz{2})^2}{F}$.
\item $\ssa{h}\approx \ssa{g}_{2,2}\,,\qquad \ssb{\n{\cdot }}{h}\approx \ssb{\n{\cdot }}{g}_{2,2}$.
\end{enumerate}
\end{co}

By \pr \ref{10} c), $\ssa{f}\approx \bk_{2,2}$, so by \pr \ref{781} c),e),
$$\ssa{h}\approx  \bk_{2,2}\otimes \ssa{g}\approx \ssa{g}_{2,2}\,,\qquad \ssb{\n{\cdot }}{h}\approx \bk_{2,2}\otimes \ssb{\n{\cdot }}{g}\approx \ssb{\n{\cdot }}{g}_{2,2}\;.\qedd$$

\begin{e}\label{829}
Let $\bk:=\bc$ and $E:=\ccb{\bt}{\bc}$.
\begin{enumerate}
\item With the notation of \emph{\pr \ref{807}}, if $\varrho :=(\alpha ,\beta ,\gamma ,-1)\in A$ then $\ssa{f_\varrho }\approx_E E_{2,2}$.
\item $Card\,(\me{\ssa{f}}{f\in \f{\bzz{2}\times \bzz{2}}{E}}/\approx _{\ccc{S}})=16$.
\end{enumerate}
\end{e}

Put
$$m:=w(\alpha )\,,\qquad n:=w(\beta )\,,\qquad p:=w(\gamma )\,,$$
where $w$ denotes the winding number. By \pr \ref{716} $a_1\Rightarrow a_2$, we may assume $\alpha =z^m,\,\beta =z^n,\,\gamma =z^p$. 

a) If $n+p$ is even then the assertion follows from \pr \ref{10} c). If $n+p$ is odd then either $m+p$ or $m+n$ is even and the assertion follows again from \pr \ref{10} c).

b) follows from \pr \ref{716} a),c).\qed

{\it Remark.} Assume $\bk:=\br$ and let $E$ be the real C*-algebra $\ccb{\bt}{\bc}$ ([C1] \h 4.1.1.8 a)), $\varepsilon =-1_E$, 
$$\mae{\alpha }{\bt}{\bc}{z}{z}\,,$$
$$\mae{\beta }{\bt}{\bc}{z}{-z}\,,$$
$$\mae{\gamma }{\bt}{\bc}{z}{\bar z}\,,$$
and $\varrho :=(\alpha ,\beta ,\gamma ,\varepsilon )$. Then by Corollary \ref{809} c), $\ssa{f_\varrho }\approx \bh\otimes E$.

\begin{e}\label{11}
We put $E:=\ccb{\bt^2}{\bc}$, $\gamma :=1_E$,
$$\mae{\alpha }{\bt^2}{\bc}{(z_1,z_2)}{z_1}\,,\qquad \qquad \mae{\beta }{\bt^2}{\bc}{(z_1,z_2)}{z_2}\,,$$
and (with the notation of \emph{\pr \ref{807}}) $\varrho :=(\alpha ,\beta ,\gamma ,-1_E)\in A$.
\begin{enumerate}
\item $\ssa{f_\varrho }$ is not commutative and not $E$-C*-isomorphic to $E_{2,2}$.
\item If we put
$$\mae{\tilde x}{\bt^2}{\bc}{(z_1,z_2)}{x(z_1^2,z_2^2)}$$
for every $x\in E$ then the map
$$\mad{\ssa{f_\varrho }}{E_{2,2}}{X}{\mt{\tilde X_0+\alpha \beta \tilde X_c}{\beta \tilde X_a-\alpha \tilde X_b}{\beta \tilde X_a+\alpha \tilde X_b}{\tilde X_0-\alpha \beta \tilde X_c}}$$
is a C*-isomorphism.
\item $E_{2,2}\approx \ssa{f_\varrho }\not\approx _EE_{2,2}$.
\end{enumerate}
  \end{e}

a) By \pr \ref{807} d), $\ssa{f_\varrho }$ is not commutative. Assume $\ssa{f_\varrho }\approx _EE_{2,2}$ and let us use the notation of \pr \ref{10} a).
\begin{center}
Step 1 $\{A_a\not=0\}\subset \{A_b=0\}$
\end{center}

Assume $\{A_a\not=0\}\cap \{A_b\not=0\}\not=\emptyset $. By \pr \ref{10} a), 
$$2A_aA_b+B_aC_b+B_bC_a=0\,,\qquad B_a^*=\beta ^*C_a\,,\qquad B_b^*=-\alpha ^*C_b$$
so $B_a\not=0$ and $B_b\not=0$ on this set. We put 
$$A_a=:|A_a|e^{i\tilde A_a},\,A_b=:|A_b|e^{i\tilde A_b},\,B_a=:|B_a|e^{i\tilde B_a},\,B_b=:|B_b|e^{i\tilde B_b},$$
$$\,z_1=:e^{i\theta _1},\,z_2=:e^{i\theta _2}\,,$$
with $\tilde A_a,\,\tilde A_b,\,\tilde B_a,\,\tilde B_b\in \br$. By \pr \ref{10} a), $2\tilde A_a=\theta _2$, $2\tilde A_b=\theta _1+\pi $,
$$B_aC_b+B_bC_a=-\alpha \gamma B_aB_b^*+\beta \gamma B_bB_a^*=|B_a||B_b|(e^{i(\theta _2+\tilde B_b-\tilde B_a)}-e^{i(\theta _1+\tilde B_a-\tilde B_b)})=$$
$$=|B_a||B_b|e^{i\frac{\theta _1+\theta _2}{2}}(e^{i(\frac{\theta _2-\theta _1}{2}+\tilde B_b-\tilde B_a)}-e^{i(\frac{\theta _1-\theta _2}{2}+\tilde B_a-\tilde B_b)})=$$
$$=2|B_a||B_b|\sin(\frac{\theta _2-\theta _1}{2}+\tilde B_b-\tilde B_a)e^{i\frac{\theta _1+\theta _2+\pi }{2}}\;.$$
Since $2A_aA_b=-(B_aC_b+B_bC_a)$ there is a $k\in \bz$ with
$$\frac{\theta _2}{2}+\frac{\theta _1+\pi }{2}=\frac{\theta _1+\theta _2+\pi }{2}+(2k+1)\pi $$
which is a contradiction.

\begin{center}
Step 2 $\{A_a\not=0\}\subset \{A_c=0\}$
\end{center}

The assertion follows from Step 1 by symmetry.

\begin{center}
Step 3 $\{A_a\not=0\}=\{A_b=A_c=0\}$
\end{center}

The assertion follows from Steps 1 and 2 and from $|A_a|+|A_b|+|A_c|\not=0$.

\begin{center}
Step 4 The contradiction
\end{center}

By Step 3 and by the symmetry, the sets $\{A_a\not=0\}$, $\{A_b\not=0\}$, and $\{A_c\not=0\}$ are clopen and by $|A_a|+|A_b|+|A_c|\not=0$ their union is equal to $\bt^2$. So there is exactly one of these sets equal to $\bt^2$ which implies
$$A_a^2=z_2 \,,\quad \mbox{or}\quad A_b^2=-z_1\quad \mbox{or}\quad A_c^2=z_1z_2$$
and no one of these identities can hold.

b) is a direct verification.

c) follows from a) and b).\qed

\begin{center}
\subsection{$T:=(\bzz{2})^n$ with $n\in \bn$}
\end{center}

\begin{e}\label{h}
Assume $f$ constant and put
$$\s{s}{t}:=\proo{i=1}{n}(-1)^{s(i)t(i)}$$
for all $s,t\in T$ (where $\bzz{2}$ is identified with $\{0,1\}$) and
$$\mae{\varphi _t}{\ssa{f}}{E}{X}{\si{s\in T}\s{t}{s}X_s}$$
for all $t\in T$. Then the map
$$\mae{\varphi }{\ssa{f}}{E^{2^n}}{X}{(\varphi _tX)_{t\in T}}$$
is an $E$-C*-isomorphism.
\end{e}  

For $r,s,t\in T$,
$$t+t=0\,,\quad \s{s}{t}=\s{t}{s}\,,\quad \s{r+s}{t}=\s{r}{t}\s{s}{t}\,,$$
$$ \s{r}{s+t}=\s{r}{s}\s{r}{t}\;.$$
For $t\in T$ and $X,Y\in \ssa{f}$, by \h \ref{745} c),g),
$$\varphi _t(X^*)=\si{s\in T}\s{t}{s}(X^*)_s=\si{s\in T}\s{t}{s}(X_s)^*=(\varphi _tX)^*\,,$$
$$(\varphi _tX)(\varphi _tY)=\si{r,s\in T}\s{t}{r}\s{t}{s}X_rY_s=\si{q,r\in T}\s{t}{r}\s{t}{q-r}X_rY_{q-r}=$$
$$=\si{q,r\in T}\s{t}{q}X_rY_{q-r}=\si{q\in T}\s{t}{q}(XY)_q=\varphi _t(XY)$$
so $\varphi _t$ and $\varphi $ are $E$-C*-homomorphisms.

We have 
$$\si{t\in T}\s{0}{t}=2^n\;.$$
We want to prove
$$\si{t\in T}\s{s}{t}=0$$
for all $s\in T$, $s\not=0$, by induction with respect to $Card\me{i\in \bnn{n}}{s(i)\not=0}$. Let $i\in \bnn{n}$ with $s(i)\not=0$ and put $r:=s+e_i$,
$$T_0:=\me{t\in T}{t(i)=0}\,,\quad\quad T_1:=\me{t\in T}{t(i)=1}\;.$$
Then
$$\si{t\in T_0}\s{s}{t}=\si{t\in T_0}\s{r}{t}\,,\quad\quad \si{t\in T_1}\s{s}{t}=-\si{t\in T_1}\s{r}{t}\;.$$
But
$$\si{t\in T_0}\s{r}{t}=\si{t\in T_1}\s{r}{t}=2^{n-1}$$
if $r=0$. By the hypothesis of the induction
$$\si{t\in T_0}\s{r}{t}=\si{t\in T_1}\s{r}{t}=0$$
if $r\not=0$ (with $\bnn{n}$ replaced by $\bnn{n}\setminus \{i\}$, since $r(i)=0$). This finishes the proof by induction.

For $r\in T$ and $X\in \ssa{f}$, by the above,
$$\si{t\in T}\s{r}{t}\varphi _tX=\si{s,t\in T}\s{r}{t}\s{t}{s}X_s=\si{s,t\in T}\s{r+s}{t}X_s=$$
$$=\si{s\in T\setminus \{r\}}\si{t\in T}\s{r+s}{t}X_s+\si{t\in T}\s{0}{t}X_r=2^nX_r\;.$$
Hence $\varphi $ is bijective.\qed

\begin{e}\label{803}
Let $E:=\ccb{\bt^n}{\bc}$, denote by $z:=(z_1,z_2,\cdots,z_n)$ the points of  $\,\bt^n$, and put $z^2:=(z_1^2,z_2^2,\cdots,z_n^2)$ for every $z\in \bt^n$. We identify $(\bzz{2})^n$ with $\fr{P}(\bn_n)$ by using the bijection
$$\mad{\fr{P}(\bnn{n})}{(\bzz{2})^n}{I}{e_I}$$
and denote by 
$$I\triangle J:=(I\setminus J)\cup (J\setminus I)$$
the addition on $\fr{P}(\bnn{n})$ corresponding to this identification. We put $\lambda _I:=\pro{i\in I}z_i$ for every $I\subset \bnn{n}$ and
$$\mae{f}{\fr{P}(\bnn{n})\times \fr{P}(\bnn{n})}{\un{E}}{(I,J)}{\lambda _{I\cap J}}\;.$$
Then $f\in \f{(\bzz{2})^n}{E}$ and, if we put
$$\tilde X:=\si{I\subset \bnn{n}}\lambda _I(z)X_I(z^2)\in E$$
for every $X\in \ssa{f}$, the map
$$\mae{\varphi }{\ssa{f}}{E}{X}{\tilde X}$$
is an isomorphism of C*-algebras.
\end{e}

Let $X,Y\in \ssa{f}$. By \h \ref{745} c),g),
$$\widetilde{X^*}=\si{I\subset \bnn{n}}\lambda _I\,(X^*)_I(z^2)=\si{I\subset \bnn{n}}\lambda _I\,\overline{\lambda _I^2}\,X_I^*=\overline{\tilde X}\,,$$
$$\widetilde{XY}=\si{I\subset \bnn{n}}\lambda _I\,(XY)_I(z^2)=\si{I\subset \bnn{n}}\lambda _I\si{J\subset \bnn{n}}f(J,J\triangle I)^2X_JY_{J\triangle I}=$$
$$=\si{J,K\subset \bnn{n}}\lambda _{J\triangle K}\lambda _{J\cap K}^2X_JY_K=\si{J,K\subset \bnn{n}}\lambda _J\lambda _KX_JY_K=\tilde X\tilde Y$$
so $\varphi $ is a C*-homomorphism.

We put for $k\in \bnn{n}$, $i\in \bz^n$, and $I\subset \bnn{n}$,
$$i_k^I:=\ab{2i_k+1}{k\in I}{2i_k}{k\in \bnn{n}\setminus I}\,,\qquad\qquad i^I:=(i_1^I,i_2^I,\cdots,i_n^I)\in \bz^n$$
and
$$\ccc{G}:=\me{\si{i\in \bz^n}a_iz_1^{i_1}z_2^{i_2}\cdots z_n^{i_n}}{(a_i)_{i\in \bz^n}\in \bc^{(\bz^n)}}\;.$$
Let 
$$x:=\si{i\in \bz^n}a_iz_1^{i_1}z_2^{i_2}\cdots z_n^{i_n}\in \ccc{G}$$
and for every $I\subset \bnn{n}$ put
$$X_I:=\si{i\in \bz^n}a_{i^I}z_1^{i_1}z_2^{i_2}\cdots z_n^{i_n}\,,\qquad\quad X:=\si{I\subset \bnn{n}}(X_I\otimes 1_K)V_I\;.$$
Then $\varphi X=x$ and so $\ccc{G}\subset \varphi (\ssa{f})$. Since $\ccc{G}$ is dense in $E$, it follows that $\varphi $ is surjective.

We prove that $\varphi $ is injective by induction with respect to $n\in \bn$. The case $n=1$ was proved in Example \ref{789}. Assume the assertion holds for $n-1$. Let $X\in Ker\,\varphi $. Then
$$\si{I\subset \bnn{n}}\lambda _I(z)X_I(z^2)=0\;.$$
By replacing  $z_n$ by $-z_n$ in the above relation, we get
$$\si{I\subset \bnn{n-1}}\lambda _I(z)X_I(z^2)-\si{n\in I\subset \bnn{n}}\lambda _I(z)X_I(z^2)=0$$
and so
$$\si{I\subset \bnn{n-1}}\lambda _I(z)X_I(z^2)=\si{n\in I\subset \bnn{n}}\lambda _I(z)X_I(z^2)=0\;.$$
By the induction hypothesis, we get $X_I=0$ for all $I\subset \bnn{n}$ and so $X=0$. Thus $\varphi $ is injective and a C*-isomorphism.\qed

\begin{e}\label{19}
Let $f\in \f{(\bzz{2})^3}{E}$, put
$$a:=(0,0,1)\,,\quad b:=(0,1,0)\,,\quad c:=(0,1,1)\,,\quad s:=(1,0,0)\,,$$
and denote by $g$ the element of $\f{\bzz{2}}{E}$ defined by $g(1,1):=f(s,s)$ \emph{\pr \ref{785} a)}.
\begin{enumerate}
\item There is a family $(\alpha _i\,,\beta _i\,,\gamma _i\,,\varepsilon _i)_{i\in \bnn{7}}$ in $(\un{E})^4$ such that $f$ is given by the attached table and such that $\varepsilon _i^2=1_E$ for every $i\in \bnn{7}$ and 
$$\varepsilon _3=\varepsilon _1\varepsilon _2\,,\quad \varepsilon _5=\varepsilon _1\varepsilon _4\,,\quad \varepsilon _6=\varepsilon _2\varepsilon _4\,,\quad \varepsilon _7=\varepsilon _1\varepsilon _2\varepsilon _4\,,$$
$$\alpha _3=\varepsilon _2\varepsilon _4\alpha _1\alpha _2^*\alpha _4\alpha _6\gamma _2^*\,,\quad \alpha _5=\alpha _6\beta _1\gamma _2^*\,,\quad \alpha _7=\alpha _4\gamma _1\gamma _2^*\,,$$
$$\beta _2=\beta _1\gamma _1\gamma _2^*\,,\quad \beta _3=\varepsilon _2\alpha _4^*\alpha _6\beta _1\,,\quad \beta _4=\varepsilon _1\varepsilon _2\varepsilon _4\alpha _1\alpha _2^*\alpha _4\gamma _1\gamma _2^*\,,$$
$$\beta _5=\varepsilon _4\alpha _1\alpha _2^*\alpha _6\,,\quad \beta _6=\varepsilon _4\alpha _1\alpha _2^*\alpha _6\beta _1\gamma _2^*\,,\quad \beta _7=\varepsilon _1\varepsilon _2\varepsilon _4\alpha _1\alpha _2^*\alpha _6\,,$$
$$\gamma _3=\varepsilon _2\alpha _4\alpha _6^*\gamma _1\,,\quad \gamma _4=\varepsilon _2\varepsilon _4\alpha _2\alpha _4^*\gamma _1^*\gamma _2\,,\quad \gamma _5=\varepsilon _1\varepsilon _4\alpha _2\alpha _6^*\gamma _1\,,$$
$$\gamma _6=\varepsilon _4\alpha _2\alpha _6^*\gamma _2\,,\qquad \gamma _7=\varepsilon _1\varepsilon _2\varepsilon _4\alpha _2\alpha _4^*\beta _1\;.$$

$$\begin{array}{||c||c|c|c|c|c|c|c||}\hline \hline
f&a&b&c&s&a+s&b+s&c+s\\\hline\hline a&\beta _1\gamma _1&\gamma _1&\beta _1&\gamma _2&\beta _2&\gamma _3&\beta _3\\\hline b&\varepsilon _1\gamma _1&\varepsilon _1\alpha _1\gamma _1&\alpha _1&\gamma _4&\gamma _5&\beta _4&\beta _5\\\hline c&\varepsilon _1\beta _1&\varepsilon _1\alpha _1&\alpha _1\beta _1&\gamma _6&\gamma _7&\beta _7&\beta _6\\\hline s&\varepsilon _2\gamma _2&\varepsilon _4\gamma _4&\varepsilon _6\gamma _6&\varepsilon _2\alpha _2\gamma _2&\alpha _2&\alpha _4&\alpha _6\\\hline a+s&\varepsilon _2\beta _2&\varepsilon _5\gamma _5&\varepsilon _7\gamma _7&\varepsilon _2\alpha _2&\alpha _2\beta _2&\alpha _7&\alpha _5\\\hline b+s&\varepsilon _3\gamma _3&\varepsilon _4\gamma _4&\varepsilon _7\gamma _7&\varepsilon _4\alpha _4&\varepsilon _7\alpha _7&\varepsilon _3\alpha _3\gamma _3&\alpha _3\\\hline c+s&\varepsilon _3\beta _3&\varepsilon _5\beta _5&\varepsilon _6\beta _6&\varepsilon _6\beta _6&\varepsilon _5\alpha _5&\varepsilon _3\alpha _3&\alpha _3\beta _3\\\hline\hline
\end{array}$$

\item If $\varepsilon _1=-1_E$, $\varepsilon _2=\varepsilon _4$, $\gamma _1=1_E$, and there is an $x\in E^c$ with $x^2=\alpha _1\beta _1^*$ then there are $P_\pm \in (E\widetilde\otimes 1_K)^c\cap Pr\,\ssa{f}$ with $P_++P_-=V^f_1$ and \emph{(\h \ref{895} b))}
$$P_+\ssa{f}P_+\approx _E\ssa{g}\approx _EP_-\ssa{f}P_-\;.$$
\item If $\varepsilon _1=-1_E$, $\varepsilon _2=\varepsilon _4=\gamma _1=1_E$, and there is an $x\in E^c$ with $x^2=\alpha _1\beta _1^*$ then $\ssa{f}\approx _E\ssa{g}_{2,2}$.
\item Assume $\varepsilon _1=-1_E$, $\varepsilon _2=\varepsilon _4=\alpha _1=\beta _1=\gamma _1=1_E$, $\gamma _2=\alpha _2^*$, and $\alpha _2^4=\alpha _4^4=\alpha _6=1_E$ and put
$$\mae{\varphi _\pm }{\ssa{f}}{E_{2,2}}{X}$$
$${\mt{X_0+X_c\pm X_s\pm X_{c+s}}{X_a-X_b\pm \alpha _2^*X_{a+s}\mp \alpha _4^*X_{b+s}}{X_a+X_b\pm \alpha _2^*X_{a+s}\pm \alpha _4^*X_{b+s}}{X_0-X_c\pm X_s\mp X_{c+s}}}\;.$$ 
Then the map
$$\mad{\ssa{f}}{E_{2,2}\times E_{2,2}}{X}{(\varphi _+X,\varphi _-X)}$$
is an $E$-C*-isomorphism.
\end{enumerate}
\end{e}

a) is a long calculation.

b) and c) follow from a) and \h \ref{895} e).

d) is a long calculation using a).\qed

\begin{center}
\subsection{$T:=\bzz{n}$ with $n\in \bn$}
\end{center}
\begin{p}\label{6437}
Put $A:=\un{E}$ and for every $\alpha \in A^{n-1}$ put
$$\mae{f_\alpha }{\bzz{n}\times \bzz{n}}{A}{(p,q)}{\left(\proo{j=p}{p+q-1}\alpha _j\right)\left(\proo{k=1}{q-1}\alpha _k^*\right)}\,,$$
where $\bzz{n}$ and $\bnn{n}$ are canonically identified and $\alpha _n:=1_E$.
\begin{enumerate}
\item For every $f\in \f{\bzz{n}}{E}$ and $X\in \ssa{f}$, $X\in \ssa{f}^c$ iff $X_t\in E^c$ for all $t\in T$. In particular, $\ssa{f}$ is commutative if $E$ is commutative.
\item $f_\alpha \in \f{\bzz{n}}{E}$ for every $\alpha \in A^{n-1}$ and the map
$$\mad{A^{n-1}}{\f{\bzz{n}}{E}}{\alpha }{f_\alpha }$$
is a group isomorphism.
\item The following are equivalent for all $\alpha ,\beta \in A^{n-1}$.
\begin{enumerate}
\item $\ssa{f_\alpha }\approx _{\ccc{S}}\ssa{f_\beta }$.
\item There is a $\gamma \in A$ such that
$$\gamma ^n=\proo{j=1}{n-1}(\alpha _j\beta _j^*)\;.$$
\item There is a $\lambda \in \Lambda (\bzz{n},E)$ such that $f_\alpha =f_\beta \delta \lambda $.

\vspace{3ex}
If these equivalent conditions are fulfilled then the map
$$\mad{\ssa{f_\alpha }}{\ssa{f_\beta }}{X}{U_\lambda ^*XU_\lambda }$$
is an $\ccc{S}$-isomorphism and
$$\lambda (1)^n=\proo{j=1}{n-1}(\alpha _j\beta _j^*)=\gamma ^n\,,\quad p\in \bzz{n}\Longrightarrow \lambda (p)=\lambda (1)^p\proo{j=1}{p-1}(\alpha _j^*\beta _j)\;.$$
\end{enumerate}
\item Let $\alpha \in A^{n-1}$ and put
$$\mae{\beta }{\bnn{n-1}}{A}{j}{\ab{1}{j<n-1}{\left(\proo{k=1}{n-1}\alpha _k^*\right)^{n-1}}{j=n-1}}\;.$$
Then $\alpha $ and $\beta $ fulfill the equivalent conditions of c).
\item There is a natural bijection
$${\me{\ssa{f}}{f\in \f{\bzz{n}}{E}}/\approx _{\ccc{S}}}\longrightarrow {A/\me{x^n}{x\in A}}\;.$$
If $E:=\ccb{\bt^m}{\bc}$ for some $m\in \bn$ then 
$$Card\,(\me{\ssa{f}}{f\in \f{\bzz{n}}{E}}/\approx _{\ccc{S}})=mn\;.$$
\item Let $\alpha \in A^{n-1}$, $\beta \in A$ such that $\beta ^n=\proo{j=1}{n-1}\alpha _j$, 
$$F:=\ab{E}{\bk=\bc}{\stackrel{\circ }{E}}{\bk=\br}\,,$$
where $\stackrel{\circ }{E}$ denotes the complexification of $E$, and
$$\mae{w_k}{\ssa{f_\alpha }}{F}{X }{\sii{j=1}{n}\beta ^j\left(\proo{l=1}{j-1}\bar \alpha _l\right)e^{\frac{2\pi ijk}{n}}X_j}$$
for every $k\in \bnn{n}(=\bzz{n})$.
\begin{enumerate}
\item If $\bk=\bc$ then the map
$$\mad{\ssa{f_\alpha }}{E^n}{X}{(w_kX)}_{k\in \bzz{n}}$$
is an $E$-C*-isomorphism.
\item If $\bk=\br$ and $n$ is odd then we may take $\beta \in \br$ and the map
$$\mad{\ssa{f_\alpha }}{E\times (\stackrel{\circ }{E})^{\frac{n-1}{2}}}{X}{(w_nX,(w_kX)_{k\in \bnn{\frac{n-1}{2}}})}$$
is an $E$-C*-isomorphism.
\item If $\bk=\br$, $n$ is even, and $\proo{j=1}{n-1}\alpha _j=-1$ then the map
$$\mad{\ssa{f_\alpha }}{(\stackrel{\circ }{E})^{\frac{n}{2}}}{X}{(w_{k-1}X)_{k\in \bnn{\frac{n}{2}}}}$$
is an $E$-C*-isomorphism.
\item If $\bk=\br$, $n$ is even, and $\proo{j=1}{n-1}\alpha _j=1$, and $\beta =1$ then the map
$$\mad{\ssa{f_\alpha }}{E\times E\times (\stackrel{\circ }{E})^{\frac{n}{2}-1}}{X}{(w_nX,w_{\frac{n}{2}}X,(w_kX)_{k\in \bnn{\frac{n}{2}-1}})}$$
is an $E$-C*-isomorphism.
\item If $n$ is even then there is a $\gamma \in A$ such that $f_\alpha (\frac{n}{2},\frac{n}{2})=\gamma ^2$.
\end{enumerate}
\end{enumerate}
\end{p}\qed

\begin{e}\label{6439}
Let $E:=\ccb{\bt}{\bc}$, $r\in \bz^{n-1}$, $z:\bt\rightarrow \bc$ the canonical inclusion, and
$$\mae{f}{\bzz{n}\times \bzz{n}}{\un{E}}{(p,q)}{z^{\left(\sii{j=p}{p+q-1}r_j-\sii{j=1}{q-1}r_j\right)}}\,,$$
where $\bzz{n}$ and $\bnn{n}$ are canonically identified. Then $f\in \f{\bzz{n}}{E}$. Let further $S$ be the subgroup of $\bzz{n}$ generated by $\rho (\sii{j=1}{n-1}r_j)$, where $\rho :\bz\rightarrow \bzz{n}$ is the quotient map,
$$m:=Card\;S\,,\qquad h:=\frac{n}{m}\,,\qquad \omega :=e^{\frac{2\pi i}{n}}\,,$$
$$\mae{\sigma }{\bnn{n}}{\bz}{p}{\frac{p}{h}\sii{j=1}{n-1}r_j-m\sii{j=1}{p-1}r_j}\,,$$
and
$$\mae{\varphi _k}{\ssa{f}}{E}{X}{\sii{p=1}{n}(X_p\circ z^m)z^{\sigma (p)}\omega ^{pk}}$$
for every $k\in \bnn{h}$. Then the map
$$\mae{\varphi }{\ssa{f}}{E^h}{X}{(\varphi _kX})_{k\in \bnn{h}}$$
is an $E$-C*-isomorphism.\qed
\end{e}

The next example shows that the set $\me{\ssa{f}}{f\in \f{\bzz{n}}{\ccb{\bt}{\bc}}}$ is not reduced by restricting the Schur functions to have the form indicated in \ee{6439}.

\begin{e}\label{6444}
Let $E:=\ccb{\bt}{\bc}$ and $g\in \f{\bzz{n}}{E}$. Put
$$\mae{\varphi }{[0,2\pi [}{\br}{\theta }{\log\proo{j=1}{n-1}(g(j,1))(e^{i\theta })}\,,$$
where we take a fixed (but arbitrary) branch of $\log$. If we define
$$\mae{r}{\bnn{n-1}}{\bz}{j}{\ab{\lim\limits_{\theta \rightarrow 2\pi }\varphi (\theta )-\varphi (0)}{j=1}{0}{j\not=1}}$$
then there is a $\lambda \in \Lambda (\bzz{n},E)$ such that $g=f\delta \lambda $, where $f$ is the Schur function defined in \emph{ \ee{6439}}. In particular $\ssa{f}\approx _{\ccc{S}}\ssa{g}$.\qed
\end{e}

\begin{center}
\subsection{$T:=\bz$}
\end{center}
\begin{e}\label{819}
Let $f\in \f{\bz}{E}$.
\begin{enumerate}
\item $\ssb{\n{\cdot }}{f}\approx \ccb{\bt}{E}$.
\item If $E$ is a W*-algebra then
$$\ssb{W}{f}\approx E\bar \otimes L^\infty (\mu )\approx L^\infty (\mu ,E)\,,$$
where $\mu $ denotes the Lebesgue measure on $\bt$.
\end{enumerate}
\end{e}

By Corollary \ref{719} c) and \pr \ref{716} $a_1\Rightarrow a_2$, we may assume $f$ constant. By \pr \ref{771} c),e), we may assume $E:=\bc$. Let $\alpha :\bt\rightarrow \bc$ be the inclusion map. Then
$$\mad{l^2(\bz)}{L^2(\mu )}{\xi }{\si{n\in \bz}\xi _n\alpha ^n}$$
is an isomorphism of Hilbert spaces. If we identify these Hilbert spaces using this isomorphism then $V_1$ becomes the multiplicator operator
$$\mad{L^2(\mu )}{L^2(\mu )}{\eta }{\alpha \eta }$$
so
$$\mad{\ccc{R}(f)}{L^\infty (\mu )}{X}{\si{n\in \bz}X_n\alpha ^n}$$
is an injective, involutive algebra homomorphism. The assertion follows.\qed 

\begin{center}
\section{Clifford Algebras}
\subsection{The general case}
\end{center}
\fbox{\parbox{12 cm}{Throughout this subsection $I$ is a totally ordered set, $(T_i)_{i\in I}$ is a family of groups, and $(f_i)_{i\in I}\in \pro{i\in I}\f{T_i}{E}$. We put \vspace{-1ex}
$$\bar t:=\me{i\in I}{t_i\not=1_i}$$ \vspace{-1ex}
 for every $t\in \pro{i\in I}T_i$ (where $1_i$ denotes the neutral element of $T_i$) and \vspace{-1ex}
$$T:=\me{t\in \pro{i\in I}T_i}{\bar t\;\mbox{is finite}}\,,\qquad T':=\me{t\in T}{t^2=1}\;.$$ 
An associated $f\in \fte$ will be defined in \pr \ref{22} b).}}

$T$ is a subgroup of $\pro{i\in I}T_i$. We canonically associate to every element $t\in T$ in a bijective way the "word" $t_{i_1}t_{i_2}\cdots t_{i_n}$, where 
$$\{i_1,i_2,\cdots,i_n\}=\bar t\qquad \mbox{and}\qquad i_1<i_2<\cdots<i_n$$
and use sometimes this representation instead of $t$ ( to $1\in T$ we associate the "empty word").

\begin{p}\label{22}
\rule{1em}{0ex}
\begin{enumerate}
\item Let $t_{i_1}t_{i_2}\cdots t_{i_n}$ be a finite sequence of letters with $t_{i_j}\in T_{i_j}\setminus \{1_{i_j}\}$ for every $j\in \bnn{n}$ and use transpositions of successive letters \emph{with distinct indices} in order to bring these indices in an increasing order. If $\tau $ denotes the number of used transpositions then $(-1)^\tau $ does not depend on the manner in which this operation was done.
\item Let $s,t\in T$ and let
$$s_{i_1}s_{i_2}\cdots s_{i_m}\,,\qquad t_{i'_1}t_{i'_2}\cdots t_{i'_n}$$
be the canonically associated words of $s$ and $t$, respectively. We put for every $k\in I$, $\tilde s_k:=s_{i_j}$ if there is a $j\in \bnn{m}$ with $k=i_j$ and $\tilde s_k:=1_k$ if the above condition is not fulfilled and define $\tilde t$ in a similar way. Moreover we put \emph{(\pr \ref{704} a))}
$$f(s,t):=(-1)^\tau \pro{k\in I}f_k(\tilde s_k,\tilde t_k)\,,$$
where $\tau $ denotes the number of transpositions of successive letters \emph{with distinct indices} in the finite sequence of letters
$$s_{i_1}s_{i_2}\cdots s_{i_m}t_{i'_1}t_{i'_2}\cdots t_{i'_n}$$
in order to bring the indices in an increasing order. Then $f\in \fte$.
\item Let $I_0$ be a subset of $I$, $T_0$ the subgroup $\me{t\in T}{\bar t\subset I_0}$ of $T$, and $f_0$ the element of $\f{T_0}{E}$ defined in a similar way as $f$ was defined in b). Then $f_0=f|(T_0\times T_0)$ and the map
$$\mad{\ssb{\n{\cdot }}{f_0}}{\ssb{\n{\cdot }}{f}}{\sii{t\in T_0}{\n{\cdot }}(X_t\widetilde\otimes 1_K)V_t^{f_0}}{\sii{t\in T_0}{\n{\cdot }}(X_t\widetilde\otimes 1_K)V_t^f}$$
is an injective $E$-C**-homomorphism with image 
$$\me{X\in \ssa{f}}{(t\in T\,\&\,X_t\not=0)\Rightarrow t\in T_0}\;.$$
\end{enumerate}
\end{p}

a) We define a new total order relation on the indices of the given word by putting for all $j,k\in \bnn{n}$
$$i_j\prec i_k:\Longleftrightarrow ((i_j<i_k)\; \mbox{or}\; (i_j=i_k \;\mbox{and}\; j<k))\;.$$
Let $P$ be a sequence of transpositions of successive letters in order to bring the indices in an increasing form with respect to the new order and let $\tau '$ be the number of used transpositions. Then $\tau -\tau '$ is even and so $(-1)^\tau =(-1)^{\tau '}$. By the theory of permutations $(-1)^{\tau '}$ does not depend on $P$, which proves the assertion.

b) By a), $f$ is well-defined. Let $r,s,t\in T$ and let
$$r_{i_1}r_{i_2}\cdots r_{i_m}\,,\qquad s_{i'_1}s_{i'_2}\cdots s_{i'_n}\,,\qquad t_{i''_1}t_{i''_2}\cdots t_{i''_p}$$
be the words canonically associated to $r$, $s$, and $t$, respectively. There are $\alpha ,\beta \in \{-1,+1\}$ such that
$$f(r,s)f(rs,t)=\alpha \pro{i\in I}f(\tilde r_i,\tilde s_i)f(\widetilde{r_is_i},\tilde t_i)\,,$$
$$ f(r,st)f(s,t)=\beta \pro{i\in I}f_i(\tilde r_i,\widetilde{s_it_i})f(\tilde s_i,\tilde t_i)\;.$$
Write the finite sequence of letters
$$r_{i_1}r_{i_2}\cdots r_{i_m}s_{i'_1}s_{i'_2}\cdots s_{i'_n}t_{i''_1}t_{i''_2}\cdots t_{i''_p}$$
and use transpositions of successive letters with distinct indices in order to bring the indices in an increasing order. We can do this acting first on the letters of $r$ and $s$ only and then in a second step also on the letters of $t$. Then $\alpha =(-1)^\mu $, where $\mu $ denotes the number of all performed transpositions. For $\beta $ we may start first with the letters of $s$ and $t$ and then in a second step also with the letters of $r$. Then $\beta =(-1)^\nu $, where $\nu $ is the number of all effectuated transpositions. By a), $\alpha =(-1)^\mu =(-1)^\nu =\beta $. The rest of the proof is obvious.

c) follows from Corollary \ref{776} d).\qed

\begin{co}\label{29}
If $I:=\bnn{2}$ then for all $s,t\in T$,
$$f(s,t)=\ac{f_1(s_1,t_1)}{s_2=1_2}{f_2(s_2,t_2)}{t_1=1_1}{-f_1(s_1,t_1)f_2(s_2,t_2)}{s_2\not=1_2\,,t_1\not=1_1}\;.\qedd$$
\end{co}

\begin{p}\label{23}
Let $s,t\in T$. 
\begin{enumerate}
\item $f(s,t)=(-1)^{Card\,(\bar s\times \bar t)-Card\,(\bar s\cap \bar t)}f(t,s)$.
\item $st=ts$ iff $V_sV_t=(-1)^{Card\,(\bar s\times \bar t)-Card\,(\bar s\cap \bar t)}V_tV_s$.
\item Assume $\bar s\subset \bar t$. If $Card\,\bar s$ is even or if $Card\,\bar t$ is odd then $f(s,t)=f(t,s)$. If in addition $st=ts$ then $V_sV_t=V_tV_s$.
\item If $Card\,I$ is an odd natural number and $T$ is commutative then $V_t\in \ssa{f}^c$ for every $t\in T$ with $\bar t=I$.
\item Assume $t\in T'$. If $n:=Card\,\bar t$ and $\alpha :=\pro{i\in \bar t}f_i(t_i,t_i)$ then
$$f(t,t)=(-1)^{\frac{n(n-1)}{2}}\alpha \,,\qquad \tilde f(t)=(-1)^{{\frac{n(n-1)}{2}}}\alpha ^*\,,$$
$$(V_t)^2=(-1)^{\frac{n(n-1)}{2}}(\alpha \widetilde\otimes 1_K)V_1\,,\qquad V_t^*=(-1)^{\frac{n(n-1)}{2}}(\alpha ^*\widetilde\otimes 1_K)V_t\;.$$
\end{enumerate}
\end{p}

a) For $i\in \bar s$,
$$f(s_i,t)=\ab{(-1)^{Card\,\bar t}f(t,s_i)}{i\not\in \bar t}{(-1)^{Card\,\bar t-1}f(t,s_i)}{i\in \bar t}$$
so
$$f(s,t)=(-1)^{Card\,(\bar s\times \bar t)-Card\,(\bar s\cap \bar t)}f(t,s)\;.$$

b) By \pr \ref{674} b),
$$V_sV_t=(f(s,t)\widetilde\otimes 1_K)V_{st}\,,\qquad V_tV_s=(f(t,s)\widetilde\otimes 1_K)V_{ts}\;.$$
Thus if $st=ts$ then by a),
$$V_sV_t=((f(s,t)f(t,s)^*)\widetilde\otimes 1_K)V_tV_s=(-1)^{Card\,(\bar s\times \bar t)-Card\,(\bar s\cap \bar t)}V_tV_s\;.$$
Conversely, if this relation holds then by a),
$$V_{st}=(f(s,t)^*\widetilde\otimes 1_K)V_sV_t=
(-1)^{Card\,(\bar s\times \bar t)-Card\,(\bar s\cap \bar t)}(f(t,s)^*\widetilde\otimes 1_K)V_sV_t=$$
$$=(f(t,s)^*\widetilde\otimes 1_K)V_tV_s=V_{ts}$$
and we get $st=ts$ by \h \ref{745} a).

c) follows from a) and b).

d) follows from c) (and \pr \ref{674} d)).

e) We have
$$f(t,t)=(-1)^{(n-1)+\cdots+2+1}\alpha =(-1)^\frac{n(n-1)}{2}\alpha \;.$$
By \pr \ref{674} b),e),
$$(V_t)^2=(f(t,t)\widetilde\otimes 1_K)V_1=(-1)^{\frac{n(n-1)}{2}}(\alpha \widetilde\otimes 1_K)V_1\,,$$
$$V_t^*=\tilde f(t)V_{t^{-1}}=f(t,t)^*V_t=(-1)^{\frac{n(n-1)}{2}}(\alpha ^*\widetilde\otimes 1_K)V_t\;.\qedd$$

\begin{p}\label{839}
Let $S$ be a finite subset of $T'\setminus \{1\}$ such that $st=ts$ and $Card\;(\bar s\times \bar t)-Card\;(\bar s\cap \bar t)$ is odd for all distinct $s,t\in S$ and for every $t\in S$ let $\alpha _t,\varepsilon _t\in \un{E}$ and $X_t\in E$ be such that
$$\varepsilon _t^2=1_E\,,\qquad (V_t)^2=(\alpha _t\widetilde \otimes 1_K)V_1 \,,\qquad X_t^*=\alpha _tX_t\,,$$
$$\si{t\in S}|X_t|^2=\frac{1}{4}1_E\;.$$
\begin{enumerate}
\item 
$$P:=\frac{1}{2}V_1+\si{t\in S}((\varepsilon _tX_t)\tilde \otimes 1_K)V_t\in Pr\;\ssa{f}\,,$$ 
$$V_1 -P=\frac{1}{2}V_1 +\si{t\in S}((-\varepsilon _tX_t)\tilde \otimes 1_K)V_t\in Pr\;\ssa{f}\;.$$
\item If $s\in S$ and $\beta \in E^c$ such that $X_s=0$ and $\beta ^2=\alpha _s$ then $P$ is homotopic in $Pr\;\ssa{f}$ to 
$$\frac{1}{2}(V_1 +((\beta ^*\varepsilon _s)\widetilde \otimes 1_K)V_s)\;.$$
\end{enumerate}
\end{p}

a) By \pr \ref{23} b),e), 
$$P^*=\frac{1}{2}V_1+\si{t\in S}((\varepsilon _tX_t^*\alpha _t^*)\widetilde \otimes 1_K)V_t=
\frac{1}{2}V_1 +\si{1\in S}((\varepsilon _tX_t)\widetilde \otimes 1_K)V_t=P\,,$$
$$P^2=\frac{1}{4}V_1 +\si{t\in S}(X_t^2\widetilde \otimes 1_K)(V_t)^2+\si{t\in S}((\varepsilon _tX_t)\widetilde \otimes 1_K)V_t+$$
$$+\si{s,t\in S\atop s\not=t}((\varepsilon _s\varepsilon _tX_sX_t)\widetilde \otimes 1_K)(V_sV_t+V_tV_s)=$$
$$=\frac{1}{4}V_1 +\si{t\in S}((X_t^2\alpha _t)\widetilde \otimes 1_K)V_1 +\si{t\in S}((\varepsilon _tX_t)\widetilde \otimes 1_K)V_t=$$
$$=\frac{1}{4}V_1 +\si{t\in S}(|X_t|^2\widetilde \otimes 1_K)V_1 +\si{t\in S}((\varepsilon _tX_t)\widetilde \otimes 1_K)V_t=P\;.$$

b) Remark first that $\beta \in \un{E}$ and put
$$\mae{Y}{[0,1]}{E^c_+}{u}{(\frac{1}{4}1_E-u^2\si{t\in S}|X_t|^2)^{\frac{1}{2}}}\,,$$
$$\mae{Z}{[0,1]}{E^c}{u}{\beta ^*\varepsilon _sY(u)}\,,$$
$$\mae{Q}{[0,1]}{\ssa{f}}{u}{\frac{1}{2}V_1 +(Z(u)\widetilde \otimes 1_K)V_s+\si{t\in S\setminus \{s\}}((u\varepsilon _tX_t)\widetilde \otimes 1_K)V_t}\;.$$
For $u\in [0,1]$,
$$\alpha _sZ(u)=\beta ^2\beta ^*\varepsilon _s Y(u)=\beta \varepsilon _s Y(u)=Z(u)^*\,,$$
$$|Z(u)|^2+\si{t\in S\setminus \{s \}}|uX_t|^2=\frac{1}{4}1_E$$
so by a), $Q(u)\in Pr\,\ssa{f}$. Moreover
$$Q(0)=\frac{1}{2}(V_1 +((\beta ^*\varepsilon _s )\widetilde \otimes 1_K)V_s)\,,\qquad Q(1)=P\;.\qedd$$

\begin{co}\label{843}
Let $s,t\in T'\setminus \{1 \}$, $s\not=t$, $st=ts$, $\alpha _s,\alpha _t,\varepsilon _s,\varepsilon _t\in \un{E}$ such that
$$\varepsilon _s^2=\varepsilon _t^2=1_E\,,\qquad (V_s)^2=(\alpha _s^2\widetilde\otimes 1_K)V_1\,,\qquad (V_t)^2=(\alpha _t^2\widetilde\otimes 1_K)V_1 \,,$$
and put
$$P_s:=\frac{1}{2}(V_1 +((\varepsilon _s\alpha _s^*)\widetilde\otimes 1_K)V_s)\,, \quad P_t:=\frac{1}{2}(V_1 +((\varepsilon _t\alpha _t^*)\widetilde\otimes 1_K)V_t)\;.$$
\begin{enumerate}
\item $P_s,P_t\in Pr\,\ssa{f}$; we denote by $P_s\wedge P_t$ the infimum of $P_s$ and $P_t$ in $\ssa{f}_+$ \emph{(by b) and c) it exists)}.
\item If $V_sV_t\not=V_tV_s$ then $P_s\wedge P_t=0$.
\item If $V_sV_t=V_tV_s$ then $P_s\wedge P_t=P_sP_t\in Pr\,\ssa{f}$.
\end{enumerate}
\end{co}

a) follows from \pr \ref{25} $b\Rightarrow a$.

b) By \pr \ref{23} b), $V_sV_t=-V_tV_s$. Let $X\in \ssa{f}_+$ with $X\leq P_s$ and $X\leq P_t$. By [C1] \pr 4.2.7.1 $d\Rightarrow c$,
$$X=P_sX=\frac{1}{2}X+\frac{1}{2}((\varepsilon _s\alpha _s^*)\widetilde\otimes 1_K)V_sX\,,$$
$$X=((\varepsilon _s\alpha _s^*)\widetilde\otimes 1_K)V_sX=((\varepsilon _s\varepsilon _t\alpha _s^*\alpha _t^*)\widetilde\otimes 1_K)V_sV_tX=$$
$$=-((\varepsilon _s\varepsilon _t\alpha _s^*\alpha _t^*)\widetilde\otimes 1_K)V_tV_sX=-X$$
so $X=0$ and $P_s\wedge P_t=0$.

c) We have $P_sP_t=P_tP_s$ so $P_sP_t\in Pr\,\ssa{f}$ and $P_sP_t=P_s\wedge P_t$ by [C1] Corollary 4.2.7.4 $a\Rightarrow b\&d$.\qed

\begin{co}\label{850}
Let $m,n\in \bn$, $\bnn{m+n}\subset I$, $(\alpha _i)_{i\in \bnn{m}}\in (\un{E})^m$, and for every $i\in \bnn{m}$ let $t_i\in T'$ with $\bar t_i:=\bnn{n}\cup \{n+i\}$ and $t_it_j=t_jt_i$ for all $i,j\in \bnn{m}$. If for every $i\in \bnn{m}$,
$$(V_{t_i})^2=(\alpha _i^2\otimes 1_K)V_1 $$
then
$$\frac{1}{2}\left(V_1 +\frac{1}{\sqrt{m}}\si{i\in \bnn{m}}(\alpha _i^*\otimes 1_K)V_{t_i}\right)\in Pr\,\ssa{f}\;.$$
\end{co}

For distinct $i,j\in \bnn{m}$,
$$Card\;(\bar t_i\times \bar t_j)-Card\;(\bar t_i\cap \bar t_j)=(n+1)^2-n=n(n+1)+1$$
is odd. For every $i\in \bnn{m}$ put $X_i:=\frac{1}{2\sqrt{m}}\alpha _i^*$. Then
$$\alpha _i^2X_i=\frac{1}{2\sqrt{m}}\alpha _i=X_i^*\,,\qquad |X_i|^2=\frac{1}{4m}1_E\,,\qquad \si{i\in \bnn{m}}|X_i|^2=\frac{1}{4}1_E$$
and the assertion follows from \pr \ref{839} a).\qed

\begin{theo}\label{27}
Let $n\in \bn$ such that $\bnn{2n}$ is an ordered subset of $I$, $S:=\me{t\in T}{\bar t\subset \bnn{2n-2}}$, $g:=f|(S\times S)$, $a,b\in T$ such that $a^2=b^2=1$,
$$ \bar a=\bnn{2n-1}\,,\qquad \bar b=\bnn{2n-2}\cup \{2n\}\,,\qquad i\in \bnn{2n-2}\Longrightarrow a_i=b_i\,,$$
$\omega :\bzz{2}\times \bzz{2}\rightarrow T$ the (injective) group homomorphism defined by $\omega (1,0):=a$, $\omega (0,1):=b$, $\alpha _1:=f(a,a)$, $\alpha _2:=f(b,b)$, $\beta _1,\beta _2\in \un{E}$ such that $\alpha _1\beta _1^2+\alpha _2\beta _2^2=0$,
$$\gamma :=\frac{1}{2}(\alpha _1^*\beta _1^*\beta _2-\alpha _2^*\beta _1\beta _2^*)=\alpha _1^*\beta _1^*\beta _2=-\alpha _2^*\beta _1\beta _2^*\,,$$
$$X:=\frac{1}{2}((\beta _1 \widetilde\otimes 1_K)V_a+(\beta _2\widetilde\otimes 1_K)V_b)\,,\qquad P_+:=X^*X\,,\quad P_-:=XX^*\;.$$
We consider $\ssa{g}$ as an $E$-C**-subalgebra of $\ssa{f}$ \emph{(Corollary \ref{776} e))}.
\begin{enumerate}
\item $ab=ba$, $\gamma ^2=-\alpha _1^*\alpha _2^*$. We put $c:=ab=\omega (1,1)$.
\item $X,\,V_c,\,P_\pm \in \ssa{g}^c$.
\item We have
$$P_\pm =\frac{1}{2}(V_1\pm (\gamma \widetilde\otimes 1_K)V_c)\in Pr\,\ssa{f}\,,\qquad P_++P_-=V_1\,,\qquad P_+P_-=0,$$
$$X^2=0,\; XP_+=X,\; P_-X=X,\;P_+X=XP_-=0,\; X+X^*\in \unn{\ssa{f}}\;.$$
\item The map
$$\mad{E}{P_\pm \ssa{f}P_\pm }{x}{P_\pm (x\widetilde\otimes 1_K)P_\pm }$$
is an injective unital C**-homomorphism. We identify $E$ with its image with respect to this map and consider $P_\pm \ssa{f}P_\pm $ as an $E$-C**-algebra.
\item The map
$$\mae{\varphi _\pm }{\ssa{g}}{P_\pm \ssa{f}P_\pm }{Y}{P_\pm YP_\pm =P_\pm Y=YP_\pm }$$
is an injective unital C**-homomorphism. If $Y_1,Y_2\in \unn{\ssa{g}}$ then $\varphi _+Y_1+\varphi _-Y_2\in \unn{\ssa{f}}$.
\item The map
$$\mae{\psi }{\ssa{f}}{\ssa{f}}{Z}{(X+X^*)Z(X+X^*)}$$
is an $E$-C**-isomorphism such that
$$\psi ^{-1}=\psi \,,\quad \psi (P_+\ssa{f}P_+)=P_-\ssa{f}P_-\,,\quad \psi \circ \varphi _+=\varphi _-\,,\quad \psi \circ \varphi _-=\varphi _+\;.$$
If $Y_1,Y_2\in \ssa{g}$ then 
$$\varphi _+Y_1+\varphi _-Y_2=(\varphi _+Y_1+\varphi _-V_1)\psi (\varphi _+Y_2+\varphi _-V_1)\;.$$
\item If $p\in Pr\,\ssa{g}$ then
$$(X(\varphi _+p)^*(X(\varphi _+p))=\varphi _+p\,,\qquad (X(\varphi _+p))(X(\varphi _+p))^*=\varphi _-p\;.$$
\item Let $R$ be the subgroup $\{1,a,b,c\}$ of $T$, $h:=f|(R\times R)$, $d\in T$ such that $\bar d=\bnn{2n-2}$ and $d_i=a_i$ for every $i\in \bnn{2n-2}$, and
$$\alpha :=f(d,d)\,,\qquad \alpha ':=f_{2n-1}(2n-1,2n-1)\,,\qquad \alpha '':=f_{2n}(2n,2n)\;.$$
Then $\alpha _1=\alpha \alpha '$, $\alpha _2=\alpha \alpha ''$, $-\alpha '\alpha ''=(\alpha ^*\gamma ^*)^2$,
$$\begin{array}{||c||c|c|c||}\hline \hline
h&a&b&c\\\hline\hline a&\alpha \alpha '&\alpha &\alpha ' \\\hline b&-\alpha  &\alpha \alpha ''&-\alpha ''\\\hline c&-\alpha '&\alpha ''&-\alpha '\alpha '' \\\hline\hline
\end{array}$$
is the table of $h$, $P_\pm \in Pr\, \ssa{h}$, and the map
$$\mae{\varphi }{\ssa{h}}{E_{2,2}}{Z}{\mt{Z_0+\gamma ^*Z_c}{\alpha \alpha 'Z_a-\alpha \gamma ^*Z_b}{Z_a+\alpha '^*\gamma ^*Z_b}{Z_0-\gamma ^*Z_c}}$$
is an $E$-C**-isomorphism. In particular 
$$\varphi P_+=\mt{1_E}{0}{0}{0}\,,\qquad\qquad \varphi P_-=\mt{0}{0}{0}{1_E}$$
and $E_{2,2}$ is $E$-C**-isomorphic to an $E$-C**-subalgebra of $\ssa{f}$.
\item Assume $I=\bnn{2n}$ and $T_{2n-1}=T_{2n}=\bzz{2}$. Then $T\approx S\times \bzz{2}\times \bzz{2}$, $\varphi _\pm $ is an $E$-C*-isomorphism with inverse
$$\mad{P_\pm \ssa{f}P_\pm }{\ssa{f_0}}{Z}{2\si{u\in T_0}(Z_u\otimes 1_K)V_u}\,,$$
and $\ssa{f}\approx _E\ssa{g}_{2,2}$
\end{enumerate}
\end{theo}

a) is easy to see.

b) follows from \pr \ref{23} b).

c) follows from a) and \h \ref{895} b),h).

d) follows from \h \ref{895} c).

e) By b) and c), the map is well-defined. The assertion follows now from \h \ref{895} d),h).

f) follows from b),c), and \h \ref{895} h).

g) follows from b) and \pr \ref{20} d).

h) follows from c), d), \pr \ref{807} a), Corollary \ref{809} d), and \pr \ref{10} c).

i) follows from \h \ref{895} f).\qed

\begin{p}\label{28}
We use the notation and the hypotheses of \emph{\h \ref{27}} and assume $I:=\bnn{2}$, $T_1:=\bzz{2}$, and $T_2:=\bzz{2m}$ with $m\in \bn$.
\begin{enumerate}
\item $a=(1,0)$, $b=(0,m)$, $c=(1,m)$, $\alpha =1_E$, $\alpha '=\alpha _1=f_1(1,1)$, $\alpha ''=\alpha _2=f_2(m,m)$, and
$$P_\pm \ssa{f}P_\pm =\me{(x\widetilde\otimes 1_K)P_\pm }{x\in E}\;.$$
\item If $m=1$ then there are $\alpha ,\beta ,\gamma ,\delta \in \un{E}$ such that $f$ is given by the following table:
 $$\begin{array}{||c||c|c|c|c|c|c|c||}\hline \hline
f&(0,1)&(0,2)&(0,3)&(1,0)&(1,1)&(1,2)&(1,3)\\\hline\hline (0,1)&\alpha &\beta &\gamma &-1_E&-\alpha &-\beta&-\gamma\\\hline (0,2)&\beta &\alpha ^*\beta \gamma &\alpha^*\gamma &-1_E&-\beta&-\alpha ^*\beta \gamma &-\alpha ^*\gamma \\\hline(0,3)&\gamma &\alpha ^*\gamma &\beta ^*\gamma &-1_E&-\gamma&-\alpha ^*\gamma &-\beta^*\gamma \\\hline (1,0)&1_E&1_E&1_E&\delta &\delta &\delta &\delta \\\hline (1,1)&\alpha &\beta &\gamma &-\delta &- \alpha \delta &-\beta \delta &-\gamma \delta \\\hline (1,2)&\beta &\alpha ^*\beta \gamma &\alpha ^*\gamma &-\delta &-\beta \delta &-\alpha ^*\beta \gamma \delta &-\alpha ^*\gamma \delta \\\hline (1,3)&\gamma &\alpha ^*\gamma &\beta ^*\gamma &-\delta &-\gamma \delta &-\alpha ^*\gamma \delta &-\beta ^*\gamma \delta \\\hline\hline
\end{array}\;.$$
\item We assume $\bk:=\bc$ and $m:=1$ and put for all $j,k\in \{0,1\}$
$$\mae{\varphi _{j,k}}{\ssa{f}}{E}{Z}{Z_0+(-1)^jZ_b+i^jZ_{(k,1)}-i^jZ_{(k,3)}}\,,$$
$$\mae{\phi }{\ssa{f}}{E^4}{Z}{(\varphi _{0,0}Z,\,\varphi _{0,1}Z,\,\varphi _{1,0}Z,\,\varphi _{1,1}Z)}\;.$$
If we take $\alpha :=\beta :=\gamma :=-\delta :=\beta _1:=\beta _2:=1_E$ in b) then the map
$$\mad{\ssa{f}}{E_{2,2}\times E^4}{Z}{\left(\mt{Z_0+Z_{(1,2)}}{Z_{(1,0)}-Z_b}{Z_{(1,0)}+Z_b}{Z_0-Z_{(1,2)}},\phi Z\right)}$$
is an $E$-C**-isomorphism.
 \end{enumerate}
\end{p}

a) Use Corollary \ref{29} and \pr \ref{674} b).

b) Use \pr \ref{6437} a) and \pr \ref{22}.

c) follows from b) and \pr \ref{6437} $f_1$.\qed

\begin{center}
\subsection{A special case}
\end{center}

\fbox{\parbox{12 cm}{Throughout this subsection we denote by $S$ a totally ordered set, put $T:=(\bzz{2})^{(S)}$, and fix a map $\rho :S\rightarrow \un{E}$. We define for every $s\in S$, $f_s\in \f{\bzz{2}}{E}$ by putting $f_s(1,1)=\rho (s)\;$( \pr \ref{785} a)). Moreover we denote by $f_\rho $ the Schur function $f$ defined in \pr \ref{22} b) (with $I$ replaced by $S$) and put $\ssc{\rho }:=\ssa{f_\rho }$.}}

{\it Remark.} If $S:=\bnn{2}$ then $T=\bzz{2}\times \bzz{2}$ so $\ssc{\rho }$ is a special case of the example treated in subsection 3.2. With the notation used in the left table of \pr \ref{807} this case appears for $a:=(1,0)$ and $b:=(0.1)$ exactly when $\varepsilon =-1_E$, $\alpha =-\rho (b)$, $\beta =\rho (a)$, and $\gamma =1_E$.

\begin{lem}\label{833}
$\fr{P}_f(S)$ endowed with the composition law
$$\mad{\fr{P}_f(S)\times \fr{P}_f(S)}{\fr{P}_f(S)}{(A,B)}{A\triangle B:=(A\setminus B)\cup (B\setminus A)}$$
is a locally finite commutative group (\emph{\dd  $\;$\ref{779}}) with $\emptyset $ as neutral element and the map
$$\mad{\fr{P}_f(S)}{T}{A}{e_A}$$ 
is a group isomorphism with inverse
$$\mad{T}{\fr{P}_f(S)}{x}{\me{s\in S}{x(s)=1}}\;.$$
We identify $T$ with $\fr{P}_f(S)$ by using this isomorphism and write $s$ instead of $\{s\}$ for every $s\in S$. For $A,B\in T$,
$$f_\rho (A,B)=(-1)^\tau \pro{s\in A\cap B}\rho (s)\,,$$
where  $\tau $ is defined in \emph{\pr \ref{22} b)}.\qed
\end{lem}

\begin{p}\label{835}
Assume $S$ finite and let $F$ be an $E$-C*-algebra. Let further $(x_s)_{s\in S}$ be a family in $F$ such that for all distinct $s,t\in S$ and for every $y\in E$,
$$x_sx_t=-x_tx_s\,,\qquad x_s^2=\rho (s)1_F\,,\qquad x_s^*=\rho (s)^*x_s\,,\qquad x_sy=yx_s\;.$$
Then there is a unique $E$-C*-homomorphism $\varphi :\ssc{\rho }\rightarrow F$ such that $\varphi V_s=x_s$ for all $s\in S$. If the family $\left(\pro{s\in A}x_s\right)_{A\subset S}$ is $E$-linearly independent (resp. generates $F$ as an $E$-C*-algebra) then $\varphi $ is injective (resp. surjective).
\end{p}

Put
$$\varphi V_A:=x_{s_1}x_{s_2}\cdots x_{s_m}$$
for every $A:=\{s_1,s_2,\cdots ,s_m\}$, where $s_1<s_2<\cdots <s_m$, and
$$\mae{\varphi }{\ssc{\rho }}{F}{X}{\si{A\subset S}X_A\,\varphi V_A}\;.$$
It is easy to see that $(\varphi V_s)(\varphi V_t)=\varphi (V_sV_t)$ and $y\,\varphi V_s=(\varphi V_s)y$ for all $s,t\in S$ and $y\in E$ (\pr \ref{674} b)). Let
$$A:=\{s_1,s_2,\cdots ,s_m\}\subset S\,,\; B:=\{t_1,t_2,\cdots ,t_n\}\subset S\,,\; \{r_1,r_2,\cdots ,r_p\}:=A\triangle B\,,$$
where the letters are written in strictly increasing order. Then
$$(\varphi V_A)(\varphi V_B)=x_{s_1}x_{s_2}\cdots x_{s_m}x_{t_1}x_{t_2}\cdots x_{t_n}=f_\rho (A,B)x_{r_1}x_{r_2}\cdots x_{r_p}=$$
$$=f_\rho (A,B)\varphi V_{A\triangle B}=\varphi ((f_\rho (A,B)\widetilde\otimes 1_K)V_{A\triangle B})=\varphi (V_AV_B)\,,$$
$$(\varphi V_A)^*=x_{s_m}^*\cdots x_{s_2}^*x_{s_1}^*=(-1)^{\frac{m(m-1)}{2}}x_{s_1}^*x_{s_2}^*\cdots  x_{s_m}^*=$$
$$=(-1)^{\frac{m(m-1)}{2}}\pro{i\in \bnn{m}}\rho (s_i)^*x_{s_1}x_{s_2}\cdots x_{s_m}=
(-1)^{\frac{m(m-1)}{2}}\pro{i\in \bnn{m}}\rho (s_i)^*\varphi V_A=$$
$$=\varphi ((-1)^{\frac{m(m-1)}{2}}((\pro{i\in \bnn{m}}\rho (s_i)^*)\widetilde\otimes 1_K)V_A)=\varphi (V_A^*)$$
by \pr \ref{23} e).

For $X,Y\in \ssc{\rho }$ (by \h \ref{745} c),g)),
$$(\varphi X)(\varphi Y)=\left(\si{A\in T}X_A\varphi V_A\right)\left(\si{B\in T}Y_B\varphi V_B\right)=\si{A,B\in T}X_AY_B(\varphi V_A)(\varphi V_B)=$$
$$=\si{A,B\in T}X_AY_B\varphi (V_AV_B)=\si{A,B\in T}X_AY_Bf_\rho (A,B)\varphi V_{A\triangle B}=$$
$$=\si{A,C\in T}X_AY_{A\triangle C}f_\rho (A,A\triangle C)\varphi V_C
=\si{C\in T}\left(\si{A\in T}f_\rho (A,A\triangle C)X_AY_{A\triangle C}\right)\varphi V_C=$$
$$=\si{C\in T}(XY)_C\varphi V_C=\varphi (XY)\,,$$
$$(\varphi X)^*=\si{A\in T}X_A^*(\varphi V_A)^*=\si{A\in T}X_A^*\varphi (V_A)^*=$$
$$=\si{A\in T}\tilde f_\rho (A)^*(X^*)_A\tilde f_\rho (A)\varphi V_A=\si{A\in T}(X^*)_A\varphi V_A=\varphi (X^*)$$
(\pr \ref{23} e)) i.e. $\varphi $ is an $E$-C*-homomorphism. The uniqueness and the last assertions are obvious (by \h \ref{745} a)).\qed

\begin{p}\label{851}
Let $m,n\in \bn\cup \{0\}$, $S:=\bnn{2n}$, $S':=\bnn{2n+m}$, $K':=l^2(\fr{P}(S'))$, $(\alpha _i)_{i\in \bnn{m}}\in (\un{E})^m$,
$$\mae{\rho '}{S'}{\un{E}}{s}{\ab{\rho (s)}{s\in S}{\alpha _i^2\tilde f_\rho (S)}{s=2n+i \; \mbox{with}\; i\in \bnn{m}}}\,,$$
and $A_i:=A\cup \{2n+i\}$ for every $A\subset S$ and $i\in \bnn{m}$.
\begin{enumerate}
\item $i\in \bnn{m}\Longrightarrow \tilde f_{\rho '}(S_i)=\alpha _i^{*2}\,,\;\;(V_{S_i}^{\rho '})^2=(\alpha _i^2\otimes 1_{K'})V_\emptyset ^{\rho'}$.
\item $P:=\frac{1}{2}V_\emptyset ^{\rho '}+\frac{1}{2\sqrt{m}}\si{i\in \bnn{m}}(\alpha _i^*\otimes 1_{K'})V_{S_i}^{\rho '}\in Pr\,\ssc{\rho '}$.
\item There is a unique injective $E$-C*-homomorphism $\varphi :\ssc{\rho }\rightarrow P\ssc{\rho '}P$ such that
$$\varphi V_s^\rho =PV_s^{\rho '}P=PV_s^{\rho '}=V_s^{\rho '}P$$
for every $s\in S$.
\item If $m\in \bnn{2}$ then $\varphi $ is an $E$-C*-isomorphism.
\end{enumerate}
\end{p}

a) By \pr \ref{23} e),
$$\tilde f_{\rho '}(S_i)=(-1)^{n(2n+1)}\pro{s\in S_i}\rho '(s)^*=\left((-1)^{n(2n-1)}\pro{s\in S}\rho (s)^*\right)\alpha _i^{*2}\tilde f_\rho (S)^*=\alpha _i^{*2}\,,$$
$$(V_{S_i}^{\rho '})^2=(\alpha _i^2\otimes 1_{K'})V_\emptyset ^{\rho '}\;.$$

b) follows from a) and Corollary \ref{850}.

c) By \pr \ref{23} c), for $s\in S$, $V_s^{\rho '}V_{S_i}^{\rho '}=V_{S_i}^{\rho '}V_s^{\rho '}$ for every $i\in \bnn{m}$ so $V_s^{\rho '}P=PV_s^{\rho '}$. By b),
for distinct $s,t\in S$ (\pr \ref{23} b)),
$$(PV_s^{\rho '})(PV_t^{\rho '})=PV_s^{\rho '}V_t^{\rho '}=-PV_t^{\rho '}V_s^{\rho '}=-(PV_t^{\rho '})(PV_s^{\rho '})\,,$$
$$(PV_s^{\rho '})^2=P(V_s^{\rho '})^2=P(\rho '(s)\otimes 1_{K'})V_\emptyset ^{\rho '}=(\rho (s)\otimes 1_{K'})P\,,$$
$$(PV_s^{\rho '})^*=P(V_s^{\rho '})^*=P(\rho '(s)^*\otimes 1_{K'})V_s^{\rho '}=(\rho (s)\otimes 1_{K'})^*PV_s^{\rho '}\;.$$
By \pr \ref{835} there is a unique $E$-C*-homomorphism $\varphi :\ssc{\rho }\rightarrow P\ssc{\rho '}P$ with the given properties.

Let $X\in \ssc{\rho }$ with $\varphi X=0$. Then
$$0=\left(\si{A\subset S}(X_A\otimes 1_{K'})V_A^{\rho '}\right)P=$$
$$=\frac{1}{2}\si{A\subset S}(X_A\otimes 1_{K'})V_A^{\rho '}+\frac{1}{2\sqrt{m}}\si{i\in \bnn{m}}\si{A\subset S}(X_A\otimes 1_{K'})f_{\rho '}(A,S_i)V_{A\triangle S_i}^{\rho '}$$
and this implies $X_A=0$ for all $A\subset S$ (\h \ref{745} a)). Thus $\varphi $ is injective.

d) \vspace{-6.4ex}
\begin{center}
The case $m=1$
\end{center}
Let $Y\in P\ssc{\rho '}P$. Then (by \pr \ref{674} b))
$$Y=YP=\frac{1}{2}Y+\frac{1}{2}\si{A\subset S'}(\alpha _1^*\otimes 1_{K'})V_{S_1}^{\rho '}Y\,,$$
$$Y=\si{A\subset S}((\alpha _1^*f_{\rho '}(S_1,A)Y_A)\otimes 1_{K'})V_{S_1\triangle A}^{\rho '}+$$
$$+\si{A\subset S}(((\alpha _1^*f_{\rho '}(S_1,A_1)Y_{A_1}))\otimes 1_{K'})V_{S\triangle A}^{\rho '}$$
so
$$\left\{
\begin{array}{l}Y_A=\alpha _1^*f_{\rho '}(S_1,(S\triangle A)_1)Y_{(S\triangle A)_1}\\Y_{A_1}=\alpha _1^*f_{\rho '}(S_1,S\triangle A)Y_{S\triangle A}
\end{array}\right.$$
for every $A\subset S$. If we put
$$X:=2\si{A\subset S}(Y_A\otimes 1_K)V_A^\rho \in \ssc{\rho }$$
then
$$\varphi X=\frac{1}{2}\varphi X+\si{A\subset S}((\alpha _1^*f_{\rho '}(S_1,A)Y_A)\otimes 1_{K'})V_{S_1\triangle A}^{\rho '}=$$
$$=\si{A\subset S}(Y_A\otimes 1_{K'})V_A^{\rho '}+\si{A\subset S}((\alpha _1^*f_{\rho '}(S_1,S\triangle A)Y_{S\triangle A})\otimes 1_{K'})V_{A_1}^{\rho '}=$$
$$=\si{A\subset S}(Y_A\otimes 1_{K'})V_A^{\rho '}+\si{A\subset S}(Y_{A_1}\otimes 1_{K'})V_{A_1}^{\rho '}=Y\;.$$
Thus $\varphi $ is surjective.

\begin{center}
The case $m=2$
\end{center}

Let $Y\in P\ssc{\rho '}P$. Then
$$\left\{
\begin{array}{l}
Y=PY=\frac{1}{2}Y+\frac{1}{2\sqrt{2}}((\alpha _1^*\otimes 1_{K'})V_{S_1}^{\rho '}+(\alpha _2^*\otimes 1_{K'})V_{S_2}^{\rho '})Y\\
Y=YP=\frac{1}{2}Y+\frac{1}{2\sqrt{2}}Y((\alpha _1^*\otimes 1_{K'})V_{S_1}^{\rho '}+(\alpha _2^*\otimes 1_{K'})V_{S_2}^{\rho '})\,,
\end{array}
\right.$$
$$\sqrt{2}Y=(\alpha _1^*\otimes 1_{K'})V_{S_1}^{\rho '}Y+(\alpha _2^*\otimes 1_{K'})V_{S_2}^{\rho '}Y=(\alpha _1^*\otimes 1_{K'})YV_{S_1}^{\rho '}+(\alpha _2^*\otimes 1_{K'})YV_{S_2}^{\rho '}\;.$$

For every $B\subset S$ put
$$B_a:=B\cup \{2n+1\}\,,\qquad B_b:=B\cup \{2n+2\}\,,\qquad B_c:=B\cup \{2n+1\,,\;2n+2\}\;.$$
Then
$$V_{S_1}^{\rho '}Y=\si{B\subset S}((Y_Bf_{\rho '} (S_1,B))\otimes 1_{K'})V_{(S\triangle B)_a}^{\rho '}+\si{B\subset S}((Y_{B_a}f_{\rho '} (S_1,B_a))\otimes 1_{K'})V_{S\triangle B}^{\rho '}+$$
$$+\si{B\subset S}((Y_{B_b}f_{\rho '} (S_1,B_b))\otimes 1_{K'})V_{(S\triangle B)_c}^{\rho '}+\si{B\subset S}((Y_{B_c}f_{\rho '} (S_1,B_c))\otimes 1_{K'})V_{(S\triangle B)_b}\,,$$
$$V_{S_2}^{\rho '}Y=$$
$$=\si{B\subset S}((Y_Bf_{\rho '} (S_2,B))\otimes 1_{K'})V_{(S\triangle B)_b}^{\rho '}+\si{B\subset S}((Y_{B_a}f_{\rho '} (S_2,B_a))\otimes 1_{K'})V_{(S\triangle B)_c}^{\rho '}+$$
$$+\si{B\subset S}((Y_{B_b}f_{\rho '} (S_2,B_b))\otimes 1_{K'})V_{S\triangle B}^{\rho '}+\si{B\subset S}((Y_{B_c}f_{\rho '} (S_2,B_c))\otimes 1_{K'})V_{(S\triangle B)_a}^{\rho '}\,,$$
$$YV_{S_1}^{\rho '}=\si{B\subset S}((Y_Bf_{\rho '} (B,S_1))\otimes 1_{K'})V_{(S\triangle B)_a}^{\rho '}+\si{B\subset S}((Y_{B_a}f_{\rho '} (B_a,S_1))\otimes 1_{K'})V_{S\triangle B}^{\rho '}+$$
$$+\si{B\subset S}((Y_{B_b}f_{\rho '} (B_b,S_1))\otimes 1_{K'})V_{(S\triangle B)_c}^{\rho '}+\si{B\subset S}((Y_{B_c}f_{\rho '} (B_c,S_1))\otimes 1_{K'})V_{(S\triangle B)_b}^{\rho '}\,,$$
$$YV_{S_2}^{\rho '}=$$
$$=\si{B\subset S}((Y_Bf_{\rho '} (B,S_2))\otimes 1_{K'})V_{(S\triangle B)_b}^{\rho '}+\si{B\subset S}((Y_{B_a}f_{\rho '} (B_a,S_2))\otimes 1_{K'})V_{(S\triangle B)_c}^{\rho '}+$$
$$+\si{B\subset S}((Y_{B_b}f_{\rho '} (B_b,S_2))\otimes 1_{K'})V_{S\triangle B}^{\rho '}+\si{B\subset S}((Y_{B_c}f_{\rho '} (B_c,S_2))\otimes 1_{K'})V_{(S\triangle B)_a}^{\rho '}\,,$$
$$\sqrt{2}Y=\si{B\subset S}((\alpha _1^*Y_{B_a}f_{\rho '} (S_1,B_a)+\alpha _2^*Y_{B_b}f_{\rho '} (S_2,B_b))\otimes 1_{K'})V_{S\triangle B}^{\rho '}+$$
$$+\si{B\subset S}((\alpha _1^*Y_Bf_{\rho '} (S_1,B)+\alpha _2^*Y_{B_c}f_{\rho '} (S_2,B_c))\otimes 1_{K'})V_{(S\triangle B)_a}^{\rho '}+$$
$$+\si{B\subset S}((\alpha _1^*Y_{B_c}f_{\rho '} (S_1,B_c)+\alpha _2^*Y_{B}f_{\rho '} (S_2,B))\otimes 1_{K'})V_{(S\triangle B)_b}^{\rho '}+$$
$$+\si{B\subset S}((\alpha _1^*Y_{B_b}f_{\rho '} (S_1,B_b)+\alpha _2^*Y_{B_a}f_{\rho '} (S_2,B_a))\otimes 1_{K'})V_{(S\triangle B)_c}^{\rho '}\,,$$
$$\sqrt{2}Y=\si{B\subset S}((\alpha _1^*Y_{B_a}f_{\rho '} (B_a,S_1)+\alpha _2^*Y_{B_b}f_{\rho '} (B_b,S_2))\otimes 1_{K'})V_{S\triangle B}^{\rho '}+$$
$$+\si{B\subset S}((\alpha _1^*Y_Bf_{\rho '} (B,S_1)+\alpha _2^*Y_{B_c}f_{\rho '} (B_c,S_2))\otimes 1_{K'})V_{(S\triangle B)_a}^{\rho '}+$$
$$+\si{B\subset S}((\alpha _1^*Y_{B_c}f_{\rho '} (B_c,S_1)+\alpha _2^*Y_{B}f_{\rho '} (B,S_2))\otimes 1_{K'})V_{(S\triangle B)_b}^{\rho '}+$$
$$+\si{B\subset S}((\alpha _1^*Y_{B_b}f_{\rho '} (B_b,S_1)+\alpha _2^*Y_{B_a}f_{\rho '} (B_a,S_2))\otimes 1_{K'})V_{(S\triangle B)_c}^{\rho '}\;.$$
It follows for $B\subset S$,
$$\sqrt{2}Y_{B_a}=\alpha _1^*Y_{S\triangle B}f_{\rho '} (S\triangle B,S_1)+\alpha _2^*Y_{(S\triangle B)_c}f_{\rho '} ((S\triangle B)_c,S_2)\,,$$
$$\sqrt{2}Y_{B_b}=\alpha _1^*Y_{(S\triangle B)_c}f_{\rho '} ((S\triangle B)_c,S_1)+\alpha _2^*Y_{S\triangle B}f_{\rho '} (S\triangle B,S_2)\,,$$
$$\sqrt{2}Y_{B_c}=\alpha _1^*Y_{(S\triangle B)_b}f_{\rho '} (S_1,(S\triangle B)_b)+\alpha _2^*Y_{(S\triangle B)_a}f_{\rho '} (S_2,(S\triangle B)_a)=$$
$$=\alpha _1^*Y_{(S\triangle B)_b}f_{\rho '} ((S\triangle B)_b,S_1)+\alpha _2^*Y_{(S\triangle B)_a}f_{\rho '} ((S\triangle B)_a,S_2)\,,$$
so by \pr \ref{23} a),b), $Y_{B_c}=0$. If we put
$$X:=2\si{B\subset S}(Y_B\otimes 1_K)V_B^\rho \in \ssc{\rho }$$
then
$$\varphi X=\left(2\si{B\subset S}(Y_B\otimes 1_{K'})V_B^{\rho '}\right)P=$$
$$=\si{B\subset S}(Y_B\otimes 1_{K'})V_B^{\rho '}+\frac{1}{\sqrt{2}}\si{B\subset S}((\alpha _1^*Y_Bf_{\rho '} (B,S_1))\otimes 1_{K'})V_{S_1\triangle B}^{\rho '}+$$
$$+\frac{1}{\sqrt{2}}\si{B\subset S}((\alpha _2^*Y_Bf_{\rho '} (B,S_2))\otimes 1_{K'})V_{S_2\triangle B}^{\rho '}\,,$$
and so for $B\subset S$,
$$(\varphi X)_B=Y_B\,,\qquad (\varphi X)_{B_a}=\frac{1}{\sqrt{2}}\alpha _1^*Y_{S\triangle B}f_{\rho '} (S\triangle B,S_1)=Y_{B_a}\,,$$
$$(\varphi X)_{B_b}=\frac{1}{\sqrt{2}}\alpha _2^*Y_{S\triangle B}f_{\rho '} (S\triangle B,S_2)=Y_{B_b}\,,\qquad (\varphi X)_{B_c}=0=Y_{B_c}\;.$$
Thus $\varphi X=Y$ and $\varphi $ is surjective.\qed

{\it Remark.} 
If $m=3$ then $\varphi $ may be not surjective. 

\begin{p}\label{838}
Let $\bk:=\br$, $n\in \bn\cup \{0\}$, $S:=\bnn{2n}$, and
$$\mae{\rho '}{\bnn{2n+1}}{\un{E}}{s}{\ab{\rho (s)}{s\in S}{-\tilde f_\rho (S)}{s=2n+1}}\;.$$
Let $\overbrace{\ssc{\rho }}^\circ $ be the complexification of $\ssc{\rho }$, considered as a real $E$-C*-algebra \emph{([C1] \h 4.1.1.8 a))} by using the embedding
$$\mad{E}{\overbrace{\ssc{\rho }}^\circ }{x}{((x\otimes 1_K)V_\emptyset ^\rho ,0)}\;.$$
Then there is a unique $E$-C*-isomorphism $\varphi :\ssc{\rho '}\rightarrow \overbrace{\ssc{\rho }}^\circ $ such that $\varphi V_s^{\rho '}=(V_s^\rho ,0)$ for every $s\in S$ and 
$$\varphi V_{2n+1}^{\rho '}=(0,-(\tilde f_\rho (S)\otimes 1_K)V_S^\rho )\;.$$
\end{p}

We put
$$x_s:=\ab{(V_s^\rho ,0)}{s\in S}{(0,-(\tilde f_\rho (S)\otimes 1_K)V_S^\rho )}{s=2n+1}\;.$$
For $s\in S$, by \pr \ref{23} b),
$$x_sx_{2n+1}=(V_s^\rho ,0)(0,-(\tilde f_\rho (S)\otimes 1_K)V_S^\rho )=(0,-(\tilde f_\rho (S)\otimes 1_K)V_s^\rho V_S^\rho )=$$
$$=(0,(\tilde f_\rho (S)\otimes 1_K)V_S^\rho V_s^\rho )=(0,(\tilde f_\rho (S)\otimes 1_K)V_s^\rho )(V_s^\rho,0 )=-x_{2n+1}x_s\;.$$
By \pr \ref{674} b),e),
$$x_{2n+1}^2=(-((\tilde f_\rho (S)\otimes 1_K)V_S^\rho )^2,0)=$$
$$=(-(\tilde f_\rho (S)^2\otimes 1_K)(f_\rho (S,S)\otimes 1_K)V_\emptyset ^\rho ,0)=(\rho '(2n+1)\otimes 1_K)(V_\emptyset ^\rho ,0)\,,$$
$$x_{2n+1}^*=(0,((\tilde f_\rho (S)\otimes 1_K)V_S^\rho )^*)=$$
$$=(0,(\tilde f_\rho (S)^*\otimes 1_K)(\tilde f_\rho (S)\otimes 1_K)V_S ^\rho )=(\rho '(2n+1)^*\otimes 1_K)x_{2n+1}\,,$$
and the assertion follows from \pr \ref{835}.\qed

\begin{p}\label{857}
Let $n\in \bn\cup \{0\}$, $S:=\bnn{n}$, $S':=\bnn{n+2}$, $K':=l^2(\fr{P}(S'))$, $\alpha _1,\alpha _2\in \un{E}$, and
$$\mae{\rho '}{S'}{\un{E}}{s}{\ac{\rho (s)}{s\in S}{\alpha _1^2}{s=n+1}{-\alpha _2^2}{s=n+2}}\;.$$
\begin{enumerate}
\item There is a unique $E$-C*-isomorphism $\varphi : \ssc{\rho '}\rightarrow \ssc{\rho }_{2,2}$ such that
$$\varphi V_s^{\rho '}=\mt{V_s^\rho }{0}{0}{-V_s^\rho }$$
for every $s\in S$ and
$$\varphi V_{n+1}^{\rho '}=(\alpha _1\otimes 1_K)\mt{0}{V_\emptyset ^\rho }{V_\emptyset ^\rho }{0},
\; \varphi V_{n+2}^{\rho '}
=(\alpha _2\otimes 1_K)\mt{0}{-V_\emptyset ^\rho }{V_\emptyset ^\rho }{0}\;.$$
\item $$\varphi \frac{1}{2} (V_\emptyset ^{\rho '}+((\alpha _1^*\alpha _2^*)\otimes 1_{K'})V_{\{n+1,\,n+2\}}^{\rho '})=\mt{V_\emptyset ^\rho }{0}{0}{0}\,,$$
$$\varphi \frac{1}{2}(V_\emptyset ^{\rho '}-((\alpha _1^*\alpha _2^*)\otimes 1_{K'})V_{\{n+1,\,n+2\}}^{\rho '})=\mt{0}{0}{0}{V_\emptyset ^\rho }\;.$$
\end{enumerate}
\end{p}

a) Put 
$$x_s:=\mt{V_s^\rho }{0}{0}{-V_s^\rho }$$
for every $s\in S$ and
$$x_{n+1}:=(\alpha _1\otimes 1_K)\mt{0}{V_\emptyset ^\rho }{V_\emptyset ^\rho }{0}\,,\quad x_{n+2}:=(\alpha _2\otimes 1_K)\mt{0}{-V_\emptyset ^\rho }{V_\emptyset ^\rho }{0}\;.$$
For distinct $s,t\in S$ and $i\in \bnn{2}$,
$$x_sx_t=-x_tx_s\,,\quad x_s^2=(\rho '(s)\otimes 1_K)\mt{V_\emptyset ^\rho }{0}{0}{V_\emptyset ^\rho }\,,\quad x_s^*=(\rho '(s)\otimes 1_K)^*x_s\,,$$
$$x_sx_{n+i}=-x_{n+i}x_s\,,\qquad x_{n+i}^2=(\rho '(n+i)\otimes 1_K)\mt{V_\emptyset ^\rho }{0}{0}{V_\emptyset ^\rho }\,,$$
 $$x_{n+i}^*
=(\rho '(n+i)\otimes 1_K)^*x_{n+i}\,,\qquad 
x_{n+1}\,x_{n+2}=-x_{n+2}\,x_{n+1}\;.$$
By \pr \ref{835} there is a unique $E$-C*-homomorphism $\varphi :\ssc{\rho '}\rightarrow \ssc{\rho }_{2,2}$ satisfying the given conditions.

We put for every $A\subset S$ and $i\in \bnn{2}$
$$|A|:=Card\;A\,,\qquad A_i:=A\cup \{n+i\}\,,\qquad A_3:=A\cup \{n+1,\,n+2\}\;.$$
For $A\subset S$,
$$\varphi V_{A_1}^{\rho '}=(\alpha _1\otimes 1_K)\mt{V_A^\rho }{0}{0}{(-1)^{|A|}V_A^\rho }\mt{0}{V_\emptyset ^\rho }{V_\emptyset ^\rho }{0}=$$
$$=(\alpha _1\otimes 1_K)\mt{0}{V_A^\rho }{(-1)^{|A|}V_A^\rho }{0}\,,$$
$$\varphi V_{A_2}^{\rho '}=(\alpha _2\otimes 1_K)\mt{V_A^\rho }{0}{0}{(-1)^{|A|}V_A^\rho }\mt{0}{-V_\emptyset ^\rho }{V_\emptyset ^\rho }{0}=$$
$$=(\alpha _2\otimes 1_K)\mt{0}{-V_A^\rho }{(-1)^{|A|}V_A ^\rho }{0}\,,$$
$$\varphi V_{A_3}^{\rho '}=((\alpha _1\alpha _2)\otimes 1_K)\mt{0}{V_A^\rho }{(-1)^{|A|}V_A^\rho }{0}\mt{0}{-V_\emptyset ^\rho }{V_\emptyset ^\rho }{0}=$$
$$=((\alpha _1\alpha _2)\otimes 1_K)\mt{V_A^\rho }{0}{0}{-(-1)^{|A|}V_A^\rho }\;.$$
Then for $Y\in \ssc{\rho '}$,
$$\left\{
\begin{array}{l}
(\varphi Y)_{11}=\si{A\subset S}((Y_A+(\alpha _1\alpha _2)Y_{A_3})\otimes 1_K)V_A^\rho \\(\varphi Y)_{12}=\si{A\subset S}((\alpha _1Y_{A_1}-\alpha _2Y_{A_2})\otimes 1_K)V_A^\rho \\(\varphi Y)_{21}=\si{A\subset S}(-1)^{|A|})((\alpha _1Y_{A_1}+\alpha _2Y_{A_2})\otimes 1_K)V_A^\rho \\(\varphi Y)_{22}=\si{A\subset S}(-1)^{|A|}((Y_A-\alpha _1\alpha _2Y_{A_3})\otimes 1_K)V_A^\rho \;.
\end{array}
\right.$$
It follows from the above identities that $\varphi $ is bijective.

b) By the above, 
$$\varphi V_{\{n+1,\,n+2\}}^{\rho '}=\varphi V_{\emptyset _3}^{\rho '}
=((\alpha _1\alpha _2)\otimes 1_K)\mt{V_\emptyset ^\rho }{0}{0}{-V_\emptyset ^\rho }$$
and the assertion follows.\qed

\begin{co}\label{889}
Let $m,n\in \bn\cup \{0\}$, $S:=\bnn{n}$, $(\alpha _i)_{i\in \bnn{2m}}\in  (\un{E})^{2m}$, and 
$$\mae{\rho '}{\bnn{n+2m}}{\un{E}}{s}{\ab{\rho (s)}{s\in S}{-(-1)^i\alpha _i^2}{s=n+i}}\;.$$
Then $\ssc{\rho '}\approx _E\ssc{\rho }_{2^m,2^m}$.\qed
\end{co}

\begin{p}\label{859}
Let $\bk:=\br$, $n\in \bn\cup \{0\}$, $S:=\bnn{2n}$, $S':=\bnn{2n+2}$, $\alpha _1,\alpha _2\in \un{E}$, and
$$\mae{\rho '}{S'}{\un{E}}{s}{\ab{\rho (s)}{s\in S}{-\alpha _l^2\tilde f_\rho (S)}{s=2n+l \;with\; l\in \bnn{2}}}\;.$$
Then there is a unique $E$-C*-isomorphism $\varphi :\ssc{\rho '}\rightarrow \ssc{\rho }\otimes \bh$ such that
$$\varphi V_s^{\rho '}=\ac{V_s^\rho \otimes 1_{\bh}}{s\in S}{(((\alpha _1\tilde f_\rho (S))\otimes 1_K)V_S^\rho )\otimes i}{s=2n+1}{(((\alpha _2\tilde f_\rho (S))\otimes 1_K)V_S^\rho )\otimes j}{s=2n+2}\,,$$
where $i,j,k$ are the canonical unitaries of $\bh$.
\end{p}

Put
$$x_s:=\ac{V_s^\rho \otimes 1_{\bh}}{s\in S}{(((\alpha _1\tilde f_\rho (S))\otimes 1_K)V_S^\rho )\otimes i}{s=2n+1}{(((\alpha _2\tilde f_\rho (S))\otimes 1_K)V_S^\rho )\otimes j}{s=2n+2}\;.$$
For distinct $s,t\in S$ and $l\in \bnn{2}$, by \pr \ref{23} b),
$$x_sx_t=-x_tx_s\,,\qquad x_s^2=(\rho' (s)\otimes 1_K)(V_\emptyset ^\rho \otimes 1_{\bh})\,,\qquad x_s^*=(\rho '(s)\otimes 1_K)^*x_s\,,$$ 
$$x_sx_{2n+l}=-x_{2n+l}x_s,\,
 x_{2n+1}x_{2n+2}=(((\alpha _1\alpha _2\tilde f_\rho (S))\otimes 1_K)V_\emptyset ^\rho )\otimes k=-x_{2n+2}x_{2n+1},$$
$$(x_{2n+l})^2=(((\alpha _l^2\tilde f_\rho (S)^2)\otimes 1_K)(\tilde f_\rho (S)^*\otimes 1_K)V_\emptyset ^\rho )\otimes (-1_{\bh})=$$
$$=(\rho '(2n+l)\otimes 1_K)(V_\emptyset ^\rho \otimes 1_{\bh})\,,$$
$$(x_{2n+l})^*=(((\alpha _l^*\tilde f_\rho (S)^*)\otimes 1_K)(\tilde f_\rho (S)\otimes 1_K)V_S^\rho )\otimes -(i\; \mbox{or}\;j)=$$
$$=(\rho '(2n+l)\otimes 1_K)^*x_{2n+l}\;.$$
By \pr \ref{835} there is a unique $E$-C*-homomorphism $\varphi :\ssc{\rho '}\rightarrow \ssc{\rho }\otimes \bh$ satisfying the given conditions.

For $X\in \ssc{\rho '}$,
$$\varphi X=\left(\si{A\subset S}(X_A\otimes 1_K)V_A^\rho \right)\otimes 1_{\bh}+$$
$$+\left(\si{A\subset S}((X_{A\cup \{2n+1\}}\alpha _1\tilde f_\rho (S)f_\rho (A,S))\otimes 1_K)V_{S\triangle A}\right)\otimes i+$$
$$+\left(\si{A\subset S}((X_{A\cup \{2n+2\}}\alpha _2\tilde f_\rho (S)f_\rho (A,S))\otimes 1_K)V_{S\triangle A}^\rho \right)\otimes j+$$
$$+\left(\si{A\subset S}((X_{A\cup \{2n+1,\,2n+2\}}\alpha _1\alpha _2\tilde f_\rho (S))\otimes 1_K)V_A^\rho \right)\otimes k$$
and so $\varphi $ is bijective.\qed

\begin{p}\label{863}
Let $n\in \bn\cup \{0\}$, $S:=\bnn{2n}$, $A':=A\cup \{2n+1\}$ for every $A\subset S$,
$$\mae{\rho '}{S'}{\un{E}}{s}{\ab{\rho (s)}{s\in S}{\tilde f(S)}{s=2n+1}}\,,$$
$P_\pm :=\frac{1}{2}(V_\emptyset ^{\rho '}\pm V_{S'}^{\rho '})$, and $\theta _\pm :\cb{A\subset S}\breve E\rightarrow \cb{A\subset S'}\breve E$ defined by
$$(\theta _\pm \xi )_A:=\frac{1}{\sqrt{2}}\xi _A\,,\qquad (\theta _\pm \xi )_{A'}:=\pm \frac{1}{\sqrt{2}}f_\rho (S\triangle A,S)\xi _{S\triangle A}$$
for every $\xi \in \cb{A\subset S}\breve E$ and $A\subset S$.
\begin{enumerate}
\item $$\tilde f_{\rho '}(S')=1_E\,, \qquad (V_{S'}^{\rho '})^2=V_\emptyset ^{\rho '}\,,\qquad P_\pm \in Pr\,\ssc{\rho '}^c\,,$$
$$P_++P_-=V_\emptyset ^{\rho '}\,,\qquad V_{S'}^{\rho '}\in \ssc{\rho '}^c\,,\qquad V_{S'}^{\rho '}P_\pm =\pm P_\pm \;.$$
\item For $A\subset S$,
$$f_\rho (A,S)^*=f_{\rho '}(S',A)^*=f_{\rho '}(S',(S\triangle A)')\;.$$
\item $\theta _\pm \in \lc{E}{\cb{A\subset S}\breve E}{\cb{A\subset S'}\breve E}$ and for $\eta \in \cb{A\subset S'}\breve E$ and $A\subset S$,
$$(\theta _\pm ^*\eta )_A=\frac{1}{\sqrt{2}}(\eta _A\pm f_\rho (A,S)^*\eta _{(S\triangle A)'})=\sqrt{2}(P_\pm \eta )_A\;.$$
\item $\theta _\pm ^*\theta _\pm $ is the identity map of $\cb{A\subset S}\breve E$.
\item $\theta _\pm \theta _\pm ^*=P_\pm $.
\item For every $A\subset S$,
$$\theta _\pm V_A^\rho \theta _\pm ^*=V_A^{\rho '}P_\pm =P_\pm V_A^{\rho '}=P_\pm V_A^{\rho '}P_\pm \;.$$
\item For every closed ideal $F$ of $E$ the map
$$\mae{\varphi }{\ssc{\rho ,F}}{P_\pm \ssc{\rho ',F}P_\pm }{X}{\theta _\pm X\theta _\pm ^*}$$
is an $E$-C*-isomorphism with inverse
$$\mad {P_\pm \ssc{\rho ',F}P_\pm }{\ssc{\rho ,F}}{Y}{\theta _\pm ^*Y\theta _\pm }$$
and the map
$$\mae{\psi }{\ssc{\rho ',F}}{\ssc{\rho ,F}\times \ssc{\rho ,F}}{Y}$$
$${(\theta _+^*P_+YP_+\theta _+,\,\theta _-^*P_-YP_-\theta _-)=(\theta _+Y\theta _+,\,\theta _-^*Y\theta _-)}$$
is an $E$-C*-isomorphism.
\end{enumerate}
\end{p}

a) By \pr \ref{23} d),e), $V_{S'}^{\rho '}\in \ssc{\rho '}^c$,
$$\tilde f_{\rho '}(S')=(-1)^{n(2n+1)}\pro{s\in S'}\rho '(s)^*=(-1)^{n(2n-1)}\left(\pro{s\in S}\rho (s)^*\right)\rho '(2n+1)^*=1_E\,,$$
$$(V_{S'}^{\rho '})^*=\tilde f_{\rho '}(S')V_{S'}^{\rho '}=V_{S'}^{\rho '}\,,\qquad (V_{S'}^{\rho '})^2=\tilde f(S')^*V_\emptyset ^{\rho '}=V_\emptyset ^{\rho '}\,,$$
so
$$P_\pm \in Pr\,\ssc{\rho '}^c\,,\qquad V_{S'}^{\rho '}P_\pm =\pm P_\pm \;.$$

b) By a), \pr \ref{23} c),d), \pr \ref{22} b), and \pr \ref{704} b),
$$f_\rho (A,S)^*=f_{\rho '}(A,S)^*=f_{\rho '}(A,S')^*=$$
$$=f_{\rho '}(S',A)^*=f_{\rho '}(S',(S\triangle A)')\tilde f_{\rho '}(S')=f_{\rho '}(S',(S\triangle A)')\;.$$

c) For $\xi \in \cb{A\subset S}\breve E$,
$$\s{\theta \xi }{\eta }=\si{A\subset S}\eta _A^*\frac{1}{\sqrt{2}}\xi _A\pm \si{A\subset S}\eta _{A'}^*\frac{1}{\sqrt{2}}f_\rho (S\triangle A,S)\xi _{S\triangle A}=$$
$$=\si{A\subset S}\eta _A^*\frac{1}{\sqrt{2}}\xi _A\pm \si{A\subset S}\eta _{(S\triangle A)'}^*\frac{1}{\sqrt{2}}f_\rho (A,S)\xi _A=$$
$$=\si{A\subset S}\frac{1}{\sqrt{2}}(\eta _A\pm f_\rho (A,S)^*\eta _{(S\triangle A)'})^*\xi _A$$
so $\theta \in \lc{E}{\cb{A\subset S}\breve E}{\cb{A\subset S'}\breve E}$ and
$$(\theta ^*\eta )_A=\frac{1}{\sqrt{2}}(\eta _A\pm f_\rho (A,S)^*\eta _{(S\triangle A)'})\;.$$
By a) and b),
$$(P_\pm \eta )_A=\frac{1}{2}\eta _A\pm \frac{1}{2}f_{\rho '}(S',(S\triangle A)')\eta _{(S\triangle A)'}=$$
$$=\frac{1}{2}(\eta _A\pm f_\rho (A,S)^*\eta _{(S\triangle A)'})=\frac{1}{\sqrt{2}}(\theta _\pm ^*\eta )_A\;.$$

d) For $\xi \in \cb{A\subset S}\breve E$ and $A\subset S$, by c),
$$(\theta _\pm ^*\theta _\pm \xi )_A=\frac{1}{\sqrt{2}}((\theta \xi )_A\pm f_\rho (A,S)^*(\theta \xi )_{(S\triangle A)'})=$$
$$=\frac{1}{2}(\xi  _A+f_\rho (A,S)^*f_\rho (A,S)\xi _A)=\xi _A\;.$$

e) For $\eta \in \cb{A\subset S'}\breve E$ and $A\subset S$, by b) and c),
$$(\theta _\pm \theta _\pm ^*\eta )_A=\frac{1}{\sqrt{2}}(\theta _\pm ^*\eta )_A=(P_\pm \eta )_A\,,$$
$$(\theta _\pm \theta _\pm ^*\eta )_{A'}=\pm \frac{1}{\sqrt{2}}f_\rho (S\triangle A,S)(\theta _\pm ^*\eta )_{S\triangle A}=$$
$$=\pm \frac{1}{2}f_\rho (S\triangle A,S)(\eta _{S\triangle A}\pm f_\rho (S\triangle A,S)^*\eta _{A'})=\pm \frac{1}{2}f_\rho (S\triangle A,S)\eta _{S\triangle A}+\frac{1}{2}\eta _{A'}=$$
$$=\frac{1}{2}(\eta _{A'}\pm f_{\rho '}(S',S\triangle A)\eta _{S\triangle A})=\frac{1}{2}((V_\emptyset ^{\rho '}\eta )_{A'}\pm (V_{S'}^{\rho '}\eta )_{A'})=(P_\pm \eta )_{A'}\,,$$
so $\theta _\pm \theta _\pm ^*=P_\pm $.

f) For $\eta \in \cb{B\subset S'}\breve E$ and $B\subset S$, by a),b),c),e) and \pr \ref{22} b) (and Corollary \ref{776} e)),
$$(V_A^{\rho '}P_\pm \eta )_B=f_{\rho '}(A,A\triangle B)(P_\pm \eta )_{A\triangle B}=f_\rho (A,A\triangle B)(\theta _\pm \theta _\pm ^*\eta )_{A\triangle B}=$$
$$=\frac{1}{\sqrt{2}}f_\rho (A,A\triangle B)(\theta _\pm ^*\eta )_{A\triangle B}=\frac{1}{\sqrt{2}}(V_A^\rho \theta _\pm ^*\eta )_B=(\theta _\pm V_A^\rho \theta _\pm ^*\eta )_B\,,$$
$$(\theta _\pm V_A^\rho \theta _\pm ^*\eta )_{B'}=\pm \frac{1}{\sqrt{2}}f_\rho (S\triangle B,S)(V_A^\rho \theta _\pm ^*\eta )_{S\triangle B}=$$
$$=\pm \frac{1}{\sqrt{2}}f_\rho (S\triangle B,S)f_\rho (A,S\triangle A\triangle B)(\theta _\pm ^*\eta )_{S\triangle A\triangle B}=$$
$$=\pm f_\rho (S\triangle B,S)f_\rho (A,S\triangle A\triangle B)(P_\pm \eta )_{S\triangle A\triangle B}=\pm f_\rho (S\triangle B,S)(V_A^{\rho '}P_\pm \eta )_{S\triangle B}=$$
$$=\pm f_{\rho '}(S',S'\triangle B')(V_A^{\rho '}P_\pm \eta )_{S'\triangle B'}=\pm (V_{S'}^{\rho '}V_A^{\rho '}P_\pm \eta )_{B'}=$$
$$=\pm (V_A^{\rho '}V_{S'}^{\rho '}P_\pm \eta )_{B'}=(V_A^{\rho '}P_\pm \eta )_{B'}$$
so by a), 
$$\theta _\pm V_A^\rho \theta _\pm ^*=V_A^{\rho '}P_\pm =P_\pm V_A^{\rho '}P_\pm =P_\pm V_A^{\rho '}\;.$$

g) The assertion concerning $\varphi $ as well as the identity in the definition of $\psi $ follow from a),d),e), and f). Thus $\psi $ is a surjective $E$-C*-homomorphism. For $Y\in Ker\,\psi $,
$$\theta _+^*Y\theta _+=\theta _-^*Y\theta _-=0\,,$$
so by a) and e),
$$P_+Y=P_-Y=0$$
and we get
$$Y=P_+Y+P_-Y=0$$
i.e. $\psi $ is injective.\qed

\addtocontents{toc}{REFERENCES}
\begin{center}
REFERENCES
\end{center}

\begin{flushleft}
[C1] Corneliu Constantinescu, {\it C*-algebras.} Elsevir, 2001. \newline
[C2] Corneliu Constantinescu, {\it W*-tensor products and \ris{W^*}s.} Rev. Roumaine Math. Pures Appl., {\bf 51}: 5-6 (2006) 583-596. \newline
[C3] Corneliu Constantinescu, {\it Selfdual \ri{W^*}s and their W*-tensor products.} Rev. Roumaine Math. Pures Appl., {\bf 55}: 3 (2010) 159-196. \newline
[K] Richard V. Kadison and John R. Ringrose, {\it Fundamentals of the theory of operator algebras.} Academic Press, 1983-1986. \newline 
[L] Christopher E. Lance, {\it Hilbert C*-modules. A toolkit for operator algebraist.} Cambridge University Press, 1995. \newline
[S] Issai Schur {\it $\ddot{U} $ber die Darstellung der endlichen Gruppen durch gebrochene lineare Substitutionen.} J. Reine Angew. Math., {\bf 127} (1904) 20-50. \newline
[T] Masamichi Takesaki, {\it Theory of Operator Algebra I.} Springer, 2002. \newline
[W] N. E. Wegge-Olsen, {\it K-theory and C*-algebras.} Oxford University Press, 1993. \newline
\end{flushleft}

\addtocontents{toc}{and INDEX}

\begin{theindex}

\item Schur $E$-function for $T$, $\f{T}{E},\, \tilde f$, $\hat f\,$  (\dd$\,$\ref{703}).
\item $\Lambda (T,E),\, \hat \lambda$, $\delta \lambda $  (\dd$\,$\ref{705}).
\item $E$-module, $E$-linear (Definition \ref{32}).
\item $E$-C*-algebra, $E$-W*-algebra, $E$-C*-subalgebra, $E$-W*-subalgebra, $E$-C*-homomorphism, \\ $E$-W*-homomorphism, $E$-C*-isomorphism, $\approx _E$ (Definition \ref{808}).
\item $\fr{C}_E,\;\fr{C}_E^1$ (Definition \ref{41}).
\item adapted, $\check F$, $\fr{M}_E$ (\pr \ref{33}).
\item $\check \varphi $ (\pr \ref{34}). 
\item $\Phi_E$ (\pr \ref{42}).
\item $\widetilde\ca\,,\,\widetilde\otimes ,\,\widetilde \sum,\,\ccc{G}_{\fr{T}},\,\stackrel{\fr{T}}{\bar {\ccc{G}}},\,\sii{}{\fr{T}}$  (\dd$\,$\ref{673}).
\item $x\widetilde\otimes 1_K$ (Lemma \ref{18}).
\item $\fr{T_1},\,\fr{T_2},\,\fr{T_3}\,$  (\dd$\,$\ref{740}).
\item $u_t,\,V_t,\,V_t^f$  (\dd$\,$\ref{671}).
\item $\varphi _{s,t}\,,\,X_t$ (\dd$\,$\ref{680}).
\item $\ccc{R}(f),\,\ssa{f},\,\ssb{C}{f},\,\ssb{W}{f},\,\ssb{\n{\cdot }}{f}, \\ \ssa{f,F},\,\ssb{\n{\cdot }}{f,F}$ (\dd$\,$\ref{682}).
\item Locally finite, $\fr{S}_T$ (\dd$\,$\ref{779}).
\item $U_\lambda $ (\dd$\,$\ref{720}).
\item $\ccc{S}$-isomorphism, $\approx _{\ccc{S}}$ (\pr \ref{716} a)).
\item $\ssa{F}$, $V_t^F$ (Definition \ref{46}).
\item $\ssa{\varphi }$ (\pr \ref{47}).
\item $\bt$ (\dd$\,$\ref{787}).
\item $w(x)$, winding number of  $x$ (\dd$\,$\ref{787'}).
\item $f_\rho $, $\ssc{\rho }$ Box of subsection 4.2.
\item $\fr{P}_f$, $\triangle $ (\dd$\,$\ref{833}).
\end{theindex}

\begin{flushright}
Corneliu Constantinescu\\
Bodenacherstr. 53\\
CH 8121 Benglen\\
e-mail: constant@math.ethz.ch
\end{flushright}
\end{document}